Министерство науки и высшего образования Российской Федерации
Федеральное государственное автономное
образовательное учреждение высшего образования
«Московский физико-технический институт
(национальный исследовательский университет)»

А. В. Гасников, Е. В. Гасникова

# МОДЕЛИ
# РАВНОВЕСНОГО РАСПРЕДЕЛЕНИЯ ТРАНСПОРТНЫХ ПОТОКОВ В БОЛЬШИХ СЕТЯХ

Учебное пособие

Р е ц е н з е н т ы :

Директор Института системного анализа ФИЦ ИУ РАН
академик РАН *Ю. С. Попков*

Заведующий лабораторией Дальневосточного федерального университета
(г. Владивосток) профессор, д.ф.-м.н. *Е. А. Нурминский*

Учебное пособие основано на семестровом курсе лекций, читаемых студентам шестого курса физтех-школы прикладной математики и информатики для ФУПМ МФТИ и пятого курса физтех-школы аэрокосмических технологий для ФАКИ МФТИ с 2008 года. Пособие посвящено моделям равновесного распределения транспортных потоков по путям и их многостадийным вариантам.

В первой главе описываются элементы теории популяционных игр загрузки (потенциальных игр) на примере задач поиска равновесного распределения транспортных потоков. Потенциальность рассматриваемых игр позволяет сводить поиск равновесий (Нэша) к задачам оптимизации. При весьма естественных условиях эти задачи оказываются выпуклыми, и появляется возможность их эффективно решать. Численные способы решения таких задач описаны во второй главе.

Предназначено для студентов старших курсов, аспирантов и преподавателей МФТИ.

# Оглавление

# Введение

## *Актуальность*

Настоящее пособие посвящено разработке новых подходов к построению многостадийных моделей транспортных потоков и эффективных численных методов поиска равновесий в таких моделях. Начиная с 50-х годов XX века вопросам поиска равновесий в транспортных сетях стали уделять большое внимание в связи с ростом городов и необходимостью соответствующего транспортного планирования. В 1955 г. появилась первая модель равновесного распределения потоков по путям: BMW-*модель*, также называемая *моделью Бэкмана* [75]. В этой модели при заданных корреспонденциях (потоках из источников в стоки) решалась задача поиска равновесного распределения этих корреспонденций по путям, исходя из принципа Вардропа, т. е. исходя из предположения о том, что каждый пользователь транспортной сети рационален и выбирает кратчайший маршрут следования. Таким образом, поиск равновесия в такой модели сводился к поиску равновесия Нэша в популяционной игре [145] (популяций столько, сколько корреспонденций). Поскольку в модели предполагалось, что время прохождения ребра есть функция от величины потока только по этому ребру, то получившаяся игра была игрой загрузки, следовательно, потенциальной. Последнее означает, что поиск равновесия сводится к решению задачи оптимизации. Получившуюся задачу выпуклой оптимизации решали с помощью метода условного градиента [99]. Описанная модель и численный метод и по настоящее время используются в подавляющем большинстве продуктов транспортного моделирования для описания блока равновесного распределения потоков по путям [61, 140, 142, 148]. Однако в работе [133] было указано на ряд существенных недостатков модели Бэкмана и предложена альтернативная модель, которую авторы назвали *моделью стабильной динамики*.

Отмеченные выше модели равновесного распределения потоков по путям могут быть использованы при решении различных задач долгосрочного планирования. Например: имея заданный бюджет, нужно решить, на каких участках графа транспортной сети стоит увеличить полосность дороги, а на каких построить новые дороги. Заданы несколько сценариев, нужно отобрать лучший. Задачу можно решить, найдя равновесные распределения потоков, отвечающие каждому из сценариев, и сравнивая найденные решения, например по критерию суммарного времени, потерянного в пути всеми пользователями сети в данном равновесии. При значительных изменениях графа транспортной сети необходимо в приведенную выше цепочку рассуждений включать дополнительный контур, связанный с тем, что изменения приведут не только к перераспределению потоков на путях, но и поменяют корреспонденции. Таким образом, корреспонденции также необходимо моделировать. В 60-е годы XX века появилось сразу несколько различных моделей для расчета матрицы корреспонденций исходя из информации о численностях районов и числе рабочих мест в них. Наибольшую популярность приобрела энтропийная модель расчета матрицы корреспонденций [18]. В этой модели поиск матрицы корреспонденций сводился к решению задачи энтропийно-линейного программирования.

К сожалению, при этом в энтропийную модель явным образом входит информация о матрице затрат на кратчайших путях по всевозможным парам районов. Возникает «порочный круг»: чтобы посчитать эту матрицу затрат, нужно сначала найти равновесное распределение потоков по путям, а чтобы найти последнее, необходимо знать матрицу корреспонденций, которая рассчитывается по матрице затрат. На практике отмеченную проблему решали методом простых итераций. Как-то «разумно» задавали начальную матрицу корреспонденций, по ней считали распределение потоков по путям, по этому распределению считали матрицу затрат, на основе которой пересчитывали матрицу корреспонденций, и процесс повторялся. Повторялся он до тех пор, пока не выполнялся критерий останова. К сожалению, до сих пор для описанной процедуры неизвестно никаких гарантий ее сходимости и тем более оценок скорости сходимости [140].

Описанный выше подход можно назвать двухстадийной моделью транспортных потоков, потому что модель состоит из прогонки двух разных блоков. В действительности, в реальных приложениях, число блоков обычно равно трем–четырем [140]. В частности, как правило, всегда включают блок расщепления потоков по типу передвижения (например, личный и общественный транспорт) – этот блок описывается моделью, аналогичной модели Бэкмана. В математическом плане это уточнение несущественно. Все основные имеющиеся тут сложности хорошо демонстрирует уже двухстадийная модель. Отметим также, что часто в приложениях вместо модели Бэкмана используется её «энтропийно-регуляризованный» вариант, который отражает ограниченную рациональность пользователей транспортной сети [70, 145, 148]. Равновесие в такой модели часто называют *стохастическим равновесием*.

Резюмируем написанное выше. До настоящего момента в учебной литературе, насколько нам известно, не существовало строгого научного обоснования используемого повсеместно на практике (и зашитого во все современные пакеты транспортного моделирования) способа формирования многостадийных моделей транспортных потоков. Не существовало также никаких гарантий сходимости численных методов, используемых для поиска равновесий в многостадийных моделях. В используемых сейчас повсеместно многостадийных моделях в качестве основных блоков фигурируют блоки с моделями типа Бэкмана, а не более современные блоки стабильной динамики. Таким образом, актуальной является задача обоснования современной многостадийной модели и разработка эффективных численных методов поиска (стохастического) равновесия в такой модели.

## *Цели и задачи*

Многие законы природы могут быть записаны в форме *вариационных принципов*, или *экстремальных принципов*. В моделировании транспортных потоков это также имеет место. Однако на текущий момент с помощью вариационных принципов описываются только отдельные блоки многостадийной транспортной модели, и эволюционный вывод вариационных принципов имеется только для блоков с моделью Бэкмана в основе. Одной из целей данного пособия является эволюционный вывод всех блоков многостадийной транспортной модели (прежде всего, речь идет о блоке расчета матрицы корреспонденций) и получение (с помощью эволюционного вывода) вариационного принципа для описания равновесия в многостадийной модели.

Целью также является разработка *алгебры* над блоками-моделями (каждый блок описывается своим вариационным принципом), которая позволит, как в конструкторе, собирать (формируя общий вариационный принцип) сколь угодно сложные модели из небольшого числа базисных элементов конструктора (блоков).

Описанный выше формализм приводит в итоге к решению задач выпуклой оптимизации в пространствах огромных размеров, которые имеют довольно специальную иерархическую (многоуровневую) структуру функционала задачи. Чтобы подчеркнуть нетривиальность таких задач, отметим, что переменные, по которым необходимо оптимизировать, – это, в частности, компоненты вектора распределения потоков по путям. Для графа в виде двумерной квадратной решетки (манхэттенская сеть) с числом вершин порядка нескольких десятков тысяч (это число соответствует транспортному графу Москвы) такой вектор с большим запасом нельзя загрузить в память любого современного суперкомпьютера, не говоря уже о том, чтобы как-то работать с такими векторами.

Важной целью пособия является разработка (с теоретическими гарантиями) эффективных численных методов, способных за несколько часов на персональном компьютере с хорошей точностью (и с высокой вероятностью) найти равновесие в многостадийной модели транспортных потоков крупного мегаполиса.

В частности, целью является разработка *алгебры* над численными методами, используемыми для расчета отдельных блоков многостадийной модели, которая позволит, как в конструкторе, собирать итоговую эффективную численную процедуру (для поиска равновесия в многостадийной модели) с помощью правильного чередования / комбинации работы алгоритмов, используемых для отдельных блоков.

## *Методы исследования*

В основе предложенного в пособии эволюционного формализма обоснования многостадийной транспортной модели лежит часто используемая в популяционной теории игр марковская logit-динамика, отражающая ограниченную рациональность агентов (водителей) [145]. Новым элементом является понимание этой динамики как модели стохастической химической кинетики с унарными реакциями и рассмотрение сразу нескольких разных типов таких унарных реакций, происходящих с разной (по порядку величины) интенсивностью и отвечающих разномасштабным процессам, протекающим в транспортной системе. Например, для двухстадийной модели динамика, отвечающая формированию корреспонденций, идет, по терминологии А. Н. Тихонова, *в медленном времени* (годы), а динамика, отвечающая распределению потоков по путям, – *в быстром времени* (дни). Тогда с некоторыми оговорками функционал в вариационном принципе с точностью до множителя и аддитивной константы можно, с одной стороны, понимать как функционал Санова (действие), отвечающий за концентрацию стационарной (инвариантной) меры введенной марковской динамики, а с другой – как функционал Ляпунова–Больцмана для кинетической динамики, полученной при (каноническом) скейлинге (по числу агентов) введённой марковской динамики.

## *Общие замечания*

Наш интерес к данной проблематике был инициирован общением с А. А. Шананиным, Ю. С. Попковым, Е. А. Нурминским и Ю. Е. Нестеровым. Данное пособие было написано в промежутке 2010–2016 гг. Однако мы не спешили с его публикацией, поскольку хотели убедиться, что полученные результаты пройдут некоторую апробацию. Результаты, собранные в разделе 2.2 гл. 2 отражают то новое, что было сделано с 2016 года. Отметим также, что в данном пособии, написанном в основном по статьям авторов, были исправлены различные опечатки, которые удавалось обнаруживать в опубликованных статьях.

Для лучшего понимания второй главы пособия рекомендуется предварительно ознакомиться с пособием [21].

Работа над пособием была поддержана грантом РФФИ 18-29-03071 мк, а также Министерством науки и высшего образования Российской Федерации (госзадание) № 075-003377-20-03, номер проекта 0714-2020-0005.

# Глава 1. Эволюционный вывод равновесных моделей распределения транспортных потоков в больших сетях

## 1.1. Трехстадийная модель равновесного распределения транспортных потоков

### 1.1.1. Введение

Одной из основных задач последнего времени, остро стоящих в Москве и ряде других крупных городов России (Санкт-Петербург, Пермь, Владивосток, Иркутск, Калининград и др.) является разработка транспортной модели города, позволяющей решать задачи долгосрочного планирования (развития) транспортной инфраструктуры города [140]. В частности, ожидается, что разработка такой модели поможет ответить на вопросы: какой из проектов дорожного строительства оптимален, где пропуская способность дороги недостаточна, как изменится транспортная ситуация, если построить в этом месте торговый центр (жилой район, стадион), как правильно определять маршруты и расписание движения общественного транспорта, какой эффект даст выделение полос для общественного транспорта и т. п.

Уже имеется программное обеспечение, позволяющее частично решать указанные выше задачи. Однако имеется много вопросов к тому, какие модели и алгоритмы используются в большинстве программных продуктах. Например, неочевидным элементом почти всех этих продуктов является использование в качестве одного из блоков модели равновесного распределения транспортных потоков Бэкмана (1955) [34, 61, 75, 141, 148]. Эта во многом хорошая модель тем не менее имеет довольно много недостатков (см. [133]). Например, калибровка такой модели требует знания функций затрат на ребрах графа транспортной сети (эти функции связывают время в пути по ребру с величиной транспортного потока по этому ребру), причем сама модель оказывается довольно чувствительной к выбору этих функций, которые в модели Бэкмана, как правило, предполагаются выпуклыми, монотонно возрастающими. Достаточно сказать, что в случае наличия платных дорог для расчета оптимальных плат за проезд требуется вычислять, например, производные этих неизвестных функций [144, 145]. Существование таких функций в модели

Бэкмана является одним из основных предположений и одновременно одним из самых слабых мест.

Реальные данные показывают (см. рис. 1.1.1, полученный В. А. Данилкиным по данным ЦОДД в 2012 г.), что предположение о классе функций затрат не выполняется. Но даже если предположить, что такая зависимость все же существует,[1] то по-прежнему остается другая проблема: как калибровать модель Бэкмана, то есть откуда брать эти зависимости. Не получится ли переобучения у создаваемой нами модели? То есть не получится ли так, что распоряжаясь большим произволом при калибровке по обучающей выборке (историческим данным) мы «переподгоним» модель: исторические данные за счет большого числа подкручиваемых параметров мы действительно можем хорошо научиться описывать, но использовать такую модель для планирования будет опасно, поскольку не будет контроля переобучения. Обычным средством борьбы с переобучением в этом месте является параметризация функций затрат, например в классе BPR-функций [34, 61, 75, 140, 141, 148]. К сожалению, каким бы то ни было научное обоснование, почему именно такая параметризация используется, нам неизвестно.

Другим, неочевидным элементом этих программных продуктов являются используемые вычислительные алгоритмы: контроль их робастности [78] к неточности (неполноте) данных, ошибкам округления (поскольку возникают задачи огромных размеров, то такие ошибки могут интенсивно накапливаться). Наконец, сама философия, использующаяся в таких продуктах при построении равновесной модели города [140], также вызывает много вопросов, о которых немного подробнее будет написано в подразделе 1.1.5. Несмотря на отмеченные выше проблемы, разработчики программного обеспечения часто находят вполне разумные инженерные компромиссы (нам известно о положительном опыте PTV VISSUM и TRANSNET), неплохо работающие на практике.

---

[1]Это допущение можно оправдать, например, тем, что, как правило, мы рассматриваем равновесные конфигурации с точки зрения пользователей сети в модели Бэкмана, которым соответствует только одна из веток – верхняя (рис. 1.1.1). В модели Бэкмана пользователи сети при принятии решений оценивают время в пути в зависимости от величины потока, которая интерпретируется как число желающих проехать по этому ребру в единицу времени. Нижняя ветка, отвечающая приблизительно линейному росту скорости с ростом потока $V \approx f/\rho_{\max}$, соответствует ситуации, когда есть узкое место, пропускная способность которого по каким-то причинам определяется не типичными характеристиками рассматриваемого ребра, а скажем, пробкой, пришедшей от впереди идущего ребра (по ходу движения). И в таких ситуациях величина потока $f$ интерпретируется не как число желающих, а как число могущих проехать. Такие ситуации просто исключаются в модели Бэкмана.

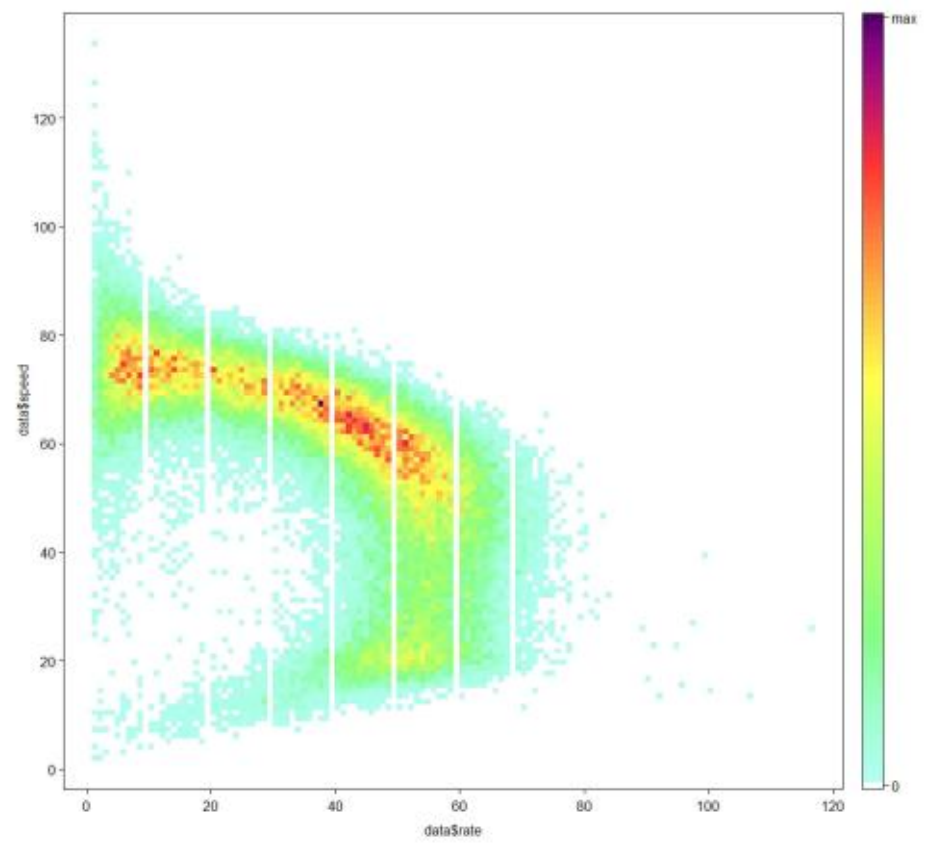

Рис. 1.1.1. По оси абсцисс – поток по двум полосам (автомобилей/мин),
а по оси ординат – скорость (км/ч)

Целью данного раздела является предложить математическую трехстадийную транспортную модель города, в которой один из блоков (модель равновесного распределения потоков) предлагается заменить с модели Бэкмана на модель стабильной динамики [130, 133]. К сожалению, целый ряд проблем, свойственных ранее известным многостадийным моделям будет присущ и модели, предлагаемой в данном разделе. Однако несколько важных недостатков, по-видимому, удалось устранить. Прежде всего, речь идет о возможности калибровки модели по реальным данным, контроле переобучения и существовании эффективного робастного вычислительного алгоритма с гарантированными (неулучшаемыми) оценками числа затраченных арифметических операций для достижения требуемой точности. Последнее обстоятельство представляется особенно важным в контексте того, как обычно используются такие модели. А именно: с помощью таких моделей просматривается множество различных сценариев. К сожалению, оптимизационные задачи вида «где и какую дорогу стоит построить при заданных бюджетных ограничениях», а также многие другие задачи, решаются непосредственным перебором различных вариантов, где для расчета каждого варианта потребуется запускать модель, меняя каждый раз что-то на входе. Кстати сказать, предложенный для расчета модели алгоритм позволяет также учитывать масштаб изменения входных данных. Если эти изменения небольшие (точнее говоря, меняется небольшое количество входных параметров), то для выполнения перерасчета по модели потребуется значительно меньше времени, чем при первом запуске.

### 1.1.2. Структура раздела и предварительные сведения

Опишем вкратце структуру раздела. В подразделе 1.1.3 описывается эволюционный способ вывода популярного на практике статического способа расчета матрицы корреспонденций (энтропийной модели). Также приводится основная идея, базирующаяся на теореме Тихонова о разделении времен [3], получения трехстадийной модели из отдельных блоков (расчет матрицы корреспонденций + равновесное расщепление потоков + равновесное распределение потоков). Точнее говоря, использование теоремы Тихонова – лишь часть идеи, которая сводит решение задачи к поиску (единственного) притягивающего положения равновесия системы в медленном времени, отвечающей за формирование корреспонденций, при подстановке в неё зависимостей времен в пути от корреспонденций (такие зависимости получаются из системы в быстром времени). Другая её часть заключается в том, что задачи поиска (единственного) притягивающего положения равновесия системы в медленном времени и поиска зависимостей времен в пути от корреспонденций с помощью некоторых вариационных принципов, о которых говорится в подразделах 1.1.3 и 1.1.4, 1.1.6–1.1.8, сводятся к задачам выпуклой оптимизации, которые можно объединить в одну общую задачу поиска седловой точки негладкой выпукло-вогнутой функции. В конце раздела указывается возможность обобщения трехстадийной модели до четырехстадийной (есть некоторые нюансы в трактовках – стоит иметь в виду, что в разных литературных источниках к таким моделям могут предъявляться немного разные требования), в которой учитываются различные типы пользователей и различные типы передвижений (моделирование идет на больших масштабах времени). В подразделе 1.1.4 описывается модель равновесного распределения потоков Бэкмана. В конце пункта приводится эволюционный способ интерпретации возникающего в этой модели равновесия Нэша–Вардропа. В подразделе 1.1.5 приводится краткий обзор многостадийных моделей, построенных на основе моделей, описанных в подразделах 1.1.3 и 1.1.4. В подразделе 1.1.6 описывается модель равновесного распределения потоков, которую мы далее будем называть *моделью стабильной динамики*. Модель стабильной динамики требует намного меньше данных для своей калибровки, наследует практически все основные «хорошие» свойства модели Бэкмана и не наследует ряд недостатков. В подразделе 1.1.7 модель стабильной динамики выводится с помощью предельного перехода из модели Бэкмана. В подразделе 1.1.8 строится обобщение модели стабильной динамики на случай, когда есть несколько способов передвижения (личный транспорт и общественный; отметим при этом, что в Москве и области более 70% пользователей сети используют общественный транспорт). Таким образом, в подразделе 1.1.8 в модель стабильной динамики органично встраивается модель равновесного расщепления потоков. В подразделе 1.1.9 энтропийная модель расчета матрицы корреспонденций из подраздела 1.1.3 объединяется с моделью из подраздела 1.1.8. В результате получается трехстадийная модель, в которой учитывается и формирование корреспонденций, и расщепление потоков, и равновесное распределение потоков по графу транспортной сети. Примечательно, что поиск равновесия в полученной трехстадийной модели (из задачи поиска седловой точки негладкой выпукло-вогнутой функции, см. подраздел 1.1.3) в итоге сводится к решению задачи негладкой выпуклой оптимизации с ограниченной константой Липшица функционала, но с неконтролируемой начальной невязкой. Важно отметить, что помимо самого решения, нужно определять и часть двойственных переменных, имеющих содержательный физический смысл. В подразделе 1.1.10 модель подраздела 1.1.9 обобщается на случай поиска стохастического равновесия, что можно проинтерпретировать как ограниченную рациональность водителей или их неполную информированность. Это допущение делает модель в подразделе 1.1.9 более приближенной к практике. С вычислительной точки зрения, сделанная модификация сводит задачу к задаче гладкой выпуклой оптимизации с немного более громоздким функционалом. Полученная задача во многом наследует все вычислительные минусы и плюсы негладкого случая. Обе задачи выпуклой оптимизации из подразделов 1.1.9, 1.1.10 требуют разработки адекватных прямодвойственных субградиентных алгоритмов решения. Заметим, что использовать методы с оракулом (сообщающим, в зависимости от своего порядка, значения функций в выбранных точках, их градиенты, и т. д.) порядка выше первого [56] не представляется возможным ввиду размеров задач. В заключительном подразделе 1.1.11 приводятся (с объяснениями) несколько практических рецептов по калибровке предложенных моделей.

На протяжении всего этого раздела (и всех последующих разделов этого пособия) мы будем активно использовать элементы выпуклого анализа и методы численного решения задач выпуклой оптимизации, ориентируясь на априорное знакомство читателя с этими дисциплинами, например в объеме книг [53, 56]. Одним из основных инструментов этого раздела будет теорема фон Неймана о минимаксе для выпукло-вогнутых функций [63]. Причем использоваться эта теорема будет не только для функций, заданных на произведении компактов, но и на неограниченных множествах. Надо лишь иметь гарантию, что максимумы и минимумы существуют (достигаются). Проблема сводится к существованию неподвижной точки у многозначного отображения. В этой области имеется большое количество результатов, с запасом покрывающих потребности данного раздела, для обоснования возможности перемены порядка взятия максимума и минимума. В частности, в этом разделе мы будем пользоваться вариантом *минимаксной теоремы*, называемой в зарубежной литературе *Sion's minimax theorem* [149], в которой предполагается компактность лишь одного из множеств, отвечающих выпуклым или вогнутым переменным, и непрерывность функции. Тем не менее далее мы будем использовать более привычное название: *минимаксная теорема фон Неймана*.

Все необходимые обозначения и ссылки будут вводиться в разделе по мере необходимости.

### 1.1.3. Энтропийная модель расчета матрицы корреспонденций

Приведем, во многом следуя книге [34] (см. также раздел 1.2 гл. 1 и приложение 1 настоящего пособия), обоснование, пожалуй, одного из самых популярных способов расчета матрицы корреспонденций, имеющего полувековую историю, – *энтропийной модели* [18, 59, 64].

Пусть в некотором городе имеется $n$ районов, $L_i > 0$ – число жителей $i$-го района, $W_j > 0$ – число работающих в $j$-м районе. При этом $N = \sum_{i=1}^{n} L_i = \sum_{j=1}^{n} W_j$ – общее число жителей города. В последующих пунктах, под $L_i \geq 0$ будет пониматься число жителей района, выезжающих в типичный день за рассматриваемый промежуток времени из $i$-го района, а под $W_j \geq 0$ – число жителей города, приезжающих на работу в $j$-й район в типичный день за рассматриваемый промежуток времени. Обычно, введенные таким образом $L_i$, $W_j$ рассчитываются через число жителей $i$-го района и число работающих в $j$-м районе с помощью более-менее универсальных (в частности, не зависящих от $i$, $j$) коэффициентов пропорциональности. Эти величины являются входными параметрами модели, т. е. они не моделируются (во всяком случае, в рамках выбранного подхода). Для долгосрочных расчетов с разрабатываемой моделью требуется иметь прогноз изменения значений этих величин.

Обозначим через $d_{ij}(t) \geq 0$ число жителей, живущих в $i$-м районе и работающих в $j$-м районе в момент времени $t$. Со временем жители могут только меняться квартирами, поэтому во все моменты времени $t \geq 0$:

$$d_{ij}(t) \geq 0,\ \sum_{j=1}^{n} d_{ij}(t) \equiv L_i,\ \sum_{i=1}^{n} d_{ij}(t) \equiv W_j,\ i, j = 1, ..., n,\ 1 \ll n^2 \ll N.^{2} \quad \textbf{(A)}$$

---

[2]Для Москвы (и других крупных мегаполисов) часто выбирают $n \sim 10^2 - 10^3$. Следовательно, корреспонденций будет $n^2 \sim 10^4 - 10^6$, и чтобы каждую из корреспонденций определить с

Опишем основной стимул к обмену: работать далеко от дома плохо из-за транспортных издержек. Будем считать, что эффективной функцией затрат [34] будет $R(T) = \beta T/2$, где $T > 0$ – время в пути от дома до работы, а $\beta > 0$ – настраиваемый параметр модели (который также можно проинтерпретировать и даже оценить).

Теперь опишем саму динамику (детали см. в [25, 34]). Пусть в момент времени $t \geq 0$ $r$-й житель живет в $k$-м районе и работает в $m$-м, а $s$-й житель живет в $p$-м районе и работает в $q$-м. Тогда $\lambda_{k,m;\,p,q}(t)\Delta t + o(\Delta t)$ есть вероятность того, что жители с номерами $r$ и $s$ ($1 \leq r < s \leq N$) «поменяются» квартирами в промежутке времени $(t, t+\Delta t)$. Вероятность обмена местами жительства зависит только от мест проживания и работы обменивающихся:

$$\lambda_{k,m;\,p,q}(t) \equiv \lambda_{k,m;\,p,q} = \lambda N^{-1} \exp\Big( \underbrace{R(T_{km}) + R(T_{pq})}_{\substack{\text{суммарные затраты} \\ \text{до обмена}}} - \underbrace{\big(R(T_{pm}) + R(T_{kq})\big)}_{\substack{\text{суммарные затраты} \\ \text{после обмена}}} \Big) > 0,$$

где коэффициент $\lambda > 0$ характеризует интенсивность обменов. Совершенно аналогичным образом можно было рассматривать случай «обмена местами работы». Здесь стоит оговориться, что «обмены» не стоит понимать буквально – это лишь одна из возможных интерпретаций. Фактически используется так называемое *приближение среднего поля* [22], т. е. некое равноправие агентов (жителей) внутри фиксированной корреспонденции и их независимость.[3]

---

точностью 10% (относительной точностью $\varepsilon \sim 10^{-1}$) потребуется (это оценка снизу) опросить не менее $n^2/\varepsilon^2 \sim 10^6 - 10^8$ жителей города, что не представляется возможным – это обстоятельство является одной из причин (неединственной), почему стараются снизить размерность пространства параметров, считая, что матрица корреспонденций задается не $n^2$ параметрами, а только, например, $2n$ ($L$ и $W$), которые и надо определять (см. также подраздел 1.1.5, в котором параметров $2n+1$). На самом деле выбирать критерием качества *относительную точность* для каждой (в том числе очень маленькой) корреспонденции не очень разумно. Более естественно восстанавливать матрицу (вектор) корреспонденций в 1-норме. Другой причиной использования описанной ниже энтропийной модели является возможность прогнозирования с её помощью того, как будет меняться матрица корреспонденций при изменении инфраструктуры города; собственно, эта одна из тех задач, которые необходимо уметь решать для получения ответов на вопросы, приведенные в подразделе 1.1.1.

Кроме того, важно заметить, что крупный мегаполис, как правило, представляется в такого рода моделях вместе со всеми своими окрестными территориями. Скажем, для Москвы – это Московская область. Отметим также, что для Москвы и области $N \sim 10^7$.

[3]Отметим также, что конечной цели (получение задачи (1.1.1)) можно добиться разными способами. Скажем, используя формализм Л. И. Розоноэра *систем обмена и распределения ресурсов* со структурной функцией энтропии (кстати, это можно не постулировать:

Согласно эргодической теореме для марковских цепей (вне зависимости от начальной конфигурации $\left\{d_{ij}(0)\right\}_{i=1,\, j=1}^{n,\, n}$) [10, 34] предельное распределение совпадает со стационарным (инвариантным), которое можно посчитать (получается проекция прямого произведение распределений Пуассона на гипергрань некоторого многогранника):

$$\exists\ c_n > 0:\ \forall\ \left\{d_{ij}\right\}_{i=1,\, j=1}^{n,\, n} \in (\mathrm{A}),\ t \ge c_n N \ln N$$

$$P\left(d_{ij}(t) = d_{ij},\, i, j = 1,\ldots,n\right) \simeq Z^{-1} \prod_{i,j=1}^{n} \exp\left(-2R\left(T_{ij}\right)d_{ij}\right)\cdot\left(d_{ij}!\right)^{-1} \overset{def}{=} p\left(\left\{d_{ij}\right\}_{i=1,\, j=1}^{n,\, n}\right),$$

где *статсумма* $Z$ находится из условия нормировки получившейся *пуассоновской* вероятностной меры. Отметим, что стационарное распределение $p\left(\left\{d_{ij}\right\}_{i=1,\, j=1}^{n,\, n}\right)$ удовлетворяет условию детального равновесия:

$$\left(d_{km}+1\right)\left(d_{pq}+1\right) p\left(\left\{d_{11},\ldots,d_{km}+1,\ldots,d_{pq}+1,\ldots,d_{pm}-1,\ldots,d_{kq}-1,\ldots,d_{nn}\right\}\right)\lambda_{k,m;\,p,q} =$$

$$= d_{pm} d_{kq}\, p\left(\left\{d_{ij}\right\}_{i=1,\, j=1}^{n,\, n}\right)\lambda_{p,m;\,k,q}.$$

При $N \gg 1$ распределение $p\left(\left\{d_{ij}\right\}_{i=1,\, j=1}^{n,\, n}\right)$ экспоненциально сконцентрировано на множестве (**A**) в $\mathrm{O}\left(\sqrt{N}\right)$-окрестности наиболее вероятного значения $\left\{d_{ij}^{*}\right\}_{i=1,\, j=1}^{n,\, n}$, которое определяется как решение задачи энтропийно-линейного программирования [34] (подробнее о таких задачах и о том, как их решать, см. [2, 18, 22, 26, 39, 59, 64, 97, 101, 102, 109, 112, 143] и подраздел 2.3.2)[4]:

---

структурная функция появляется при работе с условием интегрируемости дифференциальной формы возможных обменов ресурсами), см., например, [11] и цитированную там литературу. При таком подходе стохастика не нужна, и вариационный принцип (максимизации энтропии при аффинных ограничениях) получается мало чувствительным к особенностям возможных превращений в системе. Другие способы получения (1.1.1) связаны с информационно-статистическими соображениями, например с принципом максимума правдоподобия [18, 59, 146].

[4]При получении этой формулы использовалась асимптотическая формула Стирлинга [18, 34, 59], то есть предполагалось, что если $d_{ij} > 0$, то $d_{ij} \gg 1$ (если все $L_i > 0$, $W_j > 0$, то и все $d_{ij} > 0$ [34, 59]) и, как следствие, целочисленностью переменных $d_{ij}$ можно пренебречь, то есть решать не *NP*-полную задачу выпуклого целочисленного программирования, а обычную задачу выпуклой оптимизации (1.1.1). Сделанное предположение: «если $d_{ij} > 0$, то $d_{ij} \gg 1$» во многом будет следовать из дальнейших рассуждений (см. также [34, 59]).

$$\ln p\left(\left\{d_{ij}\right\}_{i=1,\, j=1}^{n,\, n}\right) \sim -\sum_{i,j=1}^{n} d_{ij} \ln\left(d_{ij}/e\right) \; - \beta \sum_{i,j=1}^{n} d_{ij} T_{ij} \to \max_{\left\{d_{ij}\right\}_{i=1,\, j=1}^{n,\, n} \in (\mathrm{A})} . \quad (1.1.1)$$

Естественно принимать решение этой задачи $\left\{d_{ij}^{*}\right\}_{i=1,\, j=1}^{n,\, n}$ за равновесную конфигурацию [7]. Однако имеется проблема: $T_{ij}$ – неизвестны, и зависят от $\left\{d_{ij}\right\}$. Эту проблему мы постараемся решить в дальнейшем.

Обратим внимание, что предложенный выше вывод энтропийной модели расчета матрицы корреспонденций отличается от классического [18]. В монографии А. Дж. Вильсона [18] $\beta$ интерпретируется как множитель Лагранжа к ограничению на среднее *время в пути*: $\sum_{i,j=1}^{n} d_{ij} T_{ij} = C$ .[5] Сделано это нами для того, чтобы контролировать знак параметра $\beta > 0$ и лучше понимать его физический смысл (нам это понадобится в дальнейшем).[6] Отметим также, что характерный временной масштаб формирования корреспонденций – годы (это не совсем так для корреспонденций типа дом–торговля, дом–отдых). В то время как характерное время установления равновесных значений $T_{ij}(d)$ – недели (см. ниже). Хочется сказать, что мы здесь находимся в условиях теоремы Тихонова о разделении времен [3], точнее говоря, что можно просто подставить в (1.1.1) зависимости $T_{ij}(d)$ и решать полученную задачу. Но теорема Тихонова должна при-

---

[5]При этом остальные ограничения имеют такой же вид, а функционал $F(d) = -\sum_{i,j=1}^{n} d_{ij} \ln\left(d_{ij}/e\right)$. Тогда согласно экономической интерпретации двойственных множителей Л. В. Канторовича [34, 59]:

$$\beta(C) = \partial F\left(d(C)\right)/\partial C .$$

Из такой интерпретации иногда делают вывод о том, что $\beta$ можно понимать как цену единицы времени в пути: чем больше $C$ , тем меньше $\beta$.

[6]Отметим, что, так же как и в [18], из принципа ле Шателье–Самуэльсона [11] следует, что с ростом $\beta$ среднее время в пути $\sum_{i,j=1}^{n} d_{ij}(\beta) T_{ij}\left(d(\beta)\right)$ будет убывать. В связи с этим обстоятельством, а также исходя из соображений размерности, вполне естественно понимать под $\beta$ величину, обратную к характерному (среднему) времени в пути [18] – физическая интерпретация. Собственно, такая интерпретация параметра $\beta$ , как правило, и используется в многостадийных моделях (см., например [140], а также подраздел 1.1.5).

меняться для системы ОДУ, представляющей собой в данном случае динамику квазисредних введенной выше стохастической динамики [22, 34][7]:

$$\frac{d}{dt}c_{ij} = \sum_{k,p=1}^{n} \lambda \exp\left(\frac{\beta}{2}\left(\left[T_{ip}+T_{kj}\right]-\left[T_{ij}+T_{kp}\right]\right)\right)c_{ip}c_{kj} -$$

$$-\sum_{k,p=1}^{n} \lambda \exp\left(\frac{\beta}{2}\left(\left[T_{ij}+T_{kp}\right]-\left[T_{ip}+T_{kj}\right]\right)\right)c_{ij}c_{kp},\ c_{ij}=d_{ij}/N\ ;$$

$$\varepsilon\frac{d}{dt}T_{ij} = \left[\text{сложный оператор, зависящий от } c\right],\ \varepsilon \sim 10^{-2}-10^{-3} \ll 1 .$$

Для обоснования возможности применения здесь такого рода результата, как теорема Тихонова, требуется много усилий: во-первых, изначально введенные динамики – стохастические (см. также подразделы 1.1.4, 1.1.7, 1.1.8 для $T_{ij}$), поэтому в конечном итоге нужно все обосновывать именно для них, во-вторых, нам интересна асимптотика по времени первой системы – в медленном времени (то есть нельзя ограничиться классическим случаем: ограниченного отрезка времени), в-третьих, сложный характер оператора, стоящего в правой части второй системы – в быстром времени, не позволяет явно его выписать. Тем не менее процедура построения этого оператора, которая, в свою очередь, предполагает переход к пределу (по $\mu \to 0+$ и по числу пользователей транспортной сети (аналог $N(\to\infty)$), см. подраздел 1.1.7), может быть при необходимости получена из того, что будет далее приведено в подразделах 1.1.4, 1.1.7, 1.1.8. Мы не будем здесь подробно описывать, как можно бороться с указанными проблемами, заметим лишь, что ввиду «хороших» свойств зависимостей $T_{ij}(c)$ (см. подразделы 1.1.7, 1.1.8), получающихся из приравнивания нулю левой части системы в быстром времени, функция[8]

$$H(c) = \sum_{i,j=1}^{n} c_{ij}\ln\left(c_{ij}\right) + \beta\Phi(c),\ \text{(следует сравнить с (1.1.1)),}$$

---

[7]Точнее говоря, из теоремы Т. Куртца [95] следует (с некоторыми, довольно общими, оговорками относительно зависимостей $T_{ij}(d/N)$): если предположить, что делается такой предельный переход, что (в общем случае пределы существуют неравномерно по времени, но в нашем случае равномерно):

$$\exists\ c_{ij}(0) = \lim_{N\to\infty} d_{ij}(0)/N\ , \text{то}\ \forall\ t \ge 0\ \exists\ c_{ij}(t) \overset{п.н.}{=} \lim_{N\to\infty} d_{ij}(t)/N$$

где функции $c_{ij}(t)$ (неслучайные) удовлетворяют выписанной системе ОДУ.

[8]Полученная функция будет выпуклой. Далее это будет проясняться.

где $\partial\Phi(c)/\partial c_{ij} = T_{ij}(c)$ (то, что такая функция $\Phi$ существует, показано ниже), будет функцией Ляпунова и одновременно функцией Санова (действием), то есть функцией, которая с точностью до знака и аддитивной постоянной характеризует экспоненциальную концентрацию стационарной меры введенной марковской динамики. Это довольно общий факт, справедливый для широкого класса систем [22, 34]. Поиск равновесной конфигурации $\left\{c_{ij}^*\right\}_{i=1,\, j=1}^{n,\, n}$ представляет собой решение задачи:

$$\min\left\{H(c): c \ge 0,\ \sum_{j=1}^{n} c_{ij} = l_i,\ \sum_{i=1}^{n} c_{ij} = \mathrm{w}_j\right\},$$

где $l = L/N$, $\mathrm{w} = W/N$, то есть (1). По $\left\{c_{ij}^*\right\}_{i=1,\, j=1}^{n,\, n}$ ($\left\{d_{ij}^*\right\}_{i=1,\, j=1}^{n,\, n} = N\left\{c_{ij}^*\right\}_{i=1,\, j=1}^{n,\, n}$) уже можно будет определить и равновесные $T_{ij} = T_{ij}\left(c^*\right)$.

В дальнейшем нам будет удобно привести задачу (1.1.1) к следующему виду (при помощи метода множителей Лагранжа [53], теоремы фон Неймана о минимаксе [63] и перенормировки $d := d/N$):

$$\min_{\lambda^L, \lambda^W} \max_{\sum_{i,j=1}^{n} d_{ij} = 1,\, d_{ij} \ge 0} \left[-\sum_{i,j=1}^{n} d_{ij} \ln d_{ij} - \beta \sum_{i=1,\, j=1}^{n,\, n} d_{ij} T_{ij} + \right.$$
$$\left. + \sum_{i=1}^{n} \lambda_i^L \left(l_i - \sum_{j=1}^{n} d_{ij}\right) + \sum_{i=1}^{n} \lambda_j^W \left(\mathrm{w}_j - \sum_{i=1}^{n} d_{ij}\right)\right], \tag{1.1.2}$$

где, напомним, $l = L/N$, $\mathrm{w} = W/N$.[9] Вместо того, чтобы подставлять сюда зависимость $T_{ij}(d)$, мы используем следующий трюк. В подразде-

---

[9]Обратим внимание, что мы берем максимум в (1.1.2) при дополнительном ограничении $\sum_{i,j=1}^{n} d_{ij} = 1$, которого изначально не было. Однако легко показать, что это ограничение является следствием системы ограничений (А), точнее следствием одновременно двух подсистем в (А), отвечающих $l$ и $w$. Как следствие, можно считать, что $\sum_{i=1}^{n} \lambda_i^L = 0$, $\sum_{j=1}^{n} \lambda_j^W = 0$. Отметим также, что задача (1.1.2) может быть упрощена, поскольку внутренний максимум явно находится, однако мы отложим соответствующие выкладки до подраздела 1.1.9. Здесь же мы рассмотрим случай, когда мы заносим в функционал с помощью множителей Лагранжа только часть ограничений (соответствующих $l$ и $w$). В таком случае, вычисляя внутренний максимум по $d$, получим соответственно задачи:

$$\min_{\lambda^W}\left[\sum_{i=1}^{n} l_i \ln\left[\sum_{j=1}^{n} \exp\left(-\lambda_j^W - \beta T_{ij}\right)\right] + \sum_{j=1}^{n} \lambda_j^W \mathrm{w}_j\right];\ \min_{\lambda^L}\left[\sum_{j=1}^{n} \mathrm{w}_j \ln\left[\sum_{i=1}^{n} \exp\left(-\lambda_i^L - \beta T_{ij}\right)\right] + \sum_{i=1}^{n} \lambda_i^L l_i\right],$$

лах 1.1.6–1.1.8 задача поиска зависимости $T_{ij}(d)$ будет сведена к задаче вида

$$\min_{t\geq\bar{t}}\left\{\beta\left\langle\bar{f},t-\bar{t}\right\rangle-\beta\sum_{i=1,\ j=1}^{n,\ n}d_{ij}T_{ij}(t)\right\}\overset{def}{=}-\Phi(d),$$

где $\bar{t}$ и $\bar{f}$ – известные векторы (входные параметры модели) временных затрат в пути на ребрах графа транспортной сети и максимальных пропускных способностей ребер графа транспортной сети, $T_{ij}(t)$ – длина кратчайшего пути из района $i$ в район $j$ (на графе транспортной сети имеется много вершин – намного больше числа районов, но мы считаем, что в каждом районе есть только одна вершина транспортного графа, являющаяся источником / стоком для пользователей сети; именно между такими представителями вершин районов $i$ и $j$ берется кратчайшие расстояние), если веса ребер графа транспортной сети задаются вектором $t$. Не стоит путать $T_{ij}(t)$ с искомой зависимостью $T_{ij}(d)$ – это разные функциональные зависимости. Решив указанную задачу $t^*(d)$ минимизации, считая $d$ параметрами, мы получили бы искомую зависимость $T_{ij}(d):=T_{ij}(t^*(d))$. Но явно это нельзя сделать в типичных ситуациях, да и не нужно, потому что в итоге все равно необходимо работать с потенциалом $\Phi(d)$. Поэтому предлагается (это предложение было сделано Ю. Е. Нестеровым в конце 2012 г.) ввести в задачу (1.1.1) подзадачу и добавить слагаемое $\min\limits_{\lambda^L,\lambda^W}\max\limits_{\sum_{i,j=1}^{n}d_{ij}=1,\ d_{ij}\geq 0}\min\limits_{t\geq\bar{t}}\left[\cdot+\beta\left\langle\bar{f},t-\bar{t}\right\rangle\right]$. Решение такой задачи сразу даст все, что нужно. Это следует из формулы Демьянова–Данскина [53]. Подраздел 1.1.9 посвящен упрощению только что сформулированной конструкции. Более полное обоснование и дальнейшее развитие можно отследить по работам [5, 20, 24], а также в разделах 1.3, 1.4 этой главы.

Подобно тому, как мы рассматривали в этом пункте трудовые корреспонденции (в утренние и вечерние часы более 70% корреспон-

---

где соответственно

$$d_{ij}=l_i\exp\left(-\lambda_j^W-\beta T_{ij}\right)\left[\sum_{k=1}^{n}\exp\left(-\lambda_k^W-\beta T_{ik}\right)\right]^{-1};\quad d_{ij}=w_j\exp\left(-\lambda_i^L-\beta T_{ij}\right)\left[\sum_{k=1}^{n}\exp\left(-\lambda_k^L-\beta T_{kj}\right)\right]^{-1}.$$

Отсюда можно усмотреть интерпретацию двойственных множителей как соответствующих *потенциалов притяжения / отталкивания районов* [18]. К этому мы еще вернемся в замечании 1.1.1. В принципе далее удобнее было бы работать с одной из этих задач, а не с (1.1.2), однако для сохранения симметрии, мы оставляем задачу в форме (1.1.2).

денций по Москве и области именно такие), можно рассматривать перемещения, например из дома к местам учебы, отдыха, в магазины и т. п. (по-хорошему еще надо было учитывать перемещения типа работа–магазин–детский сад–дом) – рассмотрение всего этого вкупе приведет также к задаче (1.1.1). Только будет больше типов корреспонденций $d$: помимо пары районов, еще нужно будет учитывать тип корреспонденции [18, гл. 2]. Все это следует из того, что инвариантной мерой динамики с несколькими типами корреспонденций по-прежнему будет прямое произведение пуассоновских мер. Важно отметить при этом, что в таком контексте (равновесные) $T_{ij}$ будут определяться парами районов, а не типом передвижения, то есть с точки зрения последующего изложения это означает, что ничего по сути не поменяется. Другое дело, когда мы рассматриваем разного типа пользователей транспортной сети, например: имеющих личный автомобиль и не имеющих личный автомобиль. Первые могут им воспользоваться, равно как и общественным транспортом, а вторые нет. И на рассматриваемых масштабах времени пользователи могут менять свой тип. То есть время в пути может для разных типов пользователей быть различным [18, гл. 2]. Считая, подобно тому как мы делали раньше, что желание пользователей корреспонденции[10] $(i, j)$ сети сменить свой тип (вероятность в единицу времени) есть

$$\tilde{\lambda}\exp\Big(\underbrace{\tilde{R}\left(T_{ij}^{old}\right)}_{\substack{\text{суммарные затраты}\\ \text{до смены типа}}} - \underbrace{\tilde{R}\left(T_{ij}^{new}\right)}_{\substack{\text{суммарные затраты}\\ \text{после смены типа}}}\Big), \text{ где } \tilde{R}(T) = \beta T,$$

и учитывая в «обменах» тип пользователя (будет больше типов корреспонденций $d$, но «меняются местами работы» только пользователи одного типа), можно показать, что все это вкупе приведет также к задаче типа (1.1.1), и все последующие рассуждения распространяются на этот случай. Но с оговоркой, что в подразделе 1.1.8 расщепление потоков надо делать только для тех пользователей, для которых имеется возможность использовать личный транспорт. При этом часть пользователей (в зависимости от текущего $d$) распределяется только по сети общественного транспорта. Несложно понять, что это ничего принципиально не изменит. Фактически в этом абзаце мы описали, как сделать из трехстадийной модели полноценную четырехстадийную модель, которой обычно и пользуются на практике. Чтобы не делать изложение излишне громоздким, далее мы не учитываем нюансы, описанные в этом абзаце.

---

[10]Здесь, конечно, надо учитывать не только трудовые миграции, но и все остальные, поскольку для заметной части жителей Москвы и области решение о покупки автомобиля напрямую связано с желанием ездить на нем в основном только на дачу. В этом месте возникает необходимость моделирования (учета) перемещений пользователей сети не только в рамках установленного диапазона времени в течение типичного дня, но и в целом все возможные перемещения. Это обстоятельство вынуждает использовать здесь различного рода эвристики, дабы не отказываться от важного предположения: *в рамках установленного диапазона времени в течение типичного дня*.

### 1.1.4. Модель равновесного распределения потоков Бэкмана

Следуя книге [34], опишем наиболее популярную на протяжении более чем полувека модель равновесного распределения потоков Бэкмана [34, 61, 64, 75, 141, 148].

Пусть транспортная сеть города представлена ориентированным графом $\Gamma = (V, E)$, где $V$ – узлы сети (вершины), $E \subset V \times V$ – дуги сети (рёбра графа). В современных моделях равновесного распределения потоков в крупном мегаполисе число узлов графа транспортной сети обычно выбирают порядка $|V| \sim 10^4 - 10^5$. Число ребер $|E|$ получается в три–четыре раза больше. Пусть $W \subseteq \{w = (i, j) : i, j \in V\}$ – множество корреспонденций, т. е. возможных пар *исходный пункт–цель поездки* ($|W|$ по порядку величины может быть $n^2$); $p = \{v_1, v_2, ..., v_m\}$ – путь из $v_1$ в $v_m$, если $(v_k, v_{k+1}) \in E$, $k = 1, ..., m-1$, $m > 1$; $P_w$ – множество путей, отвечающих корреспонденции $w \in W$, то есть если $w = (i, j)$, то $P_w$ – множество путей, начинающихся в вершине $i$ и заканчивающихся в $j$; $P = \bigcup_{w \in W} P_w$ – совокупность всех путей в сети $\Gamma$ (число «разумных» маршрутов $|P|$, которые потенциально могут использоваться, обычно растет с ростом числа узлов сети не быстрее, чем $O\left(|V|^3\right)$, однако теоретически может быть экспоненциально большим); $x_p$ [авт/час] – величина потока по пути $p$, $x = \{x_p : p \in P\}$; $f_e$ [авт/час] – величина потока по дуге $e$:

$$f_e(x) = \sum_{p \in P} \delta_{ep} x_p \text{, где } \delta_{ep} = \begin{cases} 1, & e \in p, \\ 0, & e \notin p; \end{cases}$$

$\tau_e(f_e)$ – удельные затраты на проезд по дуге $e$. Как правило, предполагают, что это – (строго) возрастающие, гладкие функции от $f_e$ (в конце этого раздела нам потребуется еще и выпуклость). Точнее говоря, под $\tau_e(f_e)$ правильнее понимать представление пользователей транспортной сети об оценке собственных затрат (обычно временных в случае личного транспорта и комфортности пути (с учетом времени в пути) в случае об-щественного транспорта) при прохождении дуги $e$, если поток желающих оказаться на этой дуге будет $f_e$.

Зависимость $\tau_e(f_e)$ можно попробовать и вывести, например, из следующих соображений [34] (вариация на тему модели Бобкова–Буслаева–Танака). Рассматривается одна полоса длины $L$ и транспортный поток, характеризующийся максимальной скоростью $v_{\max}$ и следующей зависимостью безопасного расстояния («комфортного» расстояния до впереди идущего транспортного средства) от скорости: $d(v) = l + \tilde{\tau}v + cv^2$ – динамический габарит, где $l$ – средняя длина автомобиля в «стоячей» пробке (эта длина немного больше средней «физической» длины автомобиля $\approx 6.5$ м), $\tilde{\tau}$ – время реакции водителей (эксперименты показывают, что для европейских водителей эта величина обычно равна одной секунде, для российских водителей она, как правило, не превышает полсекунды), $c$ – характеризует коэффициент трения шин о поверхность дороги, поскольку слагаемое $cv^2$ характеризует в $d(v)$ тормозной путь. Действительно, пока водитель среагирует на ситуацию он в среднем проедет (не изменяя своей скорости) путь $\tilde{\tau}v$. Когда уже реакция произошла, водитель начинает тормозить, и кинетическая энергия автомобиля $mv^2/2$ должна быть «погашена» работой силы трения на участке $h$ тормозного пути $\mu mgh$ ($\mu$ – коэффициент трения, $g$ – ускорение свободного падения). Отсюда можно найти $c = 1/(2\mu g)$. Имея зависимость $d(v)$, можно ввести зависимость $\rho(v) = 1/d(v)$, которая порождает зависимость потока от скорости $f(v) = v\rho(v)$. Используя то, что $\tau(f(v)) = L/v$ и $0 \le v \le v_{\max}$, можно явно выписать искомую зависимость $\tau(f)$. Несмотря на предложенный вывод, еще раз подчеркнем, что под $\tau_e(f_e)$ правильнее понимать представление пользователей транспортной сети об оценке собственных затрат при прохождении дуги $e$, если поток желающих оказаться на этой дуге будет $f_e$; поэтому функции $\tau_e(f_e)$, выбираемые на практике (см., например, подраздел 1.1.7), сильно отличаются от описанной в этом абзаце зависимости, приводящей на самом деле к наличию двух веток (см. в этой связи рис. 1.1.1 из подраздела 1.1.1).

Рассмотрим теперь $G_p(x)$ – затраты временные или финансовые на проезд по пути $p$. Естественно считать, что $G_p(x) = \sum_{e \in E} \tau_e(f_e(x))\delta_{ep}$. В приложениях часто требуется учитывать также затраты на прохождения вершин графа, которые могут зависеть от величин всех потоков через рассматриваемую вершину.

Пусть также известно, сколько перемещений в единицу времени $d_w$ осуществляется согласно корреспонденции $w \in W$. Тогда вектор $x$, характеризующий распределение потоков, должен лежать в допустимом множестве

$$X = \left\{ x \geq 0 : \sum_{p \in P_w} x_p = d_w,\ w \in W \right\}.$$

Это множество может иметь и более сложный вид, если дополнительно учитывать, например, конечность пропускных способностей рёбер (ограничения сверху на $f_e$).

Рассмотрим игру, в которой каждому элементу $w \in W$ соответствует свой, достаточно большой ($d_w \gg 1$), набор однотипных «игроков», осуществляющих передвижение согласно корреспонденции $w$. Чистыми стратегиями игрока служат пути, а выигрышем – величина $-G_p(x)$. Игрок «выбирает» путь следования $p \in P_w$, при этом, делая выбор, он пренебрегает тем, что от его выбора также «немного» зависят $|P_w|$ компонент вектора $x$ и, следовательно, сам выигрыш $-G_p(x)$. Можно показать, что отыскание равновесия Нэша–Вардропа $x^* \in X$ (макроописание равновесия) равносильно решению задачи нелинейной комплементарности (принцип Вардропа):

$$\text{для любых } w \in W,\ p \in P_w \text{ выполняется } x_p^* \cdot \left( G_p\left(x^*\right) - \min_{q \in P_w} G_q\left(x^*\right) \right) = 0.$$

Действительно, допустим, что реализовалось какое-то другое равновесие $\tilde{x}^* \in X$, которое не удовлетворяет этому условию. Покажем, что тогда найдется водитель, которому выгодно поменять свой маршрут следования. Действительно, тогда

$$\text{существуют такие } \tilde{w} \in W,\ \tilde{p} \in P_{\tilde{w}}, \text{ что } \tilde{x}_{\tilde{p}}^* \cdot \left( G_{\tilde{p}}\left(\tilde{x}^*\right) - \min_{q \in P_{\tilde{w}}} G_q\left(\tilde{x}^*\right) \right) > 0.$$

Каждый водитель (множество таких водителей непусто $\tilde{x}_{\tilde{p}}^* > 0$), принадлежащий корреспонденции $\tilde{w} \in W$ и использующий путь $\tilde{p} \in P_{\tilde{w}}$, действует неразумно, поскольку существует такой путь $\tilde{q} \in P_{\tilde{w}}$, $\tilde{q} \neq \tilde{p}$, что $G_{\tilde{q}}\left(\tilde{x}^*\right) = \min_{q \in P_{\tilde{w}}} G_q\left(\tilde{x}^*\right)$. Этот путь $\tilde{q}$ более выгоден, чем $\tilde{p}$. Аналогично показывается, что при $x^* \in X$ никому из водителей уже не выгодно отклоняться от своих стратегий. Но это по определению и называется *рав-* *новесием Нэша*, которое ввел в своей диссертации в конце 40-х годов XX века Джон Нэш, получивший именно за эту концепцию в 1994 г. нобелевскую премию по экономике [40]. Мы также добавляем фамилию Дж. Г. Вардропа, которой чуть позже Дж. Нэша привнес к этой концепции условие *конкурентного рынка*: игрок, принимающий решение, пренебрегает тем, что его решение сколько-нибудь значительно поменяет ситуацию на *рынке*. Когда игроков двое, трое (ситуации, рассматриваемые Дж. Нэшем), то очевидно, что так делать нельзя. Но когда игроков (водителей) десятки и сотни тысяч, вся эта конструкция неявно предполагает, что $x_p^* > 0 \Rightarrow x_p^* \gg 1$. Поэтому, не боясь сильно ошибиться, можно искать решение задачи нелинейной комплементарности, не предполагая целочисленности компонент вектора $x^* \in X$. Такая релаксация изначально целочисленной задачи заметно упрощает её с вычислительной точки зрения!

Хотя мы и смогли выписать условие равновесия в виде задачи нелинейной комплементарности, это несильно продвинуло нас в понимании того, как его находить. Пытаться честно решить задачу в таком виде – вычислительно бесперспективная задача. С другой стороны, современные вычислительные методы позволяют эффективно решать задачи выпуклой оптимизации. Постараемся свести нашу задачу к таковой.

Для этого, прежде всего, заметим, что рассматриваемая нами игра принадлежит к классу так называемых *потенциальных игр*. В нашем случае это означает, что существует такая функция

$$\Psi(x) = \sum_{e\in E} \int_0^{\sum_{p\in P}\delta_{ep}x_p} \tau_e(z)\,dz = \sum_{e\in E} \sigma_e\left(f_e(x)\right),$$

где $\sigma_e\left(f_e\right) = \int_0^{f_e(x)} \tau_e(z)\,dz$, что $\partial\Psi(x)/\partial x_p = G_p(x)$ для любого $p \in P$. Таким образом, мы имеем дело с потенциальной игрой. Оказывается, что $x^* \in X$ – это равновесие Нэша–Вардропа тогда и только тогда, когда оно доставляет минимум $\Psi(x)$ на множестве $X$. Действительно, предположим, что $x^* \in X$ – точка минимума. Тогда, в частности, для любых $w \in W$, $p,q \in P_w$ ($x_p^* > 0$) и достаточно маленького $\delta x_p > 0$ выполняется

$$-\frac{\partial\Psi(x^*)}{\partial x_p}\delta x_p + \frac{\partial\Psi(x^*)}{\partial x_q}\delta x_p \geq 0 .$$

Иначе, заменив $x^*$ на

$$\breve{x}^* = x^* + \left( \underbrace{\overbrace{0,\dots,0,-\delta x_p}^{p},0,\dots,0,\delta x_p}_{q},0,\dots,0 \right) \in X ,$$

мы пришли бы к вектору $\breve{x}^*$, доставляющему меньшее значение $\Psi(x)$ на множестве $X$ :

$$\Psi(\breve{x}^*) \approx \Psi(x^*) - \frac{\partial \Psi(x^*)}{\partial x_p} \delta x_p + \frac{\partial \Psi(x^*)}{\partial x_q} \delta x_p < \Psi(x^*).$$

Вспоминая, что $\partial\Psi(x)/\partial x_p = G_p(x)$, и учитывая, что $q$ можно выбирать произвольно из множества $P_w$, получаем:

для любых $w \in W,\ p \in P_w$, если $x_p^* > 0$, то выполняется

$$\min_{q \in P_w} G_q(x^*) \ge G_p(x^*).$$

Но это и есть по-другому записанное условие нелинейной комплементарности. Строго говоря, мы показали сейчас только то, что точка минимума $\Psi(x)$ на множестве $X$ будет равновесием Нэша–Вардропа. Аналогично рассуждая, можно показать и обратное: равновесие Нэша–Вардропа доставляет минимум $\Psi(x)$ на множестве $X$ . Этот минимум можно искать, например, с помощью метода условного градиента [27, 64] (Франка–Вульфа). Подробнее об этом методе также написано в разделе 1.5 этой главы.

**Теорема 1.1.1 [34, 61, 141, 148].** *Вектор $x^*$ будет равновесием Нэша–Вардропа тогда и только тогда, когда*

$$x \in \operatorname{Arg}\min_x \left[ \Psi(f(x)) = \sum_{e \in E} \sigma_e(f_e(x)) :\ f = \Theta x,\ x \in X \right]. \qquad (1.1.3)$$

*Если преобразование $G(\cdot)$ строго монотонное, то равновесие $x$ единственно. Если $\tau_e'(\cdot) > 0$, то равновесный вектор распределения потоков по ребрам $f$ – единственный* (*это еще не гарантирует единственность вектора распределения потоков по путям $x$* [34]).

Проинтерпретируем, следуя [25, 34, 70, 145] (см. также разделы 1.2, 1.4 этой главы), эволюционным образом полученное равновесие, попутно отвечая на вопрос (поскольку задача (1.1.4) ниже имеет единственное решение, детали см. в [25, 73]): какому из равновесий стоит отдать предпочтение в случае неединственности? Опишем *марковскую логит-динамику* (также говорят *гиббсовскую динамику*) в повторяющейся игре загрузки графа транспортной сети [145]. Пусть каждой корреспонденции отвечает $d_w M$ агентов ($M \gg 1$), $x := x/M$, $f := f/M$, $\tau_e(f_e) := \tau_e(f_e/M)$. Пусть имеется $tN$ шагов ($N \gg 1$). Пусть $k$-й агент, принадлежащий корреспонденции $w \in W$, независимо от остальных на шаге $m+1$ с вероятностью $1-\lambda/N$ выбирает путь $p^{k,m}$, который использовал на шаге $m = 0,...,tN$, а с вероятностью $\lambda/N$ ($\lambda > 0$) решает «поменять» путь и выбирает (возможно тот же самый) зашумленный кратчайший путь:

$$p^{k,m+1} = \arg\max_{q \in P_w} \left\{ -G_q\left(x^m\right) + \xi_q^{k,m+1} \right\},$$

где независимые случайные величины $\xi_q^{k,m+1}$, имеют одинаковое двойное экспоненциальное распределение, также называемое *распределением Гумбеля*[11] [70, 145]:

$$P\left(\xi_q^{k,m+1} < \zeta\right) = \exp\left\{-e^{-\zeta/T-E}\right\},\ T > 0.$$

Отметим, что

$$P\left(p^{k,m+1} = p \middle| \text{агент решил «поменять» путь}\right) = \frac{\exp\left(-G_p\left(x^m\right)/T\right)}{\sum_{q \in P_w} \exp\left(-G_q\left(x^m\right)/T\right)}.$$

Кроме того, если взять $E \approx 0.5772$ – константу Эйлера, то

$$M\left[\xi_q^{k,m+1}\right] = 0,\ D\left[\xi_q^{k,m+1}\right] = T^2\pi^2/6.$$

Такая динамика отражает ограниченную рациональность агентов (см. подраздел 1.1.5) и часто используется в популяционной теории игр [145] и теории дискретного выбора [70]. Оказывается, эта марковская динамика в пределе $N \to \infty$ превращается в марковскую логит-динамику в непрерывном времени (вырождающуюся при $T \to 0+$ в динамику наилучших ответов [145] – последующие рассуждения, в частности, формулы (1.1.4), (1.1.5), допускают переход к пределу $T \to 0+$). Марковская логит-динамика в непрерывном времени допускает два предельных перехода (обоснование перестановочности этих пределов см. в [7, 95]): $t \to \infty$, $M \to \infty$ или $M \to \infty$, $t \to \infty$. При первом порядке переходов мы сначала ($t \to \infty$) согласно эргодической теореме для марковских процессов (в нашем случае марковский процесс – модель стохастической химиче-ской кинетики с унарными реакциями в условиях детального баланса [7, 25]) приходим к финальной (= стационарной) вероятностной мере, имеющей в основе мультиномиальное распределение. С ростом числа агентов ( $M \to \infty$ ) эта мера

$$\sim \exp\left(-\frac{M}{T}\cdot\left(\Psi_T\left(x\right)+o\left(1\right)\right)\right)$$

экспоненциально концентрируется около наиболее вероятного состояния, поиск которого сводится к решению энтропийно-регуляризованной задачи (1.1.3) (как численно решать эту задачу, описано в работе [23], см. также [2]):

$$\Psi_T\left(x\right)=\Psi\left(f\left(x\right)\right)+T\sum_{w\in W}\sum_{p\in P_w}x_p\ln x_p \to \min_{\substack{f(x)=\Theta x\\ x\in X}}. \qquad (1.1.4)$$

---

[11]Распределение Гумбеля можно объяснить исходя из идемпотентного аналога центральной предельной теоремы (вместо суммы случайных величин – максимум) для независимых случайных величин с экспоненциальным и более быстро убывающим правым хвостом [52]. Распределение Гумбеля возникает в данном контексте, например, если при принятии решения водитель собирает информацию с большого числа разных (независимых) зашумленных источников, ориентируясь на худшие прогнозы по каждому из путей.

Функционал в этой задаче оптимизации с точностью до потенцирования и мультипликативных и аддитивных констант соответствует исследуемой стационарной мере, то есть это функционал Санова [7, 25]. При обратном порядке предельных переходов, мы сначала ( $M \to \infty$ ) осуществляем так называемый *канонический скейлинг* [7, 95], приводящий к детерминированной кинетической динамике, описываемой СОДУ на $x$ :

$$\frac{dx_p}{dt}=d_w\frac{\exp\left(-G_p\left(x\right)/T\right)}{\sum_{l\in P_w}\exp\left(-G_l\left(x\right)/T\right)}-x_p,\ p\in P_w,\ w\in W, \qquad (1.1.5)$$

а затем ( $t \to \infty$ ) ищем аттрактор получившейся СОДУ. Глобальным аттрактором оказывается неподвижная точка, которая определяется решением задачи (1.1.4). Более того, функционал $\Psi_T\left(x\right)$ является функцией Ляпунова полученной кинетической динамики (1.1.5) (то есть является функционалом Больцмана). Последнее утверждение – достаточно общий факт (функционал Санова, является функционалом Больцмана), верный при намного более общих условиях (см. [7] и цитированную там литературу).

Хотелось бы подчеркнуть, что рассматриваемая выше «игра» – потенциальная (это общий факт для игр загрузок [145]; Розенталь 1973, Мондерер–Шэпли 1996), поэтому из общих результатов эволюционной теории игр [85, 145] следует, что любые разумные содержательно интерпретируемые (суб-)градиентные спуски приводят к равновесию Нэша–Вардропа (или его стохастическому варианту, даваемому решением задачи (1.1.4)). В частности, в [34, 36, 38] содержательно интерпретируется быстро сходящаяся динамика, связанная с методом зеркального спуска, которая порождается имитационной логит-динамикой [145]. Хотелось бы также обратить внимание на эволюционную интерпретацию парадокса

Браесса: когда неэффективное по Парето, единственное равновесие Нэша–Вардропа в специально сконструированной транспортной сети является тем не менее эволюционно устойчивым [34, 145].

Нетривиальным является следующее наблюдение. Если рассмотреть энтропийно-регуляризованный функционал $\Psi_T(x)$ и взять предел (см. подразделы 1.1.6, 1.1.7):

$$\tau_e^{\mu}(f_e) \xrightarrow[\mu\to 0+]{} \begin{cases} \bar{t}_e, & 0 \le f_e < \bar{f}_e \\ [\bar{t}_e, \infty), & f_e = \bar{f}_e \end{cases},$$

$$d\tau_e^{\mu}(f_e)/df_e \xrightarrow[\mu\to 0+]{} 0,\ 0 \le f_e < \bar{f}_e,$$

то переход к двойственной задаче (для задачи минимизации этого функционала на множестве $X$ ) дает стохастический вариант модели стабильной динамики (см. подраздел 1.1.6), который используется в стохастическом варианте трехстадийной модели стабильной динамики (см. подраздел 1.1.10).

**Замечание 1.1.1 (*облачная модель* расчета матрицы корреспонденций [25], см. также раздел 1.2 этой главы).** В контексте написанного выше полезно отметить другой способ обоснования энтропийной модели расчета матрицы корреспонденций из подраздела 1.1.3. Предположим, что все вершины, отвечающие источникам корреспонденций, соединены ребрами с одной вспомогательной вершиной (облако № 1). Аналогично все вершины, отвечающие стокам корреспонденций, соединены ребрами с другой вспомогательной вершиной (облако № 2). Припишем всем новым ребрам постоянные веса. И проинтерпретируем веса ребер, отвечающих источникам $\lambda_i^L$ , например, как средние затраты на проживание (в единицу времени, скажем, в день) в этом источнике (районе), а веса ребер, отвечающих стокам $\lambda_j^W$ , как уровень средней заработной платы со знаком минус (в единицу времени) в этом стоке (районе), если изучаем трудовые корреспонденции. Будем следить за системой в медленном времени, то есть будем считать, что равновесное распределение потоков по путям стационарно. Поскольку речь идет о равновесном распределении потоков, то нет необходимости говорить о затратах на путях или ребрах детализированного транспортного графа, достаточно говорить только о затратах (в единицу времени), отвечающих той или иной корреспонденции. Таким образом, у нас есть взвешенный транспортный граф с одним источником (облако 1) и одним стоком (облако 2). Все вершины этого графа, кроме двух вспомогательных (облаков), соответствуют районам в модели расчета матрицы корреспонденций из подраздела 1.1.3. Все ребра этого графа имеют стационарные (не меняющиеся и не зависящие от текущих корре-спонденций) веса $\{T_{ij};\lambda_i^L;\lambda_j^W\}$. Если рассмотреть естественную в данном контексте логит-динамику с $T=1/\beta$ (здесь полезно напомнить, что согласно подразделу 1.1.3 $\beta$ обратно пропорционально средним затратам, а $T$ имеет как раз физическую размерность затрат), описанную выше, то поиск равновесия рассматриваемой макросистемы приводит (в прошкалированных переменных) к задаче, сильно похожей на задачу (1.1.2) из подраздела 1.1.3:

$$\max_{\sum_{i,j=1}^{n} d_{ij}=1,\ d_{ij}\geq 0}\left[-\sum_{i,j=1}^{n} d_{ij}\ln d_{ij} - \beta\sum_{i=1,\ j=1}^{n,\ n} d_{ij}T_{ij} - \beta\sum_{i=1}^{n}\left(\lambda_i^L\sum_{j=1}^{n} d_{ij}\right) - \beta\sum_{i=1}^{n}\left(\lambda_j^W\sum_{i=1}^{n} d_{ij}\right)\right].$$

Разница состоит в том, что здесь мы не оптимизируем по $2n$ двойственным множителям $\lambda^L$, $\lambda^W$. Более того, мы их и не интерпретируем здесь как двойственные множители, поскольку мы их ввели на этапе взвешивания ребер графа. Тем не менее значения этих переменных, как правило, неоткуда брать. Тем более, что приведенная выше (наивная) интерпретация вряд ли может всерьез рассматриваться, как способ определения этих параметров исходя из данных статистики. Более правильно понимать $\lambda_i^L$, $\lambda_j^W$ как потенциалы притяжения / отталкивания районов, включающие в себя плату за жилье и зарплату, но включающие также и многое другое, что сложно описать количественно. И здесь как раз помогает информация об источниках и стоках, содержащаяся в $2n$ уравнениях из формулы (**А**) подраздела 1.1.3. Таким образом, мы приходим ровно к той же самой задаче (1.1.2) с той лишь разницей, что мы получили дополнительную интерпретацию двойственных множителей в задаче (1.1.2). При этом двойственные множители в задаче (1.1.2) равны (с точностью до мультипликативного фактора $\beta$) введенным здесь потенциалам притяжения районов.

Несложно распространить на изложенную здесь модель написанное в подразделе 1.1.3 по поводу того, как с помощью разделения времен можно учитывать обратную связь: *перераспределение потоков по путям изменяется* (*в быстром времени*) *при изменении корреспонденций*, а также распространить на то, что написано в самом конце подраздела 1.1.3. Нам представляется такой способ рассуждения даже более привлекательным, чем способ, описанный в подразделе 1.1.3 и основанный на «обменах». И связано это с тем, что для получения равновесия в многостадийной модели, мы можем рассмотреть всего одну (общую) логит-динамику, в которой с малой интенсивностью (в медленном времени) происходят переходы, описанные в этом замечании (жители города меняют места жительства, работы), а с высокой интенсивностью (в быстром времени, изо дня в день) жители города перераспределяются по путям (в зависимости от текущих корреспонденций, подстраиваясь под корреспонденции) – это как раз и было описано непосредственно перед замечанием. Другая причина – бóльшая вариативность модели, построенной в этом замечании. Нам представляется очень плодотворной и перспективной идея перенесения имеющейся информации об исследуемой системе из обременительных законов сохранения динамики, описывающей эволюцию этой системы, в саму динамику путем введения дополнительных естественно интерпретируемых параметров. При таком подходе становится возможным, например, учитывать в моделях и рост транспортной сети. Другими словами, при таком подходе, например, можно естественным образом рассматривать также и ситуацию, когда число пользователей транспортной сетью меняется со временем (медленно).

В заключение отметим, что если штрафовать (назначать платы) за использование различных стратегий (маршрутов) по правилу

$$\bar{G}_p(x) = G_p(x) + \underbrace{\sum_{w \in W} \sum_{q \in P_w} x_q \frac{\partial G_q(x)}{\partial x_p}}_{\text{штраф}}, \quad p \in P,$$

то найдется такая функция $\bar{\Psi}(x) = \sum_{w \in W} \sum_{p \in P_w} x_p G_p(x)$, что $\partial \bar{\Psi}(x) / \partial x_p = \bar{G}_p(x)$. Поэтому из сказанного выше в этом пункте будет следовать, что возникающее в такой управляемой транспортной сети равновесие Нэша–Вардропа будет (единственным) глобально устойчивым и соответствовать социальному оптимуму в изначальной транспортной сети [144]. Чтобы в этом убедиться, достаточно (ввиду линейной связи $f = \Theta x$) проверить строгую выпуклость функции $\bar{\Psi}(x) = \sum_{w \in W} \sum_{p \in P_w} x_p G_p(x) = \sum_{e \in E} f_e(x) \tau_e(f_e(x))$, для чего достаточно выпуклости функций $\tau_e(f_e)$, $e \in E$. Все это хорошо соответствует механизму Викри–Кларка–Гроуса [69] (VCG mechanism) – штраф (плата) за использование маршрута новым пользователем равен дополнительным потерям, которые понесут из-за этого все остальные пользователи. Однако важно сделать две оговорки. Во-первых, все это хотя и можно попытаться практически осуществить (например, собирая транспортные налоги исходя из трековой информации, которую в перспективе можно будет иметь о каждом автомобиле), но механизм оказывается довольно сложным. Платы взимаются не за проезд по ребру, как хотелось бы, а именно за выбор (проезд) маршрута. Кроме того, плата является функцией состояния транспортной системы $x$, кото-рое, в отличие от $f$, не наблюдаемо. Во-вторых, взимая платы за проезд, мы, с одной стороны, приводим систему в социальный оптимум, а с другой стороны, для достижения этой цели вынуждены собирать с участников движения налоги. К сожалению, их размер может оказаться внушительным, и это уже необходимо учитывать с точки зрения расщепления участников движения по выбору типа передвижения. Относительно второй проблемы – готовых решений нам неизвестно. Это известная проблема в современном разделе теории игр: *mechanism design* [69]. А вот по первой проблеме (адаптивное) решение есть [144]:

$$\bar{\tau}_e(f_e) = \tau_e(f_e) + \underbrace{f_e \tau_e'(f_e)}_{\text{штраф}},\ e \in E.$$

Легко проверить, что это приведет к указанному выше пересчету $G_p(x) \to \bar{G}_p(x)$. Далее, если мы зафиксируем (знаем) социальный оптимум $f^{opt}$, то платы за проезд можно выбирать постоянными $f_e^{opt}\tau_e'(f_e^{opt})$, $e \in E$. Все сказанное выше об устойчивости останется в силе (без всяких дополнительных предположений о выпуклости функций $\tau_e(f_e)$, $e \in E$) с одной лишь оговоркой, что транспортная система должна поддерживаться при заданных корреспонденциях $d_w$ и функциях затрат $\tau_e(f_e)$. В противном случае, возникающее равновесие уже может не соответствовать социальному оптимуму.

### 1.1.5. Краткий обзор подходов к построению многостадийных моделей транспортных потоков, с моделью типа Бэкмана в качестве модели равновесного распределения потоков

Так называемая 4-*стадийная* модель является наиболее употребительной методологией моделирования транспортных систем городов и агломераций (см., например, [140]). В рамках данной модели производится поиск равновесия спроса и предложения на поездки. При этом рассматриваются в единой совокупности модели генерации трафика, его распределения по типам передвижения и дальнейшее распределение по маршрутам.

Методология, как уже описывалось ранее, включает последовательное выполнение четырех этапов, последние три из которых закольцовываются для получения самосогласованных результатов. Исходными данными для модели являются граф дорожной сети и сети общественного транспорта с заданными функциями издержек и других определяющих параметров, разделение города на транспортные зоны и параметры этих зон (например, количество рабочих мест или мест жительства). На выходе модели выдаются: оценка матрицы корреспонденций для каждого типа передвижения, загрузка элементов сети (например, конкретной дороги) и издержки, соответствующие данному уровню загрузки. Понятие *тип передвижения*, используемое выше, является обобщением понятия *типа транспорта* и обозначает последовательность (или просто множество) используемых типов транспорта. Например, типом передвижения может считаться поездка на общественном транспорте (неважно, автобусе или троллейбусе или на обоих поочередно с пересадками) или же использование схемы park-and-ride. Уровень детализации при этом определяется самим модельером. Обычно выбирается или модель с тремя типами передвижений (общественный и личный транспорт, пешие прогулки) или с детализацией до типа транспорта (автомобиль, наземный общественный транспорт, метро, пешие прогулки и различные комбинации, перечисленные ранее). При этом стараются не учитывать типы передвижений, которые используются очень редко.

Структурно расчет модели можно описать следующим образом:

1. На первом шаге из исходных данных о транспортных зонах получают векторы отправления и прибытия для транспортных зон, т. е. в наших обозначениях $\{L_i\}_{i=1}^n$ и $\{W_j\}_{j=1}^n$.
2. Рассчитывается первая оценка матрицы корреспонденций (обычно с помощью гравитационной модели [34, 64]).
3. Рассчитывается расщепление корреспонденций по типам передвижений.
4. Для каждого типа передвижений рассчитывается распределение потоков по маршрутам.
5. Получается оценка матрицы издержек корреспонденций и вектор загрузки сети.
6. Проверяется критерий остановки.
7. Если критерий выполнен, решение получено.
8. Если критерий не выполнен, вернуться на шаг 2 с переоцененной матрицей издержек.

Первый этап модели мы разбирать не будем, так как он является достаточно обособленным от других в том смысле, что он не входит в итерационную часть метода и практически «бесплатен» с точки зрения сложности операций. Отметим лишь, что векторы отправления и прибытия рассчитываются из параметров зон с помощью простейших регрессионных моделей.

Рассмотрим модели, лежащие в основе этого алгоритма более подробно. Данный алгоритм – итерационный, после «прогонки» очередной итерации мы возвращаемся на второй шаг схемы. На $m$-й итерации на 2-м шаге формируется $m$-я оценка матрицы корреспонденций. Для её построения необходимо знание о матрице издержек в сети (т. е. знание о стоимости проезда из каждой зоны в каждую). На всех шагах, кроме первого, данные значения получаются из предыдущей ($m$–1)-й итерации. На первой итерации используются издержки, соответствующие незагруженной сети или почти любая другая, разумная, оценка матрицы издержек. Если в сети для пользователей доступно несколько типов передвижений, то в качестве оценки издержек берутся или средние издержки (время в пути) для конкретной корреспонденции, или же (если пользователи распределены «равновесно») значения издержек, например, для личного транспорта, так как в «равновесии» издержки зачастую (исключая некоторые особые случаи) должны быть равны для различных, используемых, типов передвижений.

Как правило, оценки матрицы корреспонденций строятся согласно гравитационной или энтропийной моделям [34, 64].

Посмотрим на то, как замыкается модель при наличии различных слоев спроса и типов передвижений. После того, как был проведен первый этап, модельеру уже известно, сколько людей выезжает и въезжает в каждую транспортную зону, а также известно, какая доля этих людей совершает поездку того или иного типа, т. е. распределение поездок по слоям спроса. Для каждого слоя спроса формируется своя матрица корреспонденций, при этом используется одна и та же матрица издержек (так как используемая транспортная сеть для всех жителей одна и та же), однако параметры $\{\beta\}$ для каждого слоя спроса свои. После этого, данные матрицы корреспонденций суммируются. Полученная агрегированная матрица корреспонденций используется для проведения этапа расщепления корреспонденций (суммарных) по типам передвижений. Далее алгоритм работает только с ней (агрегированной матрицей) до начала следующей итерации, когда вся процедура повторяется. Более подробно вся процедура рассматривается ниже.

После построения очередной оценки матрицы корреспонденций, используя опять же матрицу издержек, происходит расщепление корреспонденций по типам передвижений. Для этого используются модели дискретного выбора, родственные уже описанной выше logit-choice модели. На выходе данного этапа (3-го шага в описанной схеме) получаются матрицы корреспонденций для каждого типа передвижений.

Отметим, что нет единой методологии и стандарта, какую из моделей дискретного выбора стоит использовать. Каждая из моделей имеет свои недостатки. Так, например, для модели logit-choice получаемое стохастическое равновесие даже в простейших постановках может сильно отличаться от равновесия Нэша. Также результаты данной модели очень чувствительны к тому, как задается набор альтернатив для выбора. Например, расщепление для дерева выбора (личный транспорт, автобус) и для дерева выбора (личный транспорт, зеленый автобус, красный автобус) будут разными (даже если выбор описывается для одной и той же транспортной системы и вся разница только в том, что во втором случае указан цвет автобуса). Наиболее используемыми на практике являются родственные логит-модели Nested Logit Model и Multinomial Logit Model, а также композитная модель Mixed Logit Model (см. [140]).

Важно сказать, что порядок шагов 2 и 3 зависит от того, какой слой спроса описывается, т. е., грубо говоря, от цели поездки и некоторых параметров жителей. Например, слой спроса может определяться как «люди предпенсионного возраста, совершающие поездки из дома на дачу». Иногда слои спроса определяются как цель поездки с учетом пункта отправления, например «поездка из дома на работу». Выше мы отмечали, что характерное время формирования корреспонденций – годы. Это справедливо для поездок дом–работа и обратно. Однако данный тезис кажется не совсем точным, например, для поездок за покупками. В этом случае можно предположить, что люди определяют место покупок уже после того, как будет определен тип передвижения, т. е. для данного слоя спроса шаги 2 и 3 оказывается возможным поменять местами.

Наконец, последним этапом (шаг 4) является расчет равновесного распределения потоков для каждого типа передвижения. Для личного транспорта при этом используется модель Бэкмана. Для общественного транспорта не существует единого подхода к моделированию и, как правило, используются композитные модели, включающие элементы моделей дискретного выбора. Стоит лишь отметить, что при численном расчете моделей реального города именно этап расчета равновесного распределения потоков является самым затратным по времени. На выходе этого этапа получаются векторы загрузки сети (т. е. потоки по ребрам транспортного графа) и матрица издержек (шаг 5).

Шаги 6–8 в приведенной выше схеме отвечают за создание обратной связи в модели, которая позволит учесть взаимное влияние матрицы корреспонденций и матрицы издержек. Напомним, что мы использовали матрицу издержек для расчета матрицы корреспонденций на шаге 2 и матрицу корреспонденций для получения матрицы издержек на шаге 4. Логично требовать, чтобы матрица издержек, подаваемая на вход на шаге 2 и матрица издержек, получаемая на выходе на шаге 5, если не совпадали, то были «близки» по какому-либо разумному критерию.

Критерием остановки служит равенство (с требуемой точностью) средних издержек для матрицы издержек, полученной по модели расчета матрицы корреспонденций, и средних издержек, полученных эмпирическим путем. Мы к этому еще вернемся чуть позже. В статьях [146, 155] было показано, что данный критерий остановки в нашей постановке соответствует оценке матрицы корреспонденций методом максимального правдоподобия [150] для описанной далее модели.

Пусть имеются данные о реализации некоторой случайной матрицы корреспонденций $\{R_{ij}\}$, каждый элемент которой (независимо от всех остальных) распределен по закону Пуассона с математическим ожиданием $M(R_{ij}) = d_{ij}$, причем верно параметрическое предположение о том, что матрица корреспонденций представляется следующим образом:

$$d_{ij} = \tilde{A}_i \tilde{B}_j f(C_{ij}),$$

где $f(C_{ij})$ – функция притяжения (часто выбирают $f(C_{ij}) = \exp(-\beta C_{ij})$), а матрица издержек $\{C_{ij}\}$ считается известной.[12]

Зададимся целью найти оценку $\{d_{ij}\}$ при заданной (наблюдаемой) матрице $\{R_{ij}\}$ и известной матрице $\{C_{ij}\}$ методом максимума правдоподобия (на основе теоремы Фишера) [150]. Точнее говоря, оценивать требуется не саму матрицу $\{d_{ij}\}$, а неизвестные параметры $\tilde{A}, \tilde{B}, \beta$. Для возможности применять теорему Фишера [150] и таким образом гарантировать хорошие (асимптотические) свойства полученных оценок (в 1-норме) будем считать, что число оцениваемых параметров $p = 2n + 1$ и объем выборки $N = \sum_{i,j=1}^{n} R_{ij}$ удовлетворяют следующему соотношению: $p/N \ll 1$.

---

[12]Важно заметить, что мы допускаем, следуя В. Г. Спокойному, *model misspecification* [58, 150], т. е. эти предположения неверны и истинное распределение вероятностей $\{R_{ij}\}$ не лежит в этом семействе. Тогда полученное решение по методу максимума правдоподобия можно интерпретировать, как дающее асимптотически (по $N = \sum_{i,j=1}^{n} R_{ij}$) наиболее близкое (по расстоянию Кульбака–Лейблера) распределение в этом семействе к истинному распределению. Другими словами, так полученные параметры являются «асимптотически наилучшими» оценками параметров проекции (по расстоянию Кульбака–Лейблера) истинного распределения вероятностей на выбранное нами параметрическое семейство. Детали см. в [150].

Из определения распределения Пуассона следует, что вероятность реализации корреспонденции $R_{ij}$ может быть посчитана как

$$P\left(R_{ij} \mid d_{ij}\right) = \frac{\exp\left(-d_{ij}\right) \cdot d_{ij}^{R_{ij}}}{R_{ij}!}.$$

Функция правдоподобия (вероятность того, что «выпадет» матрица $\left\{R_{ij}\right\}$, если значения параметров $\tilde{A}, \tilde{B}, \beta$ ) будет иметь следующий вид:

$$\Lambda\left(\left\{R_{ij}\right\} \middle| \tilde{A}, \tilde{B}, \beta\right) = \prod_{i,j=1}^{n} \frac{\exp\left(-d_{ij}\right) \cdot d_{ij}^{R_{ij}}}{R_{ij}!} = \prod_{i,j=1}^{n} \frac{\exp\left(-\tilde{A}_i \tilde{B}_j f\left(C_{ij}\right)\right) \cdot \left(\tilde{A}_i \tilde{B}_j f\left(C_{ij}\right)\right)^{R_{ij}}}{R_{ij}!}.$$

Нам нужно найти точку $\left(\tilde{A}, \tilde{B}, \beta\right)$, в которой достигается максимум функции правдоподобия. Как известно, точка максимума неотрицательной функции (каковой по определению является вероятность) не изменится, если решать задачу максимизации не для исходной функции, а для ее логарифма. Перейдем к логарифму функции правдоподобия:

$$\ln \Lambda = \sum_{i,j=1}^{n} \left(-\tilde{A}_i \tilde{B}_j f\left(C_{ij}\right) + R_{ij} \cdot \left(\ln \tilde{A}_i + \ln \tilde{B}_j + \ln f\left(C_{ij}\right)\right) - \ln R_{ij}!\right).$$

Получаем следующую задачу оптимизации:

$$\sum_{i,j=1}^{n} \left(-\tilde{A}_i \tilde{B}_j f\left(C_{ij}\right) + R_{ij} \cdot \left(\ln \tilde{A}_i + \ln \tilde{B}_j + \ln f\left(C_{ij}\right)\right) - \ln R_{ij}!\right) \to \max_{\tilde{A} \geq 0,\, \tilde{B} \geq 0,\, \beta}.$$

Выпишем условия оптимальности[13]:

$$\frac{\partial \ln \Lambda}{\partial \tilde{A}_i} = \sum_{j=1}^{n} \left(-\tilde{B}_j \exp\left(-\beta C_{ij}\right) + \frac{R_{ij}}{\tilde{A}_i}\right) = 0,$$

$$\frac{\partial \ln \Lambda}{\partial \tilde{B}_j} = \sum_{i=1}^{n} \left(-\tilde{A}_i \exp\left(-\beta C_{ij}\right) + \frac{R_{ij}}{\tilde{B}_j}\right) = 0,$$

$$\frac{\partial \ln \Lambda}{\partial \beta} = \sum_{i,j=1}^{n} \left(\tilde{A}_i \tilde{B}_j C_{ij} \exp\left(-\beta C_{ij}\right) - R_{ij} C_{ij}\right) = 0.$$

Мы получили $2n+1$ уравнение максимума правдоподобия:

---

[13]Легко понять, что максимум не может достигаться на границе. Если допустить, что, скажем, $L_i$ равно нулю в точке максимума, то, поскольку $\partial \ln \Lambda / \partial L_i = \infty$ в этой точке, сдвинувшись немного перпендикулярно гиперплоскости $L_i = 0$ внутрь области определения, мы увеличили бы значения функционала, то есть пришли бы к противоречию с предположением о равенстве нулю $L_i$ в точке максимума.

$$\sum_{j=1}^{n} \tilde{A}_i \tilde{B}_j \exp\left(-\beta C_{ij}\right) = \sum_{j=1}^{n} R_{ij} ,$$

$$\sum_{i=1}^{n} \tilde{A}_i \tilde{B}_j \exp\left(-\beta C_{ij}\right) = \sum_{i=1}^{n} R_{ij} ,$$

$$\sum_{i,j=1}^{n} \left( \tilde{A}_i \tilde{B}_j C_{ij} \exp\left(-\beta C_{ij}\right) - R_{ij} C_{ij} \right) = 0 .$$

Положим по определению $L_i = \sum_{j=1}^{n} R_{ij}$ и $W_j = \sum_{i=1}^{n} R_{ij}$, $\tilde{A}_i = L_i A_i$ и $\tilde{B}_j = W_j B_j$. Тогда

$$d_{ij} = A_i L_i B_j W_j f\left(C_{ij}\right),$$

где $A_i$, $B_j$ – структурные параметры (гравитационной) модели. Их рассчитывают с помощью простого итерационного алгоритма (метод балансировки = метод простых итераций [34, 39, 59, 64, 97]) по значениям $L_i$, $W_j$ следующим образом[14]:

$$A_i = \frac{1}{\sum_{j=1}^{n} B_j W_j f\left(C_{ij}\right)}, \quad B_j = \frac{1}{\sum_{i=1}^{n} A_i L_i f\left(C_{ij}\right)} . \qquad (1.1.6)$$

Последнее же уравнение системы уравнений максимума правдоподобия даст нам

$$\sum_{i,j=1}^{n} \left( d_{ij} C_{ij} - R_{ij} C_{ij} \right) = 0 \text{ или } \sum_{i,j=1}^{n} R_{ij} C_{ij} = \sum_{i,j=1}^{n} d_{ij} C_{ij} . \qquad (1.1.7)$$

**Замечание 1.1.2.** Обратим внимание, что к тем же самым соотношениям можно было прийти (при той же параметрической гипотезе), если вместо предположения: « $R_{ij}$ – независимые случайные величины, распре-деленные по законам Пуассона с математическим ожиданием $M\left(R_{ij}\right)=d_{ij}$, считать, что $\left\{R_{ij}\right\}\gg 1$ имеют мультиномиальное распределение с параметрами $\left\{d_{ij}\Big/\sum_{i,j}d_{ij}\right\}$. Такой подход представляется более естественным, чем изложенный выше. Кроме того, поскольку константа строгой выпуклости энтропии в 1-норме равна 1 [56] (неравенство Пинскера, оценивающее снизу расстояние Кульбака–Лейблера с помощью квадрата 1-нормы), то на базе асимптотического представления логарифма функции правдоподобия в окрестности точки максимума [150] можно построить доверительный интервал для разности оценок и оцениваемых параметров в (наиболее естественной) 1-норме.

---

[14]При $f\left(C_{ij}\right) = \exp\left(-\beta C_{ij}\right)$ эту модель называют также *энтропийной моделью* [64] или *моделью А. Дж. Вильсона* [18, 34]. Энтропийная модель, использованная нами ранее (см. подраздел 1.1.3), для связи матрицы корреспонденций и матрицы издержек является частным случаем гравитационной модели при указанном выборе функции притяжения. Точнее говоря, если бы мы считали, что $C_{ij} = T_{ij}$ не зависит от $\left\{d_{ij}\right\}$ и $R\left(C_{ij}\right) = f\left(C_{ij}\right)/2$, то решение задачи (1.1.2) в точности давало бы гравитационную модель. При этом $A_i$, $B_j$ выражались бы через множители Лагранжа (двойственные переменные) соответственно как $\lambda_i^L$, $\lambda_j^W$. Система (1.1.6) при этом получалась бы при подстановке решения задачи (1.1.2), зависящего от этих $2n$ неизвестных параметров, в ограничения (А) (которых тоже $2n$), точнее говоря, независимых параметров и уравнений было бы $2n$–1 [34].

Соответственно, если критерий остановки (1.1.7) не выполняется, то происходит перерасчет матрицы корреспонденций с учетом обновленной матрицы издержек.

Калибровку $\beta$ можно проводить в случае, когда имеется дополнительная информация о реальных издержках на дорогах. Для этого применяется следующий алгоритм.

Пусть $C_l^*$ – средние издержки на проезд в системе (известны, например, из опросов) для $l$-го слоя спроса [140], например для трудовых миграций.

Алгоритм [106]:

1. Рассчитываем $\beta_l^0=\dfrac{1}{C_l^*}$, $m=0$.
2. Рассчитываем $\left\{d_{ij}\left(\beta_l^m\right)\right\}_{i,j=1}^{n,n}$ – матрицу корреспонденций при $\beta_l=\beta_l^m$.
3. Пересчет 4-*стадийной модели*.
4. Рассчитываем $C_l^m$ – средние издержки на проезд, соответствующие матрице корреспонденций $\left\{d_{ij}\left(\beta_l^m\right)\right\}_{i,j=1}^{n,n}$.
5. При $m\geq 1$ проверяем условие: $\left|C_l^{m-1}-C_l^*\right|\leq\varepsilon$ (критерий остановки).
6. Рассчитываем $C_l^1$ по $\beta_l^1=\dfrac{\beta_l^0C_l^0}{C_l^*}$ (при $m=0$) и полагаем

$$\beta_l^{m+1}=\frac{\left(C_l^*-C_l^{m-1}\right)\beta_l^m+\left(C_l^m-C_l^*\right)\beta_l^{m-1}}{C_l^m-C_l^{m-1}}\quad(\text{при } m\geq 1).$$

7. Переходим на шаг 2.

Данный алгоритм (и его вариации) дает оценку коэффициента $\beta$ при известных средних издержках в сети и является наиболее часто применяемым на практике [140].

Критерий остановки выбран именно таким по следующей причине. Если критерий остановки из подраздела 1.1.5 выполняется, то, подставляя $C_l^{m-1} \simeq C_l^*$ в формулу расчета $\beta_l^{m+1}$ на шаге 6, получаем

$$\beta_l^{m+1} \simeq \frac{\left(C_l^* - C_l^*\right)\beta_l^m + \left(C_l^m - C_l^*\right)\beta_l^{m-1}}{C_l^m - C_l^*} = \beta_l^{m-1}.$$

Полученное значение параметра $\beta\left(C^*\right)$ «соответствует» средним издержкам $C^*$, наблюдаемым «в жизни». Сходимость алгоритма и монотонная зависимость $\beta(C)$ (следовательно, и взаимно-однозначное соответствие) между средними издержками $C$ и значением параметра $\beta$ были показаны в работе [96].

Стоит, однако, отметить, что при численной реализации алгоритма возникает множество проблем. В частности, при «плохом» выборе точки старта или при «плохой» калибровке параметров функций издержек возникают ситуации, в которых алгоритм требует вычисления экстремально больших чисел и превышает размер доступной памяти даже на самых современных машинах (кластерах).

Другой проблемой является отсутствие оценки скорости сходимости алгоритма. Зачастую при моделировании критерий остановки устанавливается жестко. Например, прописывается, что алгоритм калибровки $\beta$ должен быть использован не более 20 раз (или другое, разумное, по мнению модельера, количество итераций). При этом отсутствует четкий критерий качества оценок, полученных таким образом.

Еще одной проблемой описанного подхода является высокая чувствительность к параметрам модели, с одной стороны, и невозможность учесть реальные подтвержденные данные на отдельных элементах сети – с другой. Другими словами, допустим нам известны величины потоков на той или иной дороге каждый день. К сожалению, данную информацию использовать в такой модели не представляется возможным.

### 1.1.6. Модель стабильной динамики

Приведем основные положения модели стабильной динамики, следуя [34, 133]. В рамках модели предполагается, что водители действуют оппортунистически, т. е. выполнен первый принцип Вардропа. Рассмотрим ориентированный граф $\Gamma(V,E)$. В модели каждому ребру $e \in E$ ста-вятся в соответствие параметры $\overline{f}_e$ и $\overline{t}_e$. Они имеют следующую трактовку: $\overline{f}_e$ – максимальная пропускная способность ребра, $e$, $\overline{t}_e$ – минимальные временные издержки на прохождение ребра $e$. Таким образом, сама модель задается графом $\Gamma\left(V, E, \overline{f}, \overline{t}\right)$, где $\overline{f} = \left\{\overline{f}_1, ..., \overline{f}_{|E|}\right\}^T$, $\overline{t} = \left\{\overline{t}_1, ..., \overline{t}_{|E|}\right\}^T$. Пусть $f$ – вектор распределения потоков по ребрам, инициируемый равновесным распределением потоков по маршрутам, а $t$ – вектор временных издержек, соответствующий распределению $f$. Тогда, если транспортная система находится в стабильном состоянии, всегда выполняются неравенства $f \le \overline{f}$ и $t \ge \overline{t}$. При этом считается, что если поток по ребру $f_e$ меньше, чем максимальная пропускная способность ребра $\overline{f}_e$, то все автомобили в потоке двигаются с максимальной скоростью, а их временные издержки $t_e$ минимальны и равны $\overline{t}_e$. Если же поток по ребру $f_e$ становится равным пропускной способности ребра $\overline{f}_e$, то временные издержки водителей $t_e$ могут быть сколь угодно большими. Это удобно объяснить следующим образом. Допустим, на некоторое ребро $e$ стало поступать больше автомобилей, чем оно способно обслужить. Тогда на этом ребре начинает образовываться очередь (пробка). Временные издержки на прохождение ребра $t_e$ складываются из минимальных временных издержек $\overline{t}_e$ и времени, которое водитель вынужден отстоять в пробке. При этом если входящий поток автомобилей на ребро $e$ не снизится до максимально допустимого уровня (пропускной способности ребра), то очередь будет продолжать расти и система не будет находиться в стабильном состоянии. Если же в какой-то момент входящий на ребро $e$ поток снизится до уровня пропускной способности ребра, то в системе наступит равновесие. При этом пробка на ребре $e$ (если входящий поток $f_e$ будет равен $\overline{f}_e$) не будет рассасываться, т. е. временные издержки так и останутся на уровне $t_e$ ($t_e > \overline{t}_e$). Рассмотрим это на примере из статьи [133].

**Пример 1.1.1.** Рассмотрим (в рамках модели стабильной динамики) граф $\Gamma\left(V, E, \overline{t}, \overline{f}\right)$ (см. рис. 1.1.2). Пункты 1 и 2 – потокообразующая пара. При этом выполнено: $d_{12}$ – поток из 1 в 2, $\overline{t}_{\text{up}} < \overline{t}_{\text{down}}$. Если выполнено $d_{12} < \overline{f}_{\text{up}}$, то все водители будут использовать ребро *up*, причем пробка образовываться не будет, так как пропускная способность ребра больше, чем количество желающих проехать (в единицу времени) водителей. В момент, когда $d_{12} = \bar{f}_{\text{up}}$, возможности ребра up будут использоваться на пределе. Если же в какой-то момент величина $f$ станет больше $\bar{f}_{\text{up}}$, то на ребре *up* начнет образовываться пробка. Пробка будет расти до тех пор, пока издержки от использования маршрута *up* не сравняются с издержками от использования маршрута *down*. В этот момент оставшаяся часть начнет использовать маршрут *down*. Если же корреспонденция из 1 в 2 превысит суммарную пропускную способность ребер *up* и *down*, то пробки будут расти неограниченно (входящий поток на ребро будет больше, чем исходящий, соответственно, количество автомобилей в очереди будет расти постоянно), т. е. стабильное распределение в системе никогда не установится. Более подробно рассмотрение данной задачи (и модели) стоит смотреть в работе [133]. ■

Рассмотрим другой модельный пример из статьи [133].

**Пример 1.1.2 (Парадокс Браесса).** Задан граф $\Gamma\left(V, E, \bar{t}, \bar{f}\right)$ (см. рис. 1.1.3). При этом выполнено: $\bar{t}_{13} = 1$ час, $\bar{t}_{12} = 15$ минут, $\bar{t}_{23} = 30$ минут; (1,3) и (2,3) – потокообразующие пары, $d_{13} = 1500$ авт/ч и $d_{23} = 1500$ авт/ч – соответствующие корреспонденции. Пусть пропускные способности всех ребер одинаковы и равны 2000 авт/ч. Тогда равновесие будет такое: 500 авт/ч из 1 будут направляться в 3 через 2, а 1000 авт/ч будут ехать напрямую (для выезжающих из 2 никаких альтернатив нет). При этом все водители, выезжающие из 1, потратят 1 час, а водители, выезжающие из 2, потратят 45 минут. Таким образом, водителям, едущим из 2 в 3 выгодно, чтобы ребро 1–2 отсутствовало, в то время как для водителей, которые едут из 1 в 3, наличие ребра 1–2 безразлично. Другими словами, если бы мы имели власть запретить проезд по ребру 1–2, то часть водителей выиграла бы от такого решения и никто бы не проиграл. Интересно заметить, что для модели Бэкмана это не выполняется. Действительно, в модели Бэкмана, как и для модели стабильной динамики, издержки для водителей, следующих из 1 в 3, равны издержкам на ребре 1–3. Однако они монотонно возрастают от потока на данном ребре. Если бы мы запретили проезд по ребру 1–2, то поток на 1–3 увеличился бы, следовательно, возросли бы и издержки на ребре 1–3. Другими словами, улучшение ситуации для водителей, следующих из 2 в 3, привело бы к ухудшению ситуации для водителей, следующих из 1 в 3. ■

Введем ряд новых обозначений. Пусть $w = (i, j)$ – потокообразующая пара (источник – сток) графа $\Gamma\left(V, E, \bar{t}, \bar{f}\right)$, $P_w$ – множество соответствующих $w$ маршрутов, а $t$ – установившийся на графе вектор времен-ных издержек. Тогда временные издержки, соответствующие самому «быстрому» (наикратчайшему) маршруту из $P_w$ равны: $T_w(t) = \min_{p \in P_w} \sum_{e \in E} \delta_{ep} t_e$ . Функция $T_w(t)$ зависит от вектора временных издержек $t$ . Важно заметить, что функция $T_w(t)$ и её супердифференциал эффективно вычисляются, например, алгоритмом Дейкстры [50, 68] или более быстрыми современными методами [153].

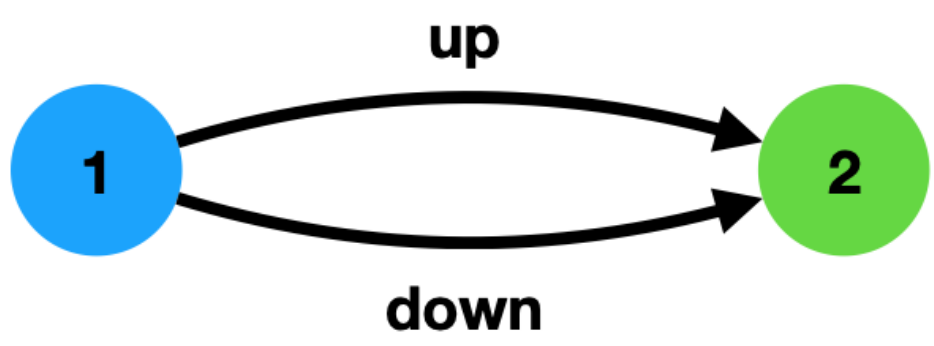

Рис. 1.1.2

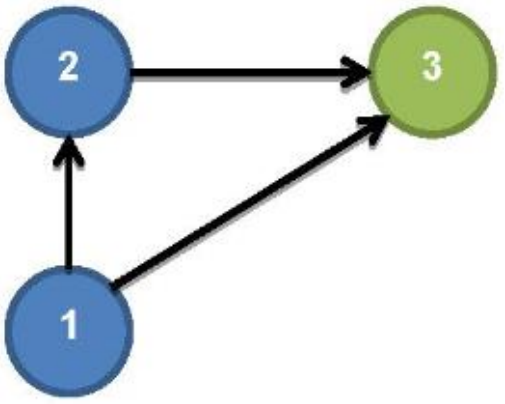

Рис. 1.1.3

Пусть корреспонденция потокообразующей пары $w$ равна $d_w$. Тогда, если $f^*$ является равновесным распределением потоков для заданного графа $\Gamma(V, E, \bar{t}, \bar{f})$, а $t^*$ – соответствующий равновесному распределению вектор временных издержек, то $f^* = \sum_{w \in W} f_w^*$ , где $f_w^*$ – вектор распределения потока, порождаемого потокообразующей парой $w$ ( $w \in OD$ ). При этом $f_w^*$ удовлетворяет

$$f_w^* \in d_w \partial T_w(t),$$

где $\partial T_w(t)$ – супердифференциал, который можно посчитать с помощью теоремы Милютина–Дубовицкого [53], используя структуру функции $T_w(t)$. Следовательно, по теореме Моро–Рокафеллара [53] $f^* \in \sum_{w \in W} d_w \partial T_w(t)$ . Остается вопрос: как искать равновесный вектор временных издержек $t^*$ .

**Теорема 1.1.2 [133].** *Распределение потоков* $f^*$ *и вектор временных издержек* $t^*$ *являются равновесными для графа* $\Gamma(V, E, \bar{t}, \bar{f})$ *заданного множества потокообразующих пар* $W$ *и соответствующих им потоков* $d_w$ *тогда и только тогда, когда* $t^*$ *является решением следующий* *задачи вогнутой оптимизации* (*максимизация количества свободно движущихся автомобилей, то есть автомобилей не стоящих в пробках*):

$$\max_{t \ge \bar{t}} \left\{ \sum_{w \in W} d_w T_w(t) - \left\langle \bar{f}, t - \bar{t} \right\rangle \right\}, \qquad (1.1.8)$$

*где* $f^* = \bar{f} - s^*$, *где* $s^*$ – (*оптимальный*) *вектор двойственных множителей для ограничений* $t \ge \bar{t}$ *в задаче* (1.1.8) *или, другими словами, решение двойственной задачи.*

**Схема доказательства.** Резюмируем далее основные постулаты модели стабильной динамики:

1. $t_e \ge \bar{t}_e$, $f_e \le \bar{f}_e$, $e \in E$;
2. $(t_e - \bar{t}_e)(\bar{f}_e - f_e) = 0$, $e \in E$;
3. $f \in \partial \sum_{w \in W} d_w T_w(t)$ (принцип Нэша–Вардропа).

Наша цель – подобрать такую задачу выпуклой (вогнутой) оптимизации, решение которой давало бы пару $(t, f)$, удовлетворяющую этим постулатам. С задачей выпуклой оптимизации намного удобнее работать (разработаны эффективные способы решения таких задач [56]), чем с описанием 1–3. Если бы мы могли ограничиться только п. 3, то тогда этот пункт можно было бы переписать следующим образом:

$$\max_{t \ge 0} \left\{ \sum_{w \in W} d_w T_w(t) - \left\langle f, t - \bar{t} \right\rangle \right\}, \qquad (1.1.9)$$

но есть еще п. 1, 2. Подсказка содержится в п. 2, который удобно понимать как условие дополняющей нежесткости [53] ограничения $t \ge \bar{t}$ (в задаче (1.1.9)) с множителями Лагранжа $s = \bar{f} - f \ge 0$. Причем нет необходимости оговаривать дополнительно, что $f \le \bar{f}$, поскольку множители Лагранжа к ограничениям вида неравенство автоматически неотрицательны. Таким образом, мы приходим к (1.1.8). Все же необходимо оговориться, что получившаяся задача (1.1.8) хотя и является довольно простой задачей (несложно заметить, что это просто задача линейного программирования (ЛП), записанная в более компактной форме; действительно: $\max_{t \ge 0;\, \{T_w\};\, C} \left\{ C - \left\langle \bar{f}, t - \bar{t} \right\rangle; T_w \le \sum_{e \in E} \delta_{ep} t_e,\ p \in P_w,\ w \in W;\ C = \sum_{w \in W} d_w T_w \right\}$), но требуется также найти двойственные множители (чтобы определить $f$). ■

Задача (1.1.8) имеет (конечное) решение тогда и только тогда, когда существует хотя бы один вектор потоков $f \le \bar{f}$, соответствующий заданным корреспонденциям. Другими словами, существует способ так организовать движение, чтобы при заданных корреспонденциях не нарушались условия $f \le \bar{f}$. В противном случае равновесного (стационарного) режима в системе не будет, и со временем весь граф превратится в одну большую пробку, которая начнет нарушать условия заданных корреспонденций, не позволяя новым пользователям приходить в сеть с той интенсивностью, с которой они этого хотят.

В отличие от модели Бэкмана, в которой при монотонно возрастающих функциях затрат на ребрах равновесие (по потокам на ребрах) единственно, в модели стабильной динамики это может быть не так. Хотя типичной ситуацией (общим положением) будет единственность равновесия, поскольку поиск стабильной конфигурации эквивалентен решению задачи ЛП, все же в определенных «вырожденных» ситуациях может возникать неединственность, а именно: если в примере 1.1.1 (см. рис. 1.1.2) $\bar{f}_{\text{ up}} = d_{12}$, то возникает вопрос: как выбрать единственное равновесие, которое реализуется? Один из ответов имеется в [133]: нужно следить за динамикой во времени $d_{12}(t)$, а именно за

$$\int_0^T \left(d_{12}(t) - \bar{f}_{up}\right) dt < \infty .$$

Другой ответ связан с прямым осуществлением для данного примера предельного перехода, описанного в следующем пункте.

По сути, с помощью этого и следующего пункта мы установим, что если

$$\tau_e^{\mu}(f_e) \xrightarrow[\mu \to 0+]{} \begin{cases} \bar{t}_e, & 0 \le f_e < \bar{f}_e, \\ [\bar{t}_e, \infty), & f_e = \bar{f}_e, \end{cases}$$

$$d\tau_e^{\mu}(f_e)/df_e \xrightarrow[\mu \to 0+]{} 0, \; 0 \le f_e < \bar{f}_e,$$

$x(\mu)$ – равновесное распределение потоков по путям в модели Бэкмана при функциях затрат на ребрах $\tau_e^{\mu}(f_e)$, то

$$\tau_e^{\mu}\left(f_e\left(x(\mu)\right)\right) \xrightarrow[\mu \to 0+]{} t_e,$$

$$f_e\left(x(\mu)\right) \xrightarrow[\mu \to 0+]{} f_e,$$

где пара $(t, f)$ – равновесие в модели стабильной динамики с тем же графом и матрицей корреспонденций, что и в модели Бэкмана, и с ребрами, характеризующимися набором $(\bar{t}, \bar{f})$ из определения функций $\tau_e^\mu(f_e)$.

Численно решать задачу (1.1.1) можно различными способами (см., например, [27]).

## 1.1.7. Связь модели стабильной динамики с моделью Бэкмана

Получим модель стабильной динамики в результате предельного перехода из модели Бэкмана. Для этого будем считать, что $\tau_e(f_e) = \bar{t}_e \cdot \left(1 - \mu \ln\left(1 - f_e / \bar{f}_e\right)\right)$.

Перепишем исходную задачу поиска равновесного распределения потоков (здесь используются обозначения: $\sigma_e^*(t_e)$ – сопряженная функция к $\sigma_e(f_e)$):

$$\min_{f,x}\left\{\sum_{e\in E}\sigma_e(f_e):\ f = \Theta x,\ x\in X\right\} = \min_{f,x}\left\{\sum_{e\in E}\max_{t_e\in \operatorname{dom}\sigma_e^*}\left[f_e t_e - \sigma_e^*(t_e)\right]:\ f = \Theta x,\ x \in X\right\} =$$

$$= \max_{t\in \operatorname{dom}\sigma^*}\left\{\min_{f,x}\left[\sum_{e\in E} f_e t_e:\ f = \Theta x,\ x\in X\right] - \sum_{e\in E}\sigma_e^*(t_e)\right\}.$$

Найдем $\sum_{e\in E}\sigma_e^*(t_e)$:

$$\sigma_e^*(t_e) = \sup_{f_e}\left(t_e\cdot f_e - \int_0^{f_e}\tau(z)\,dz\right) = \sup_{f_e}\left(t_e\cdot f_e - \int_0^{f_e}\bar{t}_e\cdot\left(1-\mu\ln\left(1 - z/\bar{f}_e\right)\right)dz\right) =$$

$$= \sup_{f_e}\left(\left(t_e - \bar{t}_e\right) f_e + \bar{t}_e\mu\int_0^{f_e}\ln\left(1 - z/\bar{f}_e\right)dz\right).$$

Из принципа Ферма [53] найдем $f_e$, доставляющее максимум:

$$\frac{\partial}{\partial f_e}\left(\left(t_e - \bar{t}_e\right) f_e + \bar{t}_e\mu\int_0^{f_e}\ln\left(1 - z/\bar{f}_e\right)dz\right) = 0 \ \Rightarrow$$

$$\Rightarrow\ \exp\left(-\frac{t_e - \bar{t}_e}{\bar{t}_e\mu}\right) = 1 - f_e/\bar{f}_e\ \Rightarrow$$

$$\Rightarrow\ f_e = \bar{f}_e\cdot\left(1 - \exp\left(-\frac{t_e - \bar{t}_e}{\bar{t}_e\mu}\right)\right).$$

Подставляем в $\sigma_e^*(t_e)$:

$$(t_e-\bar{t}_e)\,f_e+\bar{t}_e\mu\int_0^{f_e}\ln\left(1-z/\bar{f}_e\right)dz=(t_e-\bar{t}_e)\,f_e-\bar{t}_e\cdot\bar{f}_e\cdot\mu\int_1^{1-f_e/\bar{f}_e}\ln z\,dz=$$

$$=(t_e-\bar{t}_e)\,f_e-\bar{t}_e\cdot\bar{f}_e\cdot\mu\left(\left(1-\frac{f_e}{\bar{f}_e}\right)\ln\left(1-\frac{f_e}{\bar{f}_e}\right)-1+\frac{f_e}{\bar{f}_e}+1\right)=$$

$$=(t_e-\bar{t}_e)\,f_e-\left(-\left(\bar{f}_e-f_e\right)(t_e-\bar{t}_e)+f_e\cdot\bar{t}_e\cdot\mu\right)=$$

$$=(t_e-\bar{t}_e)\,\bar{f}_e-f_e\cdot\bar{t}_e\cdot\mu=(t_e-\bar{t}_e)\,\bar{f}_e-\bar{f}_e\cdot\bar{t}_e\cdot\mu\left(1-\exp\left(-\frac{t_e-\bar{t}_e}{\bar{t}_e\mu}\right)\right).$$

Возвращаясь к исходной задаче, имеем

$$\max_{t\in\operatorname{dom}\sigma^*}\left\{\min_{f,x}\left[\sum_{e\in E}f_et_e:\ f=\Theta x,\ x\in X\right]-\sum_{e\in E}\sigma_e^*(t_e)\right\}=$$

$$=\max_{t\in\operatorname{dom}\sigma^*}\left\{\sum_{w\in W}d_wT_w(t)-\left\langle\bar{f},t-\bar{t}\right\rangle+\mu\sum_{e\in E}\bar{f}_e\cdot\bar{t}_e\left(1-\exp\left(-\frac{t_e-\bar{t}_e}{\bar{t}_e\mu}\right)\right)\right\}\overset{\mu\to0+}{=}$$

$$\overset{\mu\to0+}{=}\max_{t\geq\bar{t}}\left\{\sum_{w\in W}d_wT_w(t)-\left\langle\bar{f},t-\bar{t}\right\rangle\right\}.$$

Возможность поменять местами порядок взятия максимума и минимума следует из теоремы фон Неймана о минимаксе [63].

Выбор именно функции $\tau_e(f_e)=\bar{t}_e\cdot\left(1-\mu\ln\left(1-f_e/\bar{f}_e\right)\right)$ не является определяющим. Вместо этой функции может быть взят любой другой гладкий внутренний барьер области $f_e<\bar{f}_e$, например $\tau_e(f_e)=\bar{t}_e\cdot\left(1+\mu\frac{\bar{f}_e}{\bar{f}_e-f_e}\right)$: такую функцию вполне естественно использовать при поиске равновесного распределения пользователей сети общественного транспорта [133]. Найдем $\sigma_e^*(t_e)$:

$$\sigma_e^*(t_e)=\sup_{f_e}\left(t_e\cdot f_e-\int_0^{f_e}\bar{t}_e\cdot\left(1+\mu\frac{z}{\bar{f}_e-z}\right)dz\right)=$$

$$=\sup_{f_e}\left((t_e-\bar{t}_e)\cdot f_e+\bar{t}_e\cdot\bar{f}_e\cdot\mu\cdot\ln\left(1-\frac{f_e}{\bar{f}_e}\right)+\bar{t}_e\cdot f_e\cdot\mu\right).$$

Вновь выписывая условия оптимальности первого порядка, имеем

$$\frac{\partial}{\partial f_e}\left((t_e-\overline{t}_e)\cdot f_e+\overline{t}_e\cdot\overline{f}_e\cdot\mu\cdot\ln\left(1-\frac{f_e}{\overline{f}_e}\right)+\overline{t}_e\cdot f_e\cdot\mu\right)=(t_e-\overline{t}_e)-\overline{t}_e\cdot\mu\cdot\frac{f_e}{\overline{f}_e-f_e}+\overline{t}_e\cdot\mu=0$$

$$\Rightarrow\ f_e=\overline{f}_e\cdot\left(1-\frac{\overline{t}_e\cdot\mu}{t_e-(1-\mu)\overline{t}_e}\right)$$

Подставляя в $\sigma_e^*(t_e)$, имеем

$$\sigma_e^*(t_e)=$$

$$=(t_e-\overline{t}_e)\cdot\overline{f}_e\cdot\left(1-\frac{\overline{t}_e\cdot\mu}{t_e-(1-\mu)\overline{t}_e}\right)+\overline{t}_e\cdot\overline{f}_e\cdot\mu\cdot\ln\left(\frac{\overline{t}_e\cdot\mu}{t_e-(1-\mu)\overline{t}_e}\right)+\overline{t}_e\cdot\mu\cdot\overline{f}_e\cdot\left(1-\frac{\overline{t}_e\cdot\mu}{t_e-(1-\mu)\overline{t}_e}\right)=$$

$$=(t_e-\overline{t}_e)\cdot\overline{f}_e+\overline{t}_e\cdot\overline{f}_e\cdot\mu\cdot\ln\left(\frac{\overline{t}_e\cdot\mu}{t_e-(1-\mu)\overline{t}_e}\right).$$

В итоге получаем

$$\max_{t\in\operatorname{dom}\sigma^*}\left\{\min_{f,x}\left[\sum_{e\in E}f_et_e:\ f=\Theta x,\ x\in X\right]-\sum_{e\in E}\sigma_e^*(t_e)\right\}=$$

$$=\max_{t\in\operatorname{dom}\sigma^*}\left\{\sum_{w\in W}d_wT_w(t)-\left\langle\overline{f},t-\overline{t}\right\rangle+\mu\sum_{e\in E}\overline{f}_e\cdot\overline{t}_e\cdot\ln\left(1+\frac{t_e-\overline{t}_e}{\overline{t}_e\cdot\mu}\right)\right\}\overset{\mu\to0+}{=}$$

$$=\max_{t\ge\overline{t}}\left\{\sum_{w\in W}d_wT_w(t)-\left\langle\overline{f},t-\overline{t}\right\rangle\right\}.$$

Наконец рассмотрим третий вариант выбора функции $\tau_e(f_e)$:

$$\tau_e(f_e)=\overline{t}_e\cdot\left(1+\gamma\cdot\left(\frac{f_e}{\overline{f}_e}\right)^{\frac{1}{\mu}}\right).$$

Функции такого вида называются BPR-*функциями* и наиболее часто применяются при моделировании. Так, при использовании модели Бэкмана обычно полагают $\mu=0.25$, а значение параметра $\gamma$ варьируется от 0.15 до 2 и определяется типом дороги [140]. Найдем $\sigma_e^*(t_e)$:

$$\sigma_e^*(t_e)=\sup_{f_e}\left(t_e\cdot f_e-\int_0^{f_e}\overline{t}_e\cdot\left(1+\gamma\cdot\left(\frac{z}{\overline{f}_e}\right)^{\frac{1}{\mu}}\right)dz\right)=\sup_{f_e}\left((t_e-\overline{t}_e)\cdot f_e-\overline{t}_e\cdot\gamma\cdot\int_0^{f_e}\left(\frac{z}{\overline{f}_e}\right)^{\frac{1}{\mu}}dz\right)=$$

$$=\sup_{f_e}\left((t_e-\overline{t}_e)\cdot f_e-\overline{t}_e\cdot\frac{\mu}{1+\mu}\cdot\gamma\cdot\frac{f_e^{1+\frac{1}{\mu}}}{\overline{f}_e^{\frac{1}{\mu}}}\right).$$

Аналогично рассуждая, получаем

$$\frac{\partial}{\partial f_e}\left(\left(t_e-\bar{t}_e\right)\cdot f_e-\bar{t}_e\cdot\frac{\mu}{1+\mu}\cdot\gamma\cdot\frac{f_e^{1+\frac{1}{\mu}}}{\bar{f}_e^{\frac{1}{\mu}}}\right)=t_e-\bar{t}_e-\bar{t}_e\cdot\gamma\cdot\frac{f_e^{\frac{1}{\mu}}}{\bar{f}_e^{\frac{1}{\mu}}}=0\Rightarrow f_e=\bar{f}_e\cdot\left(\frac{t_e-\bar{t}_e}{\bar{t}_e\cdot\gamma}\right)^{\mu}.$$

Тогда

$$\sigma_e^*\left(t_e\right)=\sup_{f_e}\left(\left(t_e-\bar{t}_e\right)\cdot f_e-\bar{t}_e\cdot\frac{\mu}{1+\mu}\cdot\gamma\cdot\frac{f_e^{1+\frac{1}{\mu}}}{\bar{f}_e^{\frac{1}{\mu}}}\right)=\bar{f}_e\cdot\left(\frac{t_e-\bar{t}_e}{\bar{t}_e\cdot\gamma}\right)^{\mu}\frac{\left(t_e-\bar{t}_e\right)}{1+\mu}.$$

В итоге получаем

$$\max_{t\in\operatorname{dom}\sigma^*}\left\{\sum_{w\in W}d_wT_w(t)-\sum_{e\in E}\bar{f}_e\cdot\left(\frac{t_e-\bar{t}_e}{\bar{t}_e\cdot\gamma}\right)^{\mu}\frac{\left(t_e-\bar{t}_e\right)}{1+\mu}\right\}\overset{\mu\to0+}{=}\max_{t\geq\bar{t}}\left\{\sum_{w\in W}d_wT_w(t)-\left\langle\bar{f},t-\bar{t}\right\rangle\right\}.$$

**Замечание 1.1.3.** Естественно в контексте всего написанного выше теперь задаться вопросом: а можно ли получить модель стабильной динамики эволюционным образом, то есть подобно тому, как в конце подраздела 1.1.4 была эволюционно проинтерпретирована модель Бэкмана? Действительно, рассмотрим логит-динамику (с $T\to0+$) или, скажем, просто имитационную логит-динамику [145]. Предположим, что гладкие, возрастающие, выпуклые функции затрат на всех ребрах $\tau_e\left(f_e\right)$ при $\mu\to0+$ превращаются в «ступеньки» (многозначные функции):

$$\tau_e\left(f_e\right)=\begin{cases}\bar{t}_e, & 0\leq f_e<\bar{f}_e,\\ \left[\bar{t}_e,\infty\right), & f_e=\bar{f}_e.\end{cases}$$

Согласно подразделу 1.1.4 равновесная конфигурация при таком переходе $\mu\to0+$ должна находиться из решения задачи

$$\sum_{e\in E}\int_0^{f_e}\tau_e(z)\,dz\to\min_{f=\Theta x,\,x\in X}.$$

Считая, что в равновесии не может быть $\tau_e\left(f_e\right)=\infty$ (иначе равновесие просто не достижимо, и со временем весь граф превратится в одну большую пробку, см. конец подраздела 1.1.6), можно не учитывать в интеграле вклад точек $\bar{f}_e$ (в случае попадания в промежуток интегрирования), то есть переписать задачу следующим образом:

$$\min_{f=\Theta x,\,x\in X}\sum_{e\in E}\int_0^{f_e}\left(\bar{t}_e+\delta_{\bar{f}_e}(z)\right)dz\Leftrightarrow\min_{\substack{f=\Theta x,\,x\in X\\ f\leq\bar{f}}}\sum_{e\in E}f_e\bar{t}_e,$$

где

$$\delta_{\bar{f}_e}(z) = \begin{cases} 0, & 0 \le z < \bar{f}_e, \\ \infty, & z \ge \bar{f}_e. \end{cases}$$

Мы получили задачу линейного программирования. Интересно заметить, что двойственные множители $\eta \ge 0$ к ограничениям $f \le \bar{f}$ имеют размерность («физический смысл») времени. Из решения только что выписанной задачи линейного программирования сразу может быть неясно, чему равно это время. Но если перейти к двойственной задаче, то мы получим уже рассматриваемую нами ранее задачу (1.1.9) из подраздела 1.1.6. Причем двойственные множители $\eta$ связаны с временами проезда ребер $t$ следующим образом: $\eta = t - \bar{t} \ge 0$, то есть, действительно, получают содержательную интерпретацию времени, потерянного на ребрах (дополнительно к времени проезда по свободной дороге) из-за наличия пробок. Еще одним подтверждением только что сказанному является то, что если $f_e = \bar{f}_e$, то на ребре $e$ пробка и время прохождения этого ребра $t_e \ge \bar{t}_e$ ($\eta \ge 0$), а если пробки нет, то $f_e < \bar{f}_e$ и $t_e = \bar{t}_e$ ($\eta = 0$). То есть имеют место условия дополняющей нежесткости.

Из данного замечания становится ясно, как оптимально назначать платы за проезд в модели стабильной динамики. Переведем предварительно время в деньги. Назначая платы за проезд по ребрам графа в соответствии с вектором $\eta$, получим, что равновесное распределение пользователей транспортной сети по этой сети будет описываться парой $\left(f, \bar{t}\right)$, что, очевидно, лучше, чем до введения плат $\left(f, t\right)$. Действительно, задачу поиска равновесных потоков $f$:

$$\min_{\substack{f = \Theta x,\, x \in X \\ f \le \bar{f}}} \sum_{e \in E} f_e \bar{t}_e$$

можно переписать с помощью метода множителей Лагранжа следующим образом:

$$\min_{f = \Theta x,\, x \in X} \sum_{e \in E} \left( f_e \bar{t}_e + \eta_e \cdot \left( f_e - \bar{f}_e \right) \right).$$

Другими словами, у этих двух задач одинаковые решения $f$; но $t = \bar{t} + \eta$, поэтому тот же самый вектор $f$ будет доставлять решение задаче

$$\min_{f = \Theta x,\, x \in X} \sum_{e \in E} f_e t_e .$$

Поскольку мы знаем, что $f \le \bar{f}$, то $f$ будет доставлять решение также и этой задаче

$$\min_{\substack{f=\Theta x,\, x\in X \\ f\le \bar{f}}} \sum_{e\in E} f_e t_e \ .$$

Таким образом, один и тот же вектор $f$ отвечает оптимальному (с точки зрения социального оптимума) распределению потоков в графе с ребрами, взвешенными согласно векторам $\bar{t}$ и $t$. Более того, несложно понять, что все это остается верным не только для двух векторов $\bar{t}$ и $t$, но и для целого «отрезка» векторов: $\tilde{t} = \alpha\bar{t} + (1-\alpha)t$, $\alpha \in [0,1]$. К сожалению, если не взимать платы за проезд, то всегда реализуется сценарий $\alpha = 0$, то есть вектор затрат будет $t$ (объяснению причин были посвящены подразделы 1.1.6 и 1.1.7), но если мы взимаем платы согласно вектору $\eta_\alpha = \alpha\cdot(t-\bar{t})$, то вектор реальных затрат (не учитывающих платы) будет $t-\eta_\alpha = (1-\alpha)t + \alpha\bar{t}$. Следовательно, оптимально выбирать $\alpha = 1$, то есть платы согласно $\eta = t - \bar{t}$.

Собственно, всю эту процедуру (назначения плат) можно делать адаптивно: постепенно увеличивая плату за проезд на тех ребрах, на которых наблюдаются пробки. И делать это нужно до тех пор, пока пробки не перестанут появляться.

## 1.1.8. Учет расщепления потоков по способам передвижений

Распространим модель стабильной динамики на тот случай, когда (все) пользователи (игроки) транспортной сети имеют возможность выбирать между двумя альтернативными видами транспорта: личным и общественным [9]. Соответственно, теперь у нас имеется информация о $\Gamma_л\left(V, E^л, \bar{t}^{\,л}, \bar{f}^{\,л}\right)$ – для личного транспорта и $\Gamma_o\left(V, E^o, \bar{t}^{\,o}, \bar{f}^{\,o}\right)$ – для общественного транспорта. Аналогично рассматривается общий случай. Имеет место следующий результат, получаемый аналогично [133].

**Теорема 1.1.3**. *Распределение потоков* $f^* = \left(f^л, f^o\right)$ *и вектор временных издержек* $t^* = \left(t^л, t^o\right)$ *являются равновесными для графа* $\Gamma\left(V, E, \bar{t}, \bar{f}\right)$ *заданного множества потокообразующих пар* $W$ *и соответствующих им потоков* $d_w$ *тогда и только тогда, когда* $t^*$ *является решением следующий задачи выпуклой оптимизации* (*в дальнейшем для* *нас будет удобно писать эту задачу как задачу минимизации, а не максимизации, как раньше*):

$$\min_{\substack{t^л \geq \bar{t}^л \\ t^o \geq \bar{t}^o}} \left\{ \left\langle \bar{f}^л, t^л - \bar{t}^л \right\rangle + \left\langle \bar{f}^o, t^o - \bar{t}^o \right\rangle - \sum_{w \in W} d_w \min\left\{ T_w^л\left(t^л\right), T_w^o\left(t^o\right) \right\} \right\}, \quad (1.1.10)$$

*где* $f^л = \bar{f}^л - s^л$, $f^o = \bar{f}^o - s^o$, *а* $s^л$, $s^o$ *– векторы* (*оптимальные*) *двойственных множителей для ограничений* $t^л \geq \bar{t}^л$, $t^o \geq \bar{t}^o$ *в задаче* (1.1.10) *или, другими словами, решения двойственной задачи.*

Если мы хотим учитывать возможность пересаживания пользователей сети в пути с личного транспорта на общественный (и наоборот), то выражение $\min\left\{ T_w^л\left(t^л\right), T_w^o\left(t^o\right) \right\}$ нужно будет немного изменить. Мы не будем здесь вдаваться в детали, скажем лишь, что с точки зрении всего дальнейшего это непринципиально. Более того, хотя на время работы алгоритма это и скажется (время увеличится), но некритическим образом, поскольку решение задачи о кратчайшем маршруте, которое возникает на каждом шаге субградиентного спуска, естественным образом (небольшим раздутием графа) обобщается на случай, когда есть несколько весов ребер, отвечающих разным типам транспортных средств, и в вершинах графа есть затраты на пересадку (изменение транспортного средства), при условии наличия возможности её осуществления.

Отметим, что единого подхода для моделирования общественного транспорта пока не существует. Указанный нами метод соответствует подходу, когда предполагается, что пассажиры выбирают только «оптимальные стратегии», т. е. соответствует концепции равновесия по Нэшу. Другой, альтернативный подход, который сейчас чаще используется на практике, следует концепции «стохастического равновесия». В нем предполагается, что пассажиры выбирают каждый из возможных маршрутов с некоторой вероятностью, зависящей от ожидаемых издержек в сети и ожидаемого времени ожидания соответствующего маршрута на остановочной станции. Об этом немного подробнее будет рассказано в подразделе 1.1.10.

Обратим внимание на то, что общественный транспорт все же более естественно описывать моделью Бэкмана, а не моделью стабильной динамики. Во всяком случае, к такому выводу склоняются авторы модели стабильной динамики [133]. В этой связи далее приводится смешанная модель, в которой личный транспорт описывается моделью стабильной динамики, а общественный – моделью Бэкмана (см. подраздел 1.1.7):

$$\min_{\substack{t^л \geq \bar{t}^л \\ t^o \geq \bar{t}^o \cdot (1-\mu)}} \left\{ \left\langle \bar{f}^л, t^л - \bar{t}^л \right\rangle + \left\langle \bar{f}^o, t^o - \bar{t}^o \right\rangle - \mu \sum_{e \in E} \bar{f}_e^o \cdot \bar{t}_e^o \cdot \ln\left( 1 + \frac{t_e^o - \bar{t}_e^o}{\bar{t}_e^o \cdot \mu} \right) - \sum_{w \in W} d_w \min\left\{ T_w^л\left(t^л\right), T_w^o\left(t^o\right) \right\} \right\},$$

где $f^{л} = \bar{f}^{л} - s^{л}$, $f_e^o = \bar{f}_e^o \cdot \left(1 - \frac{\bar{t}_e^o \cdot \mu}{t_e^o - (1-\mu)\bar{t}_e^o}\right)$, $s^{л}$ – (оптимальный) вектор двойственных множителей для ограничений $t^{л} \geq \bar{t}^{л}$.

Несложно распространить все проводимые в следующих пунктах рассуждения именно на такую версию модели стабильной динамики с учетом расщепления по типам передвижений.

Отметим также, что простая динамическая модель, описывающая процесс расщепления, приведена в приложении 2.

### 1.1.9. Трехстадийная модель стабильной динамики

Объединим теперь задачи (1.1.2) и (1.1.10) в одну задачу[15]:

$$\min_{\lambda^L, \lambda^W} \max_{\sum_{i,j=1}^n d_{ij}=1,\, d_{ij}\geq 0} \left[ -\sum_{i,j=1}^n d_{ij} \ln d_{ij} + \sum_{i=1}^n \lambda_i^L \left( l_i - \sum_{j=1}^n d_{ij} \right) + \sum_{j=1}^n \lambda_j^W \left( \mathrm{w}_j - \sum_{i=1}^n d_{ij} \right) + \right.$$

$$\left. + \beta \min_{\substack{t^{л} \geq \bar{t}^{л} \\ t^o \geq \bar{t}^o}} \left\{ \langle \bar{f}^{л}, t^{л} - \bar{t}^{л} \rangle + \langle \bar{f}^o, t^o - \bar{t}^o \rangle - \sum_{i,j=1}^n d_{ij} \min\left\{ T_{ij}^{л}(t^{л}), T_{ij}^o(t^o) \right\} \right\} \right] =$$

$$= \min_{\substack{t^{л} \geq \bar{t}^{л} \\ t^o \geq \bar{t}^o \\ \lambda^L, \lambda^W}} \left\{ \max_{\sum_{i,j=1}^n d_{ij}=1,\, d_{ij}\geq 0} \left[ -\sum_{i,j=1}^n d_{ij} \ln d_{ij} - \beta \sum_{i,j=1}^n d_{ij} \min\left\{ T_{ij}^{л}(t^{л}), T_{ij}^o(t^o) \right\} - \right. \right.$$

$$\left. \left. - \sum_{i,j=1}^n d_{ij} \cdot \left( \lambda_i^L + \lambda_j^W \right) \right] + \sum_{i=1}^n \lambda_i^L l_i + \sum_{j=1}^n \lambda_j^W \mathrm{w}_j + \beta \langle \bar{f}^{л}, t^{л} - \bar{t}^{л} \rangle + \beta \langle \bar{f}^o, t^o - \bar{t}^o \rangle \right\} =$$

$$= \min_{\substack{t^{л} \geq \bar{t}^{л} \\ t^o \geq \bar{t}^o \\ \lambda^L, \lambda^W}} \left\{ \ln \left( \sum_{i,j=1}^n \exp\left( -\beta \min\left\{ T_{ij}^{л}(t^{л}), T_{ij}^o(t^o) \right\} - \lambda_i^L - \lambda_j^W \right) \right) + \right.$$

$$\left. + \sum_{i=1}^n \lambda_i^L l_i + \sum_{j=1}^n \lambda_j^W \mathrm{w}_j + \beta \langle \bar{f}^{л}, t^{л} - \bar{t}^{л} \rangle + \beta \langle \bar{f}^o, t^o - \bar{t}^o \rangle \right\}, \qquad (1.1.11)$$

причем $d_{ij} = \tilde{Z}^{-1} \exp\left( -\beta \min\left\{ T_{ij}^{л}(t^{л}), T_{ij}^o(t^o) \right\} - \lambda_i^L - \lambda_j^W \right)$, где $\tilde{Z}^{-1}$ ищется из условия нормировки, $f^* = (f^{л}, f^o)$, $f^{л} = \bar{f}^{л} - s^{л}$, $f^o = \bar{f}^o - s^o$, а $s^{л}$, $s^o$ – векторы (оптимальные) двойственных множителей для ограничений $t^{л} \ge \bar{t}^{л}$, $t^{o} \ge \bar{t}^{o}$ в задаче (1.1.11). Таким образом, все, что осталось сделать, это решить задачу негладкой выпуклой оптимизации (1.1.11) прямодвойственным методом. Заметим при этом, что на переменные $t^{л}$, $t^{o}$ – ограничения сверху возникают из вполне естественных соображений. Пусть $\bar{f}$ максимально возможный поток на ребре, поскольку нас интересует оценка сверху, то можно считать $t > \bar{t}$, стало быть, $f = \bar{f}$. Представим себе самую плохую ситуацию: всё ребро стоит в пробке. Пусть $L$ – длина ребра , а $l$ – средняя длина автомобиля, $r$ – число полос. Тогда $t \le Lr/\left(l\bar{f}\right) + \bar{t}$ не может превышать нескольких часов. Задачу (1.1.11) можно решать, например, прямодвойственным универсальным композитным методом подобных треугольников [37] или методом из работ [28, 30]. Подробнее об этом будет написано в следующей главе и приложении 3.

[15]Отметим, что все приводимые далее выкладки (а также выкладки подраздела 1.1.10) можно провести, отталкиваясь не от задачи (1.1.2), а от её упрощенного варианта, описанного в сноске 9.

Отметим, что если известна информация о потоках по ряду дуг $f_e^{л} = \tilde{f}_e^{л}$, $e \in \mathrm{A}^{л}$; $f_e^{o} = \tilde{f}_e^{o}$, $e \in \mathrm{A}^{o}$, то эту информацию можно «зашить» в модель (ЗС) подобно тому, как это делается в работе [133], а именно: брать минимум по множеству $t_e^{л} \ge \bar{t}_e^{л}$, $e \in E \setminus \mathrm{A}^{л}$, $t_e^{л} \ge 0$, $e \in \mathrm{A}^{л}$, $t_e^{o} \ge \bar{t}_e^{o}$, $e \in E \setminus \mathrm{A}^{o}$, $t_e^{o} \ge 0$, $e \in \mathrm{A}^{o}$, а слагаемые $+\beta\left\langle \bar{f}^{л}, t^{л} - \bar{t}^{л} \right\rangle + \beta\left\langle \bar{f}^{o}, t^{o} - \bar{t}^{o} \right\rangle$ стоит заменить на (параметры $\bar{t}_e^{л}$, $e \in \mathrm{A}^{л}$ и $t_e^{o}$, $e \in \mathrm{A}^{o}$ – неизвестны, но они и не нужны для расчетов, поскольку входят в виде аддитивных констант в функционал):

$$+\beta \sum_{e \in E \setminus \mathrm{A}^{л}} \bar{f}_e^{л} \cdot \left(t_e^{л} - \bar{t}_e^{л}\right) + \beta \sum_{e \in \mathrm{A}^{л}} \tilde{f}_e^{л} \cdot \left(t_e^{л} - \bar{t}_e^{л}\right) + \beta \sum_{e \in E \setminus \mathrm{A}^{o}} \bar{f}_e^{o} \cdot \left(t_e^{o} - \bar{t}_e^{o}\right) + \beta \sum_{e \in \mathrm{A}^{o}} \tilde{f}_e^{o} \cdot \left(t_e^{o} - \bar{t}_e^{o}\right).$$

### 1.1.10. Стохастический вариант трехстадийной модели стабильной динамики

Рассуждая аналогично тому, как мы делали выше, можно обобщить результаты подразделе 1.1.9 на случай, когда вместо модели стабильной динамики используется её стохастический вариант (с параметром $T > 0$) [24, 34, 123] (развитие подхода этого пункта будет описано далее в разделе 1.4 этой главы):

$$\min_{\substack{t^{л} \ge \bar{t}^{л} \\ t^{o} \ge \bar{t}^{o} \\ \lambda^{L}, \lambda^{W}}} \left\{ \ln\left( \sum_{i,j=1}^{n} \exp\left( \beta T \psi_{ij}\left( \frac{t^{л}}{T}, \frac{t^{o}}{T} \right) - \lambda_i^{L} - \lambda_j^{W} \right) \right) + \right.$$

$$+\sum_{i=1}^{n}\lambda_i^L l_i+\sum_{j=1}^{n}\lambda_j^W \mathrm{w}_j+\beta\left\langle \bar{f}^{л}, t^{л}-\bar{t}^{л}\right\rangle+\beta\left\langle \bar{f}^{o}, t^{o}-\bar{t}^{o}\right\rangle\Bigg\}, \qquad (1.1.12)$$

где

$$\psi_{ij}\left(t^{л},t^{o}\right)=\ln\left(\sum_{p\in P_{(i,j)}^{Л}}\exp\left(-\sum_{e\in E}\delta_{ep}t_e^{л}\right)+\sum_{p\in P_{(i,j)}^{O}}\exp\left(-\sum_{e\in E}\delta_{ep}t_e^{o}\right)\right).$$

Причем здесь, так же как и в предыдущем пункте,

$$d_{ij}=\breve{Z}^{-1}\exp\left(\beta T\psi_{ij}\left(t^{л}/T,t^{o}/T\right)-\lambda_i^L-\lambda_j^W\right),$$

где $\breve{Z}^{-1}$ – ищется из условия нормировки, $f^*=\left(f^{л},f^{o}\right)$, $f^{л}=\bar{f}^{л}-s^{л}$, $f^{o}=\bar{f}^{o}-s^{o}$, а $s^{л}$, $s^{o}$ – векторы (оптимальные) двойственных множителей для ограничений $t^{л}\ge\bar{t}^{л}$, $t^{o}\ge\bar{t}^{o}$ в задаче (1.1.11). Таким образом, все, что осталось сделать, это решить задачу гладкой выпуклой оптимизации (1.1.12) прямодвойственным методом, например [28, 30, 31]. Мы не будем здесь подробно описывать возможные способы решения. Обратим только внимание на то, что задача вычисления значений и градиента функции $\psi\left(t^{л},t^{o}\right)$ – вычислительно сложная и, на первый взгляд, даже кажется бесперспективной из-за потенциально экспоненциально большого числа возможных маршрутов. Однако с помощью «аппарата характеристических функций на ориентированных графах» [23, 34, 123] и быстрого автоматического дифференцирования [44, 49] пересчитывать значения функции $\psi\left(t^{л},t^{o}\right)$ и ее градиента можно довольно эффективно (см. алгоритм из работ [23, 34, 123], который вырождается при $T\to 0+$ в алгоритм Беллмана–Форда [23, 34, 123]). Отметим, что при этом предельном переходе модель подраздела 1.1.10 перейдет в модель подраздела 1.1.9.

Отметим также, что всё, сказанное здесь, переносится и на функции $\psi_{ij}\left(t^{л},t^{o}\right)$ более общего вида, соответствующие различным иерархическим способам выбора (типа Nested Logit) [50, 68, 153]. Например, «практически бесплатно» можно сделать такую замену:

$$T\psi_{ij}\left(\frac{t^{л}}{T},\frac{t^{o}}{T}\right)\to\eta\ln\left(\exp\left(T\ln\left(\sum_{p\in P_{(i,j)}^{Л}}\exp\left(-\sum_{e\in E}\delta_{ep}\frac{t_e^{л}}{T}\right)\right)\Big/\eta\right)+\right.$$

$$\left.+\exp\left(T\ln\left(\sum_{p\in P_{(i,j)}^{O}}\exp\left(-\sum_{e\in E}\delta_{ep}\frac{t_e^{o}}{T}\right)\right)\Big/\eta\right)\right).$$

## 1.1.11. Калибровка модели стабильной динамики

О практическом использовании модели стабильной динамики написано, например, в [83] и немного в [133]. В этом пункте мы сконцентрируем внимание на потенциальных плюсах описанной модели с точки зрения её калибровки. Далее мы опускаем нижней индекс $e$.

Способ оценки $\bar{t}$ довольно очевиден: $\bar{t}$ определяется длиной участка (ребра) и типом ребра (в сельской местности, в городе, на шоссе). Вся эта информация обычно бывает доступной. С оценкой $\bar{f}$ немного посложнее. Пусть в конце дуги, состоящей из $r$ полос, стоит светофор, который пускает поток с этой дуги за долю времени $\chi$. Пусть $q_{\max} \approx 1800\,[\text{авт/ч}]$ – максимально возможное значение потока по одной полосе (это понятие не совсем корректное, но в первом приближении им можно пользоваться [34], и оно довольно универсально). Тогда [34] $\bar{f} \approx \chi r q_{\max}$. Тут имеется важный нюанс. Во время зеленой фазы (в зависимости от того, что это за фаза: скажем, движение прямо и налево или движение только направо) «работают» не все полосы: $\bar{f} \approx \chi_1 r_1 q_{\max} + \ldots + \chi_l r_l q_{\max}$, где $l$ – число фаз, «пускающих» поток с рассматриваемой дуги, $\chi_k$ – доля времени, отводимая фазе $k$, $r_k$ – эффективное число полос рассматриваемой дуги, задействованных на фазе $k$. Раздобыть информацию по работе светофоров, как правило, не удается. Однако существует определенные регламенты, согласно которым и устанавливаются фазы работы светофора. Поскольку все довольно типизировано, то часто бывает достаточно иметь информацию только о полосности дорог. Такая информация с 2012 года уже включается в коммерческие системы ГИС.

К сожалению, в реальных транспортных сетях в часы пик бывают пробки, которые приводят к тому, что пропускная способность ребра определяется пропускной способностью не данного ребра, а какого-то из впереди идущих (по ходу движения) ребер. Другими словами, пробка с ребра полностью «заполнила» это ребро и распространилась на ребра, входящие в это ребро. В таких ситуациях $\bar{f}$ нужно считать исходя из ограничений, накладываемых на пропускную способность впереди идущих ребер. При этом в выборе разбиения транспортного графа на ребра и в определении на этих ребрах значений $\bar{f}$, $\bar{t}$ стоит исходить из следующего: в типичной ситуации, если пробка и переходит с ребра на ребро, то желательно, чтобы это происходило без изменения пропускной способности входящего ребра. Этого не всегда можно добиться, поскольку часто приходится осуществлять разбиение без особого произвола, исходя из въездов / съездов, перекрестков. В таких случаях формула $\bar{f} \approx \chi_1 r_1 q_{\max} + \ldots + \chi_l r_l q_{\max}$ является лишь оценкой сверху для «среднего» (типичного) значения, которое надо подставлять в модель (часть этого бремени придется переносить и на $\bar{t}$, увеличивая его).

Заметим, что точно так же, как и для обычных многостадийных моделей, для калибровки предложенной модели необходимо использовать один и тот же промежуток времени (скажем, обеденные часы) каждого типового дня (например, буднего) – в зависимости от целей. Кроме того, для Москвы в утренние часы пик значения потоков существенно нестационарны, и в течение часа они, как правило, меняются сильно. Вместе с этим среднее время в пути оказывается больше часа. Таким образом, возникает вопрос: что понимается, например, под равновесным распределением потоков $f$? Ответ: «средние значения» нуждается в пояснении: средние значения не только по дням, но и по исследуемому промежутку времени (в несколько часов) внутри каждого дня. Как следствие, получаем, что выводы, сделанные по предложенной модели, нельзя, например, напрямую использовать для краткосрочного прогнозирования ситуации на дорогах или адаптивного управления светофорной сигнализацией. Таким образом, модель работает и выдает нереальные данные, хотя многие переменные модели и обозначают реальные физические параметры транспортного потока, а лишь некоторые средние (агрегированные) показатели, которые тем не менее многое могут сказать о ситуации на дорогах.

В данном разделе была предложена модель, описывающая равновесие в транспортной сети с точки зрения макромасштабов времени (месяцы). Для исследования поведения транспортного потока внутри одного дня требуются микромодели. Контекст, в котором такие модели используются, часто связан с необходимостью краткосрочного прогнозирования (на несколько часов вперед) и оптимального управления, например, светофорами. Как правило, в таких микромоделях численно исследуется начально-краевая задача для нелинейных УЧП (системы законов сохранения), в которой начальные условия имеются в наличии, а краевые условия (характеристики источников и стоков во времени, матрицы перемешивания в узлах графа транспортной сети) определяются исходя из исторической информации. Однако в последнее время в Москве все чаще можно слышать предложения о том, чтобы максимально информировать участников дорожного движения с целью лучшей маршрутизации. Такая полная информированность приводит к необходимости определять, скажем, матрицы перемешивания не из (и не только из) исторической информации, а исходя из равновесных принципов, заложенных в модели Бэкмана и стабильной динамики.

Причем нам представляется, что модель стабильной динамики подходит намного лучше для этих целей (мотивация имеется в работе [133], в которой на простом примере показывается, как возникает равновесная конфигурация модели стабильной динамики из внутредневной динамики водителей). Как отмечалось в работе [133], за исключением модельных примеров (рис. 1.1.2), очень сложной с вычислительной точки зрения представляется задача описания динамического режима функционирования модели, имеющей в своей основе модель стабильной динамики (или какую-то другую равновесную модель). Более того, пока, насколько нам известно, не было предложено ни одной более-менее обоснованной модели такого типа (иногда такие модели называют *моделями динамического равновесия* [140]). Здесь мы лишь укажем на одно, на наш взгляд, перспективное направление: объединить модель стабильной динамики с микромоделью CTM К. Даганзо [34] в современном ее варианте [43].

## 1.2. Эволюционные выводы энтропийной модели расчета матрицы корреспонденций

### 1.2.1. Введение

Настоящий раздел 1.2 развивает подраздел 1.1.3 раздела 1.1 этой главы. В этом разделе предлагается два способа объяснения популярной в урбанистике энтропийной модели расчета матрицы корреспонденций [18]. Обе предложенные модели в своей основе имеют марковский процесс, живущий в пространстве огромной размерности (говорят, что такой процесс порождает макросистему). Более точно, этот марковский процесс представляет собой ветвящийся процесс специального вида: *модель стохастической химической кинетики* [54, 95]. Более того, в обоих случаях имеет место условие детального равновесия [8, 22, 54]. Из общих результатов [8, 22, 54] отсюда сразу можно заключить, что инвариантная (стационарная) мера такого процесса будет иметь вид мультиномиального распределения, сосредоточенного на аффинном многообразии, определяемом линейными законами сохранения введенной динамики. Исходя из явления концентрации меры [80], можно ожидать концентрацию этой меры около наиболее вероятного состояния с ростом размера макросистемы. Вот такое состояние естественно принимать за равновесие изучаемой макросистемы, поскольку с большой вероятностью на больших временах мы найдем систему в малой окрестности такого состояния. Поиск равновесия сводится по теореме Санова [60] к задаче энтропийно-линейного программирования (1.2.3) [97]. Собственно, именно таким образом и планируется объяснить, почему для описания равновесной матрицы корреспонденций необходимо решить задачу ЭЛП.

В подразделе 1.2.2 на базе бинарных реакций обменного типа, популярных в различного рода физических и социально-экономических приложениях моделей стохастической химической кинетики [13, 19], будет приведен первый способ эволюционного вывода энтропийной модели расчета матрицы корреспонденций. Второй способ будет описан в подразделе 1.2.3. Он базируется на классических понятиях популяционной теории игр [145], например, таких как логит-динамика и игра загрузки (именно к такой игре сводится поиск равновесного распределения потоков по путям). Второй способ также базируется на редукции задачи расчета матрицы корреспонденций к задаче поиска равновесного распределения потоков по путям. Нетривиальным и, по-видимому, новым здесь является эволюционно-экономическая интерпретация двойственных множителей, обобщающая известную стационарную конструкцию, восходящую к Л. В. Канторовичу.

Отметим, что результаты данного раздела являются обобщением результатов работ [33, 34]. Обобщение заключается в том, что в формулировках основных результатов (теоремы 1.2.1 и 1.2.2) фигурируют точные (неулучшаемые) оценки времени сходимости изучаемых макросистем на равновесии и плотности концентрации инвариантной меры в окрестности равновесий.

Мы намеренно опускаем вопросы численного решения, возникающих задач ЭЛП. Подробнее об этом можно посмотреть, например, в [26, 97] и других частях пособия.

### 1.2.2. Вывод на основе обменов

Приведем, базируясь на работе [34], эволюционное обоснование одного из самых популярных способов расчета матрицы корреспонденций, имеющего более чем сорокалетнюю историю, – энтропийной модели [18].

Пусть в некотором городе имеется $n$ районов, $L_i > 0$ – число жителей $i$-го района, $W_j > 0$ – число работающих в $j$-м районе. При этом $N = \sum_{i=1}^{n} L_i = \sum_{j=1}^{n} W_j$, – общее число жителей города, $n^2 \ll N$. Далее под $L_i \geq 0$ будет пониматься число жителей района, выезжающих в типичный день за рассматриваемый промежуток времени из $i$-го района, $W_j > 0$ – число жителей города, приезжающих на работу в $j$-й район в типичный день за рассматриваемый промежуток времени. Обычно, так введенные, $L_i$, $W_j$ рассчитываются через число жителей $i$-го района и число работающих в $j$-м районе с помощью более-менее универсальных (в частности, не зависящих от $i$, $j$) коэффициентов пропорциональности. Эти величины являются входными параметрами модели, т. е. они не моделируются (во всяком случае, в рамках выбранного подхода). Для долгосрочных расчетов с разрабатываемой моделью требуется иметь прогноз изменения значений этих величин.

Обозначим через $d_{ij}(t) \geq 0$ – число жителей, живущих в $i$-м районе и работающих в $j$-м в момент времени $t$. Со временем жители могут только меняться квартирами, поэтому во все моменты времени $t \geq 0$:

$$d_{ij}(t) \geq 0,\ \sum_{j=1}^{n} d_{ij}(t) \equiv L_i,\ \sum_{i=1}^{n} d_{ij}(t) \equiv W_j,\ i, j = 1, ..., n.$$

Определим

$$\mathrm{A} = \left\{ d_{ij} \geq 0:\ \sum_{j=1}^{n} d_{ij} = L_i, \sum_{i=1}^{n} d_{ij} = W_j, i, j = 1, ..., n \right\}.$$

Опишем основной стимул к обмену: работать далеко от дома плохо из-за транспортных издержек. Будем считать, что эффективной функцией затрат [34] будет $R(T) = \beta T/2$, где $T > 0$ – время в пути от дома до работы (в общем случае под $T$ стоит понимать затраты, в которые может входить не только время), а $\beta > 0$ – настраиваемый параметр модели (который также можно проинтерпретировать и даже оценить, что и будет сделано ниже).

Теперь опишем саму динамику. Пусть в момент времени $t \geq 0$ $r$-й житель живет в $k$-м районе и работает в $m$-м, а $s$-й житель живет в $p$-м районе и работает в $q$-м. Тогда $\lambda_{k,m;\,p,q}(t)\Delta t + o(\Delta t)$ есть вероятность того, что жители с номерами $r$ и $s$ ($1 \leq r < s \leq N$) «поменяются» квартирами в промежутке времени $(t, t+\Delta t)$. Вероятность обмена местами жительства зависит только от мест проживания и работы обменивающихся:

$$\lambda_{k,m;\,p,q}(t) \equiv \lambda_{k,m;\,p,q} = \lambda N^{-1} \exp\Big( \underbrace{R(T_{km}) + R(T_{pq})}_{\text{суммарные затраты до обмена}} - \underbrace{\big(R(T_{pm}) + R(T_{kq})\big)}_{\text{суммарные затраты после обмена}} \Big) > 0,$$

где коэффициент $0 < \lambda = \mathrm{O}(1)$ характеризует интенсивность обменов. Совершенно аналогичным образом можно было рассматривать случай «обмена местами работы». Здесь стоит оговориться, что «обмены» не стоит понимать буквально – это лишь одна из возможных интерпретаций. Фактически используется так называемое *приближение среднего поля* [22, 95], т. е. некое равноправие агентов (жителей) внутри фиксированной корреспонденции и их независимость.

Согласно эргодической теореме для марковских цепей (вне зависимости от начальной конфигурации $\left\{d_{ij}(0)\right\}_{i=1,\ j=1}^{n,\ n}$) [10, 13, 19, 34, 54, 115, 145] предельное распределение совпадает со стационарным (инвариантным), которое можно посчитать (получается проекция прямого произведения распределений Пуассона на А):

$$\lim_{t\to\infty} P\left(d_{ij}(t) = d_{ij},\, i, j = 1,\dots,n\right) = \\ = Z^{-1}\prod_{i,j=1}^{n}\exp\left(-2R\left(T_{ij}\right)d_{ij}\right)\cdot\left(d_{ij}!\right)^{-1} \overset{def}{=} p\left(\left\{d_{ij}\right\}_{i=1,\ j=1}^{n,\ n}\right), \tag{1.2.1}$$

где $\left\{d_{ij}\right\}_{i=1,\ j=1}^{n,\ n} \in \mathrm{A}$, а «статсумма» $Z$ находится из условия нормировки получившейся «пуассоновской» вероятностной меры. Отметим, что стационарное распределение $p\left(\left\{d_{ij}\right\}_{i=1,\ j=1}^{n,\ n}\right)$ удовлетворяет условию детального равновесия [22, 145]:

$$\left(d_{km}+1\right)\left(d_{pq}+1\right)p\left(\left\{d_{11},\dots,d_{km}+1,\dots,d_{pq}+1,\dots,d_{pm}-1,\dots,d_{kq}-1,\dots,d_{nn}\right\}\right)\lambda_{k,m;\,p,q} = \\ = d_{pm}d_{kq}\,p\left(\left\{d_{ij}\right\}_{i=1,\ j=1}^{n,\ n}\right)\lambda_{p,m;\,k,q}.$$

При $N \gg 1$ распределение $p\left(\left\{d_{ij}\right\}_{i=1,\ j=1}^{n,\ n}\right)$ экспоненциально сконцентрировано на множестве А в $\mathrm{O}\left(\sqrt{N}\right)$-окрестности наиболее вероятного значения $d^* = \left\{d_{ij}^*\right\}_{i=1,\ j=1}^{n,\ n}$, которое определяется, как решение задачи энтропийно-линейного программирования [33, 34, 54]:

$$\ln p\left(\left\{d_{ij}\right\}_{i=1,\ j=1}^{n,\ n}\right) \sim -\sum_{i,j=1}^{n} d_{ij}\ln\left(d_{ij}\right) - \beta\sum_{i,j=1}^{n} d_{ij}T_{ij} \to \max_{\left\{d_{ij}\right\}_{i=1,\ j=1}^{n,\ n}\in(\mathrm{A})}. \tag{1.2.2}$$

Это следует из теоремы Санова о больших уклонениях для мультиномиального распределения [60] (на распределение (1.2.1) можно смотреть так же, как и на проекцию мультиномиального распределения на А):

$$\frac{N!}{d_{11}!\cdot\ldots\cdot d_{ij}!\cdot\ldots\cdot d_{nn}!}\,p_{11}^{d_{11}}\cdot\ldots\cdot p_{ij}^{d_{ij}}\cdot\ldots\cdot p_{nn}^{d_{nn}} = \exp\left(-N\sum_{i,j=1}^{n}\nu_{ij}\ln\left(\nu_{ij}/p_{ij}\right)+\bar{R}\right),$$

где $\nu_{ij} = d_{ij}/N$, $|\bar{R}| \le \frac{n^2}{2}(\ln N + 1)$, и формулы Тейлора (с остаточным членом второго порядка в форме Лагранжа), примененной к энтропийному функционалу в точке $d^*$, заданному на А [18].

Сформулируем более точно полученный результат (см. также приложение 1).

**Теорема 1.2.1.** *Для любого* $d(0) \in \mathrm{A}$ *существует такая константа* $c_n(d(0)) > 0$, *что для всех* $\sigma \in (0, 0.5)$, $t \ge c_n(d(0)) \ln N$, *имеет место неравенство*

$$P\left( \frac{\|d(t) - d^*\|_2}{N} \ge \frac{2\sqrt{2} + 4\sqrt{\ln(\sigma^{-1})}}{\sqrt{N}} \right) \le \sigma,$$

*где* $d(t) = \{d_{ij}(t)\}_{i=1,\, j=1}^{n,\, n} \in \mathrm{A}$.

**Схема доказательства.** Установим сначала оценку для плотности концентрации стационарной меры:

$$\lim_{t\to\infty} P\left(d_{ij}(t) = d_{ij},\, i, j = 1, \dots, n\right) = \frac{N!}{d_{11}! \cdot \ldots \cdot d_{ij}! \cdot \ldots \cdot d_{nn}!} p_{11}^{d_{11}} \cdot \ldots \cdot p_{ij}^{d_{ij}} \cdot \ldots \cdot p_{nn}^{d_{nn}}.$$

Из неравенства Хефдинга в гильбертовом пространстве [80] следует ($\varepsilon \ge \sqrt{2/N}$):

$$\lim_{t\to\infty} P\left(\|d(t) - d^*\|_2 \ge \varepsilon N\right) \le \exp\left(-\frac{1}{4N}\left(\varepsilon N - \sqrt{2N}\right)^2\right).$$

Беря в этом неравенстве

$$\varepsilon = \frac{\sqrt{2} + 2\sqrt{\ln(\sigma^{-1})}}{\sqrt{N}},$$

получим

$$\lim_{t\to\infty} P\left( \frac{\|d(t) - d^*\|_2}{N} \ge \varepsilon \right) \le \sigma.$$

Однако если не переходить к пределу по времени, а лишь обеспечить при $t \ge T_n(\varepsilon, N; d(0))$ выполнение условия

$$E\left\| \frac{d(t) - d^*}{N} \right\|_2 \le \frac{\varepsilon}{2},$$

то при $t \geq T_n\left(\varepsilon, N; d(0)\right)$:

$$P\left(\left\|d(t)-d^*\right\|_2 \geq \varepsilon N\right) \leq \exp\left(-\frac{1}{4N}\left(\varepsilon N-\left(\sqrt{2N}+\frac{\varepsilon N}{2}\right)\right)^2\right)=\exp\left(-\frac{1}{4N}\left(\frac{\varepsilon N}{2}-\sqrt{2N}\right)^2\right).$$

Беря в этом неравенстве

$$\varepsilon = \frac{2\sqrt{2}+4\sqrt{\ln\left(\sigma^{-1}\right)}}{\sqrt{N}},$$

получим при $t \geq T_n\left(\varepsilon, N; d(0)\right)$ аналогичное неравенство

$$P\left(\frac{\left\|d(t)-d^*\right\|_2}{N} \geq \varepsilon\right) \leq \sigma.$$

Отметим, что с точностью до мультипликативной константы эта оценка не может быть улучшена. Это следует из неравенства Чебышёва:

$$P\left(X \geq EX-\varepsilon\right) \geq 1-\frac{\operatorname{Var}(X)}{\varepsilon^2}$$

при

$$X=\frac{\left\|d(t)-d^*\right\|_2}{N},\ d^*=N\cdot\left(n^{-2},\dots,n^{-2}\right)^T,\ N \gg n^2.$$

Осталось оценить зависимость $T_n\left(\varepsilon, N; d(0)\right)$. Для этого воспользуемся неравенством Чигера [115]. Поставим в соответствие нашему марковскому процессу его дискретный аналог (в дискретном времени с шагом $\lambda N^{-1}$) – случайное блуждание (со скачками, соответствующими парным реакциям, введенным выше) на целочисленных точках части гиперплоскости, задаваемой множеством A. Граф, на котором происходит блуждание, будем обозначать $G=\left\langle V_G, E_G\right\rangle$. Пусть $\pi(\cdot)$ – стационарная мера этого блуждания (в нашем случае – мультиномиальная мера, экспоненциально сконцентрированая на множестве A в $O\left(\sqrt{N}\right)$-окрестности наиболее вероятного значения $d^*$), а $P=\left\|p_{ij}\right\|_{i,j=1}^{|V_G|,|V_G|}$ – матрица переходных вероятностей. Тогда, вводя константу Чигера ( $\bar{S}=V_G \setminus S$ ):

$$h(G)=\min_{S\subseteq V_G:\,\pi(S)\leq 1/2} P\left(S\to\bar{S}\,\middle|\,S\right)=\min_{S\subseteq V_G:\,\pi(S)\leq 1/2}\frac{\sum_{(i,j)\in E_G:\,i\in S,\ j\in\bar{S}}\pi(i)\,p_{ij}}{\sum_{i\in S}\pi(i)}$$

и время выхода блуждания (стартовавшего из $i \in V_G$) на стационарное распределение:

$$T(i,\varepsilon) = \mathrm{O}\left(\lambda N^{-1} \cdot h(G)^{-2}\left(\ln\left(\pi(i)^{-1}\right) + \ln\left(\varepsilon^{-1}\right)\right)\right),$$

получим для любых $i = 1,\ldots,|V_G|$, $t \ge T(i,\varepsilon)$:

$$\left\|P^t(i,\cdot) - \pi(\cdot)\right\|_2 \le \left\|P^t(i,\cdot) - \pi(\cdot)\right\|_1 = \sum_{j=1}^{n}\left|P^t(i,j) - \pi(j)\right| \le \frac{\varepsilon}{2},$$

где $P^t(i,j)$ – условная вероятность того, что, стартуя из состояния $i \in V_G$, через $t$ шагов марковский процесс окажется в состоянии $j \in V_G$.

В нашем случае можно явно геометрически описать множество $S$ в виде целочисленных точек множества А, попавших внутрь эллипсоида (сферы) размером $\mathrm{O}\left(\sqrt{N}\right)$ с $\pi(S) \simeq 1/2$ и с центром в точке $d^*$. На этом множестве достигается решение изопериметрической задачи в определении константы Чигера. Значение константы Чигера при этом будет пропорционально отношению площади поверхности этого эллипсоида к его объему, т. е. $h(G) \sim N^{-1/2}$. Используя это наблюдение, можно получить, что

$$T_n(\varepsilon, N; d(0)) \sim \ln\left(\pi(i)^{-1}\right) + \ln N.$$

Отсюда видно, что время выхода зависит от точки старта. Если ограничиться точками старта $i \in V_G$, для которых равномерно по $N \to \infty$ и $i, j = 1,\ldots,n$ имеет место неравенство $d_{ij}(0)/N \ge \varsigma_n > 0$, то

$$\ln\left(\pi(i)^{-1}\right) \sim \ln N.$$

Объединяя приведенные результаты, получаем утверждение теоремы 1.2.1. ■

Эта теорема уточняет результат работы [34], явно указывая скорость сходимости и плотность концентрации. Приведенная скорость сходимости характерна для более широкого класса моделей стохастической химической кинетики, приводящих к равновесию вида неподвижной точки [22]. Плотность концентрации оценивалась на базе конструкции работы [32].

Естественно, ввиду теоремы 1.2.1, принимать решение задачи (1.2.2) $\left\{d_{ij}^*\right\}_{i=1,\, j=1}^{n,\, n}$ за равновесную конфигурацию. Обратим внимание, что предложенный выше вывод известной энтропийной модели расчета матрицы корреспонденций отличается от классического [18]. В монографии

А. Дж. Вильсона [18] $\beta$ интерпретируется как множитель Лагранжа к ограничению на среднее «время в пути»: $\sum_{i,j=1}^{n} d_{ij}T_{ij} = C$.

**Замечание 1.2.1.** При этом остальные ограничения имеют такой же вид, а функционал имеет вид $F(d) = -\sum_{i,j=1}^{n} d_{ij}\ln(d_{ij})$. Тогда, согласно экономической интерпретации двойственных множителей Л. В. Канторовича: $\beta(C) = \partial F(d(C))/\partial C$. Из такой интерпретации иногда делают вывод о том, что $\beta$ можно понимать как цену единицы времени в пути. Чем больше $C$, тем меньше $\beta$.

Приведенный нами вывод позволяет контролировать знак параметра $\beta > 0$ и лучше понимать его физический смысл.

**Замечание 1.2.2.** Отметим, что, так же как и в [18], из принципа ле Шателье–Самуэльсона следует, что с ростом $\beta$ среднее время в пути $\sum_{i,j=1}^{n} d_{ij}(\beta)T_{ij}$ будет убывать. В связи с этим обстоятельством, а также исходя из соображений размерности, вполне естественно понимать под $\beta$ величину, обратную к характерному (среднему) времени в пути [18] – физическая интерпретация. Собственно, такая интерпретация параметра $\beta$, как правило, и используется в многостадийных моделях (см., например, [140]).

В дальнейшем, нам будет удобно привести задачу (1.2.2) к следующему виду (при помощи метода множителей Лагранжа [53], теоремы фон Неймана о минимаксе [149] и перенормировки $d := d/N$):

$$\max_{\lambda^L,\lambda^W} \min_{\sum_{i,j=1}^{n} d_{ij}=1,\, d_{ij}\ge 0} \Bigg[ \sum_{i,j=1}^{n} d_{ij}\ln d_{ij} + \beta \sum_{i=1,\, j=1}^{n,\, n} d_{ij}T_{ij} + \\ + \sum_{i=1}^{n} \lambda_i^L \left( \sum_{j=1}^{n} d_{ij} - l_i \right) + \sum_{i=1}^{n} \lambda_j^W \left( \mathrm{w}_j - \sum_{i=1}^{n} d_{ij} \right) \Bigg], \tag{1.2.3}$$

где $l = L/N$, $\mathrm{w} = W/N$.

**Замечание 1.2.3.** Используя принцип Ферма [53] несложно проверить, что решение задачи (1.2.3) можно представить в виде

$$d_{ij} = \exp(-\lambda_i^L)\exp(\lambda_j^w)\exp(-\beta T_{ij}),$$

где

$$\exp\left(\lambda_i^L\right) = \frac{1}{l_i}\sum_{j=1}^{n}\exp\left(\lambda_j^w\right)\exp\left(-\beta T_{ij}\right),$$
$$\exp\left(-\lambda_j^w\right) = \frac{1}{w_j}\sum_{i=1}^{n}\exp\left(-\lambda_i^L\right)\exp\left(-\beta T_{ij}\right). \quad (1.2.4)$$

Отсюда можно усмотреть интерпретацию двойственных множителей как соответствующих «потенциалов притяжения / отталкивания районов» [18]. К этому мы еще вернемся в подразделе 1.2.3.

**Замечание 1.2.4.** Департамент транспорта г. Москвы несколько лет назад поставил следующую задачу. Сколько человек надо обзвонить и опросить на предмет того, какой корреспонденции они принадлежат (где живут и где работают), чтобы восстановить матрицу корреспонденций с достаточной точностью и на доверительном уровне? Формализуем задачу. Предположим, что $\left\{d_{ij}^*\right\}_{i=1,\ j=1}^{n,\ n}$ – истинная (пронормированная) матрица корреспонденций, вот её и надо восстановить. В результате опросов населения получилась матрица (вектор) $r = \left\{r_{ij}\right\}_{i=1,\ j=1}^{n,\ n}$, где $r_{ij}$ – количество опрошенных респондентов, принадлежащих корреспонденции $(i, j)$, $\sum_{i,j=1}^{n} r_{ij} = N$. Задачу формализуем следующим образом: найти наименьшее $N$, чтобы

$$P_{d^*}\left(\left\|\frac{r_{ij}}{N} - d^*\right\|_2 \geq \varepsilon\right) \leq \sigma.$$

Нижний индекс $d^*$ у вероятности означает, что случайный вектор $r$ имеет мультиномиальное распределение с параметром $d^*$, т. е. считаем, что опрос проводился в идеальных условиях. Из теоремы 1.2.1 несложно усмотреть, что достаточно опросить $N = \left(4 + 8\ln\left(\sigma^{-1}\right)\right)\varepsilon^{-2}$ респондентов. Скажем, опрос ~50 000 респондентов, который и был произведен, действительно позволяет неплохо восстановить матрицу корреспонденций $d^* \approx \bar{d} \overset{def}{=} r_{ij}/N$. Однако мы привели здесь это замечание для других целей. В ряде работ (см., например, [146, 155]), исходя из данных таких опросов, также восстанавливают матрицу корреспонденций, но при другой параметрической гипотезе (меньшее число параметров):

$d_{ij} = \exp\left(-\lambda_i^L\right)\exp\left(\lambda_j^w\right)\exp\left(-\beta T_{ij}\right)$. То есть дополнительно предполагают, что $n^2$ неизвестных параметров в действительности однозначно определяются $2n$ параметрами (иногда к ним добавляют еще один параметр $\beta$), которые у нас ранее (замечание 1.2.3) интерпретировались как множители Лагранжа. Встает вопрос: как оптимально оценить эти параметры? Ответ дает теорема Фишера (современный неасимптотический вариант изложения этой теоремы можно найти в [150]) об оптимальности оценок максимального правдоподобия (ОМП). Собственно, для выборки из мультиномиального распределения оценкой максимального правдоподобия как раз и будет выборочное среднее $\bar{d}$. Для описанной параметрической модели (с $2n$ параметрами) поиск такой оценки сводится к разрешению системы (1.2.4) (замечание 1.2.3).

**Замечание 1.2.5.** Подобно тому как мы рассматривали в этом пункте трудовые корреспонденции (в утренние и вечерние часы более 70% корреспонденций по Москве и области именно такие), можно рассматривать перемещения, например, из дома к местам учебы, отдыха, в магазины и т. п. (по-хорошему еще надо было учитывать перемещения типа работа–магазин–детский сад–дом) – рассмотрение всего этого вкупе приведет также к задаче (1.2.2). Только будет больше типов корреспонденций $d$: помимо пары районов, еще нужно будет учитывать тип корреспонденции [18]. Все это следует из того, что инвариантной мерой динамики с несколькими типами корреспонденций по-прежнему будет прямое произведение пуассоновских мер. Другое дело, когда мы рассматриваем разного типа пользователей транспортной сети, например: имеющих личный автомобиль и не имеющих личный автомобиль. Первые могут им воспользоваться, равно как и общественным транспортом, а вторые нет. И на рассматриваемых масштабах времени пользователи могут менять свой тип. То есть время в пути может для разных типов пользователей быть различным [18]. Считая, подобно тому как мы делали раньше, что желание пользователей корреспонденции $(i, j)$ сети сменить свой тип (вероятность в единицу времени) есть

$$\tilde{\lambda}\exp\Big( \underbrace{\tilde{R}\left(T_{ij}^{old}\right)}_{\text{суммарные затраты до смены типа}} - \underbrace{\tilde{R}\left(T_{ij}^{new}\right)}_{\text{суммарные затраты после смены типа}} \Big), \text{ где } \tilde{R}(T) = \beta T,$$

и учитывая в «обменах» тип пользователя (будет больше типов корреспонденций $d$, но «меняются местами работы» только пользователи одного типа), можно показать, что все это вкупе приведет также к задаче типа (1.2.2).

### 1.2.3. Вывод на основе модели равновесного распределения потоков

Предварительно напомним, следуя [33, 145], эволюционный вывод модели равновесного распределения потоков [33, 34, 141, 148].

Задан ориентированный граф $\Gamma = (V, E)$ – транспортная сеть города ($V$ – узлы сети (вершины), $E \subset V \times V$ – дуги сети (ребра графа)). В графе имеются две выделенные вершины. Одна из вершин графа является источником, другая – стоком. Из источника в сток ведет много путей, которые мы будем обозначать $p \in \mathrm{P}$, $|\mathrm{P}| = m$. Обозначим через

$x_p$ – величина потока по пути $p$, $x_p \in S_m(N) = \left\{ x = \{x\}_{p\in \mathrm{P}} \geq 0 : \sum_{p\in\mathrm{P}} x_p = N \right\}$;

$y_e$ – величина потока по ребру $e \in E$: $y_e = \sum_{p\in\mathrm{P}} \delta_{ep} x_p$ ($y = \Theta x$, $\Theta = \left\{\delta_{ep}\right\}_{e\in E, p\in\mathrm{P}}$), где

$$\delta_{ep} = \begin{cases} 1, & e \in p, \\ 0, & e \notin p; \end{cases}$$

$\tau_e(y_e)$ – удельные затраты на проезд по ребру $e$ (гладкие неубывающие функции);

$G_p(x) = \sum_{e\in E} \tau_e(y_e)\delta_{ep}$ – удельные затраты на проезд по пути $p$.

Рассмотрим, следуя [145], популяционную игру, в которой имеется набор $N$ однотипных игроков (агентов). Множеством чистых стратегий каждого такого агента является $\mathrm{P}$, а выигрыш (потери со знаком минус) от использования стратегии $p \in \mathrm{P}$ определяется формулой $-G_p(x)$.

Опишем динамику поведения игроков. Пусть в момент времени $t \geq 0$ агент использует стратегию $p \in \mathrm{P}$, $\lambda_{p,q}(t)\Delta t + o(\Delta t)$ – вероятность того, что он поменяет свою стратегию на $q \in \mathrm{P}$ в промежутке времени $(t, t+\Delta t)$. Будем считать, что

$$\lambda_{p,q}(t) \equiv \lambda_{p,q} = \lambda P_q\left(\left\{G_p\left(x(t)\right)\right\}_{p\in\mathrm{P}}\right),$$

где (как и в подразделе 1.2.2) коэффициент $0 < \lambda = \mathrm{O}(1)$ характеризует интенсивность «перескоков» агентов, а

$$P_q\left(\left\{G_p\left(x(t)\right)\right\}_{p\in\mathrm{P}}\right) = P\left(q = \arg\max_{p\in\mathrm{P}}\left\{-G_p\left(x(t)\right) + \xi_p\right\}\right).$$

Если $\xi_p \equiv 0$, то получаем динамику наилучших ответов [145], если $\xi_p$ – независимые одинаково распределенные случайные величины с распределением Гумбеля [70, 157]: $P\left(\xi_p < \xi\right) = \exp\left\{-e^{-\xi/\omega - E}\right\}$, где $\omega \in \left(0, \omega_0\right]$ ($\omega_0 = O(1)$), $E \approx 0.5772$ – константа Эйлера, а $\mathrm{Var}\left[\xi_p\right] = \omega^2\pi^2/6$, то получаем логит-динамику [145, 157]:

$$P_q\left(\left\{G_p\left(x(t)\right)\right\}_{p\in P}\right) = \frac{\exp\left(-G_q\left(x(t)\right)/\omega\right)}{\sum\limits_{p\in P}\exp\left(-G_p\left(x(t)\right)/\omega\right)},$$

вырождающуюся в динамику наилучших ответов при $\omega \to 0+$. Далее мы будем считать, что задана логит-динамика. Объясняется такая динамика совершенно естественно. Каждый агент имеет какую-то картину текущего состояния загрузки системы, но либо он не имеет возможности наблюдать ее точно, либо он старается как-то спрогнозировать возможные изменения загрузок (а соответственно и затраты на путях) в будущем (либо и то и другое).

Приведем соответствующий аналог теоремы 1.2.1.

**Теорема 1.2.2.** *Для любого* $x(0) \in S_m(N)$ *существует такая константа* $c_m\left(x(0)\right) > 0$, *что для всех* $\sigma \in (0, 0.5)$, $t \geq c_m\left(x(0)\right) N \ln N$, *имеет место неравенство*

$$P\left(\left\|\frac{x(t)}{N} - x^*\right\|_2 \geq \frac{2\sqrt{2} + 4\sqrt{\ln\left(\sigma^{-1}\right)}}{\sqrt{N}}\right) \leq \sigma,\ x(t) \in S_m(N),$$

*где* (*стохастическое равновесие Нэша–Вардропа*):

$$x^* = \arg\min_{x\in S_m(1)} \Psi\left(y(x)\right) + \omega \sum_{p\in P} x_p \ln\left(x_p\right), \qquad (1.2.5)$$

$$\Psi\left(y(x)\right) = \sum_{e\in E} \int\limits_0^{y_e(x)} \tilde{\tau}_e(z)\,dz,\ \tilde{\tau}_e(z) = \tau_e(zN).$$

**Схема доказательства.** Стационарная мера описанного марковского процесса имеет (с точностью до нормирующего множителя) вид (теорема 11.5.12 [145]):

$$\frac{N!}{x_1! \cdot \ldots \cdot x_m!}\exp\left(-\frac{\Psi\left(y(x)\right)}{\omega}\right),\ x \in S_m(N).$$

Для оценок плотности концентрации достаточно заметить, что $\Psi\big(y(x)\big)$ – выпуклая функция, как композиция линейной и выпуклой функции. Поэтому оценки плотности концентрации здесь не могут быть хуже оценок в теореме 1.2.1. Рассуждения для оценки времени выхода проводятся аналогично теореме 1.2.1. ■

**Следствие (Proposition 1 [87]).** *Если* $\omega \to 0+$, *то решение задачи* (1.2.5) *сходится к*

$$x^* = \arg\min_{x\in S_m(1):\ \Theta x = y^*} \sum_{p\in P} x_p \ln\big(x_p\big), \qquad (1.2.6)$$

*где* $y^* = \arg\min_{y=\Theta x,\ x\in S_m(1)} \Psi(y)$.

Это следствие решает проблему обоснования гипотезы Бар-Гира [73]. Напомним вкратце, в чем состоит эта гипотеза. Известно (см., например, [141]), что поиск равновесного распределения потоков по путям (равновесия Нэша–Вардропа) в модели с одним источником и стоком сводится к задаче оптимизации (1.2.5) с $\omega = 0$. Хотя функционал этой задачи выпуклый, но он не всегда строго выпуклый, даже в случае, когда все функции $\tau_e\big(y_e\big)$ – строго возрастающие (тогда можно лишь говорить о единственности равновесного распределения потоков по ребрам $y^*$). Гипотеза Бар-Гира говорит, что «в жизни» с большой вероятностью реализуется то равновесие из множества равновесий, которое находится из решения задачи (1.2.6).

**Замечание 1.2.6.** Распределение потоков по путям $x$ называется *равновесием* (*Нэша–Вардропа*) в рассматриваемой популяционной игре $\left\langle \{x_p\}_{p\in P}, \{G_p(x)\}_{p\in P} \right\rangle$, если

$$\text{из } x_p > 0,\ p\in P \text{ следует } G_p(x) = \min_{q\in P} G_q(x).$$

Или, что то же самое:

$$\text{для любых } p\in P \text{ выполняется } x_p\cdot\Big(G_p(x) - \min_{q\in P} G_q(x)\Big) = 0.$$

**Замечание 1.2.7.** Если при $\omega = 0$ рассмотреть предельный случай

$$\tau_e\big(y_e\big) := \tau_e^{\mu}\big(y_e\big) \xrightarrow[\mu\to 0+]{} T_e,$$

то поиск равновесия Нэша–Вардропа просто сводится к поиску социального оптимума, что приводит в данном случае к решению транспортной задачи линейного программирования [68]. Если делать предельный переход с учетом ограничений на пропускные способности ребер

$$\tau_e\big(y_e\big) := \tau_e^{\mu}\big(f_e\big) \xrightarrow[\mu\to 0+]{} \begin{cases} T_e, & 0 \le y_e < \bar{y}_e, \\ \big[T_e, \infty\big), & y_e = \bar{y}_e, \end{cases}$$

то получится более сложная задача, которая описывает равновесие в (стохастической, если $\omega > 0$) модели стабильной динамики [33, 133].

Из рассмотренных в замечании 1.2.7 случаев выпала ситуация $\omega > 0$: $\tau_e(y_e) := \tau_e^{\mu}(y_e) \xrightarrow[\mu \to 0+]{} T_e$. Её мы сейчас отдельно и рассмотрим на примере другого эволюционного способа обоснования энтропийной модели расчета матрицы корреспонденций. По ходу обсуждения этого способа с коллегами (прежде всего с Ю. Е. Нестеровым и Ю. В. Дорном) у предложенного подхода появилось название: *облачная модель*.

Предположим, что все вершины, отвечающие источникам корреспонденций, соединены ребрами с одной вспомогательной вершиной (облако № 1). Аналогично все вершины, отвечающие стокам корреспонденций, соединены ребрами с другой вспомогательной вершиной (облако № 2). Припишем всем новым ребрам постоянные веса. И проинтерпретируем веса ребер, отвечающих источникам $\lambda_i^L$, например, как средние затраты на проживание (в единицу времени, скажем, в день) в этом источнике (районе), а веса ребер, отвечающих стокам $\lambda_j^W$, как уровень средней заработной платы (в единицу времени) в этом стоке (районе), если изучаем трудовые корреспонденции. Будем следить за системой в медленном времени, то есть будем считать, что равновесное распределение потоков по путям стационарно. Поскольку речь идет о равновесном распределении потоков, то нет необходимости говорить о затратах на путях или ребрах детализированного транспортного графа, достаточно говорить только затратах (в единицу времени), отвечающих той или иной корреспонденции. Таким образом, у нас есть взвешенный транспортный граф с одним источником (облако 1) и одним стоком (облако 2). Все вершины этого графа, кроме двух вспомогательных (облаков), соответствуют районам в модели расчета матрицы корреспонденций. Все ребра этого графа имеют стационарные (не меняющиеся и не зависящие от текущих корреспонденций) веса $\{T_{ij}; \lambda_i^L; -\lambda_j^W\}$. Если рассмотреть естественную в данном контексте логит-динамику ($x \equiv d$), описанную выше, с $\omega = 1/\beta$ (здесь полезно напомнить, что согласно замечанию 1.2.2 $\beta$ обратно пропорционально средним затратам, а $\omega$ имеет как раз физическую размерность затрат), то поиск равновесия рассматриваемой макросистемы согласно теореме 1.2.1 приводит (в прошкалированных переменных) к задаче, сильно похожей на задачу (1.2.3) из подраздела 1.2.2:

$$\min_{\sum_{i,j=1}^{n} d_{ij} = 1,\ d_{ij} \geq 0} \left[ \sum_{i,j=1}^{n} d_{ij} \ln d_{ij} + \beta \sum_{i=1,\ j=1}^{n,\ n} d_{ij} T_{ij} + \beta \sum_{i=1}^{n} \left( \lambda_i^L \sum_{j=1}^{n} d_{ij} \right) - \beta \sum_{i=1}^{n} \left( \lambda_j^W \sum_{i=1}^{n} d_{ij} \right) \right].$$

Разница состоит в том, что здесь мы не оптимизируем по $2n$ двойственным множителям $\lambda_i^L$, $\lambda_j^W$ (множителям Лагранжа). Более того, мы их и не интерпретируем здесь как двойственные множители, поскольку мы их ввели на этапе взвешивания ребер графа. Тем не менее значения этих переменных, как правило, неоткуда брать. Тем более что приведенная выше (наивная) интерпретация вряд ли может всерьез рассматриваться как способ определения этих параметров исходя из данных статистики. Более правильно понимать $\lambda_i^L$, $\lambda_j^W$ как потенциалы притяжения / отталкивания районов (см. также замечание 1.2.3), включающие в себя плату за жилье и зарплату, но включающие также и многое другое, что сложно описать количественно. И здесь как раз помогает информация об источниках и стоках, содержащаяся в $2n$ уравнениях, задающих множество A. Таким образом, мы приходим как раз к той же самой задаче (1.2.3) с той лишь разницей, что мы получили дополнительную интерпретацию двойственных множителей в задаче (1.2.3). При этом двойственные множители в задаче (1.2.3) равны (с точностью до мультипликативного фактора $\beta$) введенным здесь потенциалам притяжения районов.

Нам представляется очень плодотворной и перспективной идея перенесения имеющейся информации об исследуемой системе из обременительных законов сохранения динамики, описывающей эволюцию этой системы, в саму динамику путем введения дополнительных естественно интерпретируемых параметров. При таком подходе становится возможным, например, учитывать в моделях и рост транспортной сети. Другими словами, при таком подходе, например, можно естественным образом рассматривать также и ситуацию, когда число пользователей транспортной сетью меняется со временем (медленно).

## 1.3. Об эффективной вычислимости конкурентных равновесий в транспортно-экономических моделях

### 1.3.1. Введение

В данном разделе мы развиваем конструкцию подраздела 1.1.9 раздела 1.1 этой главы. А именно в этом разделе планируется сосредоточить внимание на транспортно-экономических моделях, объединяющих в себе, в частности, модели из недавних работ [15, 33]. Этот раздел мотивирован обоснованием существующих и созданием новых моделей транспортного планирования, включающих модели роста транспортной инфраструктуры городов, формирования матрицы корреспонденций и равновесного распределения потоков.

Имеется ориентированный транспортный граф, каждое ребро которого характеризуется неубывающей функцией затрат $\tau_e\left(f_e\right)$ на прохождение этого ребра в зависимости от потока по этому ребру. Можно еще ввести затраты на прохождения вершин графа $E$, но это ничего не добавляет с точки зрения последующих математических выкладок [15]. Часть вершин графа является источниками, часть – стоками (эти множества вершин могут пересекаться). В источниках $O$ и в стоках $D$ имеются (соответственно) пункты производства и пункты потребления. Для большей наглядности в первой половине этого раздела мы будем считать, что производится и потребляется лишь один продукт. Несложно все, что далее будет написано, обобщить на многопродуктовый рынок.

Задача разбивается на две подзадачи разного уровня [89]. На нижнем уровне, соответствующем быстрому времени, при заданных корреспонденциях $\{d_{ij}\}$ (сколько товара перевозится в единицу времени из источника $i$ в сток $j$) идет равновесное формирование способов транспортировки товаров [33]. В результате формируются функции затрат $T_{ij}\left(\{d_{ij}\}\right)$. Исходя из этих затрат на верхнем уровне, соответствующем медленному времени, решается задача поиска конкурентного равновесия [4, 104] между производителями и потребителями с учетом затрат на транспортировку. В данном случае (см., например, [33]) мы будем иметь дело с адиабатическим исключением быстрых переменных (принцип подчинения Г. Хакена) в случае стохастических динамик. Обоснование имеется в [19].

Различные частные случаи такого рода постановки задач встречались в литературе. Так, например, в классической монографии [18] рассматривается большое количество моделей верхнего уровня, связанных с расчетом матрицы корреспонденций. В монографиях [34, 61, 141], напротив, внимание сосредоточено на моделях нижнего уровня, в которых с помощью принципа Нэша–Вардропа рассчитывается $T_{ij}\left(\{d_{ij}\}\right)$. В статье [33] эти модели объединяются для создания единой многоуровневой (в транспортной науке чаще используется термин «многостадийной») модели, включающей в себя и расчет матрицы корреспонденций, и равновесное распределение потоков по путям. В препринте [15] введена терминология, которой мы будем придерживаться и в данном разделе, связанная с пунктами производства и потребления, и в отличие от [18, 33] внешняя задача в [15] больше привязана непосредственно к экономике. Но во всех этих случаях можно было обойтись (с некоторыми оговорками в слу-чае [33]) и без понятия конкурентного равновесия, поскольку получающиеся в итоге (популяционные) игры были потенциальными[16], причем имелась и эволюционная интерпретация [145]. Поиск равновесия сводился к решению задачи выпуклой оптимизации, а цены определялись из решения двойственной задачи. В препринте [135] для задачи верхнего уровня была предложена оригинальная конструкция, сводящая поиск конкурентного равновесия к поиску седловой точки (причем получившаяся игра не была потенциальной в обычном смысле). Тем не менее в [135] не рассматривалась транспортировка, т. е. не было задачи нижнего уровня.

Целью настоящего раздела является предложить такое описание задачи верхнего уровня, включающее в себя описанные выше примеры, которое сводит в итоге поиск конкурентного транспортно-экономического равновесия к поиску седловой точки в выпукло-вогнутой игре. Отметим здесь, что в общем случае поиск конкурентного равновесия сводится к решению задачи дополнительности или (при другой записи) вариационному неравенству [4, 104]. При этом известно, что в общем случае это вычислительно трудные задачи. Однако в ряде случаев экономическая специфика задачи позволяет гарантировать, что полученное вариационное неравенство – монотонно. Тогда задача становится существенно привлекательнее в вычислительном плане. В данном разделе мы рассматриваем класс задач, в которых вариационные неравенства, возникающие при поиске конкурентных равновесий, переписываются в виде седловых задач. Монотонность автоматически обеспечивается правильной выпукло-вогнутой структурой седловой задачи.

Опишем вкратце структуру раздела. В подразделе 1.3.2 описывается решение «транспортной» задачи нижнего уровня (ищется равновесное распределение потоков по путям). В подразделе 1.3.3 описывается конструкция равновесного формирования корреспонденций при заданных функциях транспортных затрат. Отметим, что в этих двух пунктах мы фактически работаем только с одним экономическим агентом «Перевозчиком» (если смотреть с точки зрения популяционной теории игр, то с агентами одного типа «Перевозчиками»). В подразделе 1.3.3 этот агент(-ы) может(-гут) производить товар, неся затраты, и его потреблять, получая выгоду. В подразделе 1.3.4 модель из подраздела 1.3.3 переносится на случай, когда помимо экономического агента «Перевозчика(-ов)» в источниках и стоках транспортного графа располагаются независимые от «Перевозчика» новые экономические агенты «Производители» и «Потребители», решающие свои задачи. В заключительном подразделе 1.3.5 мо-дели верхнего и нижнего уровня объединяются в одну общую модель, конкурентное равновесие в которой сводится к поиску седловой точки в выпукло-вогнутой игре.

[16]Вектор-функция затрат, характеризующая затраты при выборе различных стратегий как функция от распределения игроков по стратегиям, является градиентом некоторой скалярной функции.

### 1.3.2. Равновесное распределение потоков по путям

Обозначим множество пар $w=(i,j)$ источник–сток $OD$, $x_p$ – поток по пути $p$; $P_w$ – множество путей, отвечающих корреспонденции $w$, $P=\bigcup_{w\in OD}P_w$ – множество всех путей; $f_e(x)=\sum_{p\in P}\delta_{ep}x_p$ – поток на ребре $e$ (здесь и далее $x=\{x_p\}$, $f=\Theta x$), где $\delta_{ep}=\begin{cases}1, & e\in p,\\ 0, & e\notin p;\end{cases}$ $\tau_e(f_e)$ – затраты на ребре $e$ ($\tau_e'(f_e)\geq 0$); $G_p(x)=\sum_{e\in E}\tau_e(f_e(x))\delta_{ep}$ – затраты на пути $p$; $X=\left\{x\geq 0:\ \sum_{p\in P_w}x_p=d_w,\ w\in OD\right\}$ – множество допустимых распределений потоков по путям, где $d_w$ – корреспонденция, отвечающая паре $w$.

**Определение 1.3.1.** Распределение потоков по путям $x=\{x_p\}\in X$ называется *равновесием* (*Нэша–Вардропа*) в популяционной игре $\langle\{x_p\}\in X,\{G_p(x)\}\rangle$, если из $x_p>0$ ($p\in P_w$) следует $G_p(x)=\min_{q\in P_w}G_q(x)$. Или, что то же самое:

для любых $w\in OD,\ p\in P_w$ выполняется $x_p\cdot\left(G_p(x)-\min_{q\in P_w}G_q(x)\right)=0.$

**Теорема 1.3.1 (см. [33, 34, 61, 141]).** *Популяционная игра* $\langle\{x_p\}\in X,\{G_p(x)\}\rangle$ *является потенциальной. Равновесие* $x^*$ *в этой игре всегда существует, и находится из решения задачи выпуклой оптимизации*:

$$x^*\in\operatorname{Arg}\min_{x\in X}\Psi(f(x)),\qquad(1.3.1)$$

*где*

$$\Psi(f(x))=\sum_{e\in E}\int_0^{f_e(x)}\tau_e(z)dz=\sum_{e\in E}\sigma_e(f_e(x)).$$

Мы оставляем в стороне вопрос единственности равновесия (детали см., например, в [33, 34]). Отметим лишь, что при естественных условиях равновесное распределение потоков по ребрам $f^*$ единственно.

В частности, для этого достаточно, чтобы $\tau_e'(f_e) > 0$ для всех $e \in E$. Если дополнительно $f^* = \Theta x$ однозначно разрешимо относительно $x$ (в реальных транспортных сетях часто случается, что число допустимых для перевозки путей меньше числа ребер, это как раз и приводит к однозначной разрешимости), то отсюда будет следовать, что равновесное распределение потоков по путям $x^*$ единственно.

Удобно считать [15, 33], что возрастающие функции затрат $\tau_e(f_e) := \tau_e^\mu(f_e)$ зависят от параметра $\mu > 0$, причем

$$\tau_e^\mu(f_e) \xrightarrow[\mu\to 0+]{} \begin{cases} \bar{t}_e, & 0 \le f_e < \bar{f}_e, \\ [\bar{t}_e, \infty), & f_e = \bar{f}_e, \end{cases}$$

$$d\tau_e^\mu(f_e)/df_e \xrightarrow[\mu\to 0+]{} 0, \; 0 \le f_e < \bar{f}_e.$$

В таком пределе задачу выпуклой оптимизации можно переписать как задачу ЛП [133]:

$$\min_{\substack{f=\Theta x,\, x\in X \\ f\le \bar{f}}} \sum_{e\in E} f_e \bar{t}_e .$$

Такого типа транспортные задачи достаточно хорошо изучены [48, 68].

Для дальнейшего будет важно переписать задачу $\min_{x\in X} \Psi(f(x))$ через двойственную [33]:

$$\min_{x\in X} \Psi(f(x)) = \min_{f,x} \left\{ \sum_{e\in E} \sigma_e(f_e) : \; f = \Theta x, \; x \in X \right\} =$$

$$= \min_{f,x} \left\{ \sum_{e\in E} \max_{t_e \in \operatorname{dom} \sigma_e^*} \left[ f_e t_e - \sigma_e^*(t_e) \right] : \; f = \Theta x, \; x \in X \right\} =$$

$$= \max_{t \in \operatorname{dom} \sigma^*} \left\{ \min_{f,x} \left[ \sum_{e\in E} f_e t_e : \; f = \Theta x, \; x \in X \right] - \sum_{e\in E} \sigma_e^*(t_e) \right\} =$$

$$= \max_{t\ge \bar{t}} \left\{ \sum_{w\in OD} d_w T_w(t) - \left\langle \bar{f}, t - \bar{t} \right\rangle - \mu \sum_{e\in E} h\left(t_e - \bar{t}_e, \bar{t}_e, \bar{f}_e, \mu\right) \right\} =$$

$$\overset{\mu\to 0+}{=} \max_{t\ge \bar{t}} \left\{ \sum_{w\in OD} d_w T_w(t) - \left\langle \bar{f}, t - \bar{t} \right\rangle \right\}, \qquad (1.3.2)$$

где $\sigma_e^*(t_e)$ – сопряженная функция к $\sigma_e(f_e)$, $T_w(t) = \min_{p\in P_w} \sum_{e\in E} \delta_{ep} t_e$ – длина кратчайшего пути из $i$ в $j$ ($w = (i, j)$) на графе, взвешенном согласно вектору $t$, $h\left(t_e-\bar{t}_e,\bar{t}_e,\bar{f}_e,\mu\right)$ – сильно выпуклая функция по первому аргументу. При этом

$$\tau_e^{\mu}\left(f_e\left(x(\mu)\right)\right)\xrightarrow[\mu\to 0+]{}t_e,$$

где $x(\mu)$ – равновесное распределение потоков, рассчитывающееся по формуле (1.3.1), а $t=\{t_e\}$ – решение задачи (1.3.2) при естественных условиях единственное [33]. Описанный предельный переход позволяет переходить к задачам, в которых вместо функции затрат на ребрах $\tau_e(f_e)$ заданны ограничения на пропускные способности $0\le f_e\le\bar{f}_e$ и затраты $\bar{t}_e$ на прохождения ребер, когда на ребрах нет «пробок» ( $f_e<\bar{f}_e$ ) [33, 133].

Основным для дальнейшего выводом из этого всего является способ (основанный на применении теоремы Демьянова–Данскина [41, 42], как правило, в градиентном варианте ввиду единственности $t$) потенциального описания набора $T(d):=\left\{T_w\left(t(d)\right)\right\}$:

$$T(d)=\nabla_d\min_{x\in X(d)}\Psi\left(f(x)\right)=\nabla_d\max_{t\ge\bar{t}}\left\{\sum_{w\in OD}d_wT_w(t)-\left\langle\bar{f},t-\bar{t}\right\rangle-\mu\sum_{e\in E}h\left(t_e-\bar{t}_e,\bar{f}_e,\mu\right)\right\}. \tag{1.3.3}$$

В [33, 34, 145] приведены эволюционные динамики, приводящие к описанным здесь равновесиям. Отметим, однако, что если рассматривать логит-динамику [33, 145] (ограничение рациональных агентов с параметром $\tilde{\gamma}>0$ [70]), то задачу (1.3.1) необходимо будет переписать в виде (говорят, что вместо равновесия Нэша–Вардропа ищется стохастическое равновесие [33, 148]):

$$\min_{x\in X}\left\{\Psi\left(f(x)\right)+\tilde{\gamma}\sum_{w\in OD}\sum_{p\in P_w}x_p\ln\left(x_p/d_w\right)\right\}. \tag{1.3.4}$$

Это замечание понадобится нам в дальнейшем.

В заключение этого раздела отметим, что теорема 1.3.1 может быть распространена и на случай, когда затраты на ребрах $\tau_e\left(f_e;\{f_{\tilde{e}}\}\right)$ удовлетворяют условию потенциальности (частный случай – это когда $\tau_e\left(f_e;\{f_{\tilde{e}}\}\right)\equiv\tau_e(f_e)$ ) [145]:

$$\frac{\partial\tau_e\left(f_e;\{f_{\tilde{e}}\}\right)}{\partial e'}=\frac{\partial\tau_{e'}\left(f_{e'};\{f_{\tilde{e}}\}\right)}{\partial e}.$$

Тогда

$$\Psi\left(f(x)\right)=\sum_{e\in E}\int_{0}^{f_e(x)}\tau_e\left(z;\left\{\{f_{\tilde e}\}_{\tilde e\neq e}\cup\{f_e=z\}\right\}\right)dz.$$

Такого рода обобщение нужно, например, когда пропускные способности узлов зависят от потоков, пересекающих узлы. В случае транспортных потоков такими узлами являются, в частности, перекрестки. Тогда, путем раздутия исходного графа, мы, с одной стороны, «развязываем узел», сводя затраты прохождения узла по разным путям к затратам прохождения фиктивных (введенных нами) ребер, а с другой – приобретаем более общую зависимость $\tau_e\left(f_e;\{f_{\tilde e}\}\right)$.

### 1.3.3. Равновесный расчет матрицы корреспонденций

В подразделе 1.3.2 матрица корреспонденций $\{d_w\}_{w\in OD}$ была задана по постановке задачи. В данном пункте мы откажемся от этого условия, вводя в источники $O$ производство, а в стоки $D$ – потребление. Агенты «появляются» в тех пунктах производства, в которых, произведя товар, его можно с выгодой для себя реализовать в каком-нибудь из пунктов потребления. Это означает, что затраты на производство и затраты на транспортировку полностью окупаются последующей выручкой от реализации продукции в пункте потребления. Агенты, которых мы здесь считаем маленькими, будут «приходить» в систему до тех пор, пока существует цепочка (пункт производства–маршрут–пункт потребления), обеспечивающая им положительную прибыль. Важно отметить, что по ходу «наплыва» агентов транспортная сеть становится все более и более загруженной, что может сказываться на затратах на перевозку. В результате прибыль ранее пришедших агентов падает, что побуждает их перераспределяться, т. е. искать более выгодные цепочки. Возникает ряд вопросов. Например, сходится ли такая динамика (точнее семейство динамик, отражающих рациональность агентов) к равновесию? Если сходится, то единственно ли равновесие? Если равновесие единственно, то как его можно эффективно найти (описать)? Попытка ответить на эти вопросы (но, прежде всего, на последний вопрос) для достаточно широкого и важного в приложениях класса задач предпринята в последующей части раздела.

Сначала, для бóльшей наглядности, отдельно рассмотрим потенциальный случай. А именно тот случай, когда в источнике $i\in O$ производственная функция имеет вид: $\sigma_i\left(f_i\right)$, где $f_i=\sum_{k:(i,k)=e\in E}f_e=\sum_{j:(i,j)\in OD}d_{ij}$, аналогично для стоков $j\in D$ определим функции полезности со знаком ми-нус: $\sigma_j\left(f_j\right)$, $f_j = \sum_{k:(k,j)=e\in E} f_e = \sum_{i:(i,j)\in OD} d_{ij}$ . Все эти функции считаем выпуклыми. Мы обозначаем эти функции одинаковыми буквами, однако это не должно вызвать в дальнейшем путаницы ввиду характерных нижних индексов. Редуцируем рассматриваемую задачу к задаче подраздела 1.3.2. Рассмотрим новый граф с множеством вершин $O\bigcup D$, соединенных теми же ребрами, что и в изначальном графе, и с одним дополнительным фиктивным источником и одним дополнительным фиктивным стоком. Этот фиктивный источник соединим со всеми источниками $O$, аналогично фиктивный сток соединим со всеми стоками $D$. Если существует путь из источника $i$ в сток $j$ в исходном графе, то в новом графе прочертим соответствующее ребро с функцией затрат $T_{ij}\left(d\right)$. Проведем дополнительное (фиктивное) ребро, соединяющее фиктивный источник с фиктивным стоком, затраты на прохождения которого тождественный ноль. Получим в итоге ориентированный граф путей из источника в сток. Легко понять, что мы оказываемся «почти» в условиях предыдущего пункта (причем с более частным графом – с одним источником и стоком) с точностью до обозначений:

$$\left\{x_p\right\} \to \left\{d_{ij}\right\},\ \left\{\tau_e\left(f_e\right)\right\} \to \left\{\sigma_i'\left(f_i\right), T_{ij}\left(d\right), \sigma_j'\left(f_j\right)\right\}.$$

«Почти» – потому что, во-первых, затраты $T_{ij}\left(d\right)$ зависят от всего набора $\left\{d_{ij}\right\}$, а не только от $d_{ij}$, а во-вторых, неясно, что в данном случае играет роль матрицы корреспонденций (в нашем случае это матрица $1\times1$, т. е. просто число). Начнем с ответа на второй вопрос. Мы считаем, что в источниках имеется потенциальная возможность производить неограниченное количество продукта, просто в какой-то момент перестает быть выгодным что-то производить и перевозить. Для этого, собственно, и было введено нулевое ребро, поток по которому обозначим $d_0$. То есть, другими словами, мы должны считать, что $\sum_{(i,j)\in W} d_{ij} + d_0 = \bar{d}$ . Если $\bar{d}$ достаточно большое, то равновесная конфигурация не зависит от того, чему именно равно $\bar{d}$ , поскольку не требуется определять $d_0$. С первой проблемой можно разобраться, немного обобщив теорему 1.3.1. Предположим, что

$$\exists\ \Phi\left(d\right)\ – \text{выпуклая}:\ T\left(d\right) = \nabla\Phi\left(d\right). \qquad (1.3.5)$$

Тогда имеет место

**Теорема 1.3.2.** *Популяционная игра*

$$\left\langle \left\{ d_{ij}, d_0 \ge 0 \right\}, \left\{ G_{ij}(d) = \sigma_i'(f_i) + T_{ij}(d) + \sigma_j'(f_j), G_0(d) \equiv 0 \right\} \right\rangle,$$

*является потенциальной. Равновесие* $d^*$ *в этой игре всегда существует* (*если* $\sigma(\cdot)$ *– сильно выпуклые функции, то равновесие гарантировано единственно*), *и находится из решения задачи выпуклой оптимизации*:

$$d^* \in \arg\min_{d \ge 0} \tilde{\Psi}(d),$$

$$\tilde{\Psi}\left(d = \left\{ d_{ij} \right\}\right) = \sum_{i \in O} \sigma_i \left( \sum_{j:(i,j) \in OD} d_{ij} \right) + \sum_{j \in D} \sigma_j \left( \sum_{i:(i,j) \in OD} d_{ij} \right) + \Phi(d). \quad (1.3.6)$$

**Доказательство.** Выпишем условие нелинейной комплиментарности (то есть, по сути, определения равновесия Нэша в популяционной игре, заданной в условии). Для этого занумеруем все индексы $ij$ и $0$ одним индексом $k$:

$$\text{для любых } k \text{ выполняется } d_k^* \cdot \left( G_k\left(d^*\right) - \min_{k'} G_{k'}\left(d^*\right) \right) = 0.$$

Действительно, допустим, что реализовалось какое-то другое равновесие $\tilde{d}$, которое не удовлетворяет этому условию. Покажем, что тогда найдется агент, которому выгодно поменять свой выбор. Действительно, тогда

$$\text{существуют такой } \tilde{k}, \text{ что } \tilde{d}_{\tilde{k}} \cdot \left( G_{\tilde{k}}\left(\tilde{d}\right) - \min_{k'} G_{k'}\left(\tilde{d}\right) \right) > 0.$$

Каждый агент (множество таких агентов непусто $\tilde{d}_{\tilde{k}} > 0$), использующий стратегию $\tilde{k}$ действует неразумно, поскольку существует такая стратегия $\bar{k}$, $\bar{k} \ne \tilde{k}$, что $G_{\bar{k}}\left(\tilde{d}\right) = \min_{k'} G_{k'}\left(\tilde{d}\right)$. Этот стратегия $\bar{k}$ более выгодна, чем $\tilde{k}$. Аналогично показывается, что при распределении $d^*$ никому из агентов уже не выгодно отклоняться от своих стратегий.

Покажем, что рассматриваемая нами игра принадлежит к классу так называемых *потенциальных игр*. В нашем случае это означает, что существует такая функция $\tilde{\Psi}\left(d = \left\{ d_{ij} \right\}\right)$, что $\partial \tilde{\Psi}(d) / \partial d_k = G_k(d)$ для любого $k$. В нашем случае легко проверить, что этому условию удовлетворяет (по условию) функция

$$\tilde{\Psi}\left(d = \left\{ d_{ij} \right\}\right) = \sum_{i \in O} \sigma_i \left( \sum_{j:(i,j) \in OD} d_{ij} \right) + \sum_{j \in D} \sigma_j \left( \sum_{i:(i,j) \in OD} d_{ij} \right) + \Phi(d).$$

Таким образом, мы имеем дело с потенциальной игрой. Оказывается, что $d^*$ – равновесие Нэша–Вардропа тогда и только тогда, когда оно доставляет минимум $\tilde{\Psi}(d)$ на множестве $d \ge 0$. Действительно, предположим, что $d^* \ge 0$ – точка минимума. Тогда, в частности, для любых $p,q$ ($d_p^* > 0$) и достаточно маленького $\delta d_p > 0$ выполняется

$$-\frac{\partial\tilde{\Psi}(d^*)}{\partial d_p}\delta d_p + \frac{\partial\tilde{\Psi}(d^*)}{\partial d_q}\delta d_p \ge 0.$$

Иначе, заменив $d^*$ на

$$\breve{d}^* = d^* + \left( \underbrace{\overbrace{0,\ldots,0,-\delta d_p}^{p},0,\ldots,0,\delta d_p}_{q},0,\ldots,0 \right) \ge 0,$$

мы пришли бы к вектору $\breve{d}^*$, доставляющему меньшее значение $\tilde{\Psi}(d)$ на множестве $d \ge 0$:

$$\tilde{\Psi}\left(\breve{d}^*\right) \approx \tilde{\Psi}(d^*) - \frac{\partial\tilde{\Psi}(d^*)}{\partial d_p}\delta d_p + \frac{\partial\tilde{\Psi}(d^*)}{\partial d_q}\delta d_p < \tilde{\Psi}(d^*).$$

Вспоминая, что $\partial\tilde{\Psi}(d)/\partial d_p = G_p(d)$, и учитывая, что $q$ можно выбирать произвольно, получаем:

для любого $p$, если $d_p^* > 0$, то выполняется $\min_q G_q(d^*) \ge G_p(d^*)$.

Но это и есть по-другому записанное условие нелинейной комплементарности. Строго говоря, мы показали сейчас только то, что точка минимума $\tilde{\Psi}(d)$ на множестве $d \ge 0$ будет равновесием Нэша–Вардропа. Аналогично рассуждая, можно показать и обратное: равновесие Нэша–Вардропа доставляет минимум $\tilde{\Psi}(d)$ на множестве $d \ge 0$. □

Именно такая конструкция и была рассмотрена в препринте [15] (для многопродуктового рынка). Если искать стохастическое равновесие, то функционал в теореме 1.3.2 необходимо энтропийно регуляризовать[17]. Такие конструкции рассматривались (в нерегуляризованном случае), например, в работах [18, 33, 34]. Уже в этих работах можно углядеть необходимость искусственного введения потенциалов (двойственных множителей) в сами функции $\sigma$. А именно: в этих работах предполагает-ся, что все эти функции $\sigma$ – линейные с неизвестными наклонами. Тем не менее считается, что при этом известно, чему должны равняться в равновесии следующие суммы (см. также раздел 1.2 этой главы):

$$\sum_{j:(i,j)\in OD} d_{ij} = L_i, \quad \sum_{i:(i,j)\in OD} d_{ij} = W_j \quad \left(\sum_{i\in O} L_i = \sum_{j\in D} W_j = N\right). \tag{1.3.7}$$

---

[17]К сожалению, строгое обоснование (теорема 11.5.12 [145]) имеется только в случае известного (фиксированного) значения $\bar{d}$ (при этом можно считать $d_0 = 0$).

То есть имеются скрытые от нас (модельера) потенциалы [48] (параметры) $\left\{\lambda_i^L, \lambda_j^W\right\}$, которые могут быть рассчитаны исходя из дополнительной информации. Применительно к модели расчета матрицы корреспонденций [18, 33, 34] выписанные дополнительные условия (1.3.7) однозначно определяют все неизвестные потенциалы. Однако при этом вместо задачи выпуклой оптимизации мы получаем минимаксную (седловую) задачу, выпуклую по $\left\{d_{ij}\right\} \ge 0$ и вогнутую, точнее даже линейную, по потенциалам $\left\{\lambda_i^L, \lambda_j^W\right\}$:

$$\min_{\substack{\{d_{ij}\}\ge 0 \\ \sum_{(i,j)\in W} d_{ij} = N}} \max_{\left\{\lambda_i^L, \lambda_j^W\right\}} \left[ \sum_{i\in O} \lambda_i^L \cdot \left( \sum_{j:(i,j)\in OD} d_{ij} - L_i \right) + \sum_{j\in D} \lambda_j^W \cdot \left( W_j - \sum_{i:(i,j)\in OD} d_{ij} \right) + \right.$$

$$\left. + \Phi(d) + \gamma \sum_{(i,j)\in OD} d_{ij} \ln\left(d_{ij}/N\right) \right]. \tag{1.3.8}$$

Эта задача всегда имеет решение.

### 1.3.4. Сетевая модель алгоритмического рыночного поведения

В данном пункте мы предложим сетевой вариант модели поиска конкурентного равновесия из препринта [135]. Однако в контексте изложенного в конце прошлого пункта, нам будет удобнее стартовать с двухстадийной модели транспортных потоков [33], приводящей к равновесию, рассчитываемому по формуле (1.3.8).

Предположим теперь, что имеется $m$ видов товара и дополнительно имеется $q$ типов материала (количества которых можно использовать в единицу времени, ограничены вектором $b$), использующихся в производстве. В источниках располагаются производители товаров, а в стоках – потребители. Мы считаем, что любой производитель товара одновременно является и потребителем, т. е. $O \subseteq D$. Обозначим через $y$ вектор цен (руб.) на материалы; $\lambda_i^L$ – вектор цен (руб.), по которым производитель продает товары перевозчику в пункте производства *i*, $\lambda_j^W$ – вектор цен (руб.), по которым потребитель покупает товары у перевозчика в пункте потребления *j*. Опишем каждого экономического агента.

***i*-й Производитель**

$U_i \subset \mathbb{R}_+^m$ – максимальное технологическое множество (замкнутое, выпуклое);

$\alpha_i \in [0,1]$ – уровень участия;

$\chi_i(\alpha_i U_i) = \alpha_i \chi_i(U_i)$ – постоянные технологические производственные затраты (руб.) при уровне участия $\alpha_i$ (в единицу времени);

$L_i \in \alpha_i U_i$, $[L_i]_k$ – количество произведенного *k*-го продукта (в единицу времени);

$A_i$, $[A_i]_{kl}$ – количество затраченного *l*-го продукта при производстве единицы *k*-го продукта;

$c_i$, $[c_i]_k$ – затраты (руб.) на производство единицы *k*-го продукта;

$R_i$, $[R_i]_{kl}$ – количество затраченного *k*-го материала для приготовления единицы *l*-го продукта.

Описанный «Производитель» решает задачу:

$$\max_{\substack{L_i \in \alpha_i U_i \\ \alpha_i \in [0,1]}} \left\{ \left\langle \lambda_i^L, L_i \right\rangle - \chi_i(\alpha_i U_i) - \left\langle \lambda_i^W, A_i L_i \right\rangle - \left\langle c_i, L_i \right\rangle - \left\langle y, R_i L_i \right\rangle \right\} =$$

$$= \max_{L_i \in U_i} \left\{ \left( \left\langle \lambda_i^L - c^i - A_i^T \lambda_i^W - R_i^T y, L_i \right\rangle - \chi_i(U_i) \right)_+ \right\}.$$

***j*-й Потребитель**

Предположим, что каждый продукт имеет *s* различных свойств (своеобразных полезностей). Это может быть, например, содержание витаминов, белков, жиров, углеводов.

$Q_j$, $\left[Q_j\right]_{kl}$ – вклад единицы *l*-го продукта в удовлетворение *k*-го свойства;

$\sigma_j$, $\left[\sigma_j\right]_k$ – минимально допустимый уровень удовлетворения *k*-го свойства (в единицу времени);

$\beta_j \in [0,1]$ – уровень участия;

$V_j = \left\{ W_j \in \mathbb{R}_+^m : Q_j W_j \ge \sigma_j \right\}$ – допустимое множество при полном участии;

$W_j \in \beta_j V_j$, $\left[W_j\right]_k$ – количество потребленного *k*-го продукта (в единицу времени);

$\tau_j$ – постоянный доход (руб.) при полном участии (в единицу времени).

Описанный «Потребитель» решает задачу:

$$\max_{\substack{W_j \in \beta_j V_j \\ \beta_j \in [0,1]}} \left\{ \beta_j \tau_j - \left\langle \lambda_j^W, W_j \right\rangle \right\} = \max_{W_j \in V_j} \left\{ \left( \tau_j - \left\langle \lambda_j^W, W_j \right\rangle \right)_+ \right\}.$$

**Перевозчик**

Этот агент решает задачу типа (1.3.6), (1.3.8):

$$\min_{\{d_{ij}\} \ge 0} \left[ \sum_{i \in O} \left\langle \lambda_i^L, \sum_{j:(i,j) \in OD} d_{ij} \right\rangle - \sum_{j \in D} \left\langle \lambda_j^W, \sum_{i:(i,j) \in OD} d_{ij} \right\rangle + \Phi(d) + \gamma \sum_{(i,j) \in OD} \left( \sum_{k=1}^m \left[ d_{ij} \right]_k \right) \ln \left( \sum_{k=1}^m \left[ d_{ij} \right]_k \right) \right],$$

в которой корреспонденции формируются «Производителями» и «Потребителями». Мы считаем, что все товары одинаковы с точки зрения «Перевозчика», т. е. $\Phi(d) := \Phi\left( \left\{ \sum_{k=1}^m \left[ d_{ij} \right]_k \right\} \right)$ (можно рассматривать и другие варианты). Здесь и далее нам будет удобно писать энтропийную регуляризацию в виде: $\gamma \sum_{(i,j) \in OD} \left( \sum_{k=1}^m \left[ d_{ij} \right]_k \right) \ln \left( \sum_{k=1}^m \left[ d_{ij} \right]_k \right)$, т. е. опускать $N = \sum_{(i,j) \in OD} \sum_{k=1}^m \left[ d_{ij} \right]_k$. Точнее полагать $N = 1$ с той же потоковой (физической) размерностью, что и $d$, чтобы под логарифмом была безразмерная величина.[18] При естественных балансовых условиях: $\sum_{i \in O} L_i = \sum_{i \in O} A_i L_i + \sum_{i \in D} W_j$ это никак не повлияет на решение задачи.

Проблема здесь в том, что все эти три типа задач завязаны друг на друга посредством векторов цен. Выпишем, как это принято при поиске конкурентных равновесий [4, 104], все имеющиеся *законы Вальраса* (балансовые ограничения + условия дополняющей нежесткости), которые накладывают совершенно естественные ограничения на эти векторы цен:

$$\sum_{j:(i,j) \in OD} d_{ij} \le L_i, \quad \left\langle \lambda_i^L, \sum_{j:(i,j) \in OD} d_{ij} - L_i \right\rangle = 0, \quad \lambda_i^L \ge 0;$$

$$\sum_{k:(k,i) \in OD} d_{ki} \ge W_i + A_i L_i, \quad \left\langle \lambda_i^W, W_i + A_i L_i - \sum_{k:(k,i) \in OD} d_{ki} \right\rangle = 0, \quad \lambda_i^W \ge 0, \quad i \in O;$$

---

[18] В случае микроскопического обоснования такого рода вариационных принципов (см. подраздел 1.3.5, а также [145]) исходя из рассмотрения соответствующей марковской динамики нащупывания равновесия, мы должны полагать $N \gg 1$, чтобы сделать соответствующий (канонический) скейлинг и перейти к детерминированной постановке.

$$\sum_{i:(i,j)\in OD} d_{ij} \ge W_j ,\ \left\langle \lambda_j^W, W_j - \sum_{i:(i,j)\in OD} d_{ij} \right\rangle = 0 ,\ \lambda_j^W \ge 0 ,\ j \in D \setminus O ;$$

$$\sum_{i\in O} R_i L_i \le b ,\ \left\langle y, b - \sum_{i\in O} R_i L_i \right\rangle = 0 ,\ y \ge 0 .$$

**Определение 1.3.2.** Набор $\left\langle \{d_{ij}\}, \{L_i\}, \{W_j\}; y, \{\lambda_i^L\}, \{\lambda_j^W\} \right\rangle$ называется *конкурентным равновесием* (*Вальраса–Нестерова–Шихмана*), если он доставляет решения задачам всех агентов и удовлетворяет законам Вальраса.

Для того чтобы установить корректность этого определения, подобно [135], введем ***условие продуктивности***:

существуют такие $\bar{L}_i \in U_i$, $\bar{W}_j \in V_j$,

что $\sum_{i\in O} \bar{L}_i > \sum_{i\in O} A_i \bar{L}_i + \sum_{i\in D} \bar{W}_j$ и $\sum_{i\in O} R_i \bar{L}_i < b$.

**Теорема 1.3.3.** *В условиях продуктивности конкурентное равновесие существует и находится из решения правильной выпукло-вогнутой седловой задачи*:

$$\min_{\{d_{ij}\}\ge 0} \max_{\substack{\{\lambda_i^L, \lambda_j^W\}\ge 0 \\ y\ge 0}} \min_{\substack{\{L_i \in \alpha_i U_i,\ \alpha_i\in[0,1]\} \\ \{W_j \in \beta_j V_j,\ \beta_j\in[0,1]\}}} \left[ \sum_{i\in O} \left( \left\langle \lambda_i^L, \sum_{j:(i,j)\in OD} d_{ij} - L_i \right\rangle + \left\langle \lambda_i^W, \sum_{k:(k,i)\in OD} A_i L_i \right\rangle + \chi_i(\alpha_i U_i) \right) + \right.$$

$$+ \sum_{j\in D} \left( \left\langle \lambda_j^W, W_j - \sum_{i:(i,j)\in OD} d_{ij} \right\rangle - \beta_j \tau_j \right) + \left\langle y, b - \sum_{i\in O} R_i L_i \right\rangle +$$

$$\left. + \Phi(d) + \gamma \sum_{(i,j)\in OD} \left( \sum_{k=1}^m [d_{ij}]_k \right) \ln \left( \sum_{k=1}^m [d_{ij}]_k \right) \right] =$$

$$= \max_{\substack{\{\lambda_i^L, \lambda_j^W\}\ge 0 \\ y\ge 0}} \left[ \langle y, b \rangle - \sum_{i\in O} \max_{L_i \in U_i} \left\{ \left( \left\langle \lambda_i^L - c^i - A_i^T \lambda_i^W - R_i^T y, L_i \right\rangle - \chi_i(U_i) \right)_+ \right\} - \right.$$

$$- \sum_{j\in D} \max_{W_j \in V_j} \left\{ \left( \tau_j - \left\langle \lambda_j^W, W_j \right\rangle \right)_+ \right\} + \min_{\{d_{ij}\}\ge 0} \left\{ \sum_{i\in O} \left\langle \lambda_i^L, \sum_{j:(i,j)\in OD} d_{ij} \right\rangle - \sum_{j\in D} \left\langle \lambda_j^W, \sum_{i:(i,j)\in OD} d_{ij} \right\rangle + \right.$$

$$\left. \left. + \Phi(d) + \gamma \sum_{(i,j)\in OD} \left( \sum_{k=1}^m [d_{ij}]_k \right) \ln \left( \sum_{k=1}^m [d_{ij}]_k \right) \right\} \right]. \qquad (1.3.9)$$

Порядок взятия минимума и максимумов можно менять согласно *Sion's minimax theorem* [33, 41, 42, 149].

### 1.3.5. Общее конкурентное равновесие

Для того чтобы объединить модели подразделов 1.3.2, 1.3.3, 1.3.4 в одну модель рассмотрим формулы (1.3.3), (1.3.5), (1.3.8), (1.3.9). Легко понять, что формула (1.3.3) как раз и задает тот самый потенциал, существование которого (формула (1.3.5)) требуется для справедливости теоремы 1.3.2 (неявно это предполагается и в теореме 1.3.3), фактически сводящей поиск конкурентного равновесия к задаче (1.3.8), а в общем случае к (1.3.9).

**Определение 1.3.3.** Набор $\left\langle \{x_p\}, \{d_{ij}\}, \{L_i\}, \{W_j\};\ y, \{\lambda_i^L\}, \{\lambda_j^W\} \right\rangle$ называется *полным* (*общим*) *конкурентным равновесием* (*Вальраса–Нэша–Вардропа–Нестерова–Шихмана*), если $\left\langle \{d_{ij}\},\ \{L_i\},\ \{W_j\};\ y, \{\lambda_i^L\}, \{\lambda_j^W\} \right\rangle$ – конкурентное равновесие, а $\{x_p\}$ является равновесием (Нэша–Вардропа) при заданном конкурентным равновесием наборе $\{d_{ij}\}$.

**Теорема 1.3.4.** *В условиях продуктивности полное конкурентное равновесие существует и находится из решения правильной выпукло-вогнутой седловой задачи*:

$$\max_{\substack{\{\lambda_i^L, \lambda_j^W\} \ge 0 \\ y \ge 0}} \Bigg[ \langle y, b \rangle - \sum_{i \in O} \max_{L_i \in U_i} \left\{ \left( \left\langle \lambda_i^L - c^i - A_i^T \lambda_i^W - R_i^T y, L_i \right\rangle - \chi_i(U_i) \right)_+ \right\} -$$

$$- \sum_{j \in D} \max_{W_j \in V_j} \left\{ \left( \tau_j - \left\langle \lambda_j^W, W_j \right\rangle \right)_+ \right\} + \min_{\{d_{ij}\} \ge 0} \Bigg\{ \sum_{i \in O} \left\langle \lambda_i^L, \sum_{j:(i,j) \in OD} d_{ij} \right\rangle - \sum_{j \in D} \left\langle \lambda_j^W, \sum_{i:(i,j) \in OD} d_{ij} \right\rangle +$$

$$+ \max_{t \ge \bar{t}} \left\{ \sum_{(i,j) \in OD} \left( \sum_{k=1}^{m} \left[ d_{ij} \right]_k \right) T_{ij}(t) - \left\langle \bar{f}, t - \bar{t} \right\rangle - \mu \sum_{e \in E} h\left(t_e - \bar{t}_e, \bar{f}_e, \mu\right) \right\} +$$

$$+ \gamma \sum_{(i,j) \in OD} \left( \sum_{k=1}^{m} \left[ d_{ij} \right]_k \right) \ln \left( \sum_{k=1}^{m} \left[ d_{ij} \right]_k \right) \Bigg\} \Bigg]. \qquad (1.3.10)$$

Таким образом, поиск общего конкурентного равновесия также сводится к седловой задаче (если мы вынесем все маскимумы и минимумы за квадратную скобку, то получим минимаксную = седловую задачу), имеющей правильную структуру с точки зрения того, что минимум берется по переменным, по которым выражение в квадратных скобках выпукло, а максимум – по переменным, по которым выражение вогнуто. Порядок взятия всех максимумов и минимума можно менять согласно *Sion's minimax theorem*. В частности, это дает возможность явно выполнить минимизацию по $\{d_{ij}\} \ge 0$, «заплатив» за это некоторым усложнением получив-шегося в итоге функционала, который также сохранит правильные выпукло-вогнутые свойства [33].

Мы не будем здесь приводить, что получается после подстановки формулы (1.3.3) в формулы (1.3.6) и (1.3.8). Все выкладки аналогичны и даже проще. Тем не менее, ссылаясь далее на задачи (1.3.6) и (1.3.8), мы будем считать, что такая подстановка была сделана.

Такого рода задачи можно эффективно численно решать (причем содержательно интерпретируемым способом), если транспортный граф задачи нижнего уровня (поиска равновесного распределения потоков) несверхбольшой [121, 128, 129, 136]. Если же этот граф имеет, скажем, порядка $10^5$ ребер, как транспортный граф Москвы и области [34], то требуется разработка новых эффективных методов, учитывающих разреженность задачи и использующих рандомизацию. Мы не будем здесь на этом останавливаться, поскольку планируется посвятить численным методам решения таких задач больших размеров отдельную публикацию. Впрочем, некоторые возможные подходы отчасти освещены в [15, 33]. К сожалению, численный метод, предложенный в [15], не совсем корректен.

Сделаем несколько замечаний в связи с полученным результатом.

Во-первых, в рассмотренных в разделе задачах с помощью штрафных механизмов (типа платных дорог) можно добиваться того, что возникающие равновесия соответствовали бы социальному оптимуму. Для этого можно использовать VCG-механизм [69, 144], см. также раздел 1.1 этой главы.

Во-вторых, используя аппарат [15, 33, 135], несложно вычленить из выписанных задач (1.3.6), (1.3.8), (1.3.10) всевозможные цены, тарифы, длины очередей (пробок), – если делаем предельный переход $\mu \to 0+$ и т. п., – понимаемые в смысле Л. В. Канторовича как двойственные множители.

В-третьих, рассматривая два разномасштабных по времени марковских процесса нащупывания равновесной конфигурации, можно прийти к решению задач (1.3.6), (1.3.8) и, с некоторыми оговорками, (1.3.10). Например, если в быстром времени динамика перераспределения потоков по путям задается имитационной логит-динамикой [145], а в медленном времени процесс перераспределения корреспонденций (исходя из быстро подстраивающихся затрат $\{T_{ij}(d)\}$) задается просто логит-динамикой [145], то выражение в квадратных скобках (1.3.8) будет играть роль действия в теореме типа Санова, т. е. описывать экспоненциальную концентрацию инвариантной меры марковского процесса с оговоркой, что речь идет о переменных $d$ и $x$ [145]. Аналогично можно сказать и про $\tilde{\Psi}$ в теореме 1.3.2 после подстановки (1.3.3). Кроме того, эти же самые функции будут играть роль функций Ляпунова соответствующих прошкалированным (каноническим скейлингом) марковских динамик, приводящих к СОДУ Тихоновского типа [14, 33, 145]. Это также следует из общих результатов работы [7, 22] (см. также приложение в конце пособия). Отметим, что относительно нащупывания цен (потенциалов) в задачах (1.3.8) и особенно в (1.3.10) работают механизмы похожие на те, которые описаны в классической работе [48]. Другими словами, при фиксированных ценах (потенциалах) динамика соответствует классическим популяционным динамикам нащупывания равновесий [145]. Но из-за того, что потенциалы неизвестны и, в свою очередь, должны как-то параллельно подбираться, предполагается, что в медленном времени экономические агенты переоценивают эти потенциалы исходя из той обратной связи (пример имеется в [48]), которую они ожидают увидеть, и того, что они реально видят.

В-четвертых, упомянутая выше эволюционная динамика при правильной дискретизации дает разумный численный способ поиска конкурентного равновесия. В частности, упоминаемая имитационная логит-динамика при правильной дискретизации даст метод зеркального спуска / двойственных усреднений, представляющий собой метод проекции градиента с усреднением [128] и без усреднения [136], где проекция понимается в смысле «расстояния» Кульбака–Лейблера. Зеркальный спуск можно получить также из дискретизации логит-динамики, если ориентироваться не только на предыдущую итерацию, а на среднее арифметическое всех предыдущих итераций [33, 38]. В работе [38] поясняется некоторая привилегированность этих двух логит-динамик (см. также [33, 36, 70, 145] и подраздел 1.1.4 раздела 1.1 этой главы). Отметим при этом, что логит-динамики может быть проинтерпретирована так же, как и имитационная логит-динамика для потенциальной игры с энтропийно-регуляризованным потенциалом [38] (см. также подраздел 1.1.4 раздела 1.1 этой главы).

В-пятых, везде выше мы исходили из того, что есть разные масштабы времени. Из-за этого задачи подразделов 1.3.2–1.3.4 удалось завязать с помощью формул (1.3.3), (1.3.5). Однако к аналогичным выводам можно было прийти, если вместо введения разных масштабов времени ввести иерархию в принятии решений [70]. Скажем, сначала пользователь транспортной сети выбирает тип транспорта (личный / общественный), а потом маршрут [33]. Здесь особенно актуальным становятся такие модели дискретного выбора, как Nested Logit [70]. А именно, если использовать энтропийную регуляризацию только в одной из этих двух задач разного уровня (иерархии), описанных в подразделах 1.3.2, 1.3.3, то получается обычная (Multinomial) логит-модель выбора (например, в (1.3.10) мы регуляризовали только задачу верхнего уровня), но если энтропийно регу-ляризовать обе задачи, то получится двухуровневая Nested логит-модель выбора [70]. Это означает, что соответствующая Nested логит-динамика в популяционной иерархической игре приводит к равновесию, которое описывается решением задач типа (1.3.8), (1.3.10) с дополнительной энтропийной регуляризацией задачи нижнего уровня. Несложно показать, что хорошие выпукло-вогнутые свойства задач (1.3.8), (1.3.10) при этом сохраняются. Да и в вычислительном плане задача не становится принципиально сложнее, особенно если учесть конструкцию *The shortest path problem,* описанную в пятой главе монографии [116], см. также [33]. Подробнее об этом будет написано в следующем разделе.

Резюмируем полученные в разделе результаты. На конкретных семействах примеров (но тем не менее достаточно богатых в смысле встречаемости в приложениях) была продемонстрирована некоторая «алгебра» над различными конструкциями равновесия. Было продемонстрировано, как можно сочетать их друг с другом, чтобы получать все более и более содержательные задачи. Ключевым местом стал переход, связанный с формулой (1.3.3), который можно понимать как произведение (суперпозицию) транспортно-экономических моделей, и конструкция задач (1.3.8), (1.3.9), которую можно понимать как «сумму» моделей. Представляется, что в этом направлении, может возникнуть довольно интересное движение, связанное с вычленением той «минимальной алгебры операций» над моделями, с помощью которой можно было бы описывать большое семейство равновесных конфигураций, встречающихся в различных приложениях.

## 1.4. О связи моделей дискретного выбора с разномасштабными по времени популяционными играми загрузок

### 1.4.1. Введение

В работах [5, 20, 33, 34], а также в подразделе 1.1.10 раздела 1.1 этой главы, было анонсировано, что в последующем цикле публикаций будет приведен общий способ вариационного описания равновесий (стохастических равновесий) в популярных моделях распределения транспортных потоков. Также отмечалось, что планируется предложить эффективные численные методы поиска таких равновесий. В данном разделе предпринята попытка погрузить известные нам подходы к многостадийному моделированию потоков на иерархических сетях (реальных транспортных сетях или сетях принятия решение – неважно) в одну общую схему, сводящую поиск равновесия к решению многоуровневой задачи выпуклой оптимизации. В основе схемы получения вариационного принципа для описания равновесия лежит популяционная игра загрузки с соответствующими логит-динамиками (отвечающими моделям дискретного выбора [70]) пользователей на каждом уровне иерархии [145] (см. подраздел 1.4.2). Для решения описанной задачи выпуклой оптимизации в разделе изучается двойственная задача, представляющая самостоятельный интерес (см. подраздел 1.4.3). Основнымы инструментами изучения двойственной задачи являются аппарат характеристических функций на графе [23, 34, 123] и ускоренные прямодвойственные методы в композитном варианте [126, 134].

Отметим, что общность результатов раздела достигается за счет введения большого числа параметров, которые можно вырождать, стремя их к нулю или бесконечности. Игра на выборе этих параметров позволяет, например, получать различные многостадийные модели транспортных потоков [5, 15, 20, 33, 140]. Приводимые далее результаты можно обобщать и на потоки товаров в случае, когда имеется более одного наименования товара [15]. Однако в данном разделе мы не будем касаться этого обобщения. Мы также не планируем приводить конкретные примеры получения многостадийных транспортных моделей согласно изложенной в разделе общей схемы.

### 1.4.2. Постановка задачи

Рассмотрим транспортную сеть, заданную ориентированным графом $\Gamma^1 = \left\langle V^1, E^1 \right\rangle$. Часть его вершин $O^1 \subseteq V^1$ является источниками, часть – стоками $D^1 \subseteq V^1$. Множество пар источник–сток обозначим $OD^1 \subseteq O^1 \otimes D^1$. Пусть каждой паре $w^1 \in OD^1$ соответствует своя корреспонденция пользователей: $d^1_{w^1} := d^1_{w^1} \cdot M$ ($M \gg 1$), которые хотят в единицу времени перемещаться из источника в сток и которые соответствуют заданной корреспонденции $w^1$. Пусть ребра $\Gamma^1$ разделены на два типа $E^1 = \tilde{E}^1 \coprod \bar{E}^1$. Ребра типа $\tilde{E}^1$ характеризуются неубывающими функциями затрат: $\tau^1_{e^1}\left(f^1_{e^1}\right) := \tau^1_{e^1}\left(f^1_{e^1}/M\right)$. Затраты $\tau^1_{e^1}\left(f^1_{e^1}\right)$ несут те пользователи, которые используют в своем пути ребро $e^1 \in \tilde{E}^1$ в предположении, что поток пользователей по этому ребру равен $f^1_{e^1}$. Пары вершин, задающие ребра типа $\bar{E}^1$, являются, в свою очередь, парами источник–сток $OD^2$ (с корреспонденциями $d^2_{w^2} = f^1_{e^1}$, $w^2 = e^1 \in \bar{E}_1$) в транспортной сети следующего уровня: $\Gamma^2 = \left\langle V^2, E^2 \right\rangle$, ребра которой, в свою очередь, разделены на два типа: $E^2 = \tilde{E}^2 \coprod \bar{E}^2$. Ребра типа $\tilde{E}^2$ характеризуются неубывающими функциями затрат $\tau^2_{e^2}\left(f^2_{e^2}\right) := \tau^2_{e^2}\left(f^2_{e^2}/M\right)$. Затраты $\tau^2_{e^2}\left(f^2_{e^2}\right)$ несут те пользователи, которые используют в своем пути ребро $e^2 \in \tilde{E}^2$ в предположении, что поток пользователей по этому ребру равен $f^2_{e^2}$. Пары вершин, задающие ребра типа $\bar{E}^2$, являются, в свою очередь, парами источник–сток $OD^3$ (с корреспонденциями $d^3_{w^3} = f^2_{e^2}$, $w^3 = e^2 \in \bar{E}^2$) в транспортной сети более высокого уровня: $\Gamma^3 = \left\langle V^3, E^3 \right\rangle, \ldots$ и т. д. Будем считать, что всего имеется $m$ уровней: $\tilde{E}^m = E^m$. Обычно в приложениях число $m$ небольшое, равное 2–10 [5, 20, 33, 34, 70].

Каждый пользователь в графе $\Gamma^1$ выбирает путь $p^1_{w^1} \in P^1_{w^1}$ (последовательный набор проходимых пользователем ребер), соответствующий его корреспонденции $w^1 \in OD^1$ ($P^1_{w^1}$ – множество всех путей, отвечающих в $\Gamma^1$ корреспонденции $w^1$). Задав $p^1_{w^1}$ можно однозначно восстановить ребра типа $\bar{E}^1$, входящие в этот путь. На каждом из этих ребер $w^2 \in \bar{E}^1$ пользователь может выбирать свой путь $p^2_{w^2} \in P^2_{w^2}$ ($P^2_{w^2}$ – множество всех путей, отвечающих в $\Gamma^2$ корреспонденции $w^2$),… и т. д. Пусть каждый пользователь сделал свой выбор. Обозначим через $x^1_{p^1}$ величину потока пользователей по пути $p^1 \in P^1 = \coprod_{w^1 \in OD^1} P^1_{w^1}$, через $x^2_{p^2}$ – величину потока пользователей по пути $p^2 \in P^2 = \coprod_{w^2 \in OD^2} P^2_{w^2}$,… и т. д. Заметим, что

$$x^k_{p^k_{w^k}} \geq 0, \; p^k_{w^k} \in P^k_{w^k}, \; \sum_{p^k_{w^k} \in P^k_{w^k}} x^k_{p^k_{w^k}} = d^k_{w^k}, \; w^k \in OD^k, \; k = 1,\ldots,m,$$

и для компактности мы будем далее записывать

$$\left\{x^k_{p^k_{w^k}}\right\}_{p^k_{w^k} \in P^k_{w^k}} \in S_{\left|P^k_{w^k}\right|}\left(d_{w^k}\right).$$

Отметим, что здесь и везде в дальнейшем

$$w^{k+1}\left(=e^k\right) \in OD^{k+1}\left(=\bar{E}^k\right), \; d^{k+1}_{w^{k+1}} = f^k_{e^k}, \; k = 1,\ldots,m-1.$$

Введем для графа $\Gamma^k$ и множества путей $P^k$ матрицу (Кирхгофа):

$$\Theta^k = \left\|\delta_{e^k p^k}\right\|_{e^k \in E^k, p^k \in P^k}, \; \delta_{e^k p^k} = \begin{cases} 1, & e^k \in p^k, \\ 0, & e^k \notin p^k, \end{cases} \; k = 1,\ldots,m.$$

Тогда вектор потоков на ребрах $f^k$ графа $\Gamma^k$ однозначно определяется вектором потоков на путях $x^k = \left\{x^k_{p^k}\right\}_{p^k \in P^k}$:

$$f^k = \Theta^k x^k,\ k = 1,\dots,m.$$

Обозначим через

$$x = \left\{x^k\right\}_{k=1}^m,\ f = \left\{f^k\right\}_{k=1}^m,\ \Theta = \mathrm{diag}\left\{\Theta^k\right\}_{k=1}^m,$$

$$d^k = \left\{d^k_{w^k}\right\}_{w^k \in OD^k},\ X^k\left(d^k\right) = \coprod_{w^k \in OD^k} S_{\left|P^k_{w^k}\right|}\left(d^k_{w^k}\right),\ X = \coprod_{k=1}^m X^k\left(d^k\right),$$

а через

$$\breve{p}^k_{w^k} = \left(p^k_{w^k}, \left\{p^{k+1}_{w^{k+1}}\right\}_{w^{k+1} \in p^k_{w^k} \cap \bar{E}^k}, \dots\right),\ k = 1,\dots,m$$

полное описание возможного пути (в графе $\Gamma^k$ и графах следующих уровней), соответствующего корреспонденции $w^k \in OD^k$. Множество всех таких путей будем обозначать $\breve{P}^k_{w^k}$. Введем также множество путей $\breve{P}^k = \coprod_{w^k \in OD^k} \breve{P}^k_{w^k}$ и соответствующий вектор распределения потоков по этим путям $x_{\breve{P}^k}$. Определим функции затрат пользователей на пути $\breve{p}^k_{w^k}$ по индукции:

$$G^m_{\breve{p}^m_{w^m}}\left(x_{\breve{P}^m}\right) = \sum_{e^m \in \tilde{E}^m} \delta_{e^m p^m} \tau^m_{e^m}\left(f^m_{e^m}\right),$$

$$G^k_{\breve{p}^k_{w^k}}\left(x_{\breve{P}^k}\right) = \sum_{e^k \in \tilde{E}^k} \delta_{e^k \breve{p}^k_{w^k}} \tau^k_{e^k}\left(f^k_{e^k}\right) + \sum_{w^{k+1} \in \bar{E}^k} G^{k+1}_{\breve{p}^{k+1}_{w^{k+1}}}\left(x_{\breve{P}^{k+1}}\right),\ k = 1,\dots,m-1.$$

Опишем *марковскую логит-динамику* (также говорят *гиббсовскую динамику*) в повторяющейся игре загрузки графа транспортной сети [23, 145] (см. также подраздел 1.1.4 раздела 1.1 этой главы). Пусть имеется $TN$ шагов ($N \gg 1$). Каждый пользователь транспортной сети, использовавший на шаге $t$ путь $\breve{p}^1_{w^1}$, независимо от остальных, на шаге $t+1$ (все введенные новые параметры положительны):

- с вероятностью

$$\frac{\lambda^1}{N} \frac{\exp\left(-G^{t,1}_p / \gamma^1\right)}{\sum_{\tilde{p} \in \breve{P}^1_{w^1}} \exp\left(-G^{t,1}_{\tilde{p}} / \gamma^1\right)}$$

пытается изменить свой путь $\breve{p}^1_{w^1}$ на $p \in \breve{P}^1_{w^1}$, где $G^{t,1}_p = G^1_p\left(x^t_{\breve{P}^1}\right)$ – затраты на пути $p$ на шаге $t$ ($G^{0,1}_p \equiv 0$);

- равновероятно выбирает $w^2 \in p^1_{w^1} \cap \bar{E}^1$ и затем с вероятностью

$$\frac{\lambda^2 \left| p^1_{w^1} \cap \bar{E}^1 \right|}{N} \frac{\exp\left(-G^{t,2}_p / \gamma^2\right)}{\sum_{\tilde{p} \in \breve{P}^2_{w^2}} \exp\left(-G^{t,2}_{\tilde{p}} / \gamma^2\right)}$$

пытается изменить в своем пути $\breve{p}^1_{w^1} = \left( p^1_{w^1}, \left\{ \breve{p}^2_{w^2} \right\}_{w^2 \in p^1_{w^1} \cap \bar{E}^1} \right)$ участок пути $\breve{p}^2_{w^2}$, выбирая путь $p \in \breve{P}^2_{w^2}$, где $G^{t,2}_p = G^2_p\left(x^t_{\breve{P}^2}\right)$ – затраты на пути $p$ на шаге $t$ ($G^{0,2}_p \equiv 0$);

- ... и т. д.;
- с вероятностью

$$1 - \sum_{k=1}^{m} \lambda^k \Big/ N$$

решает не менять тот путь, который использовал на шаге $t$.

Такая динамика отражает ограниченную рациональность агентов (см. замечание 1.4.5 подраздела 1.4.3) и часто используется в теории дискретного выбора [70] и популяционной теории игр [145]. В основном нас будет интересовать поведение такой системы в предположении, что

$$\lambda^2/\lambda^1 \to \infty,\ \lambda^3/\lambda^2 \to \infty,\ \ldots,\ \lambda^m/\lambda^{m-1} \to \infty,\ N/\lambda^m \to \infty. \qquad (1.4.1)$$

Эта марковская динамика в пределе $N \to \infty$ превращается в марковскую динамику в непрерывном времени [95]. Далее мы, как правило, будем считать, что такой предельный переход был осуществлен.

В пределе $M \to \infty$ эта динамика (концентраций) описывается зацепляющейся системой обыкновенных дифференциальных уравнений (СОДУ):

$$\frac{dx_{\breve{p}^1_{w^1}}}{dt} = \lambda^1 \cdot \left( d^1_{w^1} \frac{\exp\left(-G^1_{\breve{p}^1_{w^1}}\left(x_{\breve{P}^1}\right) \big/ \gamma^1\right)}{\sum_{\tilde{p} \in \breve{P}^1_{w^1}} \exp\left(-G^1_{\tilde{p}}\left(x_{\breve{P}^1}\right) \big/ \gamma^1\right)} - x_{\breve{p}^1_{w^1}} \right),\ \breve{p}^1_{w^1} \in \breve{P}^1_{w^1},\ w^1 \in OD^1,$$

$$\frac{dx_{\breve{p}^2_{w^2}}}{dt} = \lambda^2 \cdot \left( d^2_{w^2} \frac{\exp\left(-G^2_{\breve{p}^2_{w^2}}\left(x_{\breve{P}^2}\right) \big/ \gamma^2\right)}{\sum_{\tilde{p} \in \breve{P}^2_{w^2}} \exp\left(-G^2_{\tilde{p}}\left(x_{\breve{P}^2}\right) \big/ \gamma^2\right)} - x_{\breve{p}^2_{w^2}} \right),\ \breve{p}^2_{w^2} \in \breve{P}^2_{w^2},\ w^2 \in OD^2 = \bar{E}^1,$$

Применяя по индукции (ввиду условия (1.4.1)) теорему Тихонова к этой СОДУ [62], можно получить описание аттрактора СОДУ – глобально устойчивой (при $T \to \infty$) неподвижной точки. Для того чтобы это сделать, введем обозначение

$$\sigma_{e^k}^k\left(f_{e^k}^k\right) = \int_0^{f_{e^k}^k} \tau_{e^k}^k(z)\,dz\,,\ k = 1,\ldots,m\,.$$

Рассмотрим задачу

$$\begin{gathered}\Psi(x,f) := \Psi^1(x) = \sum_{e^1\in\tilde{E}^1}\sigma_{e^1}^1\left(f_{e^1}^1\right) + \Psi^2(x) + \\ +\gamma^1\sum_{w^1\in OD^1}\sum_{p^1\in P_{w^1}^1} x_{p^1}^1 \ln\left(x_{p^1}^1 / d_{w^1}^1\right) \to \min_{f=\Theta x,\ x\in X},\end{gathered} \qquad (1.4.2)$$

$$\Psi^2(x) = \sum_{e^2\in\tilde{E}^2}\sigma_{e^2}^2\left(f_{e^2}^2\right) + \Psi^3(x) + \gamma^2\sum_{w^2\in\bar{E}^1}\sum_{p^2\in P_{w^2}^2} x_{p^2}^2 \ln\left(x_{p^2}^2 / d_{w^2}^2\right),\ d_{w^2}^2 = f_{w^2}^1,$$

$$\Psi^k(x) = \sum_{e^k\in\tilde{E}^k}\sigma_{e^k}^k\left(f_{e^k}^k\right) + \Psi^{k+1}(x) + \gamma^k\sum_{w^k\in\bar{E}^{k-1}}\sum_{p^k\in P_{w^k}^k} x_{p^k}^k \ln\left(x_{p^k}^k / d_{w^k}^k\right),\ d_{w^{k+1}}^{k+1} = f_{w^{k+1}}^k,$$

$$\Psi^m(x) = \sum_{e^m\in E^m}\sigma_{e^m}^m\left(f_{e^m}^m\right) + \gamma^m\sum_{w^m\in\bar{E}^{m-1}}\sum_{p^m\in P_{w^m}^m} x_{p^m}^m \ln\left(x_{p^m}^m / d_{w^m}^m\right),\ d_{w^m}^m = f_{w^m}^{m-1}.$$

Эта задача эквивалентна следующей цепочке зацепляющихся задач выпуклой (многоуровневой [119]) оптимизации:

$$\begin{gathered}\Phi^1\left(d^1\right) = \min_{\substack{f^1=\Theta^1 x^1,\ x^1\in X^1\left(d^1\right)\\ d_{e^1}^2 = f_{e^1}^1,\ e^1\in\bar{E}^1}}\left\{\sum_{e^1\in\tilde{E}^1}\sigma_{e^1}^1\left(f_{e^1}^1\right) + \Phi^2\left(d^2\right) + \right.\\ \left. +\gamma^1\sum_{w^1\in OD^1}\sum_{p^1\in P_{w^1}^1} x_{p^1}^1 \ln\left(x_{p^1}^1 / d_{w^1}^1\right)\right\},\end{gathered} \qquad (1.4.3)$$

$$\Phi^2\left(d^2\right) = \min_{\substack{f^2=\Theta^2 x^2,\ x^2\in X^2\left(d^2\right)\\ d_{e^2}^3 = f_{e^2}^2,\ e^2\in\bar{E}^2}}\left\{\sum_{e^2\in\tilde{E}^2}\sigma_{e^2}^2\left(f_{e^2}^2\right) + \Phi^3\left(d^3\right) + \gamma^2\sum_{w^2\in\bar{E}^1}\sum_{p^2\in P_{w^2}^2} x_{p^2}^2 \ln\left(x_{p^2}^2 / d_{w^2}^2\right)\right\},$$

$$\begin{gathered}\Phi^k\left(d^k\right) = \min_{\substack{f^k=\Theta^k x^k,\ x^k\in X^k\left(d^k\right)\\ d_{e^k}^{k+1} = f_{e^k}^k,\ e^k\in\bar{E}^k}}\left\{\sum_{e^k\in\tilde{E}^k}\sigma_{e^k}^k\left(f_{e^k}^k\right) + \Phi^{k+1}\left(d^{k+1}\right) + \right.\\ \left. +\gamma^k\sum_{w^k\in\bar{E}^{k-1}}\sum_{p^k\in P_{w^k}^k} x_{p^k}^k \ln\left(x_{p^k}^k / d_{w^k}^k\right)\right\}\end{gathered}$$

$$\Phi^m\left(d^m\right) = \min_{f^m=\Theta^m x^m,\ x^m\in X^m\left(d^m\right)}\left\{\sum_{e^m\in E^m}\sigma_{e^m}^m\left(f_{e^m}^m\right) + \gamma^m\sum_{w^m\in\bar{E}^{m-1}}\sum_{p^m\in P_{w^m}^m} x_{p^m}^m \ln\left(x_{p^m}^m / d_{w^m}^m\right)\right\}.$$

То, что эти задачи выпуклые, сразу может быть не очевидно. Чтобы это понять, заметим, что ограничения $x^k \in X^k\left(d^k\right)$ с помощью метода множителей Лагранжа можно убрать, добавив в функционал слагаемые

$$\sum_{w^k \in \bar{E}^{k-1}} \max_{\lambda^k_{w^k}} \left\langle \lambda^k_{w^k}, \sum_{p^k \in P^k_{w^k}} x^k_{p^k} - d^k_{w^k} \right\rangle, \ k = 1,...,m.$$

Каждое такое слагаемое есть выпуклая функция по совокупности параметров $x^k$, $d^k$ (см., например, формулу (3.1.8), с. 96 [45]). Следовательно (см., например, теорему 3.1.2, с. 92, и формулу (3.1.9), с. 96 [45]), $\Phi^m\left(d^m\right)$ – выпуклая функция, но тогда и $\Phi^k\left(d^k\right)$ – выпуклая функция, поскольку (по индукции) $\Phi^{k+1}\left(d^{k+1}\right)$ – выпуклая функция ( $k = 1,...,m-1$).

**Теорема 1.4.1.** 1. *Задачи* (1.4.2) *и* (1.4.3) *являются эквивалентными задачами выпуклой оптимизации, имеющими единственное решение.*

2. *Введенная марковская логит-динамика при* $N \to \infty$ *– эргодическая. Её финальное распределение* (*возникающее в пределе* $T \to \infty$ ) *совпадает со стационарным. В предположении* (1.4.1) *стационарное распределение экспоненциально сконцентрировано в окрестности решения задачи* (1.4.3) (*в пределе* $M \to \infty$ *стационарное распределение полностью сосредотачивается на решении задачи* (1.4.3)).

3. *Введенная марковская логит-динамика при пределах* $N \to \infty$, $M \to \infty$ *описывается СОДУ. В предположении* (1.4.1) *любая допустимая траектория СОДУ* (*соответствующая вектору корреспонденций* $d^1$) *сходится при* $T \to \infty$ *к решению задачи* (1.4.3).

**Замечание 1.4.1.** Утверждения 1, 2 теоремы 1.4.1 (кроме единственности решения) остаются верными и в предположении, что по части параметров $\gamma^k$ сделаны предельные переходы (от стохастических равновесий к равновесиям Нэша): $\gamma^k \to 0+$ (важно, что эти переходы осуществляются после предельных переходов, указанных в соответствующих пунктах теоремы 1.4.1). К этому же результату (с точки зрения того, к какой задаче оптимизации в итоге сводится поиск равновесий) приводит рассмотрение на соответствующих уровнях вместо логит-динамик имитационных логит-динамик [36, 145].

**Замечание 1.4.2.** Утверждения теоремы 1.4.1 и замечания 1.4.1 остаются верными, если на части ребер (любого уровня) сделать предельные переходы (важно, что эти переходы осуществляются после предельных переходов, указанных в соответствующих пунктах теоремы) вида (предел стабильной динамики [15, 33], см. также подраздел 1.1.7 раздела 1.1 этой главы):

$$\tau_e(f_e) := \tau_e^{\mu}(f_e) \xrightarrow[\mu\to 0+]{} \begin{cases} \bar{t}_e, & f_e < \bar{f}_e, \\ [\bar{t}_e, \infty), & f_e = \bar{f}_e, \end{cases}$$

$$d\tau_e^{\mu}(f_e)/df_e \xrightarrow[\mu\to 0+]{} 0,\ 0 \le f_e < \bar{f}_e,$$

с дополнительной оговоркой, что существует такой $x \in X$ , что условие $f = \Theta x$ совместно с $\left\{ f_e < \bar{f}_e \right\}_e$. При этом

$$\sigma_e(f_e) = \lim_{\mu\to 0+} \int_0^{f_e} \tau_e^{\mu}(z)\,dz = \begin{cases} f_e \bar{t}_e, & f_e \le \bar{f}_e, \\ \infty, & f_e > \bar{f}_e. \end{cases}$$

Величину $t_e = \lim_{\mu\to 0+} \tau_e^{\mu}\left(f_e^{\mu}\right) \ge \bar{t}_e$ можно понимать как затраты на проезд по ребру $e$ (см. также подраздел 1.4.3), а $\lim_{\mu\to 0+} \tau_e^{\mu}\left(f_e^{\mu}\right) - \bar{t}_e$ – как дополнительные затраты, приобретенные из-за наличия «пробки» на ребре $e$ [33, 130], возникшей из-за функционирования ребра на пределе пропускной способности $\bar{f}_e$. Эти дополнительные затраты в точности совпадают с множителем Лагранжа к ограничению $f_e \le \bar{f}_e$ [33, 130]. Их также можно понимать как оптимальные платы за проезд (для обычных ребер эти платы равны $f_e\, d\tau_e^{\mu}(f_e)/df_e$ [20, 69, 144]), взимаемые согласно механизму Викри–Кларка–Гроуса [69].

Если для некоторых $1 \le p \le q \le m$ имеют место равенства $\gamma^p = \ldots = \gamma^q$, то можно свернуть $\Gamma^p$, …, $\Gamma^q$ в один граф $\coprod_{k=p}^{q} \Gamma^k$. Это следует из свойств энтропии (см. свойство 3 § 4 гл. 2 [66]).

Далее мы отдельно рассмотрим специальный случай $\gamma^1 = \ldots = \gamma^m = \gamma$. В этом случае имеем граф

$$\Gamma = \coprod_{k=1}^{m} \Gamma^k = \left\langle V, E = \coprod_{k=1}^{m} \tilde{E}^k \right\rangle,$$

у которого имеет всего один уровень, а задача (1.4.1) может быть переписана следующим образом:

$$\Psi(x) = \sum_{e\in E} \sigma_e(f_e) + \gamma \sum_{w^1 \in OD^1} \sum_{p \in \tilde{P}^1_{w^1}} x_p \ln\left(x_p / d^1_{w^1}\right) \to \min_{f=\Theta x,\ x\in X}. \qquad (1.4.4)$$

**Теорема 1.4.2.** *При* $\gamma^1 = \ldots = \gamma^m = \gamma$

1) *задачи* (1.4.2) *и* (1.4.4) *являются эквивалентными задачами выпуклой оптимизации, имеющими единственное решение*;
2) *введенная марковская логит-динамика при* $N \to \infty$ *– эргодическая. Её финальное распределение* (*возникающее в пределе* $T \to \infty$) *совпадает со стационарным, которое представимо в виде* (*представление Санова*):

$$\sim \exp\left( -\frac{M}{\gamma} \cdot \left( \Psi(x) + o(1) \right) \right), \; M \gg 1.$$

*Как следствие, получаем, что стационарное распределение экспоненциально сконцентрировано в окрестности решения задачи* (1.4.4) (*в пределе* $M \to \infty$ *стационарное распределение полностью сосредотачивается на решении задачи* (1.4.4));
3) *введенная марковская логит-динамика при пределах* $N \to \infty$, $M \to \infty$ *описывается СОДУ. Функция* $\Psi(x)$ *является функцией Ляпунова этой СОДУ* (*принцип Больцмана*), *то есть убывает на траекториях СОДУ. Как следствие, любая допустимая траектория СОДУ* (*соответствующая вектору корреспонденций* $d^1$) *сходится при* $T \to \infty$ *к решению задачи* (1.4.4).

**Замечание 1.4.3.** К теореме 1.4.2 можно сделать замечания, аналогичные замечаниям 1.4.1, 1.4.2 к теореме 1.4.1.

## 1.4.3. Двойственная задача

Рассмотрим граф

$$\Gamma = \coprod_{k=1}^{m} \Gamma^k = \left\langle V, E = \coprod_{k=1}^{m} \tilde{E}^k \right\rangle.$$

Обозначим через $t_e = \tau_e(f_e)$ (здесь специально упрощаем обозначения, поскольку ввиду предыдущего раздела контекст должен восстанавливаться однозначным образом). Запишем в пространстве $t = \{t_e\}_{e \in E}$ *двойственную задачу* к (1.4.3) [23, 33, 34] (далее мы используем обозначение $\operatorname{dom} \sigma^*$ – область определения сопряженной к $\sigma$ функции):

$$\min_{f,x} \left\{ \Psi(x, f) : \; f = \Theta x, \; x \in X \right\} = -\min_{t \in \operatorname{dom} \sigma^*} \left\{ \gamma^1 \psi^1 \left( t/\gamma^1 \right) + \sum_{e \in E} \sigma_e^*(t_e) \right\}, \quad (1.4.5)$$

где

$$\sigma_e^*(t_e) = \max_{f_e}\left\{ f_e t_e - \int_0^{f_e} \tau_e(z)dz \right\},$$

$$\frac{d\sigma_e^*(t_e)}{dt_e} = \frac{d}{dt_e}\max_{f_e}\left\{ f_e t_e - \int_0^{f_e} \tau_e(z)dz \right\} = f_e:\ t_e = \tau_e(f_e),\ e \in E;$$

$$g_{p^m}^m(t) = \sum_{e^m \in \tilde{E}^m} \delta_{e^m p^m} t_{e^m} = \sum_{e^m \in E^m} \delta_{e^m p^m} t_{e^m},$$

$$g_{p^k}^k(t) = \sum_{e^k \in \tilde{E}^k} \delta_{e^k p^k} t_{e^k} - \sum_{e^k \in \bar{E}^k} \delta_{e^k p^k} \gamma^{k+1} \psi_{e^k}^{k+1}\left(t/\gamma^{k+1}\right),\ k = 1,\ldots,m-1,$$

$$\psi_{w^k}^k(t) = \ln\left( \sum_{p^k \in P_{w^k}^k} \exp\left(-g_{p^k}^k(t)\right) \right),\ k = 1,\ldots,m,$$

$$\psi^1(t) = \sum_{w^1 \in OD^1} d_{w^1}^1 \psi_{w^1}^1(t).$$

**Теорема 1.4.3.** *Имеет место сильная двойственность* (1.4.5). *Решение задачи выпуклой оптимизации* (1.4.5) $t \geq 0$ *существует и единственно. По этому решению однозначно можно восстановить решение исходной задачи* (1.4.3) (*если какой-то из* $\gamma^k \to 0+$, *то однозначность восстановления* $x$ *может потеряться*):

$$f = \Theta x = -\nabla\psi^1\left(t/\gamma^1\right),$$

$$x_{p^k}^k = d_{w^k}^k \frac{\exp\left(-g_{p^k}^k(t)/\gamma^k\right)}{\sum_{\tilde{p}^k \in P_{w^k}^k} \exp\left(-g_{\tilde{p}^k}^k(t)/\gamma^k\right)},\ p^k \in P_{w^k}^k,\ w^k \in OD^k,\ k = 1,\ldots,m. \quad (1.4.6)$$

*Верен и обратный результат. Пусть* $f = \Theta x$ *– решение задачи* (1.4.3), *тогда* $t = \left\{\tau_e(f_e)\right\}_{e \in E}$ *– единственное решение задачи* (1.4.3) (*если какой-то из* $\gamma^k \to 0+$, *то решение* $x$ *может быть не единственно, однако это никак не сказывается на возможности однозначного восстановления* $t$).

**Замечание 1.4.4.** К теореме 1.4.3 можно сделать замечания, аналогичные замечаниям 1.4.1, 1.4.2 к теореме 1.4.1. При этом оговорки, возникающие при $\gamma^k \to 0+$, частично уже были сделаны в формулировке самой теоремы. Дополним их следующим наблюдением. Слагаемое $\gamma^1\psi^1\left(t/\gamma^1\right)$ в двойственной задаче (1.4.5) имеет равномерно ограниченную константу Липшица градиента в 2-норме:

$$L_2 \le \frac{1}{\min\limits_{k=1,...,m} \gamma^k} \sum_{w^1 \in OD^1} d^1_{w^1} \max_{\breve{p}^1_{w^1} \in \breve{P}^1_{w^1}} \left|\breve{p}^1_{w^1}\right|^2 ,$$

где $\left|\breve{p}^1_{w^1}\right|$ – число ребер в пути $\breve{p}^1_{w^1}$. Эта гладкость теряется при $\gamma^k \to 0+$:

$$-\lim_{\gamma^k \to 0+} \gamma^k \psi^k_{w^k}\left(t/\gamma^k\right) = \min_{p^k \in P^k_{w^k}} g^k_{p^k}(t)$$

– длина кратчайшего пути в графе $\Gamma^k$, отвечающего корреспонденции $w^k \in OD^k$, ребра $e^k \in \bar{E}^k$ которого взвешены величинами $\gamma^{k+1}\psi^{k+1}_{e^k}\left(t/\gamma^{k+1}\right)$, которые можно понимать как «средние» затраты на $e^k \in \bar{E}^k$ (см. замечание 1.4.5). Заметим также, что в пределе (стабильной динамики) $\mu \to 0+$ (см. замечание 1.4.2) получаем

$$\sigma^*_e(t_e) = \lim_{\mu \to 0+} \max_{f_e} \left\{ f_e t_e - \int_0^{f_e} \tau^\mu_e(z)\,dz \right\} = \begin{cases} \bar{f}_e \cdot (t_e - \bar{t}_e), & t_e \ge \bar{t}_e, \\ \infty, & t_e < \bar{t}_e. \end{cases}$$

При этом $\bar{f}_e - f_e$ в точности совпадает с множителем Лагранжа к ограничению $t_e \ge \bar{t}_e$ [33, 130].

**Замечание 1.4.5.** Формулу (1.4.6) можно получить и из других соображений. Предположим, что каждый пользователь $l$ транспортной сети, использующий корреспонденцию $w^k \in OD^k$ на уровне $k$ (ребро $e^{k-1}\left(= w^k\right) \in \bar{E}^{k-1}$ на уровне $k-1$), выбирает маршрут следования $p^k \in P^k_{w^k}$ на уровне $k$, если

$$p^k = \arg\max_{q^k \in P^k_{w^k}} \left\{ -g^k_{q^k}(t) + \xi^{k,l}_{q^k} \right\},$$

где независимые случайные величины $\xi^{k,l}_{q^k}$ имеют одинаковое двойное экспоненциальное распределение, также называемое распределением Гумбеля [23, 70, 145] (см. также подраздел 1.1.4 раздела 1.1 этой главы):

$$P\left(\xi^{k,l}_{q^k} < \zeta\right) = \exp\left\{-e^{-\zeta/\gamma^k - E}\right\}.$$

Отметим также, что если взять $E \approx 0.5772$ – константа Эйлера, то

$$M\left[\xi^{k,l}_{q^k}\right] = 0, \quad D\left[\xi^{k,l}_{q^k}\right] = \left(\gamma^k\right)^2 \pi^2 / 6.$$

Распределение Гиббса (логит-распределение) (1.4.6) получается в пределе, когда число агентов на каждой корреспонденции $w^k \in OD^k$, $k = 1, ..., m$, стремится к бесконечности (случайность исчезает и описание переходит на средние величины). Полезно также в этой связи иметь в виду, что [70, 126]:

$$\gamma^k \psi^k_{w^k}\left(t/\gamma^k\right) = M_{\left\{\xi^k_{p^k}\right\}_{p^k \in P^k_{w^k}}}\left[\max_{p^k \in P^k_{w^k}}\left\{-g^k_{p^k}(t) + \xi^k_{p^k}\right\}\right].$$

Таким образом, если каждый пользователь сориентирован на вектор затрат $t$ на ребрах $E$ (одинаковый для всех пользователей) и на каждом уровне (принятия решения) пытается выбрать кратчайший путь исходя из зашумленной информации и исходя из усреднения деталей более высоких уровней (такое усреднение можно обосновывать, если, например, как в подразделе 1.4.2, ввести разный масштаб времени (частот принятия решений) на разных уровнях, а можно просто постулировать, что пользователь так действует, как это принято в моделях типа Nested Logit [70]), то такое поведение пользователей (в пределе, когда их число стремится к бесконечности) приводит к описанию распределения пользователей по путям / ребрам (1.4.6). Равновесная конфигурация характеризуется тем, что вектор $t$ породил согласно формуле (1.4.6) такой вектор $f$, что имеет место соотношение $t = \left\{\tau_e\left(f_e\right)\right\}_{e\in E}$. Поиск такого $t$ (неподвижной точки) приводит к задаче (1.4.5). Другими словами, формула (1.4.5) есть не что иное, как отражение формулы $f = -\nabla\gamma\psi^1\left(t/\gamma^1\right)$ и связи $t_e = \tau_e\left(f_e\right)$, $e \in E$. Действительно, по формуле Демьянова–Данскина–Рубинова:

$$\frac{d\sigma^*_e\left(t_e\right)}{dt_e} = \frac{d}{dt_e}\max_{f_e \geq 0}\left\{t_e f_e - \int_0^{f_e}\tau_e(z)\,dz\right\} = f_e : t_e = \tau_e\left(f_e\right).$$

В свою очередь, формула $f = -\nabla\gamma\psi^1\left(t/\gamma^1\right)$ может интерпретироваться как следствие соотношений $f = \Theta x$ и формулы распределения Гиббса (1.4.6) (logit-распределения). При такой интерпретации связь задачи (1.4.5) с logit-динамикой, порождающей стохастические равновесия, наиболее наглядна. Действительно,

$$f \xrightarrow{t_e = \tau_e\left(f_e\right)} \underbrace{t \xrightarrow{(1.4.6)} x \xrightarrow{f = \Theta x} f}_{t \xrightarrow{f = -\nabla\gamma\psi^1\left(t/\gamma^1\right)} f}.$$

Поиск неподвижной точки в этой цепочке как раз и сводится к решению задачи (1.4.5).

**Замечание 1.4.6.** Сопоставить формуле (1.4.4), теореме 1.4.2 и замечанию 1.4.3 (отвечающих случаю $\gamma^1 = ... = \gamma^m = \gamma$) вариант двойственной задачи (1.4.5) чрезвычайно просто (мы здесь опускаем соответствующие выкладки). Собственно, понять формулу (1.4.4) как раз проще не из свойств энтропии (как это было описано в подразделе 1.4.2), а с помощью обратного перехода от двойственного описания (1.4.5). Теорема 1.4.3 и замечание 1.4.5 в случае $\gamma^1 = ... = \gamma^m = \gamma$ наглядно демонстрируют отсутствие какой бы то ни было иерархии и возможность работать на одном графе с естественной интерпретацией функций затрат на путях $g_p(t)$ (без всяких «средних» оговорок).

Забегая немного вперед (подробнее см. в гл. 2), перейдем в заключение данной главы к конспективному обсуждению численных аспектов решения задачи (1.4.5). Как правило, выгоднее решать именно задачу (1.4.5), а не (1.4.3) [23]. На эту задачу удобно смотреть, как на гладкую (с липшицевым градиентом) задачу композитной оптимизации [126, 134] с евклидовой прокс-структурой (задаваемой 2-нормой). При этом даже если по ряду параметров $\gamma^k$ требуется сделать предельный переход $\gamma^k \to 0+$, то, как правило, лучше считать, что численно мы все равно решаем задачу со всеми $\gamma^k > 0$ [23]. Этого можно добиться обратным процессом: энтропийной регуляризацией прямой задачи – сглаживанием двойственной. Некоторые детали того, как именно и в каких случаях полезно сглаживать задачу (1.4.5), описаны в работе [23] (см. также [129] и раздел 2.2 гл. 2).

Композитный быстрый градиентный метод (и различные его вариации с адаптивным подбором константы Липшица градиента, универсальный метод и др. [126, 134, 144], в частности универсальный метод подобных треугольников из раздела 2.2 гл. 2 и приложения 3) обладает прямодвойственной структурой [23, 29, 31, 33, 122, 128]. Это означает, что генерируемые этим методом последовательности $\{t^i\}$ и $\{\tilde{t}^i\}$ обладают следующим свойством:

$$\gamma^1\psi^1\left(\tilde{t}^N/\gamma^1\right) + \sum_{e\in E}\sigma_e^*\left(\tilde{t}_e^N\right) -$$
$$-\min_{t\in\operatorname{dom}\sigma^*}\left\{\frac{1}{A_N}\left[\sum_{i=0}^{N} a_i\cdot\left(\gamma^1\psi^1\left(t^i/\gamma^1\right) + \left\langle\nabla\psi^1\left(t^i/\gamma^1\right), t-t^i\right\rangle\right)\right] + \sum_{e\in E}\sigma_e^*\left(t_e\right)\right\} \le \qquad (1.4.7)$$
$$\le \frac{CL_2R_2^2}{A_N},$$

где константа $C \le 10$ зависит от метода,

$$a_N \sim N\,,\; A_N = \sum_{i=0}^{N} a_i\,,\; A_N \sim N^2\,,$$

$$R_2^2 = \max\left\{\tilde{R}_2^2, \hat{R}_2^2\right\},\; \tilde{R}_2^2 = \frac{1}{2}\left\|\bar{t} - t^*\right\|_2^2,\; \hat{R}_2^2 = \frac{1}{2}\sum_{e\in E}\left(\tau_e\left(\bar{f}_e^N\right) - t_e^*\right)^2,$$

$\bar{f}^N$ определяется в теореме 1.4.4, метод стартует с $t^0 = \bar{t}$, $t^*$ – решение задачи (1.4.5).

**Теорема 1.4.4.** *Пусть задача* (1.4.5) *решается прямодвойственным методом, генерирующим последовательности* $\left\{t^i\right\}$ *и* $\left\{\tilde{t}^i\right\}$, *с оценкой скорости сходимости* (1.4.7), *тогда*

$$0 \le \left\{\gamma^1\psi^1\left(\tilde{t}^N/\gamma^1\right) + \sum_{e\in E}\sigma_e^*\left(\tilde{t}_e^N\right)\right\} + \Psi\left(\bar{x}^N, \bar{f}^N\right) \le \frac{CL_2R_2^2}{A_N},$$

*где*

$$f^i = \Theta x^i = -\nabla\psi^1\left(t^i/\gamma^1\right),\; x^i = \left\{x_{p^k}^{k,i}\right\}_{p^k\in P_{w^k}^k,\, w^k\in OD^k}^{k=1,\dots,m},$$

$$x_{p^k}^{k,i} = d_{w^k}^k\,\frac{\exp\left(-g_{p^k}^k\left(t^i\right)/\gamma^k\right)}{\sum_{\tilde{p}^k\in P_{w^k}^k}\exp\left(-g_{\tilde{p}^k}^k\left(t^i\right)/\gamma^k\right)},\; p^k \in P_{w^k}^k,\; w^k \in OD^k,\; k = 1,\dots,m,$$

$$\bar{f}^N = \frac{1}{A_N}\sum_{i=0}^{N} a_i f^i\,,\; \bar{x}^N = \frac{1}{A_N}\sum_{i=0}^{N} a_i x^i\,.$$

**Замечание 1.4.7.** В общем случае описанный выше подход представляется наиболее предпочтительным. Однако для различных специальных случаев приведенные оценки, по-видимому, можно немного улучшить [5, 29]. Подробнее об этом будет написано в следующей главе.

Приведенная теорема 1.4.4 оценивает число итераций. Но на каждой итерации необходимо считать $\nabla\psi^1\left(t/\gamma^1\right)$, а для ряда методов и $\psi^1\left(t/\gamma^1\right)$ (например, для всех адаптивных методов, настраивающихся на параметры гладкости задачи [126, 132, 134]). Подобно [34, 123, 145], можно показать (с помощью сглаженного варианта метода Форда–Беллмана), что для этого достаточно сделать $O\left(\left|O^1\right|\left|E\right|\max_{\breve{p}^1\in\breve{P}^1}\left|\breve{p}^1\right|\right)$ арифметических операций, см. также раздел 2.2 гл. 2. Однако необходимо обговорить один нюанс. Для возможности использовать сглаженный вариант метода Форда–Беллмана [34, 123, 145] необходимо предположить, что любые движе-ния по ребрам графа с учетом их ориентации являются допустимыми, т. е. множество путей, соединяющих заданные две вершины (источник и сток), – это множество всевозможных способов добраться из источника в сток по ребрам имеющегося графа с учетом их ориентации. Сделанная оговорка не сильно обременительная, поскольку нужного свойства всегда можно добиться раздутием исходного графа в несколько раз за счет введения дополнительных вершин и ребер.

В целом, хотелось бы отметить, что прием, связанный с искусственным раздутием исходного графа путем добавления новых вершин, ребер, источников, стоков, весьма полезен для ряда приложений [5, 15, 20, 33]. В частности, нередко используется введение фиктивных (с нулевыми затратами) путей-ребер, которые дают возможность ничего не делать пользователям (не перемещаться [20], не торговать [14] и т. п.), что, в свою очередь, позволяет рассматривать ситуации с нефиксированными корреспонденциями $d^1$ [14, 20]. Также популярным приемом является перенесение затрат на преодоление вершин (узлов) графа (перекрестков [20], сортировочных станций [15]) в затраты на прохождение дополнительных ребер, появившиеся при «распутывании» узлов. Но, пожалуй, наиболее важным для большинства приложений – это введение фиктивного общего источника и общего стока, соединенных дополнительными ребрами с уже имеющимися вершинами графа [5, 15, 20].

# Глава 2. Численные методы поиска равновесий в моделях распределения транспортных потоков в больших сетях

## 2.1. Неускоренные численные методы поиска равновесного распределения транспортных потоков в модели Бэкмана и модели стабильной динамики

### 2.1.1. Введение

В работах [14, 33] (см. также подразделы 1.1.4, 1.1.6, 1.1.8 раздела 1.1 гл. 1), посвященных сведению поиска равновесного распределения транспортных потоков на сетях к решению задач выпуклой оптимизации, было поставлено несколько таких задач со специальной «сетевой» структурой. Это означает, что, скажем, расчет градиента (стохастического градиента) функционала сводится к поиску кратчайших путей в графе транспортной сети. Эта специфика задач, с одной стороны, говорит о том, что в реальных приложениях размерность задачи может быть колоссально большой. Это связано с тем, что число реально используемых путей даже в планарном графе, как правило, пропорционально кубу числа вершин, а число вершин в реальных приложениях обычно не меньше тысячи, отметим, что в худших случаях число путей может расти экспоненциально с ростом числа вершин. С другой стороны, такого рода задачи имеют хорошую геометрическую интерпретацию, что позволяет эффективно снижать их размерность.

В частности, в подразделе 2.1.2 этого раздела мы описываем метод условного градиента Франка–Вульфа [72, 99, 114] поиска равновесного распределения потоков в модели Бэкмана, который на каждой итерации требует решения задачи минимизации линейной функции от потоков на ребрах в сети (число ребер порядка нескольких тысяч) на прямом произведении симплексов (симплексов столько, сколько корреспонденций). Для реальных транспортных сетей (тысяча вершин) получается задача минимизации линейной функции в пространстве размерности миллиард, поскольку она зависит от распределения потоков по путям, число которых порядка миллиарда. Ясно, что если смотреть на эту задачу формально с точки зрения оптимизации, то все сводится к полному перебору миллиарда вершин всех симплексов (причем проработка одной вершины – это расчет соответствующего скалярного произведения, то есть порядка нескольких тысяч умножений). К счастью, транспортная специфика задачи позволяет с помощью алгоритма Дейкстры и более современных подходов [74] (в том числе учитывающих «планарность» сети: A*, ALT, SHARC, Reach based routing, Highway hierarchies, Contraction hierarchies и т. п.) решать описанную задачу, делая не более десятка миллионов операций (типа умножения двух чисел с плавающей запятой), что намного быстрее. Такого рода конструкции возникают не только в связи с сетевой спецификой задачи [84], но именно для ситуаций, когда в задаче имеется сетевая структура, возможность такой редукции наиболее естественна и типична.

В подразделе 2.1.3 мы предлагаем другой способ поиска равновесия в модели Бэкмана и аналогичных моделях (в модели стабильной динамики, в промежуточных моделях). Для этого мы переходим (следуя Ю. Е. Нестерову) к двойственной задаче, в которой целевой функционал оказывается зависящим только от потоков по ребрам, а не от распределения потоков по путям. Таким образом, задача сводится к поиску равновесного распределения потоков по ребрам. При этом в ходе вычислений потоков на ребрах мы попутно (без дополнительных затрат) вычисляем порождающие их потоки по путям. Также ввиду транспортно-сетевой специфики появляется возможность содержательной интерпретации [33] (подобно интерпретации Л. В. Канторовичем цен в экономике [46]) возникающих двойственных множителей, которые в ряде приложений представляют независимый самостоятельный интерес (например, в задаче о тарифной политике грузоперевозок РЖД [14] двойственные множители – тарифы, которые и надо рассчитывать). Также сетевая структура задачи дает возможность не рассчитывать градиент целевой функции на каждой итерации заново, а пересчитывать его, используя градиент, полученный на предыдущей итерации. Грубо говоря, найдя кратчайшие пути и посчитав на их основе градиент, мы сделаем шаг по антиградиенту, немного изменив веса ребер. Ясно, что большая часть кратчайших путей при этом останется прежними, и можно специально организовать их пересчет, чтобы ускорить вычисления. Похожая философия используется в покомпонентных спусках и в современных подходах к задачам huge-scale оптимизации [125, 131]. Однако сетевая структура задачи требует переосмысления этой техники, рассчитанной изначально в основном на свойства разреженности матриц, возникающих в условии задачи.

Настоящий раздел представляет собой одну из первых попыток сочетать современные эффективные численные методы выпуклой оптимизации с сетевой структурой задачи на примере задач, пришедших из поиска равновесного распределения потоков в транспортных сетях и сетях грузовых перевозок РЖД [15, 33].

### 2.1.2. Метод Франка–Вульфа поиска равновесия в модели Бэкмана

Для удобства чтения напомним здесь наиболее популярную на протяжении более чем полувека модель равновесного распределения потоков Бэкмана [33, 34, 61, 75, 141, 148], ранее описанную нами в подразделе 1.1.4 гл. 1.

Пусть транспортная сеть города представлена ориентированным графом $\Gamma=(V,E)$, где $V$ – узлы сети (вершины), $E\subset V\times V$ – дуги сети (рёбра графа), $O\subseteq V$ – источники корреспонденций ($S=|O|$), $D\subseteq V$ – стоки. В современных моделях равновесного распределения потоков в крупном мегаполисе число узлов графа транспортной сети обычно выбирают порядка $n=|V|\sim 10^3-10^4$. Число ребер $|E|$ получается в три–четыре раза больше. Пусть $W\subseteq\{w=(i,j):i\in O, j\in D\}$ – множество корреспонденций, т. е. возможных пар *источник–сток*; $p=\{v_1,v_2,...,v_m\}$ – путь из $v_1$ в $v_m$, если $(v_k,v_{k+1})\in E$, $k=1,...,m-1$, $m>1$; $P_w$ – множество путей, отвечающих корреспонденции $w\in W$, то есть если $w=(i,j)$, то $P_w$ – множество путей, начинающихся в вершине $i$ и заканчивающихся в $j$; $P=\bigcup_{w\in W}P_w$ – совокупность всех путей в сети $\Gamma$ (число «разумных» маршрутов $|P|$, которые потенциально могут использоваться, обычно растет с ростом числа узлов сети не быстрее, чем $\mathrm{O}(n^3)$ [61, 141, 148]); $x_p$ [авт/ч] – величина потока по пути $p$, $x=\{x_p:\ p\in P\}$; $f_e$ [авт/ч] – величина потока по дуге $e$:

$$f_e(x)=\sum_{p\in P}\delta_{ep}x_p\text{, где }\delta_{ep}=\begin{cases}1, & e\in p,\\ 0, & e\notin p;\end{cases}$$

$\tau_e(f_e)$ – удельные затраты на проезд по дуге $e$. Как правило, предполагают, что это (строго) возрастающие гладкие функции от $f_e$. Точнее говоря, под $\tau_e(f_e)$ правильнее понимать представление пользователей транспортной сети об оценке собственных затрат (обычно временных в случае личного транспорта и комфортности пути (с учетом времени в пути) в случае общественного транспорта) при прохождении дуги $e$, если поток желающих проехать по этой дуге будет $f_e$.

Рассмотрим теперь $G_p(x)$ – затраты временные или финансовые на проезд по пути $p$. Естественно считать, что $G_p(x) = \sum_{e\in E} \tau_e\left(f_e(x)\right)\delta_{ep}$.

Пусть также известно, сколько перемещений в единицу времени $d_w$ осуществляется согласно корреспонденции $w \in W$. Тогда вектор $x$, характеризующий распределение потоков, должен лежать в допустимом множестве:

$$X = \left\{ x \geq 0 : \sum_{p\in P_w} x_p = d_w,\ w \in W \right\}.$$

Рассмотрим игру, в которой каждой корреспонденции $w \in W$ соответствует свой, достаточно большой ($d_w \gg 1$), набор однотипных «игроков», осуществляющих передвижение согласно корреспонденции $w$. Чистыми стратегиями игрока служат пути, а выигрышем – величина $-G_p(x)$. Игрок «выбирает» путь следования $p \in P_w$, при этом, делая выбор, он пренебрегает тем, что от его выбора также «немного» зависят $|P_w|$ компонент вектора $x$ и, следовательно, сам выигрыш $-G_p(x)$. Можно показать (см., например, [33]), что отыскание равновесия Нэша–Вардропа $x^* \in X$ (макроописание равновесия) равносильно решению задачи нелинейной комплементарности (принцип Вардропа):

$$\text{для любых } w \in W,\ p \in P_w \text{ выполняется } x_p^* \cdot \left( G_p\left(x^*\right) - \min_{q\in P_w} G_q\left(x^*\right) \right) = 0.$$

Действительно, допустим, что реализовалось какое-то другое равновесие $\tilde{x}^* \in X$, которое не удовлетворяет этому условию. Покажем, что в этом случае найдется водитель, которому выгодно поменять свой маршрут следования, тогда

$$\text{существуют такие } \tilde{w} \in W,\ \tilde{p} \in P_{\tilde{w}}, \text{ что } \tilde{x}_{\tilde{p}}^* \cdot \left( G_{\tilde{p}}\left(\tilde{x}^*\right) - \min_{q\in P_{\tilde{w}}} G_q\left(\tilde{x}^*\right) \right) > 0.$$

Каждый водитель (множество таких водителей не пусто, так как $\tilde{x}_{\tilde{p}}^* > 0$), принадлежащий корреспонденции $\tilde{w} \in W$ и использующий путь $\tilde{p} \in P_{\tilde{w}}$, действует неразумно, поскольку существует такой путь $\tilde{q} \in P_{\tilde{w}}$, $\tilde{q} \neq \tilde{p}$, что $G_{\tilde{q}}\left(\tilde{x}^*\right) = \min_{q\in P_{\tilde{w}}} G_q\left(\tilde{x}^*\right)$. Этот путь $\tilde{q}$ более выгоден, чем путь $\tilde{p}$. Аналогично показывается, что при $x^* \in X$ никому из водителей уже невыгодно отклоняться от своих стратегий.

Условие равновесия может быть переписано следующим образом [34, 61, 75, 141, 148]

$$\text{для всех } x \in X \text{ выполняется } \left\langle G\left(x^*\right), x - x^* \right\rangle \geq 0 .$$

Рассматриваемая нами игра принадлежит к классу так называемых потенциальных игр [118, 145], поскольку $\partial G_p(x)/\partial x_q = \partial G_q(x)/\partial x_p$ . Существует такая функция

$$\Psi\left(f(x)\right) = \sum_{e \in E} \int_0^{f_e(x)} \tau_e(z)\,dz = \sum_{e \in E} \sigma_e\left(f_e(x)\right),$$

где $\sigma_e(f_e) = \int_0^{f_e(x)} \tau_e(z)\,dz$ , что $\partial\Psi(x)/\partial x_p = G_p(x)$ для любого $p \in P$ . Таким образом, $x^* \in X$ – равновесие Нэша–Вардропа в этой игре тогда и только тогда, когда оно доставляет минимум $\Psi\left(f(x)\right)$ на множестве $X$ .

**Теорема 2.1.1 [33, 34, 61, 141, 148].** *Вектор* $x^*$ *будет равновесием Нэша–Вардропа тогда и только тогда, когда*

$$x \in \operatorname{Arg}\min_x \left[ \Psi\left(f(x)\right) = \sum_{e \in E} \sigma_e\left(f_e(x)\right) :\ f = \Theta x,\ x \in X \right].$$

*Если преобразование* $G(\,\cdot\,)$ *строго монотонное, то равновесие* $x$ *единственно. Если* $\tau_e'(\,\cdot\,) > 0$ *, то равновесный вектор распределения потоков по ребрам* $f$ *единственный* (*это еще не гарантирует единственность вектора распределения потоков по путям* $x$ [34]).

Итак, будем решать задачу ( $\Psi_*$ – оптимальное значение функционала):

$$\Psi(f) = \sum_{e \in E} \sigma_e(f_e) \to \min_{\substack{f = \Theta x \\ x \in X}}$$

методом условного градиента [81, 105, 108, 124] (Франка–Вульфа).

**Начальная итерация**

*Положим* $\tilde{t}_e^0 = \partial\Psi(0)/\partial f_e = \tau_e(0)$ *и рассмотрим задачу*

$$\sum_{e \in E} \tilde{t}_e^0 f_e \to \min_{\substack{f = \Theta x \\ x \in X}} .$$

*Эту задачу можно переписать, как*

$$\min_{x \in X} \sum_{e \in E} \tilde{t}_e^0 \sum_{p \in P} \delta_{ep} x_p = \sum_{w \in W} d_w \min_{p \in P_w} \left\{ \sum_{e \in E} \delta_{ep} \tilde{t}_e^0 \right\} = \sum_{w \in W} d_w T_w\left(\tilde{t}^0\right),$$

*где* $T_w\left(\tilde{t}^0\right)$ *– длина кратчайшего пути из* $i$ *в* $j$ (*где* $w=(i,j)$) *на графе, ребра которого взвешены вектором* $\tilde{t}^0=\left\{\tilde{t}_e^0\right\}_{e\in E}$. *Таким образом, выписанную задачу можно решить с учетом того, что* $n=|V|\sim|E|$*, за* $\tilde{O}(Sn)$ (*здесь и далее* $\tilde{O}(\ )=O(\ )$ *с точностью до логарифмического фактора*) *и быстрее современными вариациями алгоритма Дейкстры* [67, 74, 77, 98]. *Обозначим решение этой задачи через* $f^0$.

Можно интерпретировать ситуацию таким образом, что в начальный момент водители посчитали, что все дороги абсолютно свободны и выбрали согласно этому предположению кратчайшие пути, соответствующие их целям, и поехали по этим маршрутам (путям). На практике более равномерное распределение водителей по путям в начальный момент может оказаться более предпочтительным.

Поняв, что в действительности из-за наличия других водителей время в пути не соответствует первоначальной оценке, выраженной весами ребер $\tilde{t}_e^0$, доля $\gamma^k$ водителей (обнаруживших это и готовых что-то менять) на следующем $(k+1)$-м шаге изменит свой выбор, исходя из кратчайших путей, посчитанных по распределению водителей на предыдущем $k$-м шаге. Таким образом, возникает процедура «нащупывания» равновесия. Если выбирать специальным образом $\gamma^k$ (в частности, необходимо $\gamma^k \xrightarrow[k\to\infty]{} 0$, чтобы избежать колебания вокруг равновесия (minority game [82]), и $\sum_{k=0}^{\infty}\gamma^k=\infty$, чтобы до равновесия дойти), то система, действительно, сойдется в равновесие. Опишем теперь более формально сказанное.

**Итерации** $k=0,1,2,\ldots$

*Пусть* $f^k$ *– вектор потоков на ребрах, полученный на предыдущей итерации с номером* $k$. *Положим* $\tilde{t}_e^k=\partial\Psi\left(f^k\right)/\partial f_e=\tau_e\left(f^k\right)$ *и рассмотрим задачу*

$$\sum_{e\in E}\tilde{t}_e^k y_e \to \min_{\substack{y=\Theta x\\ x\in X}}.$$

*Так же, как и раньше, задача сводится к поиску кратчайших путей на графе, ребра которого взвешены вектором* $\tilde{t}^k=\left\{\tilde{t}_e^k\right\}_{e\in E}$.

*Обозначим решение задачи через* $y^k$. *Положим*

$$f^{k+1} = \left(1-\gamma^k\right) f^k + \gamma^k y^k, \ \gamma^k = \frac{2}{k+1}.$$

Заметим, что возникающую здесь задачу поиска кратчайших путей на графе можно попробовать решать быстрее, чем за $\tilde{O}(Sn)$. Связано это с тем, что мы уже решали на предыдущей итерации аналогичную задачу для этого же графа с близкими весами ребер [67, 74, 77, 98] (веса ребер графа с ростом $k$ меняются все слабее от шага к шагу, поскольку $\gamma^k \xrightarrow[k\to\infty]{} 0$). Тем не менее, далее в разделе мы будем считать, что одна итерация этого метода занимает $\tilde{O}(Sn)$.

Заметим также, что, решая задачи поиска кратчайших путей, мы находим (одновременно, т. е. без дополнительных затрат) не только вектор распределения потоков по ребрам $y$, но и разреженный вектор распределения потоков по путям $x$.

Строго говоря, нужно найти вектор $y^k$, а не кратчайшие пути. Чтобы получить вектор $y^k$ за $\tilde{O}(Sn)$ стоит для каждого из $S$ источников построить (например, алгоритмом Дейкстры) соответствующее дерево кратчайших путей (исходя из принципа динамического программирования: «часть кратчайшего пути сама будет кратчайшим путем» несложно понять, что получится именно дерево с корнем в рассматриваемом источнике). Это можно сделать для одного источника за $\tilde{O}(n)$. Однако, главное, правильно взвешивать ребра (их не больше $n$) такого дерева, чтобы за один проход этого дерева можно было восстановить вклад (по всем ребрам) соответствующего источника в общий вектор $y^k$. Ребро должно иметь вес, равный сумме всех проходящих через него корреспонденций с заданным источником (корнем дерева). Имея значения соответствующих корреспонденций (их также не больше $n$) за один обратный проход (то есть с листьев к корню) такого дерева, можно осуществить необходимое взвешивание (с затратами не более $O(n)$). Делается это по правилу: вес ребра равен сумме корреспонденции (возможно, равной нулю) в соответствующую вершину, в которую ребро входит, и сумме весов всех ребер (если таковые имеются), выходящих из упомянутой вершины.

**Теорема 2.1.2 [81, 105, 108, 124].** *Имеет место следующая оценка*:

$$\Psi\left(f^N\right) - \Psi_* \le \Psi\left(f^N\right) - \Psi_N \le \frac{2L_p R_p^2}{N+1}, \ f^N \in \Delta = \left\{f = \Theta x: \ x \in X\right\},$$

*где*

$$\Psi_N = \max_{k=0,\dots,N}\left\{\Psi\left(f^k\right)+\left\langle\nabla\Psi\left(f^k\right),y^k-f^k\right\rangle\right\},$$

$$R_p^2 = \max_{k=0,\dots,N}\left\|y^k-f^k\right\|_p^2 \le \max_{f,\breve{f}\in\Delta}\left\|\breve{f}-f\right\|_p^2,$$

$$L_p = \max_{\|h\|_p\le 1}\ \max_{f\in\mathrm{conv}\left(f^0,f^1,\dots,f^N\right)}\left\langle h,\mathrm{diag}\left\{\tau_e'\left(f_e\right)\right\}h\right\rangle,\ 1\le p\le\infty.$$

Из доказательства этой теоремы [81, 105, 108, 124] можно усмотреть немного более тонкий способ оценки $L_p$, в котором вместо $f\in\mathrm{conv}\left(f^0,f^1,\dots,f^N\right)$ можно брать $f\in\mathrm{conv}\left(f^0,f^1\right)\bigcup\mathrm{conv}\left(f^1,f^2\right)\bigcup\dots$ $\dots\bigcup\mathrm{conv}\left(f^{N-1},f^N\right)$. Однако для небольшого упрощения выкладок мы будем использовать приведенный в формулировке теоремы огрубленный вариант.

К сожалению, в приложениях (см. пример о расщеплении потоков на личный и общественный транспорт в подразделе 2.1.3) функции $\tau_e\left(f_e\right)$ могут иметь вертикальные асимптоты, что не позволяет равномерно по $N$ ограничить $L_p$ (даже если более тонко оценивать $L_p$). Такие случаи мы просто исключаем из рассмотрения для метода, описанного в этом разделе. Другими словами, будем считать, что функции $\tau_e\left(f_e\right)$ заданы на положительной полуоси. К таким функциям относятся, например, BPR-функции (см. подраздел 2.1.3).

Обратим внимание на то, что сам метод никак не зависит от выбора параметра $p$, от того какие получаются $R_p^2$ и $L_p$, в то время как оценка на число итераций, которые необходимо сделать для достижения заданной по функции (функционалу) точности, от этого выбора зависит. Как следствие, от этого выбора зависит и критерий останова (значение $\Psi_*$ нам априорно не известно).

Будем считать $p=2$ (сопоставимые оценки получаются и при выборе $p=\infty$):

$$L_2\left(f^0,f^1,\dots,f^N\right)=\max_{f\in\mathrm{conv}\left(f^0,f^1,\dots,f^N\right)}\ \max_{e\in E}\tau_e'\left(f_e\right)=\max_{e\in E}\tau_e'\left(\max_{k=0,\dots,N}f_e^k\right),$$

$$R_2^2=\max_{f,\breve{f}\in\Delta}\left\|\breve{f}-f\right\|_2^2.$$

Величину $R_2^2$ мы можем оценить априорно (при этом, к сожалению, получается довольно грубая оценка), т. е. можно считать её нам известной. Труднее обстоит дело с $L_2$. Далее предлагается оригинальный способ за-пуска метода Франка–Вульфа, критерий останова которого не требует априорного знания $L_2$.

Задаемся точностью $\varepsilon > 0$. Оцениваем $R_2^2$. Полагаем $L_2 = 1$ (для определенности). Запускаем метод Франка–Вульфа с $N(L_2) = 2L_2R_2^2/\varepsilon$. На каждом шаге проверяем условие (это делается за $\mathrm{O}(n)$):

$$L_2\left(f^0, f^1, ..., f^k\right) = \max_{e\in E} \tau'_e\left(\max_{l=0,...,k} f_e^l\right) \le L_2 .$$

Если на всех шагах условие выполняется, то, сделав $N(L_2)$ шагов, гарантированно получим решение с нужной точностью. Если же на каком-то шаге $k < N(L_2)$ условие нарушилось, т. е. $L_2\left(f^0, f^1, ..., f^k\right) > L_2$, то полагаем $L_2 := L_2\left(f^0, f^1, ..., f^k\right)$, пересчитываем $N(L_2)$ и переходим к следующему шагу. Таким образом, по ходу итерационного процесса мы корректируем критерий останова, оценивая необходимое число шагов по получаемой последовательности $\left\{f^k\right\}$. Специфика данной постановки, которая позволила так рассуждать, заключается в наличии явного представления

$$L_2\left(f^0, f^1, ..., f^k\right) = \max_{e\in E} \tau'_e \underbrace{\left(\max_{l=0,...,k} f_e^l\right)}_{\hat{f}_e},$$

и независимости используемого метода от выбора $L_2$ (шаг метода Франка–Вульфа $\gamma^k = 2(k+1)^{-1}$ вообще ни от каких параметров не зависит).

На практике, однако, приведенный способ работает не очень хорошо из-за использования завышенных оценок для $L_2$ и $R_2^2$. Более эффективным оказался другой способ, который использует неравенство (см. теорему 2.1.2) $\Psi\left(f^N\right) - \Psi_N \le \varepsilon$. В этом способе важно, что $\Psi_k$, $k = 0, ..., N$, автоматически рассчитываются на каждой итерации без дополнительных затрат, а $\Psi\left(f^k\right)$ может быть рассчитано на каждой итерации по известному $f^k$ за $\tilde{\mathrm{O}}(n)$. Однако нет необходимости проверять этот критерий на каждой итерации, можно это делать, например, с периодом $\alpha\tilde{\varepsilon}^{-1}$, где $\tilde{\varepsilon}$ – относительная точность по функции (скажем, $\tilde{\varepsilon} = 0.01$ – означает, что $\varepsilon = 0.01\Psi\left(f^0\right)$), $\alpha \approx 1$ подбирается эвристически, исходя из задачи. Следуя [108], можно еще немного упростить рассуждения за счет небольшого увеличения числа итераций. А именно, можно использовать оценки (при этом следует полагать $\gamma^k = 2(k+2)^{-1}$):

$$\Psi\left(f^k\right) - \Psi_* \le \left\langle \nabla\Psi\left(f^k\right), f^k - y^k \right\rangle, \ \min_{k=1,\ldots,N} \left\langle \nabla\Psi\left(f^k\right), f^k - y^k \right\rangle \le \frac{7L_2R_2^2}{N+2}.$$

Таким образом, в данном разделе был описан способ поиска равновесного распределения потоков по ребрам $f$, который за время

$$\tilde{O}\left(SnL_2R_2^2/\varepsilon\right)$$

находит такой $f^{N(\varepsilon)}$, что

$$\Psi\left(f^{N(\varepsilon)}\right) - \Psi_* \le \varepsilon.$$

Численные эксперименты на данных [158] показали высокую эффективность метода [1]. Заметим, что именно этот метод используется в большинстве коммерческих продуктах транспортного моделирования для поиска равновесного распределения потоков по путям.

### 2.1.3. Поиск равновесия в модели стабильной динамики (Нестерова–де Пальма) и смешанной модели с помощью двойственного субградиентного метода

В ряде постановок задач вместо функций затрат на ребрах $\tau_e\left(f_e\right)$ заданы ограничения на пропускные способности $f_e \le \bar{f}_e$ и затраты на прохождения свободного (не загруженного $f_e < \bar{f}_e$) ребра $\bar{t}_e$. В модели стабильной динамики это сделано для всех ребер [33, 133], а в модели грузоперевозок РЖД – только для части [15]. Согласно работе [33], такую новую модель можно получить предельным переходом из модели Бэкмана, с помощью введения внутренних штрафов в саму модель. Ранее мы уже об этом писали в разделе 1.1 гл. 1. Далее напомним основные моменты.

Будем считать, что (как и в модели Бэкмана) у всех ребер есть свои функции затрат $\tau_e^\mu\left(f_e\right)$, но для части ребер $e \in E'$ (какой именно части, зависит от задачи) осуществляется предельный переход:

$$\tau_e^\mu\left(f_e\right) \xrightarrow[\mu\to 0+]{} \begin{cases} \bar{t}_e, & 0 \le f_e < \bar{f}_e, \\ \left[\bar{t}_e, \infty\right), & f_e = \bar{f}_e, \end{cases}$$

$$d\tau_e^\mu\left(f_e\right)/df_e \xrightarrow[\mu\to 0+]{} 0, \ 0 \le f_e < \bar{f}_e.$$

Обозначив через $x(\mu)$ равновесное распределение потоков по путям в модели Бэкмана при функциях затрат на ребрах $\tau_e^\mu(f_e)$, получим, что при $e \in E'$:

$$\tau_e^\mu\left(f_e\left(x(\mu)\right)\right) \xrightarrow[\mu\to 0+]{} t_e, \quad f_e\left(x(\mu)\right) \xrightarrow[\mu\to 0+]{} f_e,$$

где пара $(t, f)$ – равновесие в модели стабильной динамики и ее вариациях [15, 33, 133] с тем же графом и матрицей корреспонденций, что и в модели Бэкмана, и с ребрами $e \in E'$, характеризующимися набором $(\bar{t}, \bar{f})$ из определения функций $\tau_e^\mu(f_e)$. Заметим, что если $t_e > \bar{t}_e$, то $t_e - \bar{t}_e$ можно интерпретировать, например, как время, потерянное в пробке на этом ребре [33, 133].

Согласно подразделу 2.1.2, равновесная конфигурация при таком переходе $\mu \to 0+$ должна находиться из решения задачи

$$\Psi(f) = \sum_{e\in E\setminus E'} \int_0^{f_e} \tau_e^\mu(z)dz + \lim_{\mu\to 0+} \sum_{e\in E'} \int_0^{f_e} \tau_e^\mu(z)dz \to \min_{f=\Theta x,\ x\in X}.$$

Считая, что в равновесии не может быть $\lim_{\mu\to 0+} \tau_e^\mu(f_e) = \infty$ (иначе равновесие просто недостижимо, и со временем весь граф превратится в одну большую пробку), можно не учитывать в интеграле вклад точек $\bar{f}_e$ (в случае попадания в промежуток интегрирования), то есть переписать задачу следующим образом:

$$\min_{f=\Theta x,\ x\in X} \left\{ \sum_{e\in E\setminus E'} \int_0^{f_e} \tau_e^\mu(z)dz + \sum_{e\in E'} \int_0^{f_e} \left(\bar{t}_e + \delta_{\bar{f}_e}(z)\right)dz \right\} \Leftrightarrow$$
$$\Leftrightarrow \min_{\substack{f=\Theta x,\ x\in X \\ f_e \le \bar{f}_e,\ e\in E'}} \left\{ \sum_{e\in E\setminus E'} \int_0^{f_e} \tau_e^\mu(z)dz + \sum_{e\in E'} f_e \bar{t}_e \right\}, \qquad (2.1.1)$$

где

$$\delta_{\bar{f}_e}(z) = \begin{cases} 0, & 0 \le z < \bar{f}_e \\ \infty, & z \ge \bar{f}_e \end{cases}, \quad e \in E'.$$

**Теорема 2.1.3 [33, 133].** *Двойственная задача к выписанной выше задаче* (2.1.1) *может быть приведена к следующему виду*:

$$\Upsilon(t) = -\sum_{w\in W} d_w T_w(t) + \left\langle \bar{f}, t - \bar{t} \right\rangle - \mu \sum_{e\in E\setminus E'} h_e^\mu(t_e) \to \min_{\substack{t_e \ge \bar{t}_e,\ e\in E' \\ t_e \in \operatorname{dom} h_e^\mu(t_e),\ e\in E\setminus E'}},$$

*где* $T_w(t)$ *– длина кратчайшего пути из* $i$ *в* $j$ ($w = (i, j) \in W$) *на графе, ребра которого взвешены вектором* $t = \{t_e\}_{e\in E}$, *а функции* $h_e^\mu(t_e)$ *гладкие и вогнутые.*

*При этом решение изначальной задачи* $f$ *можно получить из формул*

$f_e = \bar{f}_e - s_e$, $e \in E'$, где $s_e \geq 0$ *– множитель Лагранжа к ограничению*

$$t_e \geq \bar{t}_e;\ \tau_e^{\mu}\left(f_e\right) = t_e,\ e \in E \setminus E'.$$

Приведем пример модели (типа стабильной динамики) расщепления пользователей на личный и общественный транспорт [33] (см. также подраздел 1.1.8 гл. 1), в которой каждое ребро $e \in E$ изначального графа продублировано для личного (*л*) и общественного (*о*) транспорта, при этом для общественного транспорта [133]:

$$\tau_e\left(f_e^o\right) = \bar{t}_e^o \cdot \left(1 + \mu \frac{\bar{f}_e^o}{\bar{f}_e^o - f_e^o}\right),$$

а для личного транспорта был осуществлен предельный переход $\mu \to 0+$ в аналогичных формулах

$$\tau_e\left(f_e^л\right) = \bar{t}_e^л \cdot \left(1 + \mu \frac{\bar{f}_e^л}{\bar{f}_e^л - f_e^л}\right).$$

Поиск равновесного расщепления на [личный] и [общественный] транспорт приводит к следующей задаче [33]:

$$-\sum_{w \in W} d_w \min\left\{T_w^л\left(t^л\right), T_w^o\left(t^o\right)\right\} + \left\langle \bar{f}^л, t^л - \bar{t}^л \right\rangle + \left\langle \bar{f}^o, t^o - \bar{t}^o \right\rangle - \mu \sum_{e \in E} \bar{f}_e^o \cdot \bar{t}_e^o \cdot \ln\left(1 + \frac{t_e^o - \bar{t}_e^o}{\bar{t}_e^o \mu}\right) \to \min_{\substack{t^л \geq \bar{t}^л \\ t^o \geq \bar{t}^o \cdot (1-\mu)}},$$

при этом $f^л = \bar{f}^л - s^л$, где $s^л$ – вектор множителей Лагранжа для ограничений $t^л \geq \bar{t}^л$,

$$f_e^o = \bar{f}_e^o \cdot \left(1 - \frac{\bar{t}_e^o \cdot \mu}{t_e^o - (1-\mu)\bar{t}_e^o}\right).$$

Другой способ выбора функции затрат на ребрах для общественного транспорта [142]:

$$\tau_e^{\mu}\left(f_e\right) = \bar{t}_e \cdot \left(1 + \gamma \cdot \left(\frac{f_e}{\bar{f}_e}\right)^{\frac{1}{\mu}}\right).$$

Упростим немного обозначения, опуская индексы $\mu$, $o$, $л$ и используя для краткости обозначение $\sigma_e\left(f_e\right) = \int_0^{f_e} \tau_e(z)\,dz = \int_0^{f_e} \tau_e^{\mu}(z)\,dz$. Итак, согласно теореме 2.1.3 для задачи (2.1.1) можно построить следующую двойственную задачу [27, 33] (см. также разделы 1.1, 1.5 гл. 1):

$$\Upsilon(t) = \underbrace{-\sum_{w\in W} d_w T_w(t)}_{F(t)} + \sum_{e\in E} \sigma_e^*(t_e) \to \min_{\substack{t_e \ge \bar{t}_e,\, e\in E' \\ t_e \in \operatorname{dom}\sigma_e^*(t_e),\, e\in E\setminus E'}}, \qquad (2.1.2)$$

где $T_w(t) = \min_{p\in P_w} \sum_{e\in E} \delta_{ep} t_e$ – длина кратчайшего пути из $i$ в $j$ ($w=(i,j)\in W$) на графе $\Gamma$, ребра которого взвешены вектором $t=\{t_e\}_{e\in E}$. При этом решение задачи (2.1.1) $f$ можно получить из формул: $f_e = \bar{f}_e - s_e$, $e\in E'$, где $s_e \ge 0$ – множитель Лагранжа к ограничению $t_e \ge \bar{t}_e$; $\tau_e(f_e) = t_e$, $e\in E\setminus E'$. Заметим, что для ребер $e\in E'$ имеем $\sigma_e^*(t_e) = \bar{f}_e \cdot (t_e - \bar{t}_e)$, а для BPR-функций [140]:

$$\tau_e(f_e) = \bar{t}_e \cdot \left(1 + \gamma \cdot \left(\frac{f_e}{\bar{f}_e}\right)^{\frac{1}{\mu}}\right) \Rightarrow \sigma_e^*(t_e) = \bar{f}_e \cdot \left(\frac{t_e - \bar{t}_e}{\bar{t}_e \cdot \gamma}\right)^{\mu} \frac{(t_e - \bar{t}_e)}{1+\mu}.$$

В приложениях обычно выбирают $\mu = 1/4$. В этом случае приводимый ниже шаг итерационного метода (2.1.3) может быть осуществлен по явным формулам, поскольку существуют квадратурные формулы (формулы Кардано–Декарта–Эйлера–Феррари [51]) для уравнений 4-й степени.

Поиск вектора $t$ представляет самостоятельный интерес, поскольку этот вектор описывает затраты на ребрах графа транспортной сети. Решение задачи (2.1.2) дает вектор затрат $t$ в равновесии.

Для решения двойственной задачи (2.1.2) воспользуемся методом зеркального спуска (МЗС) в композитном варианте [94, 120] ($k = 0,\dots,N$, $t^0 = \bar{t}$, ограничение $t_e \in \operatorname{dom}\sigma_e^*(t_e)$, $e\in E\setminus E'$ всегда будет не активным, т. е. его можно не учитывать), который в данном случае будет представлять собой двойственный композитный метод проекции субградиента:

$$t^{k+1} = \arg \min_{\substack{t_e \ge \bar{t}_e,\, e\in E' \\ t_e \in \operatorname{dom}\sigma_e^*(t_e),\, e\in E\setminus E'}} \left\{ \gamma_k \left\{ \left\langle \partial F(t^k), t - t^k \right\rangle + \sum_{e\in E} \sigma_e^*(t_e) \right\} + \frac{1}{2}\left\| t - t^k \right\|_2^2 \right\}, \qquad (2.1.3)$$

где $\partial F(t^k)$ – произвольный элемент субдифференциала выпуклой функции $F(t^k)$ в точке $t^k$; $\gamma_k = \varepsilon / M_k^2$, $M_k = \left\|\partial F(t^k)\right\|_2$, если $E' = E$ и $M_k \equiv M$, $\max_{k=0,\dots,N} \left\|\partial F(t^k)\right\|_2 \le M$ иначе, где $\varepsilon > 0$ – желаемая точность решения задач (2.1.1) и (2.1.2), см. (2.1.6).

Положим (следует сравнить этот подход с подходом из раздела 1.6 гл. 1):

$$\bar{t}^N = \frac{1}{S_N}\sum_{k=0}^{N}\gamma_k t^k \,,\; S_N = \sum_{k=0}^{N}\gamma_k \,,$$

$$f_e^k \in -\partial_e F\left(t^k\right),\;\; \bar{f}_e^N = \frac{1}{S_N}\sum_{k=0}^{N}\gamma_k f_e^k \,,\; e \in E \setminus E' \,, \tag{2.1.4}$$

$$\bar{f}_e^N = \bar{f}_e - s_e^N \,,\; e \in E' \,, \tag{2.1.4'}$$

где $s_e^N$ есть множитель Лагранжа к ограничению $t_e \ge \bar{t}_e$ в задаче

$$\frac{1}{S_N}\left\{\sum_{k=0}^{N}\gamma_k \cdot \left\{\sum_{e\in E'}\partial_e F\left(t^k\right)\cdot\left(t_e - t_e^k\right)\right\} + S_N\sum_{e\in E'}\bar{f}_e\cdot\left(t_e - \bar{t}_e\right) + \frac{1}{2}\sum_{e\in E'}\left(t_e - \bar{t}_e\right)^2\right\} \to \min_{t_e \ge \bar{t}_e,\, e\in E'} .$$

Критерий останова метода (правое неравенство)

$$\left(0 \le\right)\Upsilon\left(\bar{t}^N\right) + \Psi\left(\bar{f}^N\right) \le \varepsilon \,. \tag{2.1.5}$$

**Теорема 2.1.4.** *Пусть*

$$\tilde{M}_N^2 = \left(\frac{1}{N+1}\sum_{k=0}^{N}M_k^{-2}\right)^{-1},$$

$$R_N^2 := \frac{1}{2}\sum_{e\in E\setminus E'}\left(\tau_e\left(\bar{f}_e^N\right) - \bar{t}_e\right)^2 + \frac{1}{2}\sum_{e\in E'}\left(\tilde{t}_e^N - \bar{t}_e\right)^2 ,$$

$$\left\{\tilde{t}_e^N\right\}_{e\in E'} = \arg\min_{\left\{t_e\right\}_{e\in E'}\ge 0}\Bigg\{\underbrace{-\sum_{w\in W} d_w T_w\left(\left\{\tau_e\left(\bar{f}_e^N\right)\right\}_{e\in E\setminus E'},\left\{t_e\right\}_{e\in E'}\right)}_{F\left(\left\{\tau_e\left(\bar{f}_e^N\right)\right\}_{e\in E\setminus E'},\left\{t_e\right\}_{e\in E'}\right)} + \sum_{e\in E'}\bar{f}_e^N\cdot\left(t_e - \bar{t}_e\right)\Bigg\}.$$

*Тогда при* $N \ge \frac{2\tilde{M}_N^2 R_N^2}{\varepsilon^2}$ *имеет место неравенство* (2.1.5) *и, как следствие,*

$$0 \le \Upsilon\left(\bar{t}^N\right) - \Upsilon_* \le \varepsilon \,,\; 0 \le \Psi\left(\bar{f}^N\right) - \Psi_* \le \varepsilon \,. \tag{2.1.6}$$

**Доказательство.** Из работ [94, 128] имеем (также в выкладках используется то, что $d\sigma_e^*\left(t_e\right)/dt_e = f_e \;\Leftrightarrow\; t_e = \tau_e\left(f_e\right)$, $e \in E \setminus E'$; $\sigma_e^*\left(t_e\right) \ge 0$, $\sigma_e^*\left(t_e^0\right) = \sigma_e^*\left(\bar{t}_e\right) = 0$, $e \in E$):

$$\Upsilon\left(\bar{t}^N\right) \le \frac{1}{S_N}\sum_{k=0}^{N}\gamma_k\Upsilon\left(t^k\right) \le \frac{1}{2S_N}\sum_{k=0}^{N}\gamma_k^2 M_k^2 +$$

$$+\frac{1}{S_N}\min_{t_e\ge\bar{t}_e,\, e\in E'}\left\{\sum_{k=0}^{N}\gamma_k\left\{F\left(t^k\right) + \left\langle\partial F\left(t^k\right), t - t^k\right\rangle\right\} + S_N\sum_{e\in E}\sigma_e^*\left(t_e\right) + \frac{1}{2}\left\|t - t^0\right\|_2^2\right\} \le \frac{\varepsilon}{2} +$$

$$+\min_t\left\{\frac{1}{S_N}\left\{\sum_{k=0}^{N}\gamma_k\left\{F\left(t^k\right)+\left\langle\partial F\left(t^k\right),t-t^k\right\rangle\right\}+S_N\sum_{e\in E}\sigma_e^*\left(t_e\right)+\frac{1}{2}\left\|t-\bar{t}\right\|_2^2\right\}+\right.$$

$$\left.+\sum_{e\in E'}s_e^N\cdot\left(\bar{t}_e-t_e\right)\right\}\le\frac{\varepsilon}{2}+\min_{\{t_e\}_{e\in E'}}\min_{\{t_e\}_{e\in E\setminus E'}}\frac{1}{S_N}\left\{\sum_{k=0}^{N}\gamma_k\left\{F\left(t^k\right)+\left\langle\partial F\left(t^k\right),t-t^k\right\rangle\right\}+\right.$$

$$\left.+\left(S_N\sum_{e\in E\setminus E'}\sigma_e^*\left(t_e\right)+\frac{1}{2}\sum_{e\in E\setminus E'}\left(t_e-\bar{t}_e\right)^2\right)+\left(S_N\sum_{e\in E'}\bar{f}_e^N\cdot\left(t_e-\bar{t}_e\right)+\frac{1}{2}\sum_{e\in E'}\left(t_e-\bar{t}_e\right)^2\right)\right\}\le$$

$$\le\frac{\varepsilon}{2}+\min_{\{t_e\}_{e\in E'}\ge 0}\frac{1}{S_N}\left\{\sum_{k=0}^{N}\gamma_k\overbrace{\underbrace{\left\{F\left(t^k\right)+\left\langle\partial F\left(t^k\right),\breve{t}^N-t^k\right\rangle\right\}}_{\le F\left(\left\{\tau_e\left(\bar{f}_e^N\right)\right\}_{e\in E\setminus E'},\{t_e\}_{e\in E'}\right)}}^{\breve{t}^N=\left(\left\{\tau_e\left(\bar{f}_e^N\right)\right\}_{e\in E\setminus E'},\{t_e\}_{e\in E'}\right)}+\right.$$

$$+\left(S_N\sum_{e\in E\setminus E'}\sigma_e^*\left(\tau_e\left(\bar{f}_e^N\right)\right)+\frac{1}{2}\sum_{e\in E\setminus E'}\left(\tau_e\left(\bar{f}_e^N\right)-\bar{t}_e\right)^2\right)+$$

$$\left.+\left(S_N\sum_{e\in E'}\bar{f}_e^N\cdot\left(t_e-\bar{t}_e\right)+\frac{1}{2}\sum_{e\in E'}\left(t_e-\bar{t}_e\right)^2\right)\right\}\le\frac{\varepsilon}{2}+\min_{\{t_e\}_{e\in E'}\ge 0}\left\{F\left(\left\{\tau_e\left(\bar{f}_e^N\right)\right\}_{e\in E\setminus E'},\{t_e\}_{e\in E'}\right)+\right.$$

$$\left.+\sum_{e\in E'}\bar{f}_e^N\cdot\left(t_e-\bar{t}_e\right)+\frac{1}{2S_N}\sum_{e\in E'}\left(t_e-\bar{t}_e\right)^2\right\}+$$

$$+\sum_{e\in E\setminus E'}\sigma_e^*\left(\tau_e\left(\bar{f}_e^N\right)\right)+\frac{1}{2S_N}\sum_{e\in E\setminus E'}\left(\tau_e\left(\bar{f}_e^N\right)-\bar{t}_e\right)^2\le$$

$$\le\frac{\varepsilon}{2}+\min_{\{t_e\}_{e\in E'}\ge 0}\left\{F\left(\left\{\tau_e\left(\bar{f}_e^N\right)\right\}_{e\in E\setminus E'},\{t_e\}_{e\in E'}\right)+\sum_{e\in E\setminus E'}\sigma_e^*\left(\tau_e\left(\bar{f}_e^N\right)\right)+\sum_{e\in E'}\bar{f}_e^N\cdot\left(t_e-\bar{t}_e\right)\right\}+$$

$$+\frac{1}{2S_N}\sum_{e\in E'}\left(\tilde{t}_e^N-\bar{t}_e\right)^2+\frac{1}{2S_N}\sum_{e\in E\setminus E'}\left(\tau_e\left(\bar{f}_e^N\right)-\bar{t}_e\right)^2=$$

$$=\frac{\varepsilon}{2}-\Psi\left(\bar{f}^N\right)+\frac{R_N^2}{S_N}=\frac{\varepsilon}{2}+\frac{\tilde{M}_N^2R_N^2}{\varepsilon\cdot(N+1)}-\Psi\left(\bar{f}^N\right)\le\varepsilon-\Psi\left(\bar{f}^N\right).\ \blacksquare$$

**Следствие 2.1.1.** *Пусть* $t^*$ *– решение задачи* (2.1.2). *Положим*

$$R^2=\frac{1}{2}\left\|t^*-t^0\right\|_2^2=\frac{1}{2}\left\|t_*-\bar{t}\right\|_2^2.$$

*Тогда при* $N=\dfrac{2\tilde{M}_N^2R^2}{\varepsilon^2}$ *справедливы неравенства*

$$\frac{1}{2}\left\|t-t^{k+1}\right\|_2^2\le 2R^2,\ k=0,\dots,N \qquad (2.1.7)$$

$$0 \le \Upsilon\left(\bar{t}^N\right) - \Upsilon_* \le \varepsilon. \qquad (2.1.8)$$

**Доказательство.** Формула (2.1.8) – стандартный результат, см., например, [120]. Формула (2.1.7) также является достаточно стандартной [128], однако далее приводится схема ее вывода. Из доказательства теоремы 3.3.1 имеем для любого $k = 0,...,N$ :

$$0 \le \frac{1}{2}\sum_{l=0}^{k}\gamma_l^2 M_l^2 + \frac{1}{2}\left\|t - t^0\right\|_2^2 - \frac{1}{2}\left\|t - t^{k+1}\right\|_2^2.$$

Отсюда следует, что

$$\frac{1}{2}\left\|t - t^{k+1}\right\|_2^2 \le R^2 + \frac{1}{2}\sum_{l=0}^{k}\varepsilon^2 M_l^{-2} \le R^2 + \frac{1}{2}\sum_{k=0}^{N}\varepsilon^2 M_k^{-2} = 2R^2. \ \blacksquare$$

**Замечание 2.1.1.** Преимуществом подхода (2.1.3), (2.1.4), (2.1.4') над подходом п. 3 работы [27] являются простота описания (отсутствие необходимости делать рестарты по неизвестным параметрам) и наличие эффективно проверяемого критерия останова (2.1.5). К недостаткам стоит отнести вхождение в оценку скорости сходимости плохо контролируемого $R_N^2$, которое может оказаться большим даже в случае, когда $E = E'$. Далее мы опишем другой способ (см. также работы [2, 124], в которых описаны близкие конструкции) восстановления решения прямой задачи (2.1.1), отличный от (2.1.4), (2.1.4'), в части (2.1.4'), который в случае $E = E'$ позволяет использовать $R^2$ вместо $R_N^2$.

Положим (следует сравнить этот подход с подходом из предыдущего раздела):

$$f_e^k \in -\partial_e F\left(t^k\right), \quad \bar{f}_e^N = \frac{1}{S_N}\sum_{k=0}^{N}\gamma_k f_e^k, \ e \in E. \qquad (2.1.9)$$

$$\tilde{R}^2 = \frac{1}{2}\sum_{e\in E'}\left(t_e^* - t_e^0\right)^2 = \frac{1}{2}\sum_{e\in E'}\left(t_e^* - \bar{t}_e\right)^2.$$

**Теорема 2.1.5.** *Пусть*

$$\tilde{R}_N^2 := \frac{1}{2}\sum_{e\in E\setminus E'}\left(\tau_e\left(\bar{f}_e^N\right) - \bar{t}_e\right)^2 + 5\tilde{R}^2. \qquad (2.1.10)$$

*Тогда при* $N \ge \dfrac{4\tilde{M}_N^2\tilde{R}_N^2}{\varepsilon^2}$ *имеют место неравенства*

$$\left|\Upsilon\left(\bar{t}^N\right) - \Upsilon_*\right| \le \varepsilon, \ \left|\Psi\left(\bar{f}^N\right) - \Psi_*\right| \le \varepsilon.$$

*Более того, также имеют место неравенства*

$$\sqrt{\sum_{e\in E'}\left(\left(\bar{f}_e^N-\bar{f}_e\right)_+\right)^2}\le\tilde{\varepsilon},\ \tilde{\varepsilon}=\varepsilon/\tilde{R}, \qquad (2.1.11)$$

$$\Psi\left(\bar{f}^N\right)-\Psi_*\le\Upsilon\left(\bar{t}^N\right)+\Psi\left(\bar{f}^N\right)\le\varepsilon,$$

*которые можно использовать для критерия останова метода* (*задавшись парой* $(\varepsilon,\tilde{\varepsilon})$).

**Доказательство.** Рассуждая аналогично доказательству теоремы 2.1.4 (см. также [2, 124]), получим

$$\Upsilon\left(\bar{t}^N\right)\le\frac{1}{S_N}\sum_{k=0}^N\gamma_k\Upsilon\left(t^k\right)\overset{(12)}{\le}\frac{1}{2S_N}\sum_{k=0}^N\gamma_k^2M_k^2+$$

$$+\frac{1}{S_N}\min_{\substack{t_e,\,e\in E\setminus E';\,t_e\ge\bar{t}_e,\,e\in E'\\ \frac{1}{2}\sum_{e\in E'}(t_e-\bar{t}_e)^2\le 5\tilde{R}^2}}\left\{\sum_{k=0}^N\gamma_k\left\{F\left(t^k\right)+\left\langle\partial F\left(t^k\right),t-t^k\right\rangle\right\}+S_N\sum_{e\in E}\sigma_e^*\left(t_e\right)\right\}+\frac{\tilde{R}_N^2}{S_N}\overset{(13)}{\le}$$

$$\overset{(13)}{\le}\frac{\varepsilon}{2}-\Psi\left(\bar{f}^N\right)-\max_{\substack{t_e\ge\bar{t}_e,\,e\in E'\\ \frac{1}{2}\sum_{e\in E'}(t_e-\bar{t}_e)^2\le 5\tilde{R}^2}}\left\{\frac{1}{S_N}\sum_{k=0}^N\gamma_k\sum_{e\in E'}\left(f_e^k-\bar{f}_e\right)\left(t_e-\bar{t}_e\right)\right\}+\frac{\tilde{R}_N^2}{S_N}\le$$

$$\le\frac{\varepsilon}{2}-\Psi\left(\bar{f}^N\right)-3\tilde{R}\sqrt{\sum_{e\in E'}\left(\left(\bar{f}_e^N-\bar{f}_e\right)_+\right)^2}+\frac{\tilde{R}_N^2}{S_N}\le\varepsilon-\Psi\left(\bar{f}^N\right)-3\tilde{R}\sqrt{\sum_{e\in E'}\left(\left(\bar{f}_e^N-\bar{f}_e\right)_+\right)^2}.$$

Неравенство (2.1.12) было получено на базе следующего соотношения:

$$\left\{\tau_e\left(\bar{f}_e^N\right)\right\}_{e\in E\setminus E'}=\arg\min_{t_e,\,e\in E\setminus E'}\left\{\frac{1}{S_N}\sum_{k=0}^N\gamma_k\left\{F\left(t^k\right)+\left\langle\partial F\left(t^k\right),t-t^k\right\rangle\right\}+\sum_{e\in E}\sigma_e^*\left(t_e\right)\right\}.$$

Неравенство (2.1.13) было получено на базе следующих соотношений:

$$\partial F\left(t^k\right)=-f^k,\ F\left(t^k\right)=-\left\langle f^k,t^k\right\rangle,$$

$$\min_{t_e}\left\{-f_e^kt_e+\sigma_e^*\left(t_e\right)\right\}=-\sigma_e\left(f_e^k\right),\ e\in E\setminus E',$$

$$-f_e^kt_e+\sigma_e^*\left(t_e\right)=-\left(f_e^k-\bar{f}_e\right)\left(t_e-\bar{t}_e\right)-f_e^k\bar{t}_e=-\left(f_e^k-\bar{f}_e\right)\left(t_e-\bar{t}_e\right)-\sigma_e\left(f_e^k\right),\ e\in E',$$

$$-\frac{1}{S_N}\sum_{k=0}^N\gamma_k\sum_{e\in E}\sigma_e\left(f_e^k\right)\le-\Psi\left(\bar{f}^N\right).$$

Таким образом,

$$\Upsilon\left(\bar{t}^N\right)+\Psi\left(\bar{f}^N\right)+3\tilde{R}\sqrt{\sum_{e\in E'}\left(\left(\bar{f}_e^N-\bar{f}_e\right)_+\right)^2}\le\varepsilon.$$

Повторяя рассуждения п. 6.11 [122] и п. 3 [2], получим искомые неравенства. ■

**Замечание 2.1.2.** О преимуществе использования формулы (2.1.9) вместо (2.1.4), (2.1.4’) в случае $E = E'$ написано в замечании 2.1.1. Сейчас отметим возникающие при таком подходе недостатки: 1) возможность нарушения ограничения $f_e \le \bar{f}_e$, $e \in E'$ в прямой задаче, 2) отсутствие левых неравенств в двойных неравенствах (2.1.5), (2.1.6).

**Замечание 2.1.3.** Формулы (2.1.4), (2.1.9) вынужденно (в случае (2.1.4)) или осознанно (в случае (2.1.9) при $e \in E'$) восстанавливают решение прямой задачи исходя из «модели» – явной формулы, связывающей прямые и двойственные переменные. Наличие таких переменных неизбежно приводит к возникновению в оценках зазора двойственности трудно контролируемых размеров решений вспомогательных задач. При этом в случае, когда наличие модели сопряжено с каким-то ограничением в прямой задаче (формулы (2.1.9) при $e \in E'$ и ограничение в прямой задаче $f_e \le \bar{f}_e$, $e \in E'$), допускается нарушение этих ограничений, которые необходимо контролировать (2.1.11). Зато по этим переменным имеется полный контроль соответствующих этим переменным частей оценок зазора двойственности (2.1.10) (см. также предыдущий раздел). Подход (2.1.4’), связанный с наличием ограничений в решаемой задаче ($t_e \ge \bar{t}_e$, $e \in E'$), также приводит к возникновению в оценках зазора двойственности трудно контролируемых размеров решений вспомогательных задач, однако уже не приводит к нарушению никаких ограничений в самой задаче и сопряженной к ней (в нашем случае исходная задача – двойственная (2.1.2), а сопряженная к ней – прямая (2.1.1)). Эти два прямодвойственных подхода дополняются другим прямодвойственным подходом, в котором шаги осуществляются «по функционалу», если не нарушены или слабо нарушены ограничения в задаче, и «по нарушенному ограничению» в противном случае. Подробнее об этом см., например, в работах [122, 137] и следующем разделе.

**Замечание 2.1.4.** Оба описанных подхода можно распространить, сохраняя вид формул восстановления и структуру рассуждений, на практически произвольные пары прямая / двойственная задача, поскольку выбранный нами пример пары взаимно сопряженных задач (2.1.1), (2.1.2) и так содержал в себе практически все основные нюансы, которые могут возникать при таких рассуждениях. Более того, вместо МЗС можно было бы использовать любой другой прямодвойственный метод. Например, композитный универсальный градиентный метод Ю. Е. Нестерова. В частности, варианты этого метода из работы [37]. Подробнее об этом будет написано в следующем разделе 2.2.

**Замечание 2.1.5.** Возможность эффективного распараллеливания предложенных методов связана с возможностью эффективного вычисле-ния самой затратной части шага описанного итерационного метода: расчет элемента субдифференциала $\partial F\left(t^k\right)$ (см. формулы (2.1.3), (2.1.4), (2.1.9), в которых этот субдифференциал используется):

$$\partial F(t) = -\sum_{i\in O}\sum_{j\in D:\,(i,j)\in W} d_{ij}\partial T_{ij}(t).$$

Вычисление $\left\{\partial T_{ij}(t)\right\}_{j\in D:\,(i,j)\in W}$ может быть осуществлено алгоритмом Дейкстры [68] (и его более современными аналогами [74]) за $\mathrm{O}(n\ln n)$. При этом под $\partial T_{ij}(t)$ можно понимать описание произвольного (если их несколько) кратчайшего путь из вершины $i$ в вершину $j$ на графе $\Gamma$, ребра которого взвешены вектором $t = \{t_e\}_{e\in E}$. Под «описанием» понимается $\left[\partial T_{ij}(t)\right]_e = 1$, если $e$ попало в кратчайший путь и $\left[\partial T_{ij}(t)\right]_e = 0$ иначе. Таким образом, вычисление $\partial F(t)$ может быть распараллелено на $S$ процессорах.

**Замечание 2.1.6.** Строго говоря, нужно найти вектор $\partial F(t)$, а не кратчайшие пути. Чтобы получить вектор $\partial F(t)$ за $\mathrm{O}(Sn\ln n)$, необходимо для каждого из $S$ источников построить (например, алгоритмом Дейкстры) соответствующее дерево кратчайших путей. Исходя из принципа динамического программирования: «часть кратчайшего пути сама будет кратчайшим путем» несложно понять, что получится именно дерево с корнем в рассматриваемом источнике. Это можно сделать для одного источника за $\mathrm{O}(n\ln n)$ [68, 74]. Однако, главное, правильно взвешивать ребра (их не больше $\mathrm{O}(n)$) такого дерева, чтобы за один проход этого дерева можно было восстановить вклад (по всем ребрам) соответствующего источника в общий вектор $\partial F(t)$. Ребро должно иметь вес, равный сумме всех проходящих через него корреспонденций с заданным источником (корнем дерева). Имея значения соответствующих корреспонденций (их также не больше $\mathrm{O}(n)$), за один обратный проход (то есть с листьев к корню) такого дерева можно осуществить необходимое взвешивание (с затратами не более $\mathrm{O}(n)$). Делается это по правилу: вес ребра равен сумме корреспонденции (возможно, равной нулю) в соответствующую вершину, в которую ребро входит, и сумме весов всех ребер (если таковые имеются), выходящих из упомянутой вершины.

Численные эксперименты подтвердили справедливость полученных теорем, однако по скорости сходимости в случае, когда $E' = \varnothing$ метод условного градиента из подраздела 2.1.2 существенно обыгрывает методы, описанные в данном подразделе. Более подробный сравнительный численный описанный метод будет приведен в следующем разделе.

## 2.2. Универсальный ускоренный метод поиска равновесий и стохастических равновесий в транспортных сетях

### 2.2.1. Введение

Многочисленные эксперименты показали, что для поиска обычного (не стохастического) равновесия в модели Бэкмана наилучшим методом (по имеющимся на данный момент представлениям) является классический метод условного градиента Франка–Вульфа (см. подраздел 2.1.2). В данном разделе будет продемонстрировано, что для всех остальных случаев[19] разумно строить двойственную задачу[20] и решать ее универсальным ускоренным градиентным методом Ю. Е. Нестерова (см., например, приложение 3). Затем, пользуясь прямодвойственностью этого метода, можно восстанавливать решение исходной (прямой задачи). Интересно отметить, что численные эксперименты (см. подраздел 2.2.7), проведенные с универсальным градиентным методом для поиска обычных (не стохастических) равновесий, показали, что время работы метода пропорционально $\sim \tilde{\varepsilon}^{-1}$, где $\tilde{\varepsilon}$ – желаемая относительная точность (по функции) решения задачи, а не $\sim \tilde{\varepsilon}^{-2}$, как это можно было ожидать, исходя из теории решения негладких задач выпуклой оптимизации [127]. В отличие от специальных примеров, собранных в [129, 132], на которых наблюдался аналогичный эффект, в данном разделе приводится, по-видимому, первый совершенно реальный такой пример с существенно негладкой функцией, имеющей огромное число всевозможных изломов.

В подразделе 2.2.2 приводятся экстремальные принципы, описывающие равновесную конфигурацию в транспортных сетях в четырех рассматриваемых случаях. Эти принципы формулируются как сепарабельные задачи выпуклой оптимизации (с композитом вида энтропии в случае поиска стохастических равновесий) при аффинных ограничениях.

---

[19] Трех: поиск стохастических равновесий в модели Бэкмана и поиск обычных и стохастических равновесий в модели стабильной динамики.

[20] Эта задача оказывается негладкой при поиске обычных (не стохастических) равновесий.

В подразделе 2.2.3 с помощью достаточно стандартной техники выпуклого анализа строятся двойственные задачи, чтобы не проектироваться на аффинные ограничения при решении прямых задач (последняя операция является дорогой).

В подразделе 2.2.4, следуя [37, 127], излагается специальный вариант универсального градиентного метода – универсальный метод подобных треугольников (УМПТ), в котором (в отличие от [132]) используется только одна проекция на каждой итерации. Особое внимание уделяется прямодвойственности этого метода. Это нашло отражение в необычной форме записи того, как сходится метод.

В подразделе 2.2.5 УМПТ применяется для решения двойственных задач, описанных в подразделе 2.2.6. Приводятся формулы восстановления решения прямой задачи. Формулы восстановления решения задачи поиска равновесия в модели стабильной динамики, исходя из последовательности, сгенерированной УМПТ при решении двойственной задачи, по-видимому, приводятся впервые.

В подразделе 2.2.5 внимание сконцентрировано на том, как восстанавливать решения прямых задач, т. е. задач подраздела 2.2.2, исходя из приближенного решения двойственных задач, т. е. задач из подраздела 2.2.3. В подразделе 2.2.6 результаты подраздела 2.2.5 пополняются исследованием того, сколько итераций будет делать УМПТ (оценка констант метода) для достижения желаемой точности решения (прямой и двойственной задачи одновременно) и какая при этом будет стоимость одной итерации. Последняя задача сводится (при поиске стохастических равновесий) к вычислению значений и градиентов характеристических функций на транспортном графе [34, 123] (для вычисления градиентов используется теория быстрого автоматического дифференцирования [44]) или к поиску кратчайших путей в транспортном графе [68] (при поиске обычных равновесий).

В подразделе 2.2.7 результаты всех предыдущих подразделов собираются вместе. Формулируется теорема, в которой описаны полученные в данном разделе (верхние) теоретические оценки скорости сходимости УМПТ во всех четырех случаях. Приводятся результаты численных экспериментов. Эти результаты сопоставляются с ранее известными. Далее в изложении мы будем следовать статье [6]. По этой же статье можно отследить и ссылки на основную литературу.

## 2.2.2. Экстремальные (вариационные) принципы для поиска равновесий в транспортных сетях

Рассмотрим транспортную сеть, которую будем представлять ориентированным графом $\langle V, E\rangle$, где $V$ – множество вершин (как правило, можно считать, что $|E|/4 \le |V| \le |E|$), а $E$ – множество ребер, $|E| = n$. Обозначим множество пар $w = (i, j)$ источник–сток через $OD$, $d_w$ – корреспонденция, отвечающая паре $w$, $x_p$ – поток по пути $p$; $P_w$ – множество путей, отвечающих корреспонденции $w$ (начинающихся в $i$ и заканчивающихся в $j$), $P = \bigcup_{w \in OD} P_w$ – множество всех путей. Затраты на прохождения ребра $e \in E$ описываются функцией $\tau_e(f_e)$, где $f_e$ – поток по ребру $e$.

Опишем марковскую логит-динамику (также говорят гиббсовскую динамику) в повторяющейся игре загрузки графа транспортной сети. Пусть каждой корреспонденции отвечает $d_w M$ агентов ($M \gg 1$), $\tau_e(f_e) := \tau_e(f_e/M)$. Пусть имеется $TN$ шагов ($N \gg 1$). Каждый агент независимо от остальных на шаге $t+1$ выбирает с вероятностью

$$\frac{\lambda}{N} \frac{\exp\left(-G_p^t/\gamma\right)}{\sum_{q \in P_w} \exp\left(-G_q^t/\gamma\right)}$$

путь $p \in P_w$, где $G_p^t$ – затраты на пути $p$ на шаге $t$ ($G_p^0 \equiv 0$), а с вероятность $1 - \lambda/N$ – путь, который использовал на шаге $t$. Такая динамика отражает ограниченную рациональность агентов и часто используется в популяционной теории игр [145] и теории дискретного выбора [70]. Оказывается, эта марковская динамика в пределе $N \to \infty$ превращается в марковскую динамику в непрерывном времени (вырождающуюся при $\gamma \to 0+$ в динамику наилучших ответов), которая в свою очередь допускает два предельных перехода (обоснование перестановочности этих пределов см. в [95]): $T \to \infty$, $M \to \infty$ или $M \to \infty$, $T \to \infty$. При первом порядке переходов мы сначала ($T \to \infty$) согласно эргодической теореме для марковских процессов (в нашем случае марковский процесс – модель стохастической химической кинетики с унарными реакциями в условиях детального баланса, см. приложение 1) приходим к финальной (=стационарной) вероятностной мере, имеющей в основе мультиномиальное распределение. С ростом числа агентов ($M \to \infty$) эта мера концентрируется около наиболее вероятного состояния, поиск которого сводится к реше-нию задачи (2.2.1) ниже. Функционал в этой задаче оптимизации с точностью до потенцирования и мультипликативных и аддитивных констант соответствует исследуемой стационарной мере, то есть это функционал Санова. При обратном порядке переходов мы сначала осуществляем так называемый канонический скейлинг, приводящий к детерминированной кинетической динамике, описываемой СОДУ на $x$, а затем ($T \to \infty$) ищем аттрактор получившейся СОДУ. Глобальным аттрактором оказывается неподвижная точка, которая определяется решением задачи (2.2.1) ниже. Более того, функционал, стоящий в (2.2.1), является функцией Ляпунова полученной кинетической динамики (то есть функционалом Больцмана). Последнее утверждение – достаточно общий факт (функционал Санова, является функционалом Больцмана), верный при намного более общих условиях, см. приложение 1.

Итак, рассматривается следующая задача поиска стохастического равновесия Нэша–Вардропа в модели Бэкмана равновесного распределения транспортных потоков по путям:

$$\sum_{e\in E}\sigma_e\left(f_e\right)+\gamma\sum_{w\in OD}\sum_{p\in P_w}x_p\ln\left(x_p/d_w\right)\to\min_{f=\Theta x,\,x\in X}, \qquad (2.2.1)$$

где $\gamma>0$; $\sigma_e\left(f_e\right)=\int\limits_0^{f_e}\tau_e\left(z\right)dz$ – выпуклые функции;

$$\Theta=\left\|\delta_{ep}\right\|_{e\in E,\,p\in P}=\left\|\Theta^{\langle p\rangle}\right\|_{p\in P},\ \delta_{ep}=\begin{cases}1, & e\in p,\\ 0, & e\notin p;\end{cases}$$

$$X=\left\{x\geq 0:\ \sum_{p\in P_w}x_p=d_w,\ w\in OD\right\}$$ – прямое произведение симплексов;

в качестве $\tau_e\left(f_e\right)$ обычно выбирают BPR-функции [140]:

$$\tau_e\left(f_e\right)=\bar{t}_e\cdot\left(1+\rho\cdot\left(f_e/\bar{f}_e\right)^4\right).$$

В пределе модели стабильной динамики

$$\tau_e^{\mu}\left(f_e\right)=\bar{t}_e\cdot\left(1+\rho\cdot\left(f_e/\bar{f}_e\right)^{1/\mu}\right)\xrightarrow[\mu\to 0+]{}\begin{cases}\bar{t}_e, & 0\leq f_e<\bar{f}_e,\\ \left[\bar{t}_e,\infty\right), & f_e=\bar{f}_e,\end{cases}$$

$$d\tau_e^{\mu}\left(f_e\right)/df_e\xrightarrow[\mu\to 0+]{}0,\ 0\leq f_e<\bar{f}_e,$$

задача перепишется как

$$\sum_{e\in E}f_e\bar{t}_e+\gamma\sum_{w\in OD}\sum_{p\in P_w}x_p\ln\left(x_p/d_w\right)\to\min_{\substack{f=\Theta x,\,x\in X\\ f\leq\bar{f}}}. \qquad (2.2.2)$$

## 2.2.3. Переход к двойственной задаче

Запишем двойственную задачу к (2.2.1) (далее мы используем обозначение $\operatorname{dom}\sigma^*$ – область определения сопряженной к $\sigma$ функции):

$$\min_{f,x}\left\{\sum_{e\in E}\sigma_e\left(f_e\right)+\gamma\sum_{w\in OD}\sum_{p\in P_w}x_p\ln\left(x_p/d_w\right):\ f=\Theta x,\ \ x\in X\right\}=$$

$$=\min_{f,x}\left\{\sum_{e\in E}\max_{t_e\in\operatorname{dom}\sigma_e^*}\left[f_et_e-\sigma_e^*\left(t_e\right)\right]+\gamma\sum_{w\in OD}\sum_{p\in P_w}x_p\ln\left(x_p/d_w\right):\ f=\Theta x,\ \ x\in X\right\}=$$

$$=\max_{t\in\operatorname{dom}\sigma^*}\left\{\min_{f,x}\left[\sum_{e\in E}f_et_e+\gamma\sum_{w\in OD}\sum_{p\in P_w}x_p\ln\left(x_p/d_w\right):\ f=\Theta x,\ \ x\in X\right]-\sum_{e\in E}\sigma_e^*\left(t_e\right)\right\}=$$

$$=-\min_{t\in\operatorname{dom}\sigma^*}\left\{\gamma\psi\left(t/\gamma\right)+\sum_{e\in E}\sigma_e^*\left(t_e\right)\right\},\qquad(2.2.3)$$

где[21]

$$\psi\left(t\right)=\sum_{w\in OD}d_w\psi_w\left(t\right),\ \psi_w\left(t\right)=\ln\left(\sum_{p\in P_w}\exp\left(-\sum_{e\in E}\delta_{ep}t_e\right)\right),$$

$$f=-\nabla\gamma\psi\left(t/\gamma\right),\ x_p=d_w\frac{\exp\left(-\frac{1}{\gamma}\sum_{e\in E}\delta_{ep}t_e\right)}{\sum_{q\in P_w}\exp\left(-\frac{1}{\gamma}\sum_{e\in E}\delta_{eq}t_e\right)},\ p\in P_w,\qquad(2.2.4)$$

для

$$\tau_e\left(f_e\right)=\bar{t}_e\cdot\left(1+\rho\cdot\left(\frac{f_e}{\bar{f}_e}\right)^{\frac{1}{\mu}}\right)$$

имеем (BPR-функция, см. подраздел 2.2.2, получается при $\mu=1/4$):

---

[21]Обратим внимание, что вектор распределения потоков по путям $x$ при поиске стохастического равновесия ($\gamma>0$) получается не разреженным в отличие от поиска обычного равновесия Нэша(–Вардропа) ($\gamma\to 0+$). Как следствие, чтобы вычислить этот вектор, требуются затраты, существенно зависящие от потенциально огромной размерности $|P|$. К счастью, в приложениях, как правило, не требуется знание этого вектора, достаточно определить вектор потоков на ребрах $f$, который, как мы увидим ниже, может быть вычислен намного эффективнее (в частности, с затратами, независящими от $|P|$).

$$\sigma_e^*(t_e)=\sup_{f_e\ge 0}\left(\left(t_e-\bar{t}_e\right)\cdot f_e-\bar{t}_e\cdot\frac{\mu}{1+\mu}\cdot\rho\cdot\frac{f_e^{1+\frac{1}{\mu}}}{\bar{f}_e^{\frac{1}{\mu}}}\right)=\bar{f}_e\cdot\left(\frac{t_e-\bar{t}_e}{\bar{t}_e\cdot\rho}\right)^{\mu}\frac{\left(t_e-\bar{t}_e\right)}{1+\mu}.$$

Собственно, формула (2.2.3) есть не что иное, как отражение формулы $f=-\nabla\gamma\psi(t/\gamma)$ и связи $t_e=\tau_e(f_e)$, $e\in E$. Действительно, по формуле Демьянова–Данскина–Рубинова:

$$\frac{d\sigma_e^*(t_e)}{dt_e}=\frac{d}{dt_e}\max_{f_e\ge 0}\left\{t_e f_e-\int_0^{f_e}\tau_e(z)dz\right\}=f_e:\ t_e=\tau_e(f_e).$$

В свою очередь, формула $f=-\nabla\gamma\psi(t/\gamma)$ может интерпретироваться как следствие соотношений $f=\Theta x$ и формулы распределения Гиббса (логит-распределения):

$$x_p=d_w\frac{\exp\left(-\frac{1}{\gamma}\sum_{e\in E}\delta_{ep}t_e\right)}{\sum_{q\in P_w}\exp\left(-\frac{1}{\gamma}\sum_{e\in E}\delta_{eq}t_e\right)},\ p\in P_w.$$

При такой интерпретации связь задачи (2.2.3) с логит-динамикой, порождающей стохастические равновесия, наиболее наглядна. Последняя формула ввиду того, что $g_p(t)=\sum_{e\in E}\delta_{ep}t_e$ – затраты на пути $p$ на графе $\langle V,E\rangle$, ребра которого взвешены $t$, есть не что иное, как отражение следующего принципа поведения (ограниченной рациональности агентов): каждый агент $k$ (пользователь транспортной сети), отвечающий корреспонденции $w\in OD$, выбирает маршрут следования $p\in P_w$, если

$$p=\arg\max_{q\in P_w}\left\{-g_q(t)+\xi_q^k\right\},$$

где независимые случайные величины $\xi_q^k$ имеют одинаковое двойное экспоненциальное распределение, также называемое распределением Гумбеля[22]:

$$P\left(\xi_q^k<\zeta\right)=\exp\left\{-e^{-\zeta/\gamma-E}\right\},\ \gamma>0.$$

[22]Распределение Гумбеля можно объяснить исходя из идемпотентного аналога центральной предельной теоремы (вместо суммы случайных величин – максимум) для независимых случайных величин с экспоненциальным и более быстро убывающим правым хвостом [157]. Распределение Гумбеля возникает в данном контексте, например, если при принятии решения водитель собирает информацию с большого числа разных (независимых) зашумленных источников, ориентируясь на худшие прогнозы по каждому из путей.

Отметим также, что если взять $E \approx 0.5772$ – константа Эйлера, то $\mathrm{M}\left[\xi_q^k\right]=0$, $D\left[\xi_q^k\right]=\gamma^2\pi^2/6$. Распределение Гиббса получается в пределе, когда число агентов на каждой корреспонденции стремится к бесконечности (случайность исчезает и описание переходит на средние величины).

Сформулируем главный вывод из написанного выше. Если каждый пользователь сориентирован на вектор затрат $t$ на ребрах $E$, одинаковый для всех пользователей, и пытается выбрать кратчайший путь, исходя из зашумленной информации, то такое поведение пользователей в пределе, когда их число стремится к бесконечности $M \to \infty$, приводит к описанию распределения пользователей по путям / ребрам (2.2.4). *Равновесная конфигурация, описываемая решением задачи* (2.2.3), *характеризуется тем, что по вектору* $t$ *вычисляется согласно формуле* (2.2.4) *такой вектор* $f = \Theta x$, *что имеет место соотношение* $t=\left\{\tau_e\left(f_e\right)\right\}_{e\in E}$.

Заметим, что при $\gamma \to 0+$ распределение водителей по путям вырождается, и все водители (агенты) будут использовать только кратчайшие пути:

$$-\lim_{\gamma\to 0+}\gamma\psi_w\left(t/\gamma\right)=\min_{p\in P_w} g_p\left(t\right), \tag{2.2.5}$$

$$-\lim_{\gamma\to 0+}\nabla\gamma\psi_w\left(t/\gamma\right)=\mathrm{conv}\left\{\left(0,1,\ldots,0,1\right)^T,\ldots,(1,1,\ldots,0,0)^T\right\}, \tag{2.2.6}$$

в правой части формулы (2.2.6) стоит субградиент (негладкой) функции длины кратчайшего пути, отвечающего корреспонденции $w$. Этот субградиент (в общем случае) есть множество, представимое в виде выпуклой оболочки векторов, отвечающих всевозможным кратчайшим путям, начинающихся в вершине $i$ и заканчивающихся в вершине $j$, где $w=\left(i,j\right)$, в рассматриваемом транспортном графе, ребра которого взвешены вектором $t$. Единицы стоят в векторах на местах, которые отвечают ребрам, входящим в данный кратчайший путь. Если кратчайший путь единственен, то субградиент превращается в обычный градиент. Далее (в УМПТ) можно выбирать произвольный элемент субградиента. От этого выбора приводимые ниже оценки не будут зависеть.

Полезно также иметь в виду, что

$$\gamma\psi_w\left(t/\gamma\right)=\mathrm{M}_{\left\{\xi_p\right\}_{p\in P_w}}\left[\max_{p\in P_w}\left\{-g_p\left(t\right)+\xi_p\right\}\right].$$

В пределе стабильной динамики $\mu \to 0+$ задача (2.2.3) вырождается в задачу

$$\gamma\psi\left(t/\gamma\right)+\left\langle \bar{f},t-\bar{t}\right\rangle \to \min_{t\geq\bar{t}}.$$

В пределе $\gamma \to 0+$ все приведенные выше формулы переходят в соответствующие формулы для классических (не стохастических) моделей Бэкмана и стабильной динамики.

## 2.2.4. Универсальный метод подобных треугольников

Ниже изложен возможный способ решения двойственной задачи (2.2.3), позволяющий при этом восстанавливать решение прямой задачи ((2.2.1) или (2.2.2)). Рассмотрим общую задачу выпуклой композитной оптимизации на множестве простой структуры:

$$F(t) = \underbrace{\Phi(t)}_{\gamma\psi(t/\gamma)} + \underbrace{h(t)}_{\sum_{e\in E}\sigma_e^*(t_e)} \to \min_{t\in Q}. \qquad (2.2.7)$$

Положим $R^2 = \frac{1}{2}\left\|t_* - y^0\right\|_2^2$, где $y^0 = \bar{t}$, $t_*$ – решение задачи (2.2.7) (если решение не единственно, то выбирается то, которое доставляет минимум $\left\|t_* - y^0\right\|_2^2$).

Пусть ( $L_\nu \le \infty$, т. е. допускается равенство бесконечности)

$$\left\|\nabla\Phi(t) - \nabla\Phi(y)\right\|_2 \le L_\nu \left\|t - y\right\|_2^\nu,\ \nu \in [0,1],\ L_0 < \infty. \qquad (2.2.8)$$

Положим

$$\varphi_0(t) = \alpha_0\left[\Phi(y^0) + \left\langle\nabla\Phi(y^0), t - y^0\right\rangle + h(t)\right] + \frac{1}{2}\left\|t - y^0\right\|_2^2,$$

$$\varphi_{k+1}(t) = \varphi_k(t) + \alpha_{k+1}\left[\Phi(y^{k+1}) + \left\langle\nabla\Phi(y^{k+1}), t - y^{k+1}\right\rangle + h(t)\right],$$

Начнем описание универсального метода подобных треугольников (УМПТ) с самой первой итерации. Положим

$$A_0 = \alpha_0 = 1/L_0^0,\ k = 0,\ j_0 = 0;\ t^0 = u^0 = \arg\min_{t\in Q}\varphi_0(t).$$

До тех пор пока

$$\Phi(t^0) > \Phi(y^0) + \left\langle\nabla\Phi(y^0), t^0 - y^0\right\rangle + \frac{L_0^{j_0}}{2}\left\|t^0 - y^0\right\|_2^2 + \frac{\alpha_0}{2A_0}\varepsilon,$$

выполнять

$$j_0 := j_0 + 1;\ L_0^{j_0} := 2^{j_0} L_0^0;\ t^0 := u^0 := \arg\min_{t\in Q}\varphi_0(t),\ (A_0 :=)\alpha_0 := \frac{1}{L_0^{j_0}}.$$

**Универсальный метод подобных треугольников**

1. $L_{k+1}^0 = L_k^{j_k}/2$, $j_{k+1} = 0$.

2. $$\begin{cases} \alpha_{k+1} := \dfrac{1}{2L_{k+1}^{j_{k+1}}} + \sqrt{\dfrac{1}{4\left(L_{k+1}^{j_{k+1}}\right)^2} + \dfrac{A_k}{L_{k+1}^{j_{k+1}}}},\ A_{k+1} := A_k + \alpha_{k+1}; \\ y^{k+1} := \dfrac{\alpha_{k+1}u^k + A_k t^k}{A_{k+1}},\ u^{k+1} := \arg\min\limits_{t\in Q} \varphi_{k+1}(t),\ t^{k+1} := \dfrac{\alpha_{k+1}u^{k+1} + A_k t^k}{A_{k+1}}. \end{cases} \qquad (*)$$

До тех пор пока

$$\Phi\left(y^{k+1}\right) + \left\langle \nabla\Phi\left(y^{k+1}\right), t^{k+1} - y^{k+1}\right\rangle + \frac{L_{k+1}^{j_{k+1}}}{2}\left\|t^{k+1} - y^{k+1}\right\|_2^2 + \frac{\alpha_{k+1}}{2A_{k+1}}\varepsilon < \Phi\left(t^{k+1}\right),$$

выполнять

$$j_{k+1} := j_{k+1} + 1;\ L_{k+1}^{j_{k+1}} = 2^{j_{k+1}} L_{k+1}^0;\ (*).$$

3. $k := k+1$ и **go to** 1.

УМПТ сходится согласно оценке (см. приложение 3):

$$A_N F\left(t^N\right) \le \min_{t\in Q}\left\{\frac{1}{2}\left\|t - y^0\right\|_2^2 + \sum_{k=0}^{N}\alpha_k\left[\Phi\left(y^k\right) + \left\langle\nabla\Phi\left(y^k\right), t - y^k\right\rangle + h(t)\right]\right\}. \qquad (2.2.9)$$

Отсюда

$$A_N F\left(t^N\right) \le \min_{t\in Q}\left\{\frac{1}{2}\left\|t - y^0\right\|_2^2 + \sum_{k=0}^{N}\alpha_k\left[\Phi\left(y^k\right) + \left\langle\nabla\Phi\left(y^k\right), t - y^k\right\rangle + h(t)\right]\right\} + \frac{\varepsilon}{2} \le$$

$$\le \frac{1}{2}\left\|t_* - y^0\right\|_2^2 + \sum_{k=0}^{N}\alpha_k\underbrace{\left[\Phi\left(y^k\right) + \left\langle\nabla\Phi\left(y^k\right), t_* - y^k\right\rangle + h(t_*)\right]}_{\le\Phi(t_*)+h(t_*)} + \frac{\varepsilon}{2} \le$$

$$\le \frac{1}{2}\left\|t_* - y^0\right\|_2^2 + \sum_{k=0}^{N}\alpha_k F(t_*) + \frac{\varepsilon}{2} = R^2 + A_N F(t_*) + \frac{\varepsilon}{2},$$

т. е.

$$F\left(t^N\right) - F(t_*) \le \frac{R^2}{A_N} + \frac{\varepsilon}{2}.$$

В приложении 3 показано, что $A_N \simeq 2R^2/\varepsilon$ при

$$N \simeq \inf_{\nu\in[0,1]}\left(\frac{L_\nu\cdot(16R)^{1+\nu}}{\varepsilon}\right)^{\frac{2}{1+3\nu}}, \qquad (2.2.10)$$

т. е. при таком $N$ справедливо неравенство

$$F\left(t^N\right) - \min_{t\in Q} F\left(t\right) \le \varepsilon .$$

При этом среднее число вычислений значения функции на одной итерации будет $\approx 4$, а градиента функции $\approx 2$.

## 2.2.5. Приложение УМПТ к поиску равновесий в транспортных сетях

Задачи, получающиеся из (2.2.3) при $\mu > 0$ и $\mu \to 0+$, можно не различать по сложности при композитном подходе к ним [37, 126]. В задаче, отвечающей $\mu \to 0+$, сепарабельный композит $h(t)$ проще сепарабельного композита задачи с $\mu > 0$, но зато дополнительно добавляется сепарабельное ограничение простой структуры[23]. В любом случае, основные затраты на каждой итерации при использовании УМПТ в композитном варианте [37] связаны с необходимостью расчета градиента гладкой функции $\gamma\psi\left(t/\gamma\right)$. Тому, как это можно эффективно делать, будет посвящен следующий подраздел 2.2.6. Здесь же отметим, что константы Гёльдера (см. (2.2.8)) градиента этой функции в 2-норме могут быть оценены следующим образом:

$$L_1 \le \frac{1}{\gamma} \sum_{w\in OD} d_w \max_{p\in P_w} \left\| \Theta^{\langle p\rangle} \right\|_2^2 \le \frac{Hd}{\gamma}, \tag{2.2.11}$$

$$L_0 \le 2 \sum_{w\in OD} d_w \max_{p\in P_w} \left\| \Theta^{\langle p\rangle} \right\|_2 \le 2\sqrt{H} d \text{ (при } \gamma \to 0+\text{)}, \tag{2.2.12}$$

где $d = \sum_{w\in OD} d_w$, $H = \max_{p\in P} \left\| \Theta^{\langle p\rangle} \right\|_2^2$, т. е. $H$ – максимальное число ребер в пути. Как правило, можно считать, что $H = \mathrm{O}\left(\sqrt{n}\right)$ – диаметр графа с манхетенской структурой, то есть для квадратной решетки с $n$ ребрами.

В данном пункте внимание будет сосредоточено на том, как с помощью формулы (2.2.9) восстанавливать решение прямой задачи ((2.2.1) или (2.2.2)), исходя из последовательности, сгенерированной УМПТ при решении двойственной задачи (2.2.3). Описанные далее способы имеют много общего с написанным ранее в подразделе 2.1.3.

Положим

---

[23] Для BPR-функций ( $\mu = 0.25$ ) все сводится к решению уравнения четвертой степени, см. подраздел 2.2.3.

$$f^k = -\nabla\gamma\psi\left(y^k/\gamma\right),\ x_p^k = d_w \frac{\exp\left(-\frac{1}{\gamma}\sum_{e\in E}\delta_{ep}y_e^k\right)}{\sum_{q\in P_w}\exp\left(-\frac{1}{\gamma}\sum_{e\in E}\delta_{eq}y_e^k\right)},\ p\in P_w,\ w\in OD,$$

$$\bar{f}^N = \frac{1}{A_N}\sum_{k=0}^{N} a_k f^k,\ \bar{x}^N = \frac{1}{A_N}\sum_{k=0}^{N}\alpha_k x^k .$$

Допустим[24]

$$\left\{\tau_e\left(\bar{f}_e^N\right)\right\}_{e\in E} = \arg\min_{t\in\operatorname{dom}\sigma^*}\left\{\frac{1}{A_N}\left[\sum_{k=0}^{N}\alpha_k\cdot\left(\gamma\psi\left(y^k/\gamma\right)+\right.\right.\right. \tag{2.2.13}$$

$$\left.\left.\left.+\left\langle\nabla\gamma\psi\left(y^k/\gamma\right),t-y^k\right\rangle\right)\right]+\sum_{e\in E}\sigma_e^*\left(t_e\right)\right\}.$$

Из (2.2.13) следует, что

$$\gamma\psi\left(t^N/\gamma\right)+\sum_{e\in E}\sigma_e^*\left(t_e^N\right)\le$$

$$\le\min_{t\in\operatorname{dom}\sigma^*}\left\{\frac{1}{A_N}\left[\sum_{k=0}^{N}\alpha_k\cdot\left(\gamma\psi\left(y^k/\gamma\right)+\left\langle\nabla\gamma\psi\left(y^k/\gamma\right),t-y^k\right\rangle\right)\right]+\sum_{e\in E}\sigma_e^*\left(t_e\right)+\frac{1}{2A_N}\left\|t-\bar{t}\right\|_2^2\right\}$$

$$+\frac{\varepsilon}{2}\le\min_{t\in\operatorname{dom}\sigma^*}\left\{\frac{1}{A_N}\left[\sum_{k=0}^{N}\alpha_k\cdot\left(\gamma\psi\left(y^k/\gamma\right)+\left\langle\nabla\gamma\psi\left(y^k/\gamma\right),t-y^k\right\rangle\right)\right]+\sum_{e\in E}\sigma_e^*\left(t_e\right)\right\}+$$

$$+\frac{1}{A_N}\underbrace{\frac{1}{2}\sum_{e\in E}\left(\tau_e\left(\bar{f}_e^N\right)-\bar{t}_e\right)^2}_{\tilde{R}^2}+\frac{\varepsilon}{2},$$

следовательно,

$$\gamma\psi\left(t^N/\gamma\right)+\sum_{e\in E}\sigma_e^*\left(t_e^N\right)-$$

$$-\min_{t\in\operatorname{dom}\sigma^*}\left\{\frac{1}{A_N}\left[\sum_{k=0}^{N}\alpha_k\cdot\left(\gamma\psi\left(y^k/\gamma\right)+\left\langle\nabla\gamma\psi\left(y^k/\gamma\right),t-y^k\right\rangle\right)\right]+\sum_{e\in E}\sigma_e^*\left(t_e\right)\right\}\le\frac{\tilde{R}^2}{A_N}+\frac{\varepsilon}{2}.$$

Учитывая, что

$$-\min_{t\in\operatorname{dom}\sigma^*}\left\{\frac{1}{A_N}\left[\sum_{k=0}^{N}\alpha_k\left\langle\nabla\gamma\psi\left(y^k/\gamma\right),t\right\rangle\right]+\sum_{e\in E}\sigma_e^*\left(t_e\right)\right\}=$$

---

[24]Минимум в (2.2.13) достигается во внутренней точке множества $\operatorname{dom}\sigma^*$, которое для BPR-функций имеет вид $t\ge\bar{t}$, и лишь при $\mu\to 0+$ минимум может выходить на границу $t_e=\bar{t}_e$.

$$= \max_{t\in\operatorname{dom}\sigma^*}\left\{\left\langle \frac{1}{A_N}\sum_{k=0}^{N}\alpha_k f^k, t\right\rangle - \sum_{e\in E}\sigma_e^*\left(t_e\right)\right\} = \sum_{e\in E}\sigma_e\left(\bar{f}_e^N\right),$$

$$-\frac{1}{A_N}\sum_{k=0}^{N}\alpha_k\cdot\left(\gamma\psi\left(y^k/\gamma\right)-\left\langle\nabla\gamma\psi\left(y^k/\gamma\right),y^k\right\rangle\right)=$$

$$=\frac{1}{A_N}\sum_{k=0}^{N}\alpha_k\cdot\left(\left\langle f^k,y^k\right\rangle+\gamma\sum_{w\in OD}\sum_{p\in P_w}x_p^k\ln\left(x_p^k/d_w\right)-\left\langle f^k,y^k\right\rangle\right)=$$

$$=\gamma\frac{1}{A_N}\sum_{k=0}^{N}\alpha_k\sum_{w\in OD}\sum_{p\in P_w}x_p^k\ln\left(x_p^k/d_w\right)\geq\gamma\sum_{w\in OD}\sum_{p\in P_w}\bar{x}_p^N\ln\left(\bar{x}_p^N/d_w\right),$$

получаем следующую оценку сверху на зазор двойственности:

$$0\leq\left\{\gamma\psi\left(t^N/\gamma\right)+\sum_{e\in E}\sigma_e^*\left(t_e^N\right)\right\}+\left\{\sum_{e\in E}\sigma_e\left(\bar{f}_e^N\right)+\gamma\sum_{w\in OD}\sum_{p\in P_w}\bar{x}_p^N\ln\left(\bar{x}_p^N/d_w\right)\right\}\leq$$

$$\leq\left\{\gamma\psi\left(t^N/\gamma\right)+\sum_{e\in E}\sigma_e^*\left(t_e^N\right)\right\}+ \qquad (2.2.14)$$

$$+\left\{\sum_{e\in E}\sigma_e\left(\bar{f}_e^N\right)+\gamma\frac{1}{A_N}\sum_{k=0}^{N}\alpha_k\sum_{w\in OD}\sum_{p\in P_w}x_p^k\ln\left(x_p^k/d_w\right)\right\}\leq\frac{\tilde{R}^2}{A_N}+\frac{\varepsilon}{2},$$

следовательно,

$$\sum_{e\in E}\sigma_e\left(\bar{f}_e^N\right)+\gamma\sum_{w\in OD}\sum_{p\in P_w}\bar{x}_p^N\ln\left(\bar{x}_p^N/d_w\right)-$$

$$-\left(\sum_{e\in E}\sigma_e\left(f_e^*\right)+\gamma\sum_{w\in OD}\sum_{p\in P_w}x_p^*\ln\left(x_p^*/d_w\right)\right)\leq\left\{\gamma\psi\left(t^N/\gamma\right)+\sum_{e\in E}\sigma_e^*\left(t_e^N\right)\right\}+ \qquad (2.2.15)$$

$$+\left\{\sum_{e\in E}\sigma_e\left(\bar{f}_e^N\right)+\gamma\sum_{w\in OD}\sum_{p\in P_w}\bar{x}_p^N\ln\left(\bar{x}_p^N/d_w\right)\right\}\leq\frac{\tilde{R}^2}{A_N}+\frac{\varepsilon}{2}$$

где $\left(f^*,x^*\right)$ – решение задачи (2.2.1).

На формулу (2.2.14) можно смотреть как на критерий останова УМПТ. А именно, ждем, когда (вычислимый с помощью быстрого автоматического дифференцирования подраздел 2.2.6) зазор двойственности станет меньше $\varepsilon$. Формула (2.2.10) (с заменой $R$ на $\tilde{R}$) содержит гарантированную оценку на число итераций УМПТ, после которого метод гарантированно остановится по критерию (2.2.14).

Для модели (стохастической) стабильной динамики ($\mu\to 0+$) необходимо немного по-другому провести рассуждения:

$$\gamma\psi\left(t^N/\gamma\right)+\left\langle \bar{f},t^N-\bar{t}\right\rangle \le$$

$$\le\min_{t\ge\bar{t}}\left\{\frac{1}{A_N}\left[\sum_{k=0}^{N}\alpha_k\cdot\left(\gamma\psi\left(y^k/\gamma\right)+\left\langle\nabla\gamma\psi\left(y^k/\gamma\right),t-y^k\right\rangle\right)\right]+\left\langle \bar{f},t-\bar{t}\right\rangle+\frac{1}{2A_N}\left\|t-\bar{t}\right\|_2^2\right\}+\frac{\varepsilon}{2}\le$$

$$\le\min_{t\ge\bar{t},\,\left\|t-\bar{t}\right\|_2^2\le 2R^2}\left\{\frac{1}{A_N}\left[\sum_{k=0}^{N}\alpha_k\cdot\left(\gamma\psi\left(y^k/\gamma\right)+\left\langle\nabla\gamma\psi\left(y^k/\gamma\right),t-y^k\right\rangle\right)\right]+\left\langle \bar{f},t-\bar{t}\right\rangle\right\}+\frac{R^2}{A_N}+\frac{\varepsilon}{2},$$

где $R^2=\frac{1}{2}\left\|t_*-y^0\right\|_2^2=\frac{1}{2}\left\|t_*-\bar{t}\right\|_2^2$. Следовательно,

$$\gamma\psi\left(t^N/\gamma\right)+\left\langle \bar{f},t^N-\bar{t}\right\rangle-$$

$$-\min_{t\ge\bar{t},\,\left\|t-\bar{t}\right\|_2^2\le 10R^2}\left\{\frac{1}{A_N}\left[\sum_{k=0}^{N}\alpha_k\cdot\left(\gamma\psi\left(y^k/\gamma\right)+\left\langle\nabla\gamma\psi\left(y^k/\gamma\right),t-y^k\right\rangle\right)\right]+\left\langle \bar{f},t-\bar{t}\right\rangle\right\}\le\frac{5R^2}{A_N}+\frac{\varepsilon}{2}.$$

Поскольку

$$-\min_{t\ge\bar{t},\,\left\|t-\bar{t}\right\|_2^2\le 10R^2}\left\{\frac{1}{A_N}\left[\sum_{k=0}^{N}\alpha_k\left\langle\nabla\gamma\psi\left(y^k/\gamma\right),t\right\rangle\right]+\left\langle \bar{f},t-\bar{t}\right\rangle\right\}=$$

$$=\max_{t\ge\bar{t},\,\left\|t-\bar{t}\right\|_2^2\le 10R^2}\left\{\left\langle\frac{1}{A_N}\sum_{k=0}^{N}\alpha_k f^k,t\right\rangle-\left\langle \bar{f},t-\bar{t}\right\rangle\right\}=$$

$$=\max_{t\ge\bar{t},\,\left\|t-\bar{t}\right\|_2^2\le 10R^2}\left\{\left\langle \bar{f}^N-\bar{f},t-\bar{t}\right\rangle\right\}+\left\langle \bar{f}^N,\bar{t}\right\rangle\ge$$

$$\ge\left\langle \bar{f}^N,\bar{t}\right\rangle+3R\left\|\left(\bar{f}^N-\bar{f}\right)_+\right\|_2,$$

то

$$\begin{gathered}\gamma\psi\left(t^N/\gamma\right)+\left\langle \bar{f},t^N-\bar{t}\right\rangle+\left\langle \bar{f}^N,\bar{t}\right\rangle+\\+\gamma\sum_{w\in OD}\sum_{p\in P_w}\bar{x}_p^N\ln\left(\bar{x}_p^N/d_w\right)+3R\left\|\left(\bar{f}^N-\bar{f}\right)_+\right\|_2\le\\\le\gamma\psi\left(t^N/\gamma\right)+\left\langle \bar{f},t^N-\bar{t}\right\rangle+\left\langle \bar{f}^N,\bar{t}\right\rangle+\\+\gamma\frac{1}{A_N}\sum_{k=0}^{N}\alpha_k\sum_{w\in OD}\sum_{p\in P_w}x_p^k\ln\left(x_p^k/d_w\right)+3R\left\|\left(\bar{f}^N-\bar{f}\right)_+\right\|_2\le\frac{5R^2}{A_N}+\frac{\varepsilon}{2}.\end{gathered}\tag{2.2.16}$$

Повторяя рассуждения п. 3 [2] и п. 6.11 [122] из (2.2.16), получим

$$\left\langle \bar{f}^N,\bar{t}\right\rangle+\gamma\sum_{w\in OD}\sum_{p\in P_w}\bar{x}_p^N\ln\left(\bar{x}_p^N/d_w\right)-\left(\left\langle f^*,\bar{t}\right\rangle+\gamma\sum_{w\in OD}\sum_{p\in P_w}x_p^*\ln\left(x_p^*/d_w\right)\right)\le$$

$$\le\gamma\psi\left(t^N/\gamma\right)+\left\langle \bar{f},t^N-\bar{t}\right\rangle+\left\langle \bar{f}^N,\bar{t}\right\rangle+\gamma\sum_{w\in OD}\sum_{p\in P_w}\bar{x}_p^N\ln\left(\bar{x}_p^N/d_w\right)\le\frac{5R^2}{A_N}+\frac{\varepsilon}{2}, \tag{2.2.17}$$

$$\left\| \left( \bar{f}^N - \bar{f} \right)_+ \right\|_2 \le \frac{5R}{A_N} + \frac{\varepsilon}{2R}, \qquad (2.2.18)$$

где $\left( f^*, x^* \right)$ – решение задачи (2.2.2).

На формулу (2.2.16) можно смотреть как на критерий останова УМПТ. А именно, ждем, когда (вычислимая) левая часть (2.2.16) станет меньше $\varepsilon$. Формула (2.2.10) (с заменой $R$ на $\sqrt{5}R$) содержит гарантированную оценку на число итераций УМПТ, после которого метод гарантированно остановится по критерию (2.2.16). Заметим, что критерии останова (2.2.14) и (2.2.16) эффективно вычислимы, несмотря на энтропийный член, потенциально содержащий огромное число слагаемых [35].

Далее для краткости будем единым образом обозначать через $\bar{R}$ либо $\tilde{R}$, либо $\sqrt{5}R$. Расшифровывать $\bar{R}$ нужно будет в зависимости от контекста. Если $\mu \to 0+$, то $\bar{R} = \sqrt{5}R$, иначе $\bar{R} = \tilde{R}$.

Все приведенные в этом пункте формулы выдерживают предельные переходы $\gamma \to 0+$, $\mu \to 0+$ (см. подраздел 2.2.2).

### 2.2.6. Вычисление (суб-)градиентов в задаче поиска равновесий в транспортных сетях

Приведем, следуя [34, 123], сглаженный идемпотентный аналог метода Форда–Беллмана, позволяющий эффективно рассчитывать значение характеристической функции $\gamma\psi\left(t/\gamma\right)$. Для этого предположим, что любые движения по ребрам графа с учетом их ориентации являются допустимыми, т. е. множество путей, соединяющих заданные две вершины (источник и сток), – это множество всевозможных способов добраться из источника в сток по ребрам имеющегося графа с учетом их ориентации. Этого всегда можно добиться раздутием исходного графа в несколько раз за счет введения дополнительных вершин и ребер. Такое раздутие заведомо можно сделать за $\mathrm{O}(n)$. Отметим, что при этом в качестве путей будут присутствовать в том числе и самопересекающиеся маршруты. Однако можно показать, что вклад таких «нефизических» путей в итоговую равновесную конфигурацию будет пренебрежимо мал.

Будем считать, как и прежде, что число ребер в любом пути не больше $H = \mathrm{O}\left(\sqrt{n}\right)$. Введем классы путей: $P_{ij}^l$ – множество всех путей из $i$ в $j$, состоящих ровно из $l$ ребер, $\tilde{P}_{ij}^l$ – множество всех путей из $i$ в $j$, состоящих из не более чем $l$ ребер. Зафиксируем источник (вершину) $i \in V$ и введем следующие функции для $j \in V$, $l = 1, ..., H$:

$$\begin{cases} a_{ij}^{l}(t) = \gamma\psi_{P_{ij}^{l}}(t/\gamma) = \gamma \ln\left( \sum_{p\in P_{ij}^{l}} \exp\left( -\sum_{e\in E} \delta_{ep} t_e \Big/ \gamma \right) \right), \\ b_{ij}^{l}(t) = \gamma\psi_{\tilde{P}_{ij}^{l}}(t/\gamma) = \gamma \ln\left( \sum_{p\in \tilde{P}_{ij}^{l}} \exp\left( -\sum_{e\in E} \delta_{ep} t_e \Big/ \gamma \right) \right). \end{cases}$$

Некоторые из этих функций могут быть равны $-\infty$. Это означает, что соответствующее множество маршрутов пустое. Данные функции можно вычислять рекурсивным образом:

$$a_{ij}^{1}(t) = b_{ij}^{1}(t) = \begin{cases} -t_e, \ e = (i \to j) \in E \\ -\infty, \ e = (i \to j) \notin E \end{cases}$$

$$\begin{cases} a_{ij}^{l+1}(t) = \gamma \ln\left( \sum_{k:\, e=(k\to j)\in E} \exp\left( \left( a_{ik}^{l}(t) - t_e \right) / \gamma \right) \right), \\ b_{ij}^{l+1}(t) = \gamma \ln\left( \exp\left( b_{ij}^{l}(t)/\gamma \right) + \exp\left( a_{ij}^{l+1}(t)/\gamma \right) \right), \end{cases} \quad j \in V, \ \ l = 1,\dots,H-1. \quad (2.2.19)$$

На каждом шаге $l$ необходимо сделать $\mathrm{O}(n)$ арифметических операций. Следовательно, для вычисления $\gamma\psi(t/\gamma)$ необходимо сделать $\mathrm{O}(SHn)$ арифметических операций, где $S = |O|$ – число источников, как правило, можно считать $S \ll n$. Причем вычисление функции $\gamma\psi(t/\gamma)$ (и ее градиента) может быть распараллелено на $S$ процессорах. При $\gamma \to 0+$ процедура вырождается в известный метод Форда–Беллмана [68] (динамическое программирование). В процедуре Форда–Беллмана требуется посчитать $H$ -степень матрицы $A = \left\| a_{ij} \right\|_{i,j\in V}$:

$$a_{ij} = t_e, \ e = (i \to j) \in E;$$

$$a_{ij} = \infty, \ e = (i \to j) \notin E,$$

в идемпотентной математике (вместо обычного поля используется тропическое полуполе [113] со следующими операциями: сложение $a \oplus b = \min\{a,b\}$, произведение $a \otimes b = a + b$). Учитывая, что

$$a \oplus b = \min\{a,b\} = -\lim_{\gamma\to 0+} \gamma \ln\left( \exp(-a/\gamma) + \exp(-b/\gamma) \right),$$

$$a \otimes b = a + b = -\lim_{\gamma\to 0+} \gamma \ln\left( \exp(-a/\gamma) \cdot \exp(-b/\gamma) \right),$$

можно посчитать обычную (над обычным полем) $H$ -степень матрицы $A^{\gamma} = \left\| a_{ij}^{\gamma} \right\|_{i,j\in V}$:

$$a_{ij}^{\gamma} = e^{-t_e/\gamma}, \; e = (i \to j) \in E,$$
$$a_{ij}^{\gamma} = 0, \; e = (i \to j) \notin E$$

и применить поэлементно к полученной матрице $-\gamma \ln(\,\cdot\,)$. В пределе $\gamma \to 0+$ получим метод Форда–Беллмана. Однако если не делать предельный переход, то получается нужный нам сглаженный вариант этого алгоритма с такой же временной сложностью. Некоторый аналог этого сглаженного варианта, по сути, и был описан выше.

Используя быстрое автоматическое дифференцирования [44], опишем способ вычисления градиента функции $\gamma\psi^i(t/\gamma) = \sum_{j:\, w=(i,j)\in OD} d_w \gamma \psi_w(t/\gamma)$:

$$\frac{\partial \psi^i}{\partial b_{ij}^H} = d_{ij}, \; \frac{\partial \psi^i}{\partial a_{ij}^H} = 0,$$

$$\begin{cases} \dfrac{\partial \psi^i}{\partial b_{ij}^l} = \dfrac{\partial \psi^i}{\partial b_{ij}^{l+1}} \dfrac{\partial b_{ij}^{l+1}}{\partial b_{ij}^l}, \\ \dfrac{\partial \psi^i}{\partial a_{ij}^l} = \dfrac{\partial \psi^i}{\partial b_{ij}^l} \dfrac{\partial b_{ij}^l}{\partial a_{ij}^l} + \sum\limits_{k:\, e=(j\to k)\in E} \dfrac{\partial \psi^i}{\partial a_{ik}^{l+1}} \dfrac{\partial a_{ik}^{l+1}}{\partial a_{ij}^l}, \end{cases} \quad j \in V, \; l = H-1,\ldots,1, \quad (2.2.20)$$

$$\frac{\partial \psi^i}{\partial t_e} = \sum_{l=0}^{H-1} \sum_{j\in V} \frac{\partial \psi^i}{\partial a_{ij}^{l+1}} \frac{\partial a_{ij}^{l+1}}{\partial t_e} = \sum_{l=0}^{H-1} \frac{\partial \psi^i}{\partial a_{ij'}^{l+1}} \frac{\partial a_{ij'}^{l+1}}{\partial t_e}, \; e = (k, j'). \qquad (2.2.21)$$

Частные производные

$$\frac{\partial b_{ij}^{l+1}}{\partial b_{ij}^l}, \; \frac{\partial b_{ij}^l}{\partial a_{ij}^l}, \; \frac{\partial a_{ik}^{l+1}}{\partial a_{ik}^l}, \; \frac{\partial a_{ij}^{l+1}}{\partial t_e}$$

могут быть явно вычислены из системы (2.2.19) за $O(1)$ каждая. Поясним примером. Пусть, например, для некоторых $l$ и $j$ имеет место

$$a_{ij}^{l+1}(t) = \gamma \ln\left(\exp\left(\left(a_{ik_1}^l(t) - t_{e_1}\right)/\gamma\right) + \exp\left(\left(a_{ik_2}^l(t) - t_{e_2}\right)/\gamma\right)\right).$$

Тогда[25]

[25]Заметим, что возникающие при расчете градиента отношения экспонент стоит сразу же приводить (для большей вычислительной устойчивости) к дробям с числителем, равным 1:

$$\frac{e^a}{e^a + e^b + \ldots} = \frac{1}{1 + e^{b-a} + \ldots}.$$

Аналогично с выражениями вида

$$\ln(\exp(a) + \exp(b)) = c + \ln(\exp(a-c) + \exp(b-c)), \; c = \max\{a, b\}.$$

$$\frac{\partial a_{ij}^{l+1}}{\partial a_{ik_1}^{l}} = \frac{\exp\left(\left(a_{ik_1}^{l}(t) - t_{e_1}\right)/\gamma\right)}{\exp\left(\left(a_{ik_1}^{l}(t) - t_{e_1}\right)/\gamma\right) + \exp\left(\left(a_{ik_2}^{l}(t) - t_{e_2}\right)/\gamma\right)},$$

$$\frac{\partial a_{ij}^{l+1}}{\partial t_{e_1}} = -\frac{\exp\left(\left(a_{ik_1}^{l}(t) - t_{e_1}\right)/\gamma\right)}{\exp\left(\left(a_{ik_1}^{l}(t) - t_{e_1}\right)/\gamma\right) + \exp\left(\left(a_{ik_2}^{l}(t) - t_{e_2}\right)/\gamma\right)}.$$

Остальные частные производные являются неизвестными. Несложно понять, что система (2.2.20) может быть последовательно разрешена за $\mathrm{O}(Hn)$. Фактически тут происходит процесс прохождения графа вычислений (2.2.19) в обратном направлении, с той лишь разницей, что обратный проход получается приблизительно в 4 раза дороже. К сожалению, в отличие от прямого процесса (2.2.19), теперь, чтобы вычислить градиент, необходимо хранить в памяти промежуточные вычисления, см. формулу (2.2.21). То есть для вычисления градиента $\gamma\psi^i(t/\gamma)$ требуется время $\mathrm{O}(Hn)$ и память $\mathrm{O}(Hn)$. Таким образом, для вычисления градиента $\gamma\psi(t/\gamma) = \sum_{i\in O}\gamma\psi^i(t/\gamma)$ требуется время $\mathrm{O}(SHn)$ и память $\mathrm{O}(SHn)$. Так же, как и при вычислении значения функции $\gamma\psi(t/\gamma)$, описанная выше схема вычисления градиента может быть распараллелена на $S = |O|$ процессорах.

Заметим, что расчет $\nabla\gamma\psi(t/\gamma)$, реализованный на языке программирования MatLab 8 на ноутбуке с тактовой частотой 1.9 ГГц, для манхетенского графа с числом ребер $n \simeq 10^3$, занимал порядка одной минуты, в то время как использование прямого дифференцирования занимало порядка 100 минут.

В связи с написанным выше, важно отметить, что в пределе $\gamma \to 0+$ оценка $\mathrm{O}(SHn)$ сложности вычисления градиента $\gamma\psi(t/\gamma)$ может быть улучшена до оценки $\mathrm{O}(Sn\ln n)$ вычисления субградиента (2.2.5), (2.2.6) [68]. Действительно, для каждого из $S$ источников можно построить (например, алгоритмом Дейкстры [68]) соответствующее дерево кратчайших путей. Исходя из принципа динамического программирования: «часть кратчайшего пути сама будет кратчайшим путем» несложно понять, что получится именно дерево с корнем в рассматриваемом источнике. Это можно сделать для одного источника за $\mathrm{O}(n\ln n)$ [68]. Однако, главное, правильно взвешивать ребра (их $n$) такого дерева, чтобы за один проход этого дерева можно было восстановить вклад (по всем ребрам) соответствующего источника в субградиент. Ребро должно иметь вес, равный сумме всех проходящих через него корреспонденций с заданным источником (корнем дерева). Имея значения соответствующих корреспонденций (их не больше $O(n)$) за один обратный проход (то есть с листьев к корню) такого дерева можно осуществить необходимое взвешивание (с затратами не более $O(n)$). Делается это по правилу: вес ребра равен сумме корреспонденции (возможно, равной нулю) в соответствующую вершину, в которую ребро входит, и сумме весов всех ребер (если таковые имеются), выходящих из упомянутой вершины.

### 2.2.7. Сопоставление оценок, численные эксперименты

Используя наработки предыдущих подразделов, в частности, формулы (2.2.10) –(2.2.12), (2.2.15), (2.2.17), (2.2.18), можно резюмировать полученные результаты в виде следующей теоремы.

**Теорема 2.2.1.** *УМПТ из подраздела* 2.2.4 *с критериями останова* (2.2.14), (2.2.16) *гарантированно остановится, достигнув желаемой точности* $\varepsilon$ *, сделав арифметических операций не больше, чем приведено в соответствующем поле табл.* 1.

Таблица 1

| Время работы | $\gamma > 0$ | $\gamma \to 0+$ |
|---|---|---|
| $\mu \to 0+, \mu > 0$ – неважно | $O\left(SHn \cdot \sqrt{\dfrac{Hd\bar{R}^2}{\gamma\varepsilon}}\right)$ (2.2.22) | $O\left(Sn\ln n \cdot \dfrac{Hd^2\bar{R}^2}{\varepsilon^2}\right)$ (2.2.23) |

Обратим внимание, что оценки (2.2.22), (2.2.23) не зависят от (потенциально экспоненциально большого) числа путей $|P|$. В частности, $|P| \gg 2^{\sqrt{n}}$, для манхетенских сетей. Также обратим внимание, что оценка (2.2.22) весьма чувствительна к предельному переходу $\gamma \to 0+$.

Важно также заметить, что обычные (не стохастические) равновесия можно искать с помощью искусственного введения энтропийной регуляризации [23]. При этом, чтобы с точностью $\varepsilon > 0$ по функции решить исходную задачу (2.2.1) или (2.2.2) с $\gamma = 0$, можно действовать следующим образом. Выбрать

$$\gamma_* = \frac{\varepsilon}{2\sum_{w\in OD} d_w \ln|P_w|}$$

и решать регуляризованную задачу (2.2.1) или (2.2.2) с $\gamma = \gamma_*$, но с точностью $\varepsilon/2$ [21, 23]. В этом случае оценка (22) из табл. 1 примет вид

$$\mathrm{O}\left( SHn \cdot \frac{d}{\varepsilon} \sqrt{H\bar{R}^2 \omega} \right), \qquad (2.2.24)$$

где $\omega = \max_{w \in OD} \ln |P_w|$ (в худшем случае можно ожидать $\omega \sim \sqrt{n}$). Сопоставляя оценки (2.2.23), (2.2.24), можно прийти к выводу, что введенная регуляризация (приводящая к сглаживанию двойственной задачи (2.2.3) [129]) оправдана. Однако не стоит забывать, что все приведенные оценки – это верхние оценки. Численные эксперименты показали, что эти оценки не всегда точны, особенно в части формулы (2.2.23).

Прежде всего заметим, что оценки общего числа арифметических операций из табл. 1 можно понимать как [стоимость итерации] × [число итераций]. А число итераций $N$, с некоторой натяжкой, можно просто считать пропорциональным относительной точности $\tilde{\varepsilon}$ в некоторой степени[26] $N \sim \tilde{\varepsilon}^{-\beta}$. Именно такая (относительная) точность, как правило, интересна на практике. Из теоретических верхних оценок (2.2.22) – (2.2.24) напрашивается вывод, что при $\gamma \gg \gamma_*$ имеем $N \sim \tilde{\varepsilon}^{-1/2}$, при $\gamma \simeq \gamma_*$ имеем $N \sim \tilde{\varepsilon}^{-1}$, а при $\gamma \to 0+$ имеем $N \sim \tilde{\varepsilon}^{-2}$. На самом деле, конечно, опущенные числовые множители тут также играют важную роль. Однако еще более важно то, что численные эксперименты [159] не подтвердили оценку $N \sim \tilde{\varepsilon}^{-2}$ при $\gamma \to 0+$. Более того, наблюдалась совсем другая зависимость $N \sim \tilde{\varepsilon}^{-1}$ [159], точнее говоря, наблюдалась зависимость $N \sim C_1 + C_2 \tilde{\varepsilon}^{-1}$ (см. рис. 1, 2).

Учитывая, что стоимость итерации заметно меньше в случае $\gamma \to 0+$ ($SHn \to Sn \ln n$, см. подраздел 2.2.6), то вывод об оправданности регуляризации приходится поставить под сомнение. Для графа г. Анахайм (Anaheim) [158] с $n \simeq 10^3$, $S \sim 40$ и $|OD| \sim 1.5 \cdot 10^3$ УМПТ, реализованный на ноутбуке с тактовой частотой 1.9 ГГц на языке программирования Python 2.4 [159], находил с относительной точностью $\tilde{\varepsilon} \simeq 0.01$ равновесие в модели Бэкмана ($\mu > 0, \gamma \to 0+$) приблизительно за 7 минут, в модели стабильной динамики ($\mu > 0, \gamma \to 0+$) за 10 минут (см. рис. 1, 2). Для соответствующей регуляризованной модели – заметно дольше [160].

---

[26] $\tilde{\varepsilon} = 1\% = 0.01$ – означает, что начальную невязку по функции или зазору двойственности необходимо уменьшить в 100 раз.

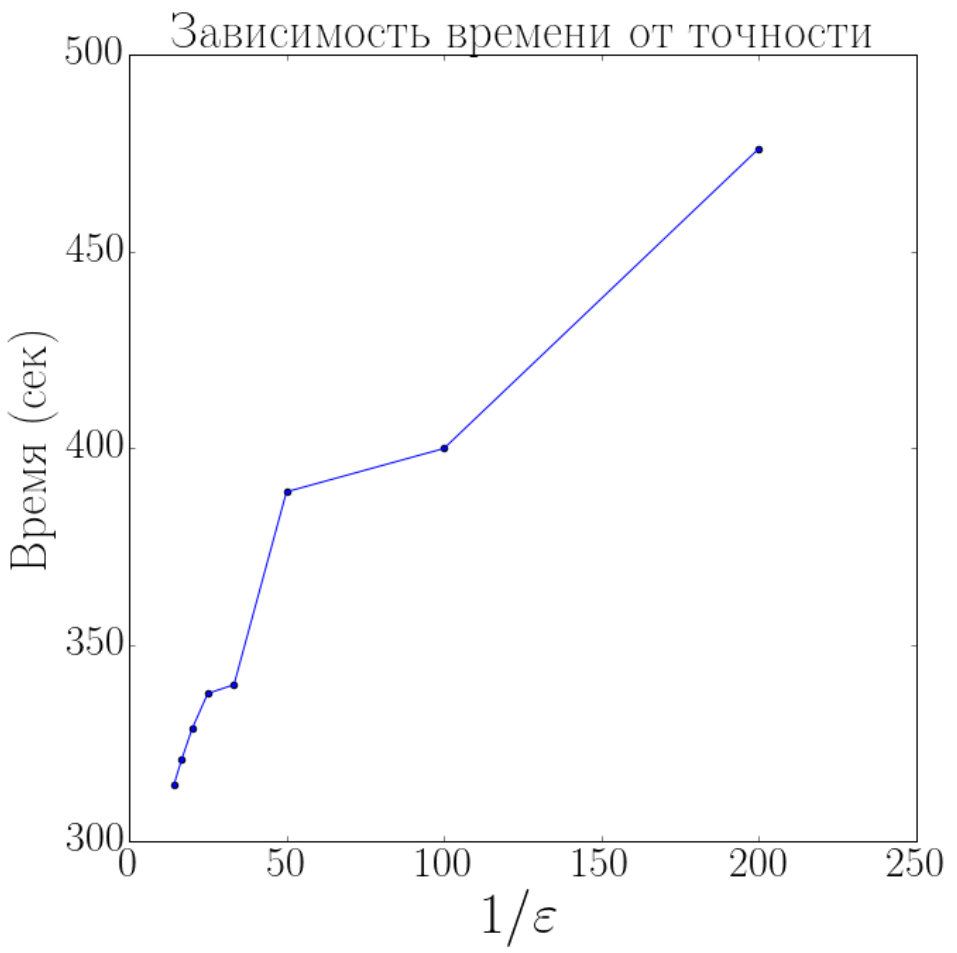

Рис. 1. Модель Бэкмана ( $\mu = 0.25$ ), $\gamma \to 0+$ ( $\tilde{\varepsilon}$ обозначено здесь через $\varepsilon$ )

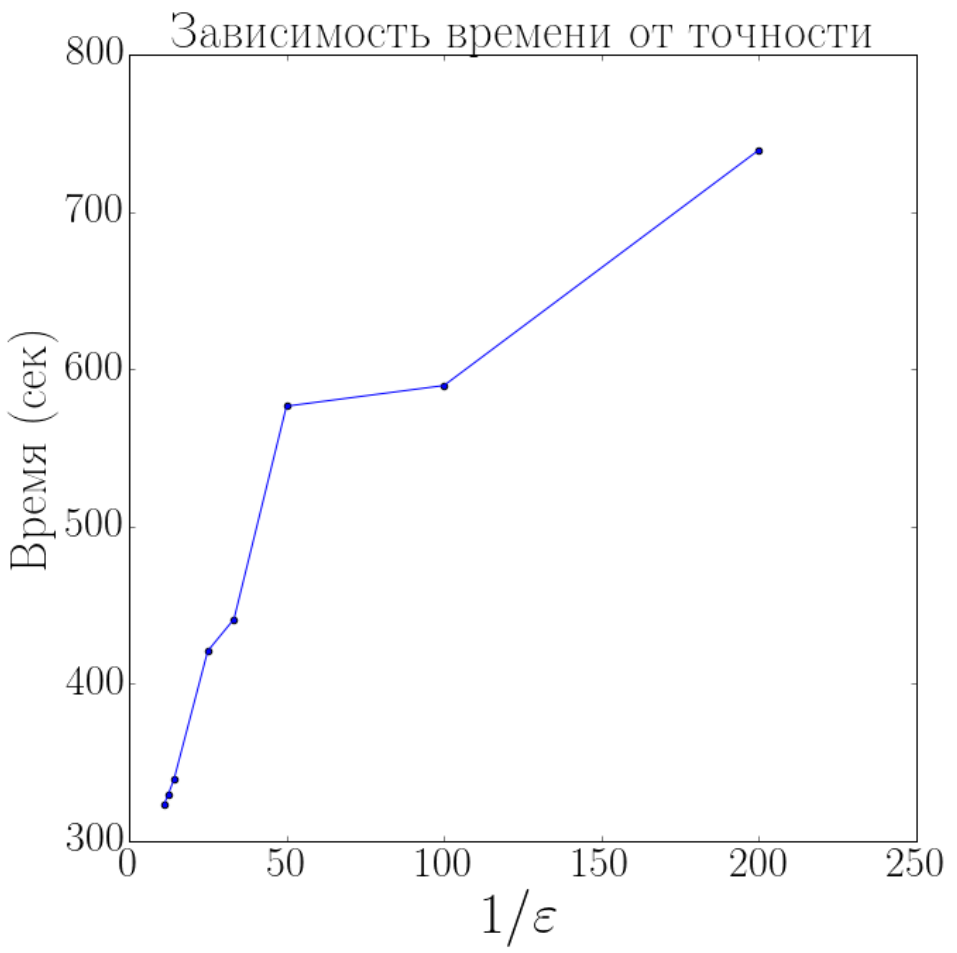

Рис. 2. Модель стабильной динамики ( $\mu \to 0+$ ),
$\gamma \to 0+$ ( $\tilde{\varepsilon}$ обозначено здесь через $\varepsilon$ )

На рис. 3, следуя [160], приведено сравнение того, как сходятся алгоритмы, описанные в подразделах 2.1.2, 2.1.3, 2.2.4–2.2.5, на примере задачи поиска равновесного распределения потоков по путям в модели

Бэкмана по г. Анахайм [158]. Отметим, что стоимость итерации у этих алгоритмов практически одинакова, поэтому рисунок демонстрирует также сравнительный анализ алгоритмов и во времени. Из рисунка видно, что наилучшие методы – это метод условного градиента (лучший), описанный в подразделе 2.1.2, и УПМТ (немного хуже), описанный в подразделах 2.2.4–2.2.5. Но в отличие от метода условного градиента УМПТ можно применять и для поиска равновесий в модели стабильной динамики, в смешанных моделях и для поиска стохастических равновесий в этих моделях. В следующем разделе мы увидим, что УМПТ можно применять и для поиска равновесий в многостадийных моделях.

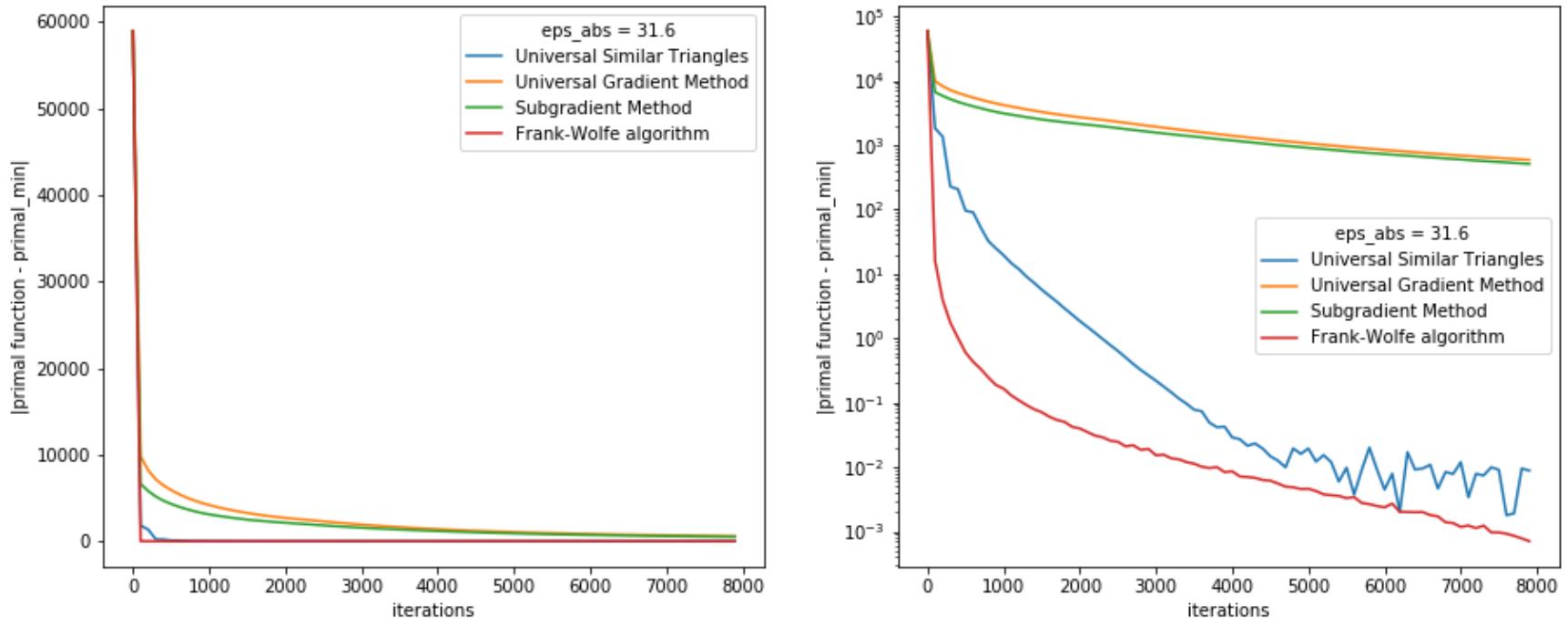

Рис. 3. Сравнение алгоритмов для модели Бэкмана $(\mu = 0.25)$, $\gamma \to 0+$

# 2.3. Поиск равновесий в многостадийных транспортных моделях

## 2.3.1. Введение

Поиск (стохастических) равновесий в многостадийных моделях транспортных потоков приводит к решению следующей седловой задачи с правильной (выпукло-вогнутой) структурой [5, 20, 24, 33] (см. также подраздел 1.1.10 раздела 1.1 гл. 1, а также разделы 1.3, 1.4 гл. 1):

$$\min_{\substack{\sum_{j=1}^{n} x_{ij} = L_i,\ \sum_{i=1}^{n} x_{ij} = W_j \\ x_{ij} \geq 0,\ i,j=1,\dots,n}} \max_{y \in Q} \left\{ \sum_{i,j=1}^{n} x_{ij} \ln x_{ij} + \sum_{i,j=1}^{n} c_{ij}(y) x_{ij} + g(y) \right\} =$$

$$= \max_{y \in Q} \max_{\lambda, \mu \in \mathbb{R}^n} \left\{ \langle \lambda, L \rangle + \langle \mu, W \rangle - \sum_{i,j=1}^{n} \exp\left( -c_{ij}(y) - 1 + \lambda_i + \mu_j \right) + g(y) \right\}, \qquad (2.3.1)$$

где $g(y)$ и $c_{ij}(y)\geq 0$ – вогнутые гладкие функции (если ищутся не стохастические равновесия, то $c_{ij}(y)$ могут быть негладкими), $Q$ – множество простой структуры, например,

$$Q=\{y:\ y\geq\bar{y}\}.$$

Легко понять, что система балансовых ограничений в (2.3.1) либо несовместна $\sum_{i=1}^{n}L_i\neq\sum_{j=1}^{n}W_j$, либо вырождена (имеет не полный ранг). В последнем случае это приводит к тому, что двойственные переменные $(\lambda,\mu)$ определены с точностью до произвольной постоянной $C$:

$$(\lambda+Ce,\mu-Ce),\ e=\underbrace{(1,\ldots,1)}_{n}.$$

Задачу (2.3.1) также можно переписать следующим образом (не ограничивая общности, считаем $\sum_{i=1}^{n}L_i=\sum_{j=1}^{n}W_j=1$):

$$\min_{\substack{\sum_{j=1}^{n}x_{ij}=L_i,\ \sum_{i=1}^{n}x_{ij}=W_j\\ x_{ij}\geq 0,\ i,j=1,\ldots,n;\ \sum_{i,j=1}^{n,n}x_{ij}=1}}\ \max_{y\in Q}\left\{\sum_{i,j=1}^{n}x_{ij}\ln x_{ij}+\sum_{i,j=1}^{n}c_{ij}(y)x_{ij}+g(y)\right\}=$$

$$=\max_{y\in Q}\max_{\lambda,\mu\in\mathbb{R}^n}\left\{\langle\lambda,L\rangle+\langle\mu,W\rangle-\ln\left(\sum_{i,j=1}^{n}\exp\left(-c_{ij}(y)+\lambda_i+\mu_j\right)\right)+g(y)\right\}=$$

$$=-\min_{y\in Q}f(y),\qquad(2.3.2)$$

где выпуклая функция $f(y)$ определяется как

$$f(y)=\max_{\substack{\sum_{j=1}^{n}x_{ij}=L_i,\ \sum_{i=1}^{n}x_{ij}=W_j\\ x_{ij}\geq 0,\ i,j=1,\ldots,n;\ \sum_{i,j=1}^{n,n}x_{ij}=1}}\left\{-\sum_{i,j=1}^{n}x_{ij}\ln x_{ij}-\sum_{i,j=1}^{n}c_{ij}(y)x_{ij}-g(y)\right\}=$$

$$=\min_{\lambda,\mu}\left\{\ln\left(\sum_{i,j=1}^{n}\exp\left(-c_{ij}(y)+\lambda_i+\mu_j\right)\right)-\langle\lambda,L\rangle-\langle\mu,W\rangle-g(y)\right\}.\qquad(2.3.3)$$

Поскольку мы добавили в ограничения условие $\sum_{i,j=1}^{n,n}x_{ij}=1$, являющееся следствием балансовых уравнений, то это привело к тому, что двойствен-ные переменные $(\lambda, \mu)$ определены с точностью до двух произвольных постоянных $C_\lambda$, $C_\mu$: $(\lambda + C_\lambda e, \mu + C_\mu e)$.

В данном разделе мы покажем, как можно решать задачу (2.3.2).

Заметим также, что расчет градиента $\nabla f(y)$ (в ряде транспортных приложений вогнутые функции $c_{ij}(y)$ – негладкие, тогда вместо градиентов стоит понимать суперградиенты $c_{ij}(y)$ и субградиент $f(y)$) осуществляется по следующей формуле (Демьянова–Данскина–Рубинова, см., например, [26, 33]):

$$\nabla f(y) = -\frac{\sum_{i,j=1}^{n} \exp\left(-c_{ij}(y) + \lambda_i^* + \mu_j^*\right) \nabla c_{ij}(y)}{\sum_{i,j=1}^{n} \exp\left(-c_{ij}(y) + \lambda_i^* + \mu_j^*\right)} - \nabla g(y) = \tag{2.3.4}$$

$$= -\sum_{i,j=1}^{n} x_{ij}\left(\lambda^*, \mu^*\right) \nabla c_{ij}(y) - \nabla g(y),$$

где $(\lambda^*, \mu^*)$ – решение задачи (2.3.3), не важно, какое именно, градиент $\nabla f(y)$ от выбора $C_\lambda$, $C_\mu$ (см. выше) не зависит. В данном разделе мы ограничимся изучением только полноградиентных методов для задачи (2.3.2), т. е. не будем рассматривать, например, рандомизацию при вычислении градиента по формуле (2.3.4). Планируется отдельно исследовать вопрос о возможности ускорения вычислений за счет введения рандомизации для внешней задачи. На данный момент нам представляется (см. формулу (2.3.8) в подразделе 2.3.2), что это может принести дивиденды только в случае, когда вспомогательная задача расчета $\nabla c_{ij}(y)$ достаточно сложная. Тут требуется много оговорок, в частности, в большинстве приложений умение рассчитывать $\nabla c_{ij}(y)$ для конкретной пары $(i, j)$ без дополнительных затрат позволяет заодно рассчитать и все $\nabla c_{ij}(y)$, $j = 1, ..., n$. Также отдельно планируется исследовать вопрос о том, какие подходы и насколько хорошо допускают распараллеливание. Вопрос о целесообразности рандомизации оказывается завязанным и на вопрос о возможности распараллеливания.

В основе подхода подраздела 2.3.2 лежит сочетание метода балансировки для решения внутренней задачи оптимизации по двойственным множителям и универсального метода с неточным оракулом для внешней задачи (2.3.2). В подразделе 2.3.3 делаются заключительные замечания.

## 2.3.2. Универсальный градиентный метод с неточным оракулом

Начнем с важного наблюдения о том, что внутренняя задача максимизации в (2.3.2) по $(\lambda,\mu)$ может быть явно решена по $\mu$ при фиксированном $\lambda$, и наоборот (это верно для задач (2.3.1) и (2.3.2), и приводит к одним и тем же формулам). Собственно, таким образом, получается метод балансировки расчета матрицы корреспонденций по энтропийной модели, см., например, [26] (тесно связанный с методом Синхорна [76, 86]), как метод простой итерации для явно выписываемых условий экстремума (принципа Ферма): $\lambda = \Lambda(\mu)$, $\mu = \mathrm{M}(\lambda)$.

Метод балансировки имеет вид ( $[\lambda]_0 = [\mu]_0 = 0$ ) [143]:

$$[\lambda_i]_{k+1} = -\ln\left(\frac{1}{L_i}\sum_{j=1}^{n}\exp\left(-c_{ij}(y) - 1 + [\mu_j]_k\right)\right),$$

$$[\mu_j]_{k+1} = -\ln\left(\frac{1}{W_j}\sum_{i=1}^{n}\exp\left(-c_{ij}(y) - 1 + [\lambda_i]_k\right)\right)$$

или

$$[\mu_j]_{k+1} = -\ln\left(\frac{1}{W_j}\sum_{i=1}^{n}\exp\left(-c_{ij} - 1 + [\lambda_i]_{k+1}\right)\right).$$

В этих формулах «–1» в экспоненте для метода (2.3.2) (в отличие от метода (2.3.1)) можно не писать, поскольку двойственные множители определяются неоднозначным образом с бóльшим произволом для задачи (2.3.2) (см. выше), достаточным для справедливости этого замечания.

Оператор $(\lambda,\mu) \to (\Lambda(\mu), \mathrm{M}(\lambda))$ является сжимающим в метрике Гильберта $\rho$ [100, 151]. Это означает, что после $N \sim \ln(\sigma^{-1})$ итераций метода балансировки можно получить такие $(\lambda_N, \mu_N)$, что ( $\{(\lambda_*(y), \mu_*(y))\}$ – двумерное аффинное множество решений, см. подраздел 2.3.1):

$$\rho\left((\lambda_N, \mu_N); \{(\lambda_*(y), \mu_*(y))\}\right) \le \sigma. \qquad (2.3.5)$$

Причем на практике наблюдается очень быстрая сходимость, т. е. коэффициент пропорциональности небольшой [26]. Таким образом, мы можем приближенно решить внутреннюю задачу.

Далее предлагается воспользоваться прямодвойственным (эта важно, поскольку нужно восстанавливать двойственные переменные) универ-сальным ускоренным (быстрым) градиентным методом [132] (можно использовать и универсальный метод подобных треугольников (УМПТ) [37], описанный в разделе 2.2 гл. 2 и в приложении 3) для решения внешней задачи оптимизации по $y$. К сожалению, в формулировке (2.3.1) (в отличие от формулировки (2.3.2)) кроме того, что внешняя задача гладкая (при условии гладкости $c_{ij}\left(y\right)$ [23, 27]), больше ничего о ней сказать нельзя (константа Липшица градиента не ограничена). Поэтому и по ряду других причин, о которых будет сказано далее, было отдано предпочтение универсальному методу, оптимально адаптивно настраивающемуся на гладкость функционала $f\left(y\right)$ на текущем участке пребывания итерационного процесса.[27] Однако нам потребуется использовать этот метод в варианте с неточным оракулом, выдающим градиент [28].

Напомним (см. подраздел 2.3.1), что мы решаем задачу (2.3.2), представимую в виде (здесь в max представлении $x = x$, $\bar{Q} = \left\{ x_{ij} \geq 0,\ i, j = 1,\ldots,n:\ \sum_{j=1}^{n} x_{ij} = L_i,\ \sum_{i=1}^{n} x_{ij} = W_j \right\}$, а в min представлении $x = \left(\lambda, \mu\right)$, $\bar{Q} = \mathbb{R}^{2n}$, см. формулу (2.3.3)):

$$f\left(y\right) = \max_{x \in \bar{Q}} \Psi\left(x, y\right) = \min_{x \in \bar{Q}} \Phi\left(x, y\right) \to \min_{y \in Q}.$$

Далее везде будем предполагать, что $y \in Q$.

**Определение 2.3.1 (см. гл. 4 [90]).** $\left(\delta, L\right)$*-оракул выдает* (*на запрос, в котором указывается только одна точка* $y$) *такие* $\left(F\left(y\right), G\left(y\right)\right)$*, что и для любых* $y, y' \in Q$:

---

[27]Бытует мнение, что любой универсальный метод должен чем-то платить за свою универсальность, и в этой связи возникает много вопросов, в частности: насколько дорога эта плата? В принципе, в статье [132] (см. также приложение 3) довольно подробно проясняется этот момент. Тем не менее мы повторим здесь соображения из [132]. Действительно, плата за универсальность есть. Универсальный метод из работы [132] может сделать где-то в четыре раза больше обращений к оракулу (что можно понимать как увеличения числа итераций в четыре раза) для задачи с более-менее одинаковой константой Липшица градиента во всей области (где довелось пройти итерационному процессу), по сравнению с обычным быстрым градиентным методом [129]. Тем не менее замечательная особенность универсального метода не только в том, что он настраивается на гладкость задачи и применим к любым задачам, но и в том, что (в отличие от подавляющего большинства методов) этот метод локально настраивается на гладкость функционала. И для сильно неоднородных функционалов типично, что универсальный метод делает заметно меньше итераций, чем, скажем, быстрый градиентный метод (плата за это уже учтена в отмеченном выше потенциально возможном увеличении числа итераций в четыре раза в худшем случае). Примеры, поясняющие сказанное, имеются в работе [132] и приложении 3.

$$0 \le f(y') - F(y) - \langle G(y), y' - y \rangle \le \frac{L}{2}\|y' - y\|^2 + \delta .$$

Из определения 2.3.1 сразу следует, что для любого $x \in Q$:

$$F(x) \le f(x) \le F(x) + \delta$$

и для любых $x, y \in Q$:

$$f(y) \ge f(y) - \langle G(x), y - x \rangle - \delta .$$

Из последнего свойства получаем, что определение $(\delta, L)$-оракула можно понимать как обобщение на гладкие задачи классического понятия негладкой выпуклой оптимизации: $\delta$-субградиента (см. п. 5 § 1 гл. 5 [57]). В приводимом далее утверждении в первой его части следует сохранить обозначения для задачи (2.3.2), (2.3.3); а во второй части утверждения следует обозначить $x = (\lambda, \mu)$, $y = y$ для задачи (2.3.2), (2.3.3). Таким образом, на задачу (2.3.2), (2.3.3) можно посмотреть с двух разных ракурсов, однако второй ракурс менее привлекателен ввиду необходимости рассмотрения ограниченных множеств $\bar{Q}$, что в интересующих нас приложениях место не имеет.

**Утверждение 2.3.1.** *Если* $\psi(y) = \max_{x \in \bar{Q}} \Psi(x, y)$, *где* $\Psi(x, y)$ *– выпуклая по* $y$ *и вогнутая по* $x$ *функция, и найден такой* $\tilde{x} \in \bar{Q}$, *что*

$$\psi(y) - \Psi(\tilde{x}, y) \le \delta ,$$

*то субградиент* $\partial_y \Psi(\tilde{x}, y)$ *есть* $\delta$*-субградиент функции* $\psi(y)$ *в точке* $y$.

*Если* $\varphi(y) = \min_{x \in \bar{Q}} \Phi(x, y)$, *где* $\Phi(x, y)$ *– выпуклая по совокупности переменных функция, и найден такой* $\tilde{x} \in \bar{Q}$, *что*

$$\max_{z \in \bar{Q}} \langle \Phi_x(\tilde{x}, y), \tilde{x} - z \rangle \le \delta ,$$

*то*

$$\Phi(\tilde{x}, y) - \varphi(y) \le \delta$$

*и субградиент* $\Phi_y(\tilde{x}, y) = \partial_y \Phi(\tilde{x}, y)$ *есть* $\delta$*-субградиент функции* $\varphi(y)$ *в точке* $y$.

**Доказательство.** Ограничимся доказательством только второй части этого утверждения. Доказательство первой части см. на с. 124 (лемма 13) книги [57].

Из выпуклости $\Phi(x, y)$ по совокупности переменных имеем

$$\Phi(x',y') \ge \Phi(x,y) + \langle \Phi_x(x,y), x'-x \rangle + \langle \Phi_y(x,y), y'-y \rangle. \qquad (2.3.6)$$

Определим зависимость $x(y)$ из соотношения

$$\varphi(y) = \min_{x \in \bar{Q}} \Phi(x,y) = \Phi(x(y), y).$$

Заметим, что $\Phi(\tilde{x}, y) \ge \varphi(y)$. Положим в (2.3.6) $x' = x(y')$, $x = \tilde{x}$. Тогда

$$\varphi(y') = \Phi(x',y') \ge \Phi(\tilde{x},y) + \langle \Phi_x(\tilde{x},y), x'-\tilde{x} \rangle + \langle \Phi_y(\tilde{x},y), y'-y \rangle \ge$$

$$\ge \varphi(y) + \langle \Phi_x(\tilde{x},y), x(y')-\tilde{x} \rangle + \langle \Phi_y(\tilde{x},y), y'-y \rangle \ge \varphi(y) + \langle \Phi_y(\tilde{x},y), y'-y \rangle - \delta.$$

В последней формуле мы использовали, что

$$\langle \Phi_x(\tilde{x},y), \tilde{x} - x(y') \rangle \le \delta.$$

В свою очередь, из выпуклости $\Phi(x,y)$ по $x$ (для всех допустимых $y$) имеем

$$\Phi(\tilde{x},y) - \Phi(x(y'),y) \le \langle \Phi_x(\tilde{x},y), \tilde{x} - x(y') \rangle.$$

Беря в этой формуле $y' = y$ и воспользовавшись определением $x(y)$, получаем, что

$$\Phi(\tilde{x},y) - \varphi(y) \le \langle \Phi_x(\tilde{x},y), \tilde{x} - x(y) \rangle. \ \blacksquare$$

Однако не хочется довольствоваться возможностью находить только $\delta$-субградиент (из утверждения 2.3.1 эта возможность очевидна), поскольку в определенных ситуациях явно можно рассчитывать на некоторую гладкость итоговой (внешней) задачи (2.3.2). Понятие $(\delta, L)$-оракула в некотором смысле налагает наиболее слабые условия на возможные неточности в вычислении функции и градиента, при которых можно рассчитывать, что скорость сходимости метода, учитывающего гладкость (липшицевость градиента функционала) задачи, сильно не пострадает (см. теорему 2.3.1 ниже).

На первый взгляд, может показаться, что применимость описанной концепции $(\delta, L)$-оракула к задаче (2.3.1) следует из следующего результата (см. п. 4.2.2 [90]).

**Утверждение 2.3.2.** *Пусть подзадача энтропийно-линейного программирования* (*ЭЛП*) *в* (2.3.2) *решена* (*по функции*) *с точностью* $\delta$, *т. е. найден такой* $\tilde{x}(c)$, *удовлетворяющий балансовым ограничениям, что*

$$\sum_{i,j=1}^{n} \tilde{x}_{ij}(c) \ln \tilde{x}_{ij}(c) + \sum_{i,j=1}^{n} c_{ij} \tilde{x}_{ij}(c) - \min_{\substack{\sum_{j=1}^{n} x_{ij} = L_i, \ \sum_{i=1}^{n} x_{ij} = W_j \\ i,j=1,\dots,n}} \left\{ \sum_{i,j=1}^{n} x_{ij} \ln x_{ij} + \sum_{i,j=1}^{n} c_{ij} x_{ij} \right\} \le \delta.$$

*Тогда для функции*

$$\bar{f}(c) = -\min_{\substack{\sum_{j=1}^{n} x_{ij} = L_i,\ \sum_{i=1}^{n} x_{ij} = W_j \\ i,j=1,\dots,n}} \left\{ \sum_{i,j=1}^{n} x_{ij} \ln x_{ij} + \sum_{i,j=1}^{n} c_{ij} x_{ij} \right\}$$

*набор*

$$-\left( \sum_{i,j=1}^{n} \tilde{x}_{ij}(c) \ln \tilde{x}_{ij}(c) + \sum_{i,j=1}^{n} c_{ij} \tilde{x}_{ij}(c), \left\{ \tilde{x}_{ij}(c) \right\}_{i,j=1}^{n,n} \right)$$

*является* $\left(\delta, 2 \cdot \max\limits_{i,j=1,\dots,n} c_{ij}\right)$*-оракулом.*

К сожалению, большинство методов (в том числе метод балансировки) не удовлетворяют одному пункту утверждения 2.3.2, а именно, они выдают вектор $\tilde{x}$, который лишь приближенно удовлетворяет балансовым ограничениям (в утверждении требование точного удовлетворения балансовых ограничений является существенным и не может быть как-то равнозначно релаксировано). Связанно это с тем, что для задачи ЭЛП, когда ограничений намного меньше числа прямых переменных, обычно решается двойственная задача, по которой восстанавливается решение прямой задачи [26, 71]. Как следствие, приобретается невязка и в ограничениях. Собственно, в представлении градиента функционала по формуле (2.3.4) имеются два способа. Первый – через двойственные множители $(\lambda, \mu)$, второй – через решение прямой задачи $x$. Функционал прямой задачи сильно выпуклый по $x$, поскольку энтропия 1-сильно выпуклая функция в 1-норме [129]. Поэтому сходимость в решении прямой задачи по функции обеспечивает сходимость и по аргументу, что и означает возможность определения с хорошей точностью градиента по формуле (2.3.4) через $x$. Другая ситуация возникает, если смотреть на двойственную задачу к задаче ЭЛП (в приводимом далее утверждении следует обозначить $x = (\lambda, \mu)$, $y = y$ для задачи (2.3.2), (2.3.3)).

**Утверждение 2.3.3.** *Пусть* $\varphi(y) = \min\limits_{x \in Q} \Phi(x, y)$, *где* $\Phi(x, y)$ *– такая достаточно гладкая, выпуклая по совокупности переменных функция, что*[28]

$$\left\| \nabla\Phi(x', y') - \nabla\Phi(x, y) \right\|_2 \le L \left\| (x', y') - (x, y) \right\|_2 .$$

---

[28]Это утверждение имеет достаточно простую геометрическую интерпретацию. Проекция надграфика выпуклой функции будет выпуклым множеством, то есть, в свою очередь, надграфиком некоторой выпуклой функции. Кривизна границы у полученного при проектировании множества будет не больше, чем была у исходного множества. Это следует из того, что проектирование – сжимающий оператор.

*Пусть для произвольного* $y \in Q$ (*считаем, что множество* $Q$ *содержит внутри себя шар радиуса более* $\sqrt{2\delta/L}$) *можно найти такой* $\tilde{x} = \tilde{x}(y) \in \bar{Q}$, *что*

$$\max_{z \in \bar{Q}} \left\langle \Phi_x(\tilde{x}, y), \tilde{x} - z \right\rangle \le \delta .$$

*Тогда*

$$\Phi(\tilde{x}, y) - \varphi(y) \le \delta ,$$
$$\left\| \nabla\varphi(y') - \nabla\varphi(y) \right\|_2 \le L \left\| y' - y \right\|_2 ,$$

*и* $\left(\Phi(\tilde{x}, y) - 2\delta, \Phi_y(\tilde{x}, y)\right)$ *будет* $(6\delta, 2L)$*-оракулом для* $\varphi(y)$ *на выпуклом множестве, полученном из множества* $Q$ *отступанием от границы* $\partial Q$ *во внутрь* $Q$ *на расстояние* $\sqrt{2\delta/L}$ (*по условию это множество не пусто*).

**Доказательство.** По условию задачи имеем при всех допустимых значениях аргументов $\Phi$ :

$$\lambda_{\max}\left( \left\| \begin{matrix} \Phi_{xx} & \Phi_{xy} \\ \Phi_{yx} & \Phi_{yy} \end{matrix} \right\| \right) = \sup_{\|h\|_2 \le 1} \left\langle h, \left\| \begin{matrix} \Phi_{xx} & \Phi_{xy} \\ \Phi_{yx} & \Phi_{yy} \end{matrix} \right\| h \right\rangle \le L . \qquad (2.3.7)$$

Заметим, что также по условию при всех допустимых значениях аргументов $\Phi$ :

$$\left\| \begin{matrix} \Phi_{xx} & \Phi_{xy} \\ \Phi_{yx} & \Phi_{yy} \end{matrix} \right\| \succ 0,\ \Phi_{xx} \succ 0,\ \Phi_{yy} \succ 0,\ \Phi_{yx} = \Phi_{xy}^T,\ \Phi_{xx} = \Phi_{xx}^T,\ \Phi_{yy} = \Phi_{yy}^T .$$

Для упрощения последующих рассуждений (в частности, чтобы не работать с псевдообратными матрицами) будем, считать, что матрица $\Phi_{xx} \succ 0$ положительно определена (исходя из условий, гарантировать можно лишь неотрицательную определенность). Также будем считать (в интересующем нас приложении к задаче (2.3.2) это имеет место), что зависимость $x(y)$, определяемая из соотношения

$$\varphi(y) = \min_{x \in \bar{Q}} \Phi(x, y) = \Phi(x(y), y)$$

однозначным образом, и удовлетворяет соотношению

$$\Phi_x(x(y), y) \underset{y}{\equiv} 0 ,$$

из которого имеем

$$\Phi_{xx}\left(x(y),y\right)\left\|\frac{\partial x(y)}{\partial y}\right\|+\Phi_{xy}\left(x(y),y\right)\underset{y}{\equiv}0,$$

т. е.

$$\left\|\partial x/\partial y\right\|=\left\|\partial x_i/\partial y_j\right\|=-\Phi_{xx}^{-1}\Phi_{xy}.$$

Поскольку $\varphi(y)=\Phi\left(x(y),y\right)$, то

$$\varphi_{yy}=\left\|\partial x/\partial y\right\|^T\Phi_{xx}\left\|\partial x/\partial y\right\|+\left\|\partial x/\partial y\right\|^T\Phi_{xy}+\Phi_{yx}\left\|\partial x/\partial y\right\|+\Phi_{yy}=\Phi_{yy}-\Phi_{yx}\Phi_{xx}^{-1}\Phi_{xy}.$$

С учетом этой формулы и из формулы дополнения по Шуру [156] получаем

$$\begin{Vmatrix}\Phi_{xx} & \Phi_{xy}\\ \Phi_{yx} & \Phi_{yy}\end{Vmatrix}=\begin{Vmatrix}E_x & 0\\ \Phi_{yx}\Phi_{xx}^{-1} & E_y\end{Vmatrix}\begin{Vmatrix}\Phi_{xx} & 0\\ 0 & \varphi_{yy}\end{Vmatrix}\begin{Vmatrix}E_x & \Phi_{xx}^{-1}\Phi_{xy}\\ 0 & E_y\end{Vmatrix},$$

где $E_x$, $E_y$ – единичные матрицы соответствующих размеров. Поскольку

$$\begin{Vmatrix}E_x & 0\\ \Phi_{yx}\Phi_{xx}^{-1} & E_y\end{Vmatrix}=\begin{Vmatrix}E_x & \Phi_{xx}^{-1}\Phi_{xy}\\ 0 & E_y\end{Vmatrix}^T$$

и эти матрицы полного ранга, то из (2.3.7) имеем, что

$$\sup_{\|h\|_2\le 1}\left\langle h,\varphi_{yy}\,h\right\rangle=\lambda_{\max}\left(\varphi_{yy}\right)\le\max\left\{\lambda_{\max}\left(\Phi_{xx}\right),\lambda_{\max}\left(\varphi_{yy}\right)\right\}=\lambda_{\max}\left(\begin{Vmatrix}\Phi_{xx} & \Phi_{xy}\\ \Phi_{yx} & \Phi_{yy}\end{Vmatrix}\right)\le L.$$

Таким образом, установлено, что

$$\varphi(y)\le\varphi(x)+\left\langle\nabla\varphi(x),y-x\right\rangle+\frac{L}{2}\|y-x\|_2^2.$$

Согласно утверждению 2.3.1

$$\varphi(y)\ge\varphi(x)+\left\langle\Phi_y\left(\tilde{x},y\right),y-x\right\rangle-\delta.$$

Далее проведем рассуждения, аналогичные рассуждениям на с. 107 (и немного отлично от с. 115) диссертации [90]. Вычитая из первого неравенства второе, получим

$$\left\langle\Phi_y\left(\tilde{x},y\right)-\nabla\varphi(x),y-x\right\rangle\le\frac{L}{2}\|y-x\|_2^2+\delta.$$

Положим ($t>0$):

$$y-x=t\frac{\Phi_y\left(\tilde{x},y\right)-\nabla\varphi(x)}{\left\|\Phi_y\left(\tilde{x},y\right)-\nabla\varphi(x)\right\|_2},$$

получим

$$\left\| \Phi_y(\tilde{x}, y) - \nabla\varphi(x) \right\|_2 \le \frac{Lt}{2} + \frac{\delta}{t}.$$

Минимизируя правую часть неравенства по $t > 0$, получим (при $t = \sqrt{2\delta/L}$), что

$$\left\| \Phi_y(\tilde{x}, y) - \nabla\varphi(x) \right\|_2 \le \sqrt{2\delta L}.$$

Отсюда и из утверждения 2.3.1 имеем, что

$$\varphi(y) \le \varphi(x) + \langle \nabla\varphi(x), y - x \rangle + \frac{L}{2}\|y - x\|_2^2 \le$$

$$\le \varphi(x) + \langle \Phi_y(\tilde{x}, y), y - x \rangle + \sqrt{2\delta L}\|y - x\|_2 + \frac{L}{2}\|y - x\|_2^2 \le$$

$$\le \Phi(\tilde{x}, y) - 2\delta + \langle \Phi_y(\tilde{x}, y), y - x \rangle + \sqrt{2\delta L}\|y - x\|_2 + \frac{L}{2}\|y - x\|_2^2 + 2\delta \le$$

$$\le \Phi(\tilde{x}, y) - 2\delta + \langle \Phi_y(\tilde{x}, y), y - x \rangle + \frac{2L}{2}\|y - x\|_2^2 + 6\delta.$$

С учетом того, что (см. утверждение 2.3.1)

$$\varphi(y) \ge \varphi(x) + \langle \Phi_y(\tilde{x}, y), y - x \rangle - \delta \ge \Phi_y(\tilde{x}, y) - 2\delta + \langle \Phi_y(\tilde{x}, y), y - x \rangle,$$

из определения 2.3.1 получаем доказываемое утверждение. ■

Это утверждение позволяет установить гладкость задачи (2.3.2), (2.3.3), если гладким оказывается $c_{ij}(y)$.

К сожалению, практическое применение утверждения 2.3.3 натыкается на следующие сложности:

1) необходимости отступать от границы множества $Q$ во внутрь на $\sqrt{2\delta/L}$,

2) корректности записи (см. доказательство утверждения 2.3.3)

$$\left\| \partial x_i / \partial y_j \right\| = -\Phi_{xx}^{-1}\Phi_{xy},$$

3) необходимости предположения о компактности множества $\bar{Q}$, иначе невозможно будет добиться выполнения условия

$$\max_{z \in \bar{Q}} \langle \Phi_x(\tilde{x}, y), \tilde{x} - z \rangle \le \delta.$$

Сложность 1), как правило, на практике преодолима за счет возможности доопределения функционала задачи с сохранением всех свойств на $\sqrt{2\delta/L}$-окрестность множества $Q$ (заметим, что доопределение часто не требуется, поскольку функционал и так задан «с запасом»).

Например, для рассматриваемых нами транспортных приложений с $Q=\{y:\ y\geq\bar{y}\}$ это возможно [5, 20, 24, 33]. Сложность 2) часто вообще не возникает (разве что оговорка о существовании $\Phi_{xx}^{-1}$, впрочем, приведенные выше рассуждения можно провести, сохранив все результаты в идентичном виде, так что эта оговорка будет не нужна), поскольку $\bar{Q}$ совпадает со всем (двойственным) пространством. А вот сложность 3), действительно, портит дело. К сожалению, простых теоретически обоснованных способов борьбы с этой сложностью мы пока не знаем. Тем не менее полезно заметить, что в действительности нужно гарантировать выполнение (см. доказательство утверждения 2.3.1)

$$\left\langle \Phi_x\left(\tilde{x}(y),y\right),\tilde{x}(y)-x(y')\right\rangle\leq\delta,$$

где точки $y$ и $y'$ близки, поскольку возникают на соседних итерациях внешнего метода. С учетом ожидаемой «близости» $\tilde{x}=\tilde{x}(y)$ и $x(y)$ мы можем заменить в этом критерии настоящее множество $\bar{Q}$, которое, как правило, совпадает со всем пространством, на шар конечного радиуса. Более детальные исследования (для задачи (2.3.2), (2.3.3)) и практические эксперименты показывают, что для выполнения приведенного выше условия достаточно обеспечить для внутреннего итерационного процесса $\{x_k\}\to x(y)$ условия

$$\left\|\Phi_x\left(x_k,y\right)\right\|_2\left\|x_k\right\|_2\leq\delta/2,\ \left\|\Phi_x\left(x_k,y\right)\right\|_2\leq\delta.$$

Соответствующее $x_k=\left(\lambda_k,\mu_k\right)$ порождает нужное $\tilde{x}(y)=x_k$. С учетом специфики рассматриваемой нами задачи (2.3.2), (2.3.3) имеем следующий критерий (возвращаемся к обозначениям (2.3.2), (2.3.3)):

$$\left\|Ax\left(\lambda_k,\mu_k\right)-b\right\|_2\left\|\left(\lambda_k,\mu_k\right)\right\|_2\leq\delta/2,\left\|Ax\left(\lambda_k,\mu_k\right)-b\right\|_2\leq\delta,$$

где $x\left(\lambda_k,\mu_k\right)$ определяется в формуле (2.3.4), а введённая линейная система балансовых уравнений $Ax=b$ есть общая запись аффинных (транспортных) ограничений:

$$\sum_{j=1}^{n}x_{ij}=L_i,\sum_{i=1}^{n}x_{ij}=W_j,\ i,j=1,\dots,n.$$

В связи со сказанным выше заметим, что (это следует из оценки (2.3.5)) метод балансировки обеспечивает сходимость и по аргументу, что для других методов (без введения регуляризации) решения двойственной задачи, вообще говоря, нельзя гарантировать. Это свойство наряду с линейной скоростью сходимости метода (со скоростью геометрической про-грессии) позволяет надеяться, что выбранный критерий является достаточно точным (точнее не слишком грубым).

Принципиально важно для гладкого случая ($c_{ij}(y)$ – функции с липшицевым градиентом), как это будет следовать из дальнейших оценок (см. теорему 2.3.1), не просто уметь решать двойственную задачу, т. е. находить $(\lambda, \mu)$ так, чтобы была сходимость по аргументу, а делать это так, чтобы сложность решения задачи зависела от точности ее решения логарифмическим образом. Выше мы отмечали, что это имеет место для метода балансировки. Также это имеет место и для быстрых методов, применяемых к регуляризованной двойственной задачи. При фиксации параметра регуляризации, исходя из итоговой желаемой точности, быстрые градиентные методы (для сильно выпуклых функций) решают регуляризованную двойственную задачу так, что зависимость сложности от точности ее решения логарифмическая.

Хочется, чтобы при решении внешней задачи в (2.3.2), т. е. задачи

$$\min_{y \in Q} f(y),$$

можно было не задумываться ни о какой гладкости. Если она есть, то метод бы это хорошо учитывал, не требуя знания констант Липшица градиента, если ее нет, то метод также бы работал оптимальным (для негладкого случая) образом. Именно таким свойством и обладает универсальный метод [132], работающий и в концепции неточного оракула [28, 111] (см. определение 2.3.1).

Заметим [132], что можно погрузить задачу с гёльдеровым градиентом ($\nu \in [0,1]$):

$$\left\| \nabla f(y') - \nabla f(y) \right\|_* \le L_\nu \left\| y' - y \right\|^\nu$$

(в том числе и негладкую задачу с ограниченной нормой разности субградиентов при $\nu = 0$) в класс гладких задач с оракулом, характеризующимся точностью $\delta$ и (см. также приложение 3)

$$L = L_\nu \left[ \frac{L_\nu (1-\nu)}{2\delta(1+\nu)} \right]^{\frac{1-\nu}{1+\nu}}.$$

Это позволяет даже в случае, когда можно рассчитывать только на $\delta$-субградиент[29] (с ограниченной нормой субградиента (разности субградиен-тов), причем какой именно константой ограниченной, методу знать не обязательно), все равно работать в концепции $(\delta, L)$-оракула.

[29]На $\delta$-субградиент всегда можно рассчитывать согласно утверждению 2.3.1. Причем, как уже отмечалось раньше, для получения $\delta$-субградиента не нужна сходимость по аргументу для вспомогательной задачи.

Итак, у нас есть внешняя задача (2.3.2):

$$\min_{y\in Q} f(y),$$

для которой обращение к $(\delta, L)$-оракулу за значением функции и градиента стоит $\sim \ln\left(\delta^{-1}\right)$. Насколько быстро мы можем решить такую задачу, т. е. при каком $N(\varepsilon)$ можно гарантировать, что

$$f\left(y_{N(\varepsilon)}\right) - \min_{y\in Q} f(y) \le \varepsilon\,?$$

Ответ можно получить из следующего результата.

**Теорема 2.3.1 (см. [28, 29, 120, 132]).** *Существует однопараметрическое семейство универсальных градиентных методов* (*параметр* $p\in[0,1]$), *не получающих на вход, кроме* $p$, *больше никаких параметров* (*в частности, не использующих значения* $L_\nu$ и $R$ – «*расстояние*» *от точки старта до решения, априорно не известное*), *которое приводит к следующей оценке на требуемое число итераций*:

$$N_p(\varepsilon) = \mathrm{O}\left(\inf_{\nu\in[0,1]}\left(\frac{L_\nu R^{1+\nu}}{\varepsilon}\right)^{\frac{2}{1+2p\nu+\nu}}\right),$$

*если* $\delta \le \mathrm{O}\left(\varepsilon / N_p(\varepsilon)^p\right)$.

В работе [111] удалось предложить такой метод, который и параметр $p$ адаптивно выбирает.

Из теоремы 2.3.1 можно заключить, что если мы рассчитываем на некоторую гладкость $f(y)$, то стоит выбирать значение параметра $p=1$, при этом общие трудозатраты машинного времени будут

$$\mathrm{O}\left(N_1(\varepsilon)\left(T\ln\left(\varepsilon^{-1}\right)+\tilde{T}\right)\right), \tag{2.3.8}$$

где $\tilde{T}$ – время вычисления (суб-)градиента функционала (в основном это вычисления $\left\{\nabla c_{ij}(y)\right\}_{i,j=1}^{n,n}$ [23, 27]), $T$ – время решения вспомогательной задачи методом балансировки с относительной точностью 1%. Численные эксперименты показывают, что на одном современном ноутбуке при $n\sim 10^2$ время $T\approx 1$ с [28], что сопоставимо со временем $\tilde{T}$ для таких $n$ [27].

Заметим, что УМПТ из подразделов 2.2.5–2.2.6 как раз и является универсальным методом, отвечающим $p = 1$.

Выгода от описанной выше конструкции по сравнению с обычным способом решения исходной задачи минимизации (2.3.2), (2.3.3) сразу по совокупности всех переменных (см., например, [24]) заключается в гарантированном неувеличении константы Липшица градиента в оценке необходимого числа итераций (см. утверждение 2.3.3) и ожидаемом уменьшении в этой же оценке «расстояния» от точки старта до (неизвестного априори) решения. Выгода здесь вполне может достигать одного порядка и более. При этом можно ожидать лишь незначительного увеличения стоимости одной итерации. Причем стоит иметь в виду, что при оптимизации сразу по всем переменным требуется рассчитывать градиент функционала по большему числу переменных, чем в описанном выше подходе, что также играет нам на пользу. В конечном итоге, сокращение числа итераций заметно превалирует над небольшим увеличением стоимости одной итерации.

Резюмируем ключевой результат этого раздела (и всего раздела) следующим образом.

*Для решения задачи* (2.3.2) *предлагается использовать универсальный метод из работы* [132] (*а точнее его модификацию из* [28, 111]). *Если рассчитываем на гладкость*[30] $f(y)$, *то полагаем в методе* $p = 1$. *Если на гладкость рассчитывать не приходится*[31], *то полагаем* $p = 0$. *В обоих случаях, кроме априорной подсказки относительно параметра* $p$, *методу больше ничего от нас знать не надо!*

### 2.3.3. Заключительные замечания

В приложениях часто возникают задачи, имеющие следующий вид (см., например, [31, 79]):

$$f(x) = \Phi(x, y(x)) \to \min_x, \qquad (2.3.9)$$

при этом $y(x)$ и $\nabla y(x)$ могут быть получены из решения отдельной подзадачи лишь приближенно. Довольно типично, что существует способ, который выдает $\varepsilon$-приближенное решение за время, зависящее от $\varepsilon$, логарифмическим образом $\sim \ln\left(\varepsilon^{-1}\right)$. В данном разделе намечен общий способ решения таких задач. Его наиболее важной отличительной чертой является адаптивность (самонастраиваемость), т. е. методу на вход не надо подавать никаких констант Липшица (более того, метод будет работать и в негладком случае). Метод сам настраивается локально на оптимальную гладкость функции. Это свойство метода делает его привилегированным, поскольку в реальных приложениях, чтобы что-то знать о свойствах $f(x)$, нужно что-то знать о свойствах зависимости $y(x)$, а это часто недоступно по постановке задачи или, при попытке оценить, приводит к сильно завышенным оценкам.

[30]В этом случае как раз существенна логарифмическая сложность приближенного вычисления градиента и значения функции $f(y)$ от точности, обеспеченная методом балансировки.

[31]В этом случае точность решения вспомогательной задачи расчета $\delta$ -субградиента можно завязать на желаемую точность решения задачи (2.3.2) $\varepsilon$ по формуле $\delta = O(\varepsilon)$ (см. теорему 2.3.1 при $\nu = 0$) с константой порядка 1.

В задаче (2.3.9) важно уметь эффективно пересчитывать значения $y(x)$, а не рассчитывать их каждый раз заново (на каждой итерации внешнего цикла). Поясним сказанное. Предположим, что мы уже как-то посчитали, скажем, $y(x)$, решив, например, с какой-то точностью соответствующую задачу оптимизации. Тогда для вычисления $y(x+\Delta x)$ (на следующей итерации внешнего цикла) у нас будет хорошее начальное приближение $y(x)$. А, как известно (см., например, [31]), расстояние от точки старта до решения (не в случае сходимости метода со скоростью геометрической прогрессии или быстрей) существенным образом определяет время работы алгоритма. Тем не менее известные нам приложения (см., [31, 79]) пока как раз всецело соответствуют сходимости процедуры поиска $y(x)$ со скоростью геометрической прогрессии. Связано это с тем, что если расчет $y(x)$ с точностью $\varepsilon$ осуществляется за $O\left(\ln\left(\varepsilon^{-1}\right)\right)$ операций, то для внешней задачи можно выбирать самый быстрый метод (а стало быть, и самый требовательный к точности), и с точностью до того, что стоит под логарифмом, общая трудоемкость будет рассчитываться по формуле, аналогичной формуле (2.3.8). Как правило, такое сочетание оказывается недоминируемым. Здесь мы ограничимся ссылкой на пример 4 и последующий текст из работы [31] и общим тезисом, который пока неплохо подтверждался на практике (и частично в теории [102]):

*если есть возможность в задаче оптимизации (в седловой задаче) явно прооптимизировать по части переменных (или как-то эффективно это сделать с хорошей точностью), то, как правило, это и надо сделать, и строить итерационный метод исходя из этого.*

В реальных транспортных приложениях [20, 23, 27, 33] достаточно сложным является расчет $c_{ij}(y)$ и их градиентов (особенно при поиске стохастических равновесий). Тем не менее эти задачи имеют вполне четкую привязку к решению некоторых задач на графах типа поиска кратчайших путей (см. раздел 2.1 этой главы). Так же, как и в предыдущем абзаце для внутренней задачи, для внешней задачи можно не рассчитывать $c_{ij}(y)$ и их градиенты каждый раз заново, а пересчитывать; также можно допускать неточность в их вычислении (и ненулевую вероятность ошибки), надеясь на ускорение (благо метод работает в концепции неточного оракула, природа которой не принципиальна, см. [90]). Также здесь оказываются полезными идеи БАД (быстрого автоматического дифференцирования [44, 49]), которые позволяют практически за то же время, что занимает вычисление самих функций, вычислять их градиенты.

# Приложение 1. Теория макросистем с точки зрения стохастической химической кинетики

## 1.1. Введение

В данном приложении мы постараемся пояснить важность формализма стохастической химической кинетики в изучении равновесий макросистем. Под макросистемой мы понимаем систему большого количества случайно взаимодействующих агентов. Мы будем интересоваться поведением такой системы на больших временах, в частности, необходимыми и достаточными условиями, при которых на больших временах такая система сходится к равновесию. При этом под равновесием будем понимать такое макросостояние, в малой окрестности которого концентрируется стационарная мера (динамику считаем марковской). Имеется большое количество конкретных примеров макросистем [13, 17, 18, 19, 32, 65, 93, 115, 117, 138, 152], встречающихся в самых разнообразных предметных областях (физика, биология, экономика, транспортное моделирование и т. д.). Однако в подавляющем большинстве случаев при описании макросистемы не описывается ее эволюция (считается, что система уже находится в стационарном состоянии), не исследуется концентрация в окрестности равновесия. В данном приложении будет изучаться эволюция макросистем и концентрация возникающих при этом стационарных мер. В отличие от подавляющего большинства работ, мы будем допускать, что макросистема характеризуется вектором чисел заполнений, размерность которого может зависеть от числа агентов.

Структура приложения следующая. В разделе 1.2 мы приведем три игрушечных примера макросистем, изучение равновесий в которых может быть единообразно осуществлено. Пункт г) теоремы из раздела 1.3 является оригинальным и представляющим наибольший интерес в этом приложении.

## 1.2. Примеры макросистем

Далее мы рассматриваем три хорошо известных примера макросистем. Однако способ преподнесения материала представляется оригинальным. Также отметим, что приводимые ниже количественные оценки мы ранее не встречали в литературе. Эти оценки получаются исходя из формализма, описанного в разделе 1.3, применимого к большому количеству макросистем.

**Пример 1 (кинетика социального неравенства и предельные формы).** В некотором городе живет $N \gg 1$ жителей (четное число). В начальный момент у каждого жителя имеется по $\bar{s}$ монеток. Каждый день жители случайно разбиваются на пары. В каждой паре жители скидываются по монетке (если один или оба участника банкроты, то банкрот не скидывается, в то время как небанкрот в любом случае обязан скинуть монетку). Далее в каждой паре случайно разыгрывается победитель, который и забирает «призовой фонд». Обозначим через $c_s(t)$ – долю жителей города, у которых ровно $s$, $s = 0,...,\bar{s}N$, монеток на $t$-й день. Имеют место следующие оценки:

$$\exists\ a > 0:\ \forall\ \sigma > 0, t \ge a(\bar{s})\ln N \to P\left( \left\| c(t) - c^* \right\|_2 \ge \frac{2\sqrt{2} + 4\sqrt{\ln\left(\sigma^{-1}\right)}}{\sqrt{N}} \right) \le \sigma,$$

$$\exists\ b, D > 0:\ \forall\ \sigma > 0, t \ge b(\bar{s})\ln N \to P\left( \left\| c(t) - c^* \right\|_1 \ge D\sqrt{\frac{\ln^2 N + \ln\left(\sigma^{-1}\right)}{N}} \right) \le \sigma,$$

где $c_s^* \simeq C\exp(-s/\bar{s})$, а $C \simeq 1/\bar{s}$ находится из условия нормировки $\sum_{s=0}^{\bar{s}N} C\exp(-s/\bar{s}) = 1$. Таким образом, кривая (предельная форма [17], равновесие макросистемы [93]) $C\exp(-s/\bar{s})$ характеризует распределение населения по богатству на больших временах. Этот результат восходит к работам итальянского экономиста Вильфредо Парето (см., например, [93]), пытавшегося объяснить расслоение населения по богатству. Насколько нам известно, ранее в таком контексте не приводились оценки скорости сходимости и плотности концентрации.

Для объяснения этого результата полезно рассмотреть схожий процесс, в котором каждой паре жителей приписан свой (независимый) «пуассоновский будильник» (звонки происходят в случайные моменты времени, соответствующие скачкам пуассоновского процесса; параметр интенсивность этого пуассоновского процесса называют интенсивностью / параметром будильника). Все будильники «приготовлены» одинаково: у всех у них одна и та же интенсивность $\lambda N^{-1}$. Далее следует погрузить задачу в модель стохастической химической кинетики с бинарными реакциями и воспользоваться результатом, приведенным в А.3. Наиболее технически сложными моментами является оценка с помощью неравенства Чигера mixing time $\sim \ln N$ [115] и получение поправки под корнем $\ln^2 N$.

К такой же предельной форме можно было бы прийти и по-другому (подход Булгакова–Маслова «разбрасывание червонцев в варьете»). В некотором городе живет $N \gg 1$ жителей (изначально банкротов). Каждый день одному из жителей (случайно выбранному в этот день) дается одна монетка. Тогда $\left\{c_s\left(\bar{s}N\right)\right\}_{s=0}^{\bar{s}N} \xrightarrow[N\to\infty]{} \left\{C\exp\left(-s/\bar{s}\right)\right\}_{s=0}^{\bar{s}N}$.

**Пример 2 (обезьянка и печатная машинка; закон Ципфа–Мандельброта).** На печатной машинке $n+1$ символ, один из символов пробел. Обезьянка на каждом шаге случайно (независимо и равновероятно) нажимает один из символов. Прожив долгую жизнь, обезьянка сгенерировала текст огромной длины. По этому тексту составили словарь. Этот словарь упорядочили по частоте встречаемости слова (слова – это всевозможные наборы букв без пробелов, которые хотя бы раз встречались в тексте обезьянки между какими-то двумя пробелами). Так, на первом месте в словаре поставили самое часто встречаемое слово, на второе поставили второе по частоте встречаемости и т. д. Номер слова в таком словаре называется рангом и обозначается буквой $r$. Предельная форма кривой, описывающей распределение частот встречаемости слов от рангов, имеет вид

$$\text{частота}\left(r\right) \simeq \frac{C}{\left(r+B\right)^{\alpha}},$$

где

$$\alpha = \frac{\log\left(n+1\right)}{\log\left(n\right)},\ B = \frac{n}{n-1},\ C = \frac{n^{\alpha-1}}{\left(n-1\right)^{\alpha}}.$$

Такой вывод закона Ципфа (с поправкой) был одним из двух, предложенных Бенуа Мандельбротом [152]. В работе [65] предложена довольно общая схема, приводящая к закону Ципфа, в которую можно погрузить и обезьянку с печатной машинкой. А именно, предположим, что динамика (порождения слов в большом тексте) такова, что вероятность того, что в тексте из $N \gg 1$ слов $x_1$ (первое по порядку слово в ранговом словаре) встречалось $N_1$ раз, $x_2$ (второе по порядку слово в ранговом словаре) встречалось $N_2$ раза и т. д., есть

$$\sim \frac{N!}{N_1!N_2!\ldots}\exp\left(-\eta\sum_{k\in\mathbb{N}} N_k E_k\right),$$

где $E_k$ – число букв в слове с рангом $k$. Часто считают, что $\eta = 0$, но зато динамика такова, что число слов и число букв становятся асимптотически (по размеру текста) связанными (закон больших чисел). Таким образом, к закону сохранения $\sum_{k\in\mathbb{N}} N_k = N$ добавляется приближенный закон сохранения $\sum_{k\in\mathbb{N}} N_k E_k \simeq \bar{E}N$ ($\bar{E}$ – среднее число букв в слове). Поиск предельной формы приводит к задаче (воспользовались формулой Стирлинга и методом множителей Лагранжа):

$$\sum_{k\in\mathbb{N}} \left\{ N_k \ln N_k + \lambda E_k N_k \right\} \to \min_{\substack{N_k \ge 0 \\ \sum_{k\in\mathbb{N}} N_k = N}},$$

где $\lambda$ либо равняется $\eta$, либо является множителем Лагранжа к ограничению $\sum_{k\in\mathbb{N}} N_k E_k \simeq \bar{E}N$. Стоит отметить, что к аналогичной задаче (с $E_k = k$) приводит поиск предельной формы в модели примера 1. Решение нашей задачи дает $N_k = \exp(-\mu - \lambda E_k)$, где $\mu$ – множитель Лагранжа к ограничению $\sum_{k\in\mathbb{N}} N_k = N$. Далее, считают, что $r(E)$ – число различных используемых слов с числом букв, не большим $E$, приближенно представимо в виде $r(E) \simeq a^E$. Тогда $N_k \sim k^{-\gamma}$, где $\gamma = \lambda / \ln a$.

**Пример 3 (теорема Гордона–Ньюэлла и PageRank [32]).** Имеется $N \gg 1$ пользователей, которые случайно (независимо) блуждают в непрерывном времени по ориентированному графу (на $m$ вершинах) с эргодической инфинитезимальной матрицей $\Lambda$. Назовем вектор $p$ (из единичного симплекса) PageRank, если $\Lambda p = 0$. Обозначим через $n_i(t)$ – число пользователей на $i$-й странице в момент времени $t \ge 0$. Несложно показать (теорема Гордона–Ньюэлла [147]), что $n(t)$ асимптотически имеет мультиномиальное распределение с вектором параметров PageRank $p$, т. е.

$$\lim_{t\to\infty} P\left(n(t) = n\right) = \frac{N!}{n_1! \cdot \ldots \cdot n_m!} p_1^{n_1} \cdot \ldots \cdot p_m^{n_m}.$$

Следовательно (неравенство Хефдинга в гильбертовом пространстве в форме [25]),

$$\lim_{t\to\infty} P\left( \left\| \frac{n(t)}{N} - p \right\|_2 \ge \frac{2\sqrt{2} + 4\sqrt{\ln(\sigma^{-1})}}{\sqrt{N}} \right) \le \sigma.$$

Этот же результат можно получить, рассмотрев соответствующую систему унарных химических реакций. Переход одного из пользователей из вершины $i$ в вершину $j$ означает превращение одной молекулы вещества $i$ в одну молекулу вещества $j$, $n_i(t)$ – число молекул $i$-го типа в момент времени $n_i(t)$. Каждое ребро графа соответствует определенной реакции (превращению). Интенсивность реакций определяется матрицей $\Lambda$ и числом молекул, вступающих в реакцию (закон действующих масс). Условие $\Lambda p = 0$ в точности соответствует условию унитарности в стохастической химической кинетике (см. раздел 1.3).

Отметим, что если воспользоваться теоремой Санова о больших уклонениях для мультиномиального распределения [60], то получим

$$\frac{N!}{n_1!\cdot\ldots\cdot n_m!} p_1^{n_1}\cdot\ldots\cdot p_m^{n_m} = \exp\left(-N\sum_{i=1}^{m}\nu_i \ln\left(\nu_i/p_i\right)+R\right),$$

где $\nu_i = n_i/N$, $|R| \le \frac{m}{2}(\ln N + 1)$. Однако последующее применение неравенства Пинскера не дает нам равномерной по $m$ оценки в 1-норме. Как и ожидалось, выписанная в 2-норме оценка (правая часть неравенства под вероятностью) и так полученная оценка в 1-норме будут отличаться по порядку приблизительно в $\sqrt{m}$ раз, что соответствует типичному (по Б. С. Кашину) соотношению между 1 и 2 нормами. Слово «типично» здесь отвечает, грубо говоря, за ситуацию, когда компоненты вектора одного порядка. При этом важно отметить, что для многих приложений, где возникают предельные конфигурации (кривые), описывающиеся вектором с огромным числом компонент, имеет место быстрый закон убывания этих компонент (см., примеры 1, 2 и [117, 138]), т. е. такая ситуация «нетипична» и можно ожидать лучшие оценки концентрации в 1-норме (см. пример 1).

# 1.3. Изучение динамики макросистемы с точки зрения стохастической химической кинетики

Предположим, что некоторая макросистема может находиться в различных состояниях, характеризуемых вектором $n$ с неотрицательными целочисленными компонентами. Будем считать, что в системе происходят случайные превращения (химические реакции).

Пусть $n \to n-\alpha+\beta$, $(\alpha,\beta)\in J$ – всевозможные типы реакций, где $\alpha$ и $\beta$ – векторы с неотрицательными целочисленными компонентами. Введем интенсивность реакции:

$$\lambda_{(\alpha,\beta)}(n)=\lambda_{(\alpha,\beta)}(n\to n-\alpha+\beta)=N^{1-\sum_i \alpha_i} K_\beta^\alpha \prod_{i:\,\alpha_i>0} n_i\cdot\ldots\cdot(n_i-\alpha_i+1),$$

где $K_\beta^\alpha \geq 0$ – константа реакции; при этом $\sum_{i=1}^m n_i(t)\equiv N\gg 1$. Другими словами, $\lambda_{(\alpha,\beta)}(n)$ – вероятность осуществления в единицу времени перехода $n\to n-\alpha+\beta$. Здесь не предполагается, что число состояний $m=\dim n$ и число реакций $|J|$ не зависят от числа агентов $N$. Тем не менее если ничего неизвестно про равновесную конфигурацию $c^*$ (типа быстрого убывания компонент этого вектора), то дополнительно предполагается, что $m\ll N$ – это нужно для обоснования возможности применения формулы Стирлинга при получении вариационного принципа (максимума энтропии) для описания равновесия макросистемы $c^*$ (в концентрационной форме). Однако часто априорно можно предполагать (апостериорно проверив), что компоненты вектора $c^*$ убывают быстро, тогда это условие можно отбросить. Так, например, обстоит дело в примере 1 из раздела 1.2. Описанный выше марковский процесс считается неразложимым. Имеет место

**Теорема. а)** $\langle \mu, n(t)\rangle \equiv \langle \mu, n(0)\rangle$ (inv) $\Leftrightarrow$ *вектор* $\mu$ *ортогонален каждому вектору семейства* $\{\alpha-\beta\}_{(\alpha,\beta)\in J}$. *Здесь* $\langle\,\cdot,\cdot\,\rangle$ *– обычно евклидово скалярное произведение.*

**б)** *Если существует* $\lim\limits_{N\to\infty} n(0)/N=c(0), K_\beta^\alpha := K_\beta^\alpha(n/N)$, $m$ *и* $|J|$ *не зависят от* $N$, *то для любого* $t>0$ *с вероятностью* 1 *существует* $\lim\limits_{N\to\infty} n(t)/N=c(t)$, *где* $c(t)$ *– неслучайная вектор-функция, удовлетворяющая СОДУ Гульдберга–Вааге*:

$$\frac{dc_i}{dt}=\sum_{(\alpha,\beta)\in J}(\beta_i-\alpha_i)K_\beta^\alpha(c)\prod_j c_j^{\alpha_j}\,. \qquad \text{(ГВ)}$$

*Гиперплоскость* (inv) (*с очевидной заменой* $n\Rightarrow c$) *инвариантна относительно этой динамики.*[32] *Более того, случайный процесс* $n(t)/N$ *слабо сходится при* $N\to\infty$ *к* $c(t)$ *на любом конечном отрезке времени.*

**в)** *Пусть выполняется условие унитарности* (*очевидно, что* $\xi$, *удовлетворяющий условию* (*U*), *– неподвижная точка в* (*ГВ*)):

[32]Можно сказать, что «жизнь» нелинейной динамической системы определяется линейными законами сохранения, унаследованными ею при скейлинге (каноническом). Этот тезис имеет, по-видимому, более широкое применение [16].

$$\exists\ \xi>0:\ \forall\ \beta\rightarrow\sum_{\alpha:(\alpha,\beta)\in J}K_\beta^\alpha\prod_j\left(\xi_j\right)^{\alpha_j}=\sum_{\alpha:(\alpha,\beta)\in J}K_\alpha^\beta\prod_j\left(\xi_j\right)^{\beta_j}. \qquad \text{(U)}$$

*Тогда неотрицательный ортант* $\mathbb{R}_+^m$ *расслаивается гиперплоскостями* (inv), *так что в каждой гиперплоскости* (inv) *уравнение* (*U*) (*положительно*) *разрешимо, притом единственным образом. Следовательно, существует, притом единственная, неподвижная точка* $c^*\in(\mathrm{inv})$ *у системы* (*ГВ*), *являющаяся при этом глобальным аттрактором. Система* (*ГВ*) *имеет функцию Ляпунова* $KL(c,\xi)=\sum_{i=1}^m c_i\ln\left(c_i/\xi_i\right)$.

*Стационарное распределение описанного марковского процесса имеет носителем множество* (inv) *и* (*с точностью до нормирующего множителя*) *имеет вид*

$$\frac{N!}{n_1!\cdot\ldots\cdot n_m!}\left(\xi_1\right)^{n_1}\cdot\ldots\cdot\left(\xi_m\right)^{n_m}\sim\exp\left(-N\cdot KL(c,\xi)\right),$$

*где* $\xi$ – *произвольное решение* (*U*), *неважно, какое именно* (*от этого, конечно, будет зависеть нормирующий множитель, но это ни на чем не сказывается*). *При этом условие унитарности* (*U*) *является обобщением условия детального равновесия*[33] (*баланса*)[34]:

$$\exists\ \xi>0:\ \forall\ (\alpha,\beta)\in J\rightarrow K_\beta^\alpha\prod_j\left(\xi_j\right)^{\alpha_j}=K_\alpha^\beta\prod_j\left(\xi_j\right)^{\beta_j},$$

*принимающего такой вид для мультиномиальной стационарной меры.*

*Существует такая зависимость* $a(m,c(0))$ (*во многих приложениях можно убрать второй аргумент в этой зависимости*), *что*

$$\exists a=a(m,c(0)):\forall\sigma>0,\,t\geq a\ln N\rightarrow P\left(\left\|\frac{n(t)}{N}-c^*\right\|_2\geq\frac{2\sqrt{2}+4\sqrt{\ln\left(\sigma^{-1}\right)}}{\sqrt{N}}\right)\leq\sigma,$$

---

[33]В терминах примера 3 условие унитарности просто означает, что в равновесии для любой вершины имеет место баланс числа пользователей, входящих в единицу времени в эту вершину, с числом пользователей, выходящих в единицу времени из этой вершины. В то время как условие детального равновесия означает, что в равновесии для любой пары вершин число пользователей, переходящих в единицу времени из одной вершину в другую, равно числу пользователей, переходящих в обратном направлении. Понятно, что второе условие является частным случаем первого.

[34]Много интересных примеров макросистем, в которых $K_\beta^\alpha:=K_\beta^\alpha\left(n/N\right)$, и имеет место детальный баланс, собрано в книге [13].

*где* (*принцип максимума энтропии, Больцман–Джейнс*):

$$c^* = \arg\max_{c\in(\mathrm{inv})}\left(-\sum_i c_i \ln\left(c_i/\xi_i\right)\right) = \arg\min_{c\in(\mathrm{inv})} KL\left(c,\xi\right),$$

*а* $\xi$ *– произвольное решение* (*U*), *причем* $c^*$ *определяется единственным образом, т. е. не зависит от выбора* $\xi$. *Геометрически* $c^*$ *– это KL-проекция произвольного* $\xi$, *удовлетворяющего* (*U*), *на гиперплоскость* (*inv*), *соответствующую начальным данным* $c(0)$. *Независимость этой проекции от выбора* $\xi$ *из* (*U*) *просто означает, что кривая* (*U*) *проходит KL-перпендикулярно через множество* (*inv*).

**г)** *Верно и обратное утверждение, то есть условие* (*U*) *не только достаточное для того, чтобы равновесие находилось из приведенной выше задачи энтропийно-линейного программирования, но и, с небольшой оговоркой* (*для почти всех* $c(0)$), *необходимое. Также верно и более общее утверждение, связывающее понимание энтропии в смысле Больцмана* (*функция Ляпунова прошкалированной кинетической динамики*) *и Санова* (*функционал действия в неравенствах больших уклонений для стационарной меры*): *если стационарная мера асимптотически представима в виде* $\sim \exp\left(-N\cdot V(c)\right)$, *то* $V(c)$ *– функция Ляпунова* (*ГВ*).

Результаты п. а) взяты из [8], п. б) из [54, 95], п. в) из [8, 54], п. г) из [22]. Выписанная оценка скорости сходимости и плотности концентрации в п. в) ранее нам не встречалась. Также отметим, что ранее результаты п. г) были получены при дополнительном предположении $K_\beta^\alpha\left(n/N\right) \equiv K_\beta^\alpha$.

К п. г) можно сделать следующее оригинальное пояснение. Обозначим через $h(c)$ вектор-функцию, стоящую в правой части СОДУ (ГВ). Тогда (см., например, [19, 95]) при $N \gg 1$ по теореме Т. Куртца $n(t)/N$ будет $O\left(\log N/\sqrt{N}\right)$ близко к $x_t = x(t)$ – решению стохастической системы дифференциальных уравнений (с начальным условием $x_0 = c(0)$):

$$dx_t = h(x_t)dt + \sqrt{\frac{g(x_t)}{N}}dW_t,$$

где функция $g(x_t) > 0$ рассчитывается по набору реакций и константам реакций (которые могут быть не постоянны и зависеть от концентраций), $W_t$ – стандартный винеровский процесс. Стационарная мера

$m(x)=\lim_{t\to\infty} p(t,x)$ этого однородного марковского процесса удовлетворяет уравнению

$$\frac{1}{2N}\nabla^2\left(g(x)m(x)\right)-\operatorname{div}\left(h(x)m(x)\right)=0,$$

поскольку плотность распределения $p(t,x)$ процесса $x_t$ подчиняется уравнению Колмогорова–Фоккера–Планка:

$$\frac{\partial p(t,x)}{\partial t}=-\operatorname{div}\left(h(x)\,p(t,x)\right)+\frac{1}{2N}\nabla^2\left(g(x)\,p(t,x)\right).$$

Здесь использовалось обозначение: $\nabla^2 f(x)=\sum_{i,j}\partial^2 f(x)/\partial x_i\partial x_j$ . Если известно, что

$$m(x)\simeq \mathrm{const}\cdot\exp\left(-N\cdot V(x)\right),$$

то из уравнения на $m(x)$ имеем

$$N\langle h,\nabla V\rangle-\operatorname{div} h-\frac{1}{2}\langle\nabla g,1\rangle\langle\nabla V,1\rangle+\frac{1}{2N}V\nabla^2 g-\frac{1}{2}g\nabla^2 V+\frac{N}{2}g\langle\nabla V,1\rangle^2\simeq 0,$$

следовательно,

$$\langle h,\nabla V\rangle\simeq-\frac{1}{2}g\langle\nabla V,1\rangle^2+\mathrm{O}\left(\frac{1}{N}\right)\underset{N\to\infty}{=}-\frac{1}{2}g\langle\nabla V,1\rangle^2\le 0.$$

Эта выкладка поясняет, почему функция $V(c)$ может быть функцией Ляпунова системы (ГВ) $dc/dt=h(c)$. Более того, модель стохастической химической кинетики здесь может быть заменена и более общими шкалирующимися марковскими моделями [95].

В заключение отметим, что много различных примеров макросистем собрано в гл. 6 книги [12].

# Приложение 2. Эволюционный вывод простейшей модели бимодального расщепления спроса на городские передвижения

## 2.1. Введение

Расщепление спроса на ежедневные трудовые передвижения, удовлетворяемого соответственно поездкам на личных автомобилях и на общественном транспорте (*Modal Split*), базируется на стандартной (и вполне правдоподобной) гипотезе, согласно которой горожанин выбирает способ передвижения по критерию минимума обобщенной цены поездки (*Generalized Cost*), складывающейся из двух компонент:

– денежных затрат на совершение поездки (*Monetary Cost*),

– «немонетарных затрат» (*Non-Monetary Cost*), исчисляемых как произведение времени поездки на цену единицы времени горожанина (*Value of Time*), вообще говоря, сугубо индивидуальную для каждого.

Денежные затраты на совершение поездки трактуются в практике транспортного планирования как «траты из кармана» (*Out of Pocket Price*). Для общественного транспорта это плата за проезд на тех или иных маршрутах и видах транспорта. Для автомобильной поездки – траты на моторное топливо и парковочные платежи; условно-постоянные компоненты затрат (амортизация транспортного средства, налоги и т. п.) в учет обычно не берутся, так как они слабо влияют на ежедневные решения горожанина по выбору способа совершения поездки.

Способы оценки времени поездки на одно и то же расстояние различаются для общественного транспорта и личного автомобиля принципиальным образом. Для общественного транспорта это время можно считать фиксированным: применительно к внеуличным видами транспорта это очевидный факт, для наземного транспорта, работающего в общем потоке транспортных средств, неизбежные задержки в движении, обусловленные плотностью трафика, учтены в маршрутных расписаниях. Для автомобильных поездок ключевым обстоятельством является зависимость времени поездки от плотности трафика, то есть в конечном итоге, от количества претендентов на пользование улично-дорожной сетью.

Здесь, следуя М. Я. Блинкину, предложена эволюционная модель поведения жителей простейшего города (с достаточно большим числом жителей), состоящего всего из двух условных районов: селитебного (спального) и делового (рабочего). У каждого горожанина-автомобилиста для совершения ежедневной трудовой поездки есть возможность воспользоваться как личным автомобилем, так и общей для всех линией (маршрутом) общественного транспорта. В отличие от многих других работ по данной тематике (см., например, [140]), в данном приложении мы сделаем акцент не только на описание равновесного расщепления горожан по выбору способа передвижения (на личном транспорте или общественном), но и на том, как в реальном времени может происходить процесс «нащупывания» этого равновесия. Таким образом, исследуется вопрос устойчивости равновесия.

## 2.2. Простейшая эволюционная модель расщепления

Итак, рассматривается город, состоящий из двух районов: спального и рабочего. Каждый день жители города (все они живут в спальном районе, а работают, естественно, в рабочем) ездят на работу. Каждый из них имеет личный автомобиль. Кроме того, в городе есть развитая сеть общественного транспорта. Таким образом, каждый житель имеет две альтернативные возможности для ежедневных трудовых поездок: личный автомобиль и общественный транспорт.

Ежедневные потери пользователей личного транспорта, оценивающих единицу (минуту) своего времени в $p \geq 1$ рублей, могут быть рассчитаны следующим образом:

$$A_p(x) = a + pT(x),$$

где $a > 0$ характеризует постоянные затраты (цена топлива и т. п.), $x \in [0,1]$ – доля жителей города, использующих личный автомобиль, $T(x)$ – функция, характеризующая то, как пользователи транспортной сети оценивают свои временные затраты. Обычно (см., например, [140]) $T(x)$ выбирают вида BPR-функций[35], т. е. $T(x) = T_0 + \gamma x^4$.

Ежедневные потери пользователей общественного транспорта, оценивающих единицу (минуту) своего времени в $p \geq 1$ рублей, могут быть рассчитаны следующим образом:

$$B_p(x) \equiv b_1 + pb_2,$$

[35]Функции данного типа были введены на основе обширных эмпирических обследований американским Bureau of Public Roads (BPR) и много лет применяются в работах по транспортному планированию и теории транспортного потока [88].

где $b_1 > 0$ характеризует постоянные затраты (цена билета и т. п.), а $b_2 > 0$ можно понимать как время, потерянное в пути. В отличие от личного транспорта, для общественного транспорта считается, что $b_2$ не зависит от $x$ (метро, электропоезда, выделенные полосы и т. п.).

Будем считать, что жителей в городе много. Они расслоены по тому, во сколько рублей каждый из них оценивает единицу своего времени (потерянного в пути). Обычно это колеблется от 1 руб./мин до 30 руб./мин (это максимальное значение будем обозначать далее через $p_{\max}$). Введем зависимость $x(p)$ – доля жителей города, оценивающих одну минуту своего времени не меньше, чем в $p$ рублей. Эта зависимость естественным образом восстанавливается из закона распределения населения по доходу Ципфа–Парето (см. раздел 2.4). Обычно эту зависимость считают степенной $x(p) = p^{-\eta}$, где $\eta$ выбирают из диапазона 1–2.

Наложим теперь физически правдоподобные ограничения. Во-первых, будем считать, что денежная цена автомобильной поездки (постоянные затраты) дороже, чем в случае общественного транспорта

$$b_1 < a. \qquad (2.1)$$

Это предположение не так очевидно, как это представляется на первый взгляд. К примеру, в Москве до введения платной парковки в 2013 году *Out of Pocket Price* автомобильной поездки была заметно ниже цены пересадочной поездки «метро + трамвай». Понятно, что в таких условиях общественный транспорт предпочитало в основном «безлошадное» население. Ограничение (2.1) автоматически выполняется немедленно после введения платной парковки даже по самому щадящему тарифу.

Во-вторых, будем считать, что «на автомобиле всегда быстрее»:

$$T(1) < k * b_2. \qquad (2.2)$$

Для справедливости этого предположения принципиально наличие поправочного коэффициента $k$, связывающего времена, потерянные на личном и общественном транспорте (важно считать эти времена не равноценными). Понятно, что в условиях затора, практически неизбежного при $x$, близких к единице, время поездки на метро заведомо будет меньше, чем на автомобиле $T(1) > b_2$. Вопрос, однако, заключается в том, что время, проведенное соответственно в вагоне общественного транспорта и в собственном автомобиле, трудно считать равноценным. Другими словами, в данной работе под $b_2$ понимается время, потерянное в пути на общественном транспорте, приведенное к масштабу времени, потерянному на личном транспорте. Такое приведение увеличивает «физическое» время ввиду различной комфортности перемещений.

Добавим для аккуратности, что для горожан с наиболее высокой ценой времени обобщенная цена автомобильной поездки всегда ниже, чем на общественном транспорте:

$$a + VT^{\max} * T(1) < b_1 + VT^{\max} * k * b_2. \tag{2.3}$$

Предположим, наконец, что для горожан с наиболее низкой ценой времени дело обстоит прямо противоположным образом: даже при самой низкой загрузке улично-дорожной сети обобщенная цена поездки на общественном транспорте ниже, чем у автомобильной поездки:

$$a + VT^{\min} * T(0) > b_1 + VT^{\min} * k * b_2. \tag{2.4}$$

Определим зависимость $p(x)$, как корень уравнения $A_p(x)=B_p$, т. е.

$$p(x) = \frac{a - b_1}{b_2 - T(x)}.$$

При условиях (2.1) – (2.4) выписанная формула корректно определяет монотонную гладкую зависимость $p(x) \in \left[VT^{\min}, VT^{\max}\right] = \left[1, p_{\max}\right]$ при $x \in [0,1]$.

Представим себе такую динамику (повторяющуюся изо дня в день). Каждый житель в $(k+1)$-й день смотрит на то, какая доля жителей $x^k$ использовала личный автомобиль в $k$-й день. Считаем, что такая информация (статистика) по вчерашнему дню общедоступна (например, благодаря каким-нибудь интернет сервисам, скажем, Яндекс.Пробки). Исходя из этой информации каждый житель, оценивающий минуту своего времени в $p$ рублей, оценивает (экстраполируя ситуацию вчерашнего дня на день сегодняшний, за неимением точной информации о $x^{k+1}$) свои затраты от двух возможных альтернатив: $A_p\left(x^k\right)$ – личный автомобиль и $B_p$ – общественный транспорт. Мы считаем всех жителей рациональными, поэтому из двух альтернатив, каждый житель выбирает ту, которая приносит ему наименьшие затраты. Таким образом, происходит формирование $x^{k+1}$.

Из описанного выше ясно, что жители города в $(k+1)$-й день, оценивающие единицу своего времени в $p\left(x^k\right) < p \le p_{\max}$ рублей, предпочтут в этот день личный автомобиль, а жители, оценивающие единицу своего времени в $1 \le p < p\left(x^k\right)$ рублей, предпочтут в этот день обще-ственный транспорт. Таким образом, в $(k+1)$-й день доля $x^{k+1} =$ $= x(p(x^k))$ жителей города использует (выберет) личный автомобиль.

Для того чтобы сформулировать основной результат, сделаем одно упрощающее предположение, которое позволит представить этот результат в более наглядной форме. Будем считать, что в зависимости $x(p) = p^{-\eta}$ параметр $\eta = 1$. Тогда если

$$4\gamma < a - b_1, \tag{2.5}$$

то $x(p(\cdot))$ – сжимающее преобразование отрезка $[0,1]$ в себя. Легко понять, что условия (2.1) – (2.5) совместны.

**Теорема 2.1.** *При* $\eta = 1$ *в условиях* (2.1) – (2.5) *описанная выше динамика сходится* (*вне зависимости от* $x^0$) *со скоростью геометрической прогрессии к неподвижной точке*: $x^k \xrightarrow[k\to\infty]{} x^*$, *которая определяется как единственная точка пересечения графиков* $x(p)$ *и* $p(x)$ *на плоскости* $(p,x)$.

**Доказательство.** Сформулированное в теореме 2.1 утверждение сразу следует из принципа неподвижной точки для сжимающих операторов [47]. Однако намного полезнее представляется продемонстрировать доказательство (рис. 1), проясняющим, как происходит «нащупывание» равновесия.

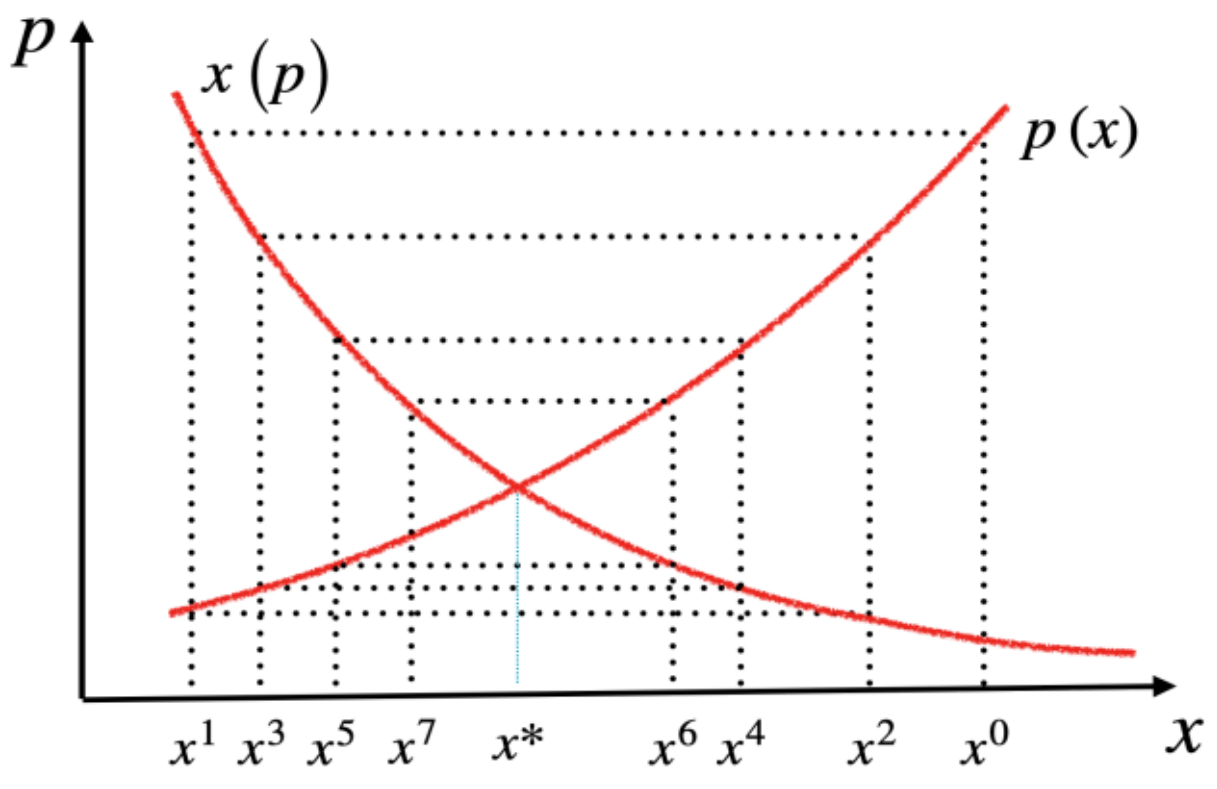

Рис. 1

Этот рисунок можно проинтерпретировать следующим образом. В начальный момент (в условный первый день) $x^0$ доля жителей города воспользовались личным транспортом. На следующий день всем известно $x^0$, исходя из этого числа, каждый оценивает, каким образом ему сегодня добираться до работы (и обратно). В результате личным транспортом воспользуется $x^1 = x\left(p\left(x^0\right)\right)$ доля жителей города. Аналогично на следующий день личным транспортом воспользуется $x^2 = x\left(p\left(x^1\right)\right)$, на следующий $x^3 = x\left(p\left(x^2\right)\right)$ и т. д. Из рис. 1 видно, что такой итерационный процесс будет сходящимся к корню уравнения $x = x\left(p\left(x\right)\right)$, что можно понимать как точку пересечения графиков функции $x\left(p\right)$ и $p\left(x\right)$. Легко проверить, что если на рис. 1 выбрать $x^0$ левее $x^*$, то сходимость также будет иметь место. □

## 2.3. Численные эксперименты

Моделировался город численностью 1000 человек. Были выбраны следующие параметры $\gamma = 2$ минуты, $T_0 = 70$ минут, $b_2 = 75$ минут, $b_1 = 50$ рублей, $a = 60$ рублей. При таких параметрах процесс «нащупывания» равновесия сходился в среднем (случайно выбиралась точка $x^0 \in \left[0,1\right]$) за 3–4 дня (итерации). Были также проведены численные эксперименты и при других значениях параметров. При самых неблагоприятных значениях, используемых в численных экспериментах, сходимость была в среднем за семь дней. Таким образом, получено модельное подтверждение известного из опыта факта, что выход (нащупывание) равновесия осуществляется крайне быстро (недели всегда хватает).

## 2.4. Закон Ципфа–Парето и процесс Юла

Далее мы постараемся пояснить возникновение закона Ципфа. Для этого мы сделаем следующее упрощающее предположение: каждый человек оценивает единицу своего времени в сумму, которую он зарабатывает в единицу времени. Упрощая еще больше, предположим, что богатство каждого жителя прямо пропорционально тому, сколько он зарабатывает в единицу времени (считаем, что все работают одинаковое время – у всех нормированный рабочий день). Таким образом, нужно показать степенной характер распределения населения по доходу. Для этого рассмотрим следующую игрушечную модель.

В некотором городе живет неограниченно много жителей (изначально банкротов), которые могут участвовать в «освоении» монеток.

База индукции: сначала выбирается один житель, он получает одну монетку. Шаг индукции: в $(k+1)$-й день выбирается очередной новый житель (отличный от $k$ уже выбранных), он получает возможность участвовать в разыгрывании монеток, монетка с вероятностью $\alpha<1$ равновероятно отдается одному из этих $k+1$ жителей, а с вероятностью $1-\alpha$ эта монетка отдается одному из $k$ старых жителей с вероятностью, пропорциональной тому, сколько у него уже есть монеток, т. е. по принципу «деньги к деньгам» (в моделях роста Интернета этот принцип называют *preferential attachment*).

Обычно эту стохастическую динамику изучают в приближении среднего поля. В данном контексте это означает, что $n_s(t)\simeq E\left[n_s(t)\right]$, где $n_s(t)$ – количество жителей, у которых ровно $s$ монеток на $t$-й день. Далее выписывают на $n_s(t)$ систему зацепляющихся обыкновенных дифференциальных уравнений. Автомодельное притягивающее решение этой системы ищут в виде $n_s(t)\sim x_s^* t$. После разрешения соответствующих уравнений получают степенной закон для зависимости $x_s^* \sim s^{-(3+\alpha/(1-\alpha))}$. Такой подход применительно к моделям роста Интернета (и изучения степенного закона для распределения степеней вершин) довольно часто сейчас встречается (в том числе и в учебной литературе). В частности, этот подход описан в обзоре [117]. По-сути, в этом приложении нами был описан процесс, возникающий в работах 20-х годов XX века по популяционной генетике, получивший названия процесса Юла (см., например, обзор [138]).

Описанная модель также восходит к работе конца XIX века В. Парето, в которой была предпринята попытка объяснить социальное неравенство и к работе Г. Ципфа конца 40-х годов XX века, в которой была отмечена важность степенных законов в природе. Эти законы для большей популярности иногда преподносят, как принцип Парето или принцип 80/20 (80% результатов проистекают всего лишь из 20% причин) – такие пропорции отвечают $x_s^* \sim s^{-2.1}$ [138]. Приведем примеры (не совсем, правда, точные): 80% научных результатов получили 20% ученых, 80% пива выпили 20% людей и т. п. Сейчас много исследований во всем мире посвящено изучению возникновения в самых разных приложениях степенных законов (распределение городов по населению, коммерческих компаний по капитализации, автомобильных пробок по длине). Особенно бурный рост возник в связи с изучением роста больших сетей (экономических, социальных, интернета), см., например, книги и работы M. O. Jackson'a, в частности [107].

# Приложение 3. Универсальный метод подобных треугольников для задач стохастической композитной оптимизации

## 3.1. Введение

Быстрые (ускоренные) градиентные методы были предложены в 1983 г. в кандидатской диссертации Ю. Е. Нестерова. Известно, см., например, [55, 127], что данный класс методов является наилучшим для задач гладкой выпуклой оптимизации. В данном приложении описывается один современный представитель такого метода (предложенный в апреле 2016 г. Ю. Е. Нестеровым [37, 127]): метод подобных треугольников (МПТ), изложенный в разделе 3.2. Особенностью МПТ является необходимость выполнения всего одного «проектирования» на каждой итерации. Как следствие, предложенный метод оказался заметно проще многих своих аналогов. В разделе 3.3 МПТ распространяется на задачи сильной выпуклой композитной оптимизации. В разделе 3.4 предлагается универсальный вариант МПТ для задач выпуклой и сильно выпуклой композитной оптимизации. Далее в изложении мы следуем работам [37, 127]. Ссылки на литературу можно почерпнуть из этих двух источников.

## 3.2. Метод подобных треугольников для задач композитной оптимизации

Рассматривается задача выпуклой композитной оптимизации

$$F(x) = f(x) + h(x) \to \min_{x \in Q}. \tag{3.1}$$

Положим $R^2 = V\left(x_*, y^0\right)$, где прокс-расстояние (расстояние Брэгмана) определяется формулой

$$V(x,z) = d(x) - d(z) - \left\langle \nabla d(z), x - z \right\rangle;$$

прокс-функция $d(x) \geq 0$ ($d\left(y^0\right) = 0$) считается сильно выпуклой относительно выбранной нормы $\| \ \|$, с константой сильной выпуклости $\geq 1$; $x_*$ – решение задачи (3.1) (если решение не единственно, то выбирается то, которое доставляет минимум $V\left(x_*, y^0\right)$).

**Предположение 3.1.** Для любых $x, y \in Q$ имеет место неравенство

$$\left\|\nabla f(y) - \nabla f(x)\right\|_* \le L\left\|y - x\right\|.$$

Опишем вариант быстрого градиентного метода для задачи (3.1) с одной «проекцией», который мы далее будем называть *метод подобных треугольников* (МПТ).

Положим

$$\varphi_0(x) = V\left(x, y^0\right) + \alpha_0\left[f\left(y^0\right) + \left\langle \nabla f\left(y^0\right), x - y^0\right\rangle + h(x)\right],$$

$$\varphi_{k+1}(x) = \varphi_k(x) + \alpha_{k+1}\left[f\left(y^{k+1}\right) + \left\langle \nabla f\left(y^{k+1}\right), x - y^{k+1}\right\rangle + h(x)\right], \quad (3.2)$$

$$A_k = \sum_{i=0}^{k} \alpha_i,\ \alpha_0 = L^{-1},\ A_k = \alpha_k^2 L,\ k = 0,1,2,\ldots,\ x^0 = u^0 = \arg\min_{x\in Q} \varphi_0(x). \quad (3.3)$$

**<u>Метод подобных треугольников</u>**

$$y^{k+1} = \frac{\alpha_{k+1}u^k + A_k x^k}{A_{k+1}},$$

$$u^{k+1} = \arg\min_{x\in Q} \varphi_{k+1}(x), \quad (3.4)$$

$$x^{k+1} = \frac{\alpha_{k+1}u^{k+1} + A_k x^k}{A_{k+1}}.$$

**Лемма 3.1 (см. [91]).** *Последовательность* $\{\alpha_k\}$, *определяемую формулой* (3.3), *можно задавать рекуррентно*:

$$\alpha_{k+1} = \frac{1}{2L} + \sqrt{\frac{1}{4L^2} + \alpha_k^2}.$$

*При этом*

$$A_k \ge \frac{(k+1)^2}{4L}.$$

**Лемма 3.2.** *Пусть справедливо предположение* 3.1. *Тогда для любого* $k = 0,1,2,\ldots$ *имеет место неравенство*

$$A_k F\left(x^k\right) \le \varphi_k^* = \min_{x\in Q} \varphi_k(x) = \varphi_k\left(u^k\right). \quad (3.5)$$

**Доказательство.** Проведем по индукции. При $k = 0$ формула (3.5) следует из того, что для любого $x \in Q$

$$F(x) = f(x) + h(x) \le f\left(y^0\right) + \left\langle \nabla f\left(y^0\right), x - y^0\right\rangle + \frac{L}{2}\left\|x - y^0\right\|^2 + h(x) \le$$

$$\le LV\left(x, y^0\right) + f\left(y^0\right) + \left\langle \nabla f\left(y^0\right), x - y^0\right\rangle + h(x) = \frac{1}{A_0}\varphi_0(x).$$

Последнее неравенство следует из того, что для любых $x, z \in Q$

$$V(x,z) \geq \frac{1}{2}\|x-z\|^2 .$$

Это неравенство, в свою очередь, следует из $\geq 1$-сильной выпуклости $d(x)$ в норме $\| \ \|$.

Итак, пусть формула (3.5) установлена для $k$, покажем, что тогда она будет справедлива и для $k+1$. По определению (3.2)

$$\varphi_{k+1}^* = \min_{x \in Q} \varphi_{k+1}(x) = \varphi_{k+1}\left(u^{k+1}\right) =$$

$$= \varphi_k\left(u^{k+1}\right) + \alpha_{k+1}\left[ f\left(y^{k+1}\right) + \left\langle \nabla f\left(y^{k+1}\right), u^{k+1} - y^{k+1} \right\rangle + h\left(u^{k+1}\right) \right]. \qquad (3.6)$$

Поскольку по предположению индукции $A_k F\left(x^k\right) \leq \varphi_k\left(u^k\right)$ и $\varphi_{k+1}(x)$ – сильно выпуклая в $\| \ \|$-норме функция с константой $\geq 1$ (это следует из аналогичного свойства функции $V\left(x, y^0\right)$, что, в свою очередь, следует из аналогичного свойства функции $d(x)$), то

$$\varphi_k\left(u^{k+1}\right) \geq \varphi_k\left(u^k\right) + \frac{1}{2}\left\|u^{k+1} - u^k\right\|^2 \geq A_k \cdot \left( f\left(x^k\right) + h\left(x^k\right)\right) + \frac{1}{2}\left\|u^{k+1} - u^k\right\|^2 .$$

Из выпуклости $f(x)$ отсюда имеем

$$\varphi_k\left(u^{k+1}\right) \geq A_k f\left(y^{k+1}\right) + \left\langle \nabla f\left(y^{k+1}\right), A_k \cdot \left(x^k - y^{k+1}\right)\right\rangle + A_k h\left(x^k\right) + \frac{1}{2}\left\|u^{k+1} - u^k\right\|^2. \qquad (3.7)$$

Подставляя (3.7) в (3.6), получим

$$\varphi_{k+1}^* \geq A_{k+1} \cdot \underbrace{\left( \frac{A_k}{A_{k+1}} h\left(x^k\right) + \frac{\alpha_{k+1}}{A_{k+1}} h\left(u^{k+1}\right) \right)}_{\geq h\left(x^{k+1}\right)} + A_{k+1} f\left(y^k\right) +$$

$$+ \left\langle \nabla f\left(y^{k+1}\right), \underbrace{\alpha_{k+1} \cdot \left(u^{k+1} - y^{k+1}\right) + A_k \cdot \left(x^k - y^{k+1}\right)}_{= A_{k+1} \cdot \left(x^{k+1} - y^{k+1}\right)} \right\rangle + \underbrace{\frac{1}{2}\left\|u^{k+1} - u^k\right\|^2}_{= \frac{A_{k+1}^2}{2\alpha_{k+1}^2}\left\|x^{k+1} - y^{k+1}\right\|^2} . \qquad (3.8)$$

Исходя из выпуклости функции $h(x)$ и описания МПТ (3.4), формулу (3.8) можно переписать следующим образом:

$$\varphi_{k+1}^* \geq A_{k+1}\left[ f\left(y^{k+1}\right) + \left\langle \nabla f\left(y^{k+1}\right), x^{k+1} - y^{k+1} \right\rangle + \frac{A_{k+1}}{2\alpha_{k+1}^2}\left\|x^{k+1} - y^{k+1}\right\|^2 + h\left(x^{k+1}\right) \right]. \qquad (3.9)$$

Из предположения 3.1 следует, что если $A_{k+1}/\alpha_{k+1}^2 \ge L$, то

$$f\left(y^{k+1}\right)+\left\langle \nabla f\left(y^{k+1}\right), x^{k+1}-y^{k+1}\right\rangle+\frac{A_{k+1}}{2\alpha_{k+1}^2}\left\|x^{k+1}-y^{k+1}\right\|^2 \ge f\left(x^{k+1}\right). \quad (3.10)$$

Согласно (3.3) $A_{k+1}/\alpha_{k+1}^2 = L$, поэтому формула (3.10) имеет место. С помощью формулы (3.10) формулу (3.9) можно переписать следующим образом:

$$\varphi_{k+1}^* \ge A_{k+1}\left[f\left(x^{k+1}\right)+h\left(x^{k+1}\right)\right]=A_{k+1}F\left(x^{k+1}\right).$$

Таким образом, шаг индукции установлен. Следовательно, лемма 3.2 доказана. ■

Из лемм 3.1, 3.2 получаем следующий результат, означающий, что МПТ сходится так же, как обычный быстрый градиентный метод, т. е. МПТ сходится оптимальным образом для рассматриваемого класса задач.

**Теорема 3.1.** *Пусть справедливо предположение* 3.1. *Тогда МПТ* (3.2) – (3.4) *для задачи* (3.1) *сходится согласно оценке*

$$F\left(x^N\right)-\min_{x\in Q} F(x) \le \frac{4LR^2}{(N+1)^2}. \quad (3.11)$$

**Доказательство.** Из леммы 3.2 следует, что (в третьем неравенстве используется выпуклость функции $f(x)$):

$$A_N F\left(x^N\right) \le \min_{x\in Q}\left\{V\left(x,y^0\right)+\sum_{k=0}^{N}\alpha_k\left[f\left(y^k\right)+\left\langle\nabla f\left(y^k\right),x-y^k\right\rangle+h(x)\right]\right\} \le$$

$$\le V\left(x_*,y^0\right)+\sum_{k=0}^{N}\alpha_k\underbrace{\left[f\left(y^k\right)+\left\langle\nabla f\left(y^k\right),x_*-y^k\right\rangle+h\left(x_*\right)\right]}_{\le f(x_*)+h(x_*)} \le$$

$$\le V\left(x_*,y^0\right)+\sum_{k=0}^{N}\alpha_k F\left(x_*\right)=R^2+A_N F\left(x_*\right). \quad (3.12)$$

Заметим, что из второго неравенства следует, что если решение задачи (3.1) $x_*$ не единственно, то можно выбирать то, которое доставляет минимум $V\left(x_*,y^0\right)$. Именно таким образом возникает $R^2$ в оценке (3.12). Для доказательства теоремы осталось подставить нижнюю оценку на $A_N$ из леммы 1 в формулу (3.12). ■

**Замечание 3.1.** В действительности в формуле (3.12) содержится более сильный результат, чем в формуле (3.11). А именно, формула (12) еще означает, что МПТ – прямодвойственный метод, см. подраздел 2.2.6 гл. 2. Мы не будем здесь подробно на этом останавливаться, отметим лишь, что это свойство позволяет получать эффективные критерии оста-нова для МПТ. Критерий останова позволяет не делать предписанного формулой (3.11) числа итераций и останавливаться раньше (по достижению желаемой точности). Это замечание можно распространить и на все последующее изложение.

МПТ получил такое название, поскольку на базе трех точек $u^k$, $u^{k+1}$, $x^k$ можно построить треугольник, а точки $y^{k+1}$ и $x^{k+1}$ лежат на сторонах этого треугольника ($y^{k+1}$ лежит на стороне $u^k x^k$, а $x^{k+1}$ – на стороне $u^{k+1}x^k$), причем прямая, проходящая через точки $y^{k+1}$ и $x^{k+1}$, параллельна прямой, проходящей через точки $u^k$ и $u^{k+1}$.

Из описанной геометрической интерпретации получаются следующие результаты, распространимые и на все последующее изложение.

**Следствие 3.1.** *Для любого* $k = 0,1,2,...$ *имеют место неравенства*

$$\left\|u^k - x_*\right\|^2 \le 2V\left(x_*, y^0\right),$$

$$\max\left\{\left\|x^k - x_*\right\|^2, \left\|y^k - x_*\right\|^2\right\} \le 4V\left(x_*, y^0\right) + 2\left\|x^0 - y^0\right\|^2 = \tilde{R}^2 .$$

**Доказательство.** Второе неравенство следует из первого, описанной выше геометрической интерпретации и неравенства $\left\|a+b\right\|^2 \le 2\left\|a\right\|^2 + 2\left\|b\right\|^2$.

Докажем первое неравенство. Для этого воспользуемся леммой 1 и тем, что $\varphi_k\left(x\right)$ – сильно выпуклая в $\|\ \|$-норме функция с константой $\ge 1$. Для любого $x \in Q$

$$A_k F\left(x^k\right) + \frac{1}{2}\left\|x - u^k\right\|^2 \le \varphi_k^* + \frac{1}{2}\left\|x - u^k\right\|^2 \le \varphi_k\left(x\right) \le$$

$$\le \underbrace{\sum_{i=0}^{k} \alpha_i \left[ f\left(y^i\right) + \left\langle \nabla f\left(y^i\right), x - y^i \right\rangle + h\left(x\right)\right]}_{\le A_k F(x)} + V\left(x, y^0\right) \le A_k F\left(x\right) + V\left(x, y^0\right).$$

Выбирая $x = x_*$ и используя то, что $F\left(x^k\right) \ge F\left(x_*\right)$, получим первое неравенство следствия 3.1. ■

Следствие 3.1 играет важную роль в случае неограниченных множеств $Q$, на которых нельзя равномерно ограничить константу $L$, поскольку гарантирует, что как бы «плохо» себя не вела функция $F\left(x\right)$ вне шара конечного радиуса $\tilde{R}$ (зависящего от качества начального приближения) с центром в решении $x_*$, это никак не скажется на скорости сходимости метода, поскольку итерационный процесс никогда не выйдет за пределы этого шара.

**Следствие 3.2.** *Пусть* $h(x) \equiv 0$ (*т. е.* $F(x) = f(x)$) *и* $\nabla f(x_*) = 0$. *Тогда*

$$\max\left\{F\left(x^N\right), F\left(y^N\right), F\left(z^N\right)\right\} - \min_{x \in Q} F(x) \le \frac{L\tilde{R}^2}{N^2}.$$

## 3.3. Метод подобных треугольников для сильно выпуклых задач композитной оптимизации

В данном пункте будем считать, что $f(x)$ в задаче (3.1) обладает следующим свойством.

**Предположение 3.2.** $f(x)$ – $\mu$-сильно выпуклая функция в норме $\| \ \|$, т. е. для любых $x, y \in Q$ имеет место неравенство

$$f(y) + \langle \nabla f(y), x - y \rangle + \frac{\mu}{2}\|x - y\|^2 \le f(x). \qquad (3.13)$$

Введем (в евклидовом случае $\tilde{\omega}_n = 1$)

$$\tilde{\omega}_n = \sup_{x, y \in Q} \frac{2V(x, y)}{\|y - x\|^2}.$$

**Замечание 3.2.** Из работы [110] следует, что в сильно выпуклом случае (когда сильно выпукла гладкая часть функционала $f(x)$, как в предположении 2) естественно выбирать именно евклидову норму и прокс-структуру, т. е. в большинстве случаев можно считать $\tilde{\omega}_n = 1$. Поясним это примером из [110]. Число обусловленности (отношение константы Липшица градиента $L_1(f)$ к константе сильной выпуклости $\mu_1(f)$), например, для квадратичных функций

$$f(x) = \frac{1}{2}x^T Ax - b^T x, \ x \in \mathbb{R}^n,$$

посчитанное, скажем, в 1-норме, не может быть меньше $n$. В то время как число обусловленности, посчитанное в евклидовой норме, может при этом равняться 1. Действительно, пусть $\xi = (\xi_1, ...., \xi_n)$, где $\xi_k$ – независимые одинаково распределенные случайные величины $P(\xi_k = 1/n) = P(\xi_k = -1/n) = 1/2$. Тогда, учитывая, что $\|\xi\|_1 \equiv 1$,

$$\mu_1(f) \le E_\xi\left[\xi^T A\xi\right] = \frac{1}{n^2}\operatorname{tr}(A) \le \frac{1}{n}\max_{i,j=1,\dots,n}\left|A_{ij}\right| = \frac{1}{n}L_1(f).$$

Кроме того, из замечания 3.2 следует также и то, что если $\| \; \| = \| \; \|_1$, то $\tilde{\omega}_n \geq n$. Это обстоятельство также говорит в пользу выбора евклидовой прокс-структуры. Тем не менее, для общности далее мы будем допускать, что используется прокс-структура отличная от евклидовой.

Перепишем формулы (3.2), (3.3) следующим образом ( $\tilde{\mu} = \mu / \tilde{\omega}_n$ ):

$$\varphi_0(x) = V(x, y^0) + \alpha_0 \left[ f(y^0) + \left\langle \nabla f(y^0), x - y^0 \right\rangle + \tilde{\mu} V(x, y^0) + h(x) \right],$$

$$\varphi_{k+1}(x) = \varphi_k(x) + \alpha_{k+1} \left[ f(y^{k+1}) + \left\langle \nabla f(y^{k+1}), x - y^{k+1} \right\rangle + \tilde{\mu} V(x, y^k) + h(x) \right], \qquad (3.14)$$

$$A_k = \sum_{i=0}^{k} \alpha_i, \; \alpha_0 = L^{-1}, \; A_{k+1} \cdot (1 + A_k \tilde{\mu}) = \alpha_{k+1}^2 L, \; k = 0, 1, 2, \ldots \qquad (3.15)$$

Сам метод по-прежнему будет иметь вид (3.4) с $x^0 = u^0 = \arg\min_{x \in Q} \varphi_0(x)$.

**Лемма 3 (см. [92]).** *Последовательность* $\{\alpha_k\}$, *определяемую формулой* (3.15), *можно задавать рекуррентно*:

$$\alpha_{k+1} = \frac{1 + A_k \tilde{\mu}}{2L} + \sqrt{\frac{1 + A_k \tilde{\mu}}{4L^2} + \frac{A_k \cdot (1 + A_k \tilde{\mu})}{L}}, \; A_{k+1} = A_k + \alpha_{k+1}. \qquad (3.16)$$

*При этом*

$$A_k \geq \frac{1}{L} \left( 1 + \frac{1}{2} \sqrt{\frac{\tilde{\mu}}{L}} \right)^{2k} \geq \frac{1}{L} \exp\left( \frac{k}{2} \sqrt{\frac{\tilde{\mu}}{L}} \right).$$

**Лемма 3.4.** *Пусть справедливы предположения* 3.1, 3.2. *Тогда для любого* $k = 0, 1, 2, \ldots$ *имеет место неравенство*

$$A_k F(x^k) \leq \varphi_k^* = \min_{x \in Q} \varphi_k(x) = \varphi_k(u^k).$$

**Доказательство.** Доказательство аналогично доказательству леммы 3.2. В основе лежит неравенство

$$V(x, y^k) \geq \frac{1}{2} \| x - y^k \|^2,$$

с помощью которого ключевая формула (3.8) перепишется следующим образом:

$$\varphi_{k+1}^* \geq A_{k+1} \cdot \left( \frac{A_k}{A_{k+1}} h(x^k) + \frac{\alpha_{k+1}}{A_{k+1}} h(u^{k+1}) \right) + A_{k+1} f(y^k) +$$

$$+ \left\langle \; \nabla f(y^{k+1}), \alpha_{k+1} \cdot (u^{k+1} - y^{k+1}) + A_k \cdot (x^k - y^{k+1}) \; \right\rangle + \frac{(1 + A_k \tilde{\mu})}{2} \| u^{k+1} - u^k \|^2.$$

Отличие от формулы (3.8) в следующем:

$$\frac{1}{2}\left\|u^{k+1}-u^k\right\|^2 \to \frac{(1+A_k\tilde{\mu})}{2}\left\|u^{k+1}-u^k\right\|^2.$$

Рассуждая дальше точно так же, как при доказательстве леммы 3.2, получим

$$\varphi_{k+1}^{*} \geq A_{k+1}\Big[f\left(y^{k+1}\right)+\left\langle\nabla f\left(y^{k+1}\right),x^{k+1}-y^{k+1}\right\rangle+ \\ +\frac{A_{k+1}\cdot(1+A_k\tilde{\mu})}{2\alpha_{k+1}^2}\left\|x^{k+1}-y^{k+1}\right\|^2+h\left(x^{k+1}\right)\Big]. \qquad (3.17)$$

Из предположения 3.1 следует, что если $A_{k+1}\cdot(1+A_k\tilde{\mu})/\alpha_{k+1}^2 \geq L$, то

$$f\left(y^{k+1}\right)+\left\langle\nabla f\left(y^{k+1}\right),x^{k+1}-y^{k+1}\right\rangle+\frac{A_{k+1}\cdot(1+A_k\tilde{\mu})}{2\alpha_{k+1}^2}\left\|x^{k+1}-y^{k+1}\right\|^2 \geq f\left(x^{k+1}\right). \qquad (3.18)$$

Согласно (3.15) $A_{k+1}\cdot(1+A_k\tilde{\mu})/\alpha_{k+1}^2 = L$, поэтому формула (3.18) имеет место. С помощью формулы (3.18) формулу (3.17) можно переписать в виде

$$\varphi_{k+1}^{*} \geq A_{k+1}\left[f\left(x^{k+1}\right)+h\left(x^{k+1}\right)\right]=A_{k+1}F\left(x^{k+1}\right). \blacksquare$$

Из лемм 3.3, 3.4 получаем следующий результат, означающий, что МПТ в сильно выпуклом случае сходится как обычный быстрый градиентный метод (с двумя проекциями), т. е. МПТ сходится оптимальным образом для рассматриваемого класса задач.

**Теорема 3.2.** *Пусть справедливы предположения* 3.1, 3.2. *Тогда МПТ* (3.14), (3.15), (3.4) *для задачи* (3.1) *сходится согласно оценке*

$$F\left(x^N\right)-\min_{x\in Q}F(x) \leq LR^2\exp\left(-\frac{N}{2}\sqrt{\frac{\tilde{\mu}}{L}}\right). \qquad (3.19)$$

**Доказательство.** Из леммы 3.4 следует, что (в третьем неравенстве используется то, что $\tilde{\mu}=\mu/\tilde{\omega}_n$ и сильная выпуклость функции $f(x)$, см. формулу (3.13) предположения 3.2):

$$A_N F\left(x^N\right) \leq \min_{x\in Q}\left\{V\left(x,y^0\right)+\sum_{k=0}^{N}\alpha_k\left[f\left(y^k\right)+\left\langle\nabla f\left(y^k\right),x-y^k\right\rangle+\tilde{\mu}V\left(x,y^k\right)+h(x)\right]\right\} \leq$$

$$\leq V\left(x_*,y^0\right)+\sum_{k=0}^{N}\alpha_k\underbrace{\left[f\left(y^k\right)+\left\langle\nabla f\left(y^k\right),x_*-y^k\right\rangle+\frac{\mu}{2}\left\|x_*-y^k\right\|^2+h(x_*)\right]}_{\leq f(x_*)+h(x_*)} \leq$$

$$\leq V\left(x_*,y^0\right)+\sum_{k=0}^{N}\alpha_k F(x_*)=R^2+A_N F(x_*). \qquad (3.20)$$

Для того чтобы получить оценку (3.19), осталось подставить нижнюю оценку на $A_N$ из леммы 3.3 в формулу (3.20). ■

В действительности, ранее установлено более сильное утверждение.

**Теорема 3.3.** *Пусть справедливы предположения* 3.1, 3.2. *Тогда МПТ* (3.14), (3.16), (3.4) *для задачи* (3.1) *сходится согласно оценке*

$$F\left(x^N\right)-\min_{x\in Q}F\left(x\right)\le\min\left\{\frac{4LR^2}{\left(N+1\right)^2},LR^2\exp\left(-\frac{N}{2}\sqrt{\frac{\mu}{L\tilde{\omega}_n}}\right)\right\}. \qquad (3.21)$$

Теорема 3.3 означает, что МПТ (3.14), (3.16), (3.4) непрерывен по параметру $\mu$. К сожалению, при этом в (3.16) явно входит данный параметр $\mu$. Если значение этого параметра неизвестно, то с помощью рестартов (см., например, [110]) можно получить оценку (3.21), увеличив константы не более чем в 4 раза, т. е. число вычислений градиента $\nabla f\left(x\right)$ (обычно именно это является самым затратным в шаге), необходимых для достижения заданной точности, увеличится не более чем в 4 раза.

Между сильно выпуклым и просто выпуклым случаями имеется глубокая связь, позволяющая, например, получить оценку (3.19) с помощью оценки (3.11), и наоборот. Другими словами, имея эффективные алгоритмы решения выпуклых / сильно выпуклых задач, можно предложить на их базе алгоритмы решения сильно выпуклых / выпуклых задач. Покажем это (приводимые далее конструкции давно и хорошо известны, мы здесь их приводим для полноты изложения).

Введем семейство $\mu$-сильно выпуклых в норме $\|\ \|$ задач ($\mu>0$):

$$F^{\mu}\left(x\right)=F\left(x\right)+\mu V\left(x,y^0\right)\to\min_{x\in Q}. \qquad (3.22)$$

**Теорема 3.4.** *Пусть*

$$\mu\le\frac{\varepsilon}{2V\left(x_*,y^0\right)}=\frac{\varepsilon}{2R^2}, \qquad (3.23)$$

*и удалось найти* $\varepsilon/2$*-решение задачи* (3.22), *т. е. нашелся такой* $x^N\in Q$, *что*

$$F^{\mu}\left(x^N\right)-F_*^{\mu}\le\varepsilon/2.$$

*Тогда*

$$F\left(x^N\right)-\min_{x\in Q}F\left(x\right)=F\left(x^N\right)-F_*\le\varepsilon.$$

**Доказательство.** Действительно,

$$F\left(x^N\right)-F_*\le F^\mu\left(x^N\right)-F_*\le F^\mu\left(x^N\right)-F_*^\mu+\varepsilon/2\le\varepsilon .$$

Здесь использовалось определение $F_*^\mu$ и формула (3.23):

$$F_*^\mu=\min_{x\in Q}\left\{F(x)+\gamma V\left(x,y^0\right)\right\}\le F\left(x_*\right)+\gamma V\left(x_*,y^0\right)\le F_*+\varepsilon/2 . \blacksquare$$

Приведем в некотором смысле обратную конструкцию.

**Теорема 3.5.** *Пусть функция* $F(x)$ – $\mu$*-сильно выпуклая в норме* $\|\ \|$. *Пусть точка* $x^{\bar N}\left(y^0\right)$ *выдается МПТ* (3.2) – (3.4), *стартующим с точки* $y^0$, *после*

$$\bar N=\sqrt{\frac{8L\omega_n}{\mu}} \tag{3.24}$$

*итераций, где*

$$\omega_n=\sup_{x\in Q}\frac{2V\left(x,y^0\right)}{\left\|x-y^0\right\|^2}.$$

*Положим*

$$\left[x^{\bar N}\left(y^0\right)\right]^1=x^{\bar N}\left(y^0\right)$$

*и определим по индукции*

$$\left[x^{\bar N}\left(y^0\right)\right]^{k+1}=x^{\bar N}\left(\left[x^{\bar N}\left(y^0\right)\right]^k\right),\ k=1,2,...$$

*При этом на* $k+1$ *перезапуске* (*рестарте*) *также корректируется прокс-функция* (*считаем, что так определенная функция корректно определена на* $Q$ *с сохранением свойства сильной выпуклости*):

$$d^{k+1}(x)=d\left(x-\left[x^{\bar N}\left(y^0\right)\right]^k+y^0\right)\ge 0,$$

*чтобы*

$$d^{k+1}\left(\left[x^{\bar N}\left(y^0\right)\right]^k\right)=0 .$$

*Тогда*

$$F\left(\left[x^{\bar N}\left(y^0\right)\right]^k\right)-F_*\le\frac{\mu\left\|y^0-x_*\right\|^2}{2^{k+1}} . \tag{3.25}$$

**Доказательство.** МПТ (3.2) – (3.4) согласно теореме 3.1 (см. формулу (3.11)) после $\bar N$ итераций выдает такой $x^{\bar N}$, что

$$\frac{\mu}{2}\left\|x^{\bar N}-x_*\right\|^2\le F\left(x^{\bar N}\right)-F_*\le\frac{4LV\left(x_*,y^0\right)}{\bar N^2} .$$

Отсюда имеем

$$\left\| x^{\bar{N}} - x_* \right\|^2 \le \frac{8LV\left(x_*, y^0\right)}{\mu \bar{N}^2} \le \frac{1}{2}\left\| y^0 - x_* \right\|^2 \frac{8L\omega_n}{\mu \bar{N}^2} .$$

Поскольку

$$\bar{N} = \sqrt{\frac{8L}{\mu}\omega_n} ,$$

то

$$\left\| x^{\bar{N}} - x_* \right\|^2 \le \frac{1}{2}\left\| y^0 - x_* \right\|^2 .$$

Повторяя эти рассуждения, по индукции получим

$$F\left(\left[x^{\bar{N}}\left(y_0\right)\right]^k\right) - F_* \le \left(\frac{1}{2}\right)^k \left\| y^0 - x_* \right\|^2 \frac{4L\omega_n}{\bar{N}^2} = \frac{\mu \left\| y^0 - x_* \right\|^2}{2^{k+1}} . \blacksquare$$

**Замечание 3.3.** Оценка (3.24), (3.25) в итоге получается похожей на оценку (3.19). Однако в подходе, описанном в теореме 3.5, всегда можно добиться, чтобы $\omega_n = O(\ln n)$, в том числе и в случае $\| \; \| = \| \; \|_1$ [110] (следует сравнить в этом случае с оценкой $\tilde{\omega}_n \ge n$, приведенной выше). Кроме того, предложенная в теореме 3.5 конструкция позволяет рассматривать более общий класс сильно выпуклых задач, в которых $f(x)$ – уже не обязательно сильно выпуклая функция (см. предположение 3.2), что позволяет избавиться в оценках числа обусловленности $L/\mu$ от ограничения, описанного в замечании 3.2. В ряде приложений эти степени свободы оказываются чрезвычайно важными. Однако, как показывают численные эксперименты (проводимые в случае евклидовой прокс-структуры и с априорно известной точной оценкой параметра $\mu$), подход, описанный в теореме 3.5, может проигрывать в скорости МПТ (3.14), (3.16), (3.4) один-два порядка, т. е. для достижения той же точности подход из теоремы 3.5 может потребовать до 100 раз больше арифметических операций. Из оценок это никак не следует. Но если метод из теоремы 3.5 работает по полученным верхним оценкам, то МПТ (3.14), (3.16), (3.4) (с критерием останова связанным с контролем малости нормы градиентного отображения) из-за отсутствия рестартов (т. е. необходимости делать предписанное число итераций) может остановиться раньше, что и происходит на практике. К тому же МПТ (3.14), (3.16), (3.4) еще и непрерывен по параметру $\mu$ .

Из вышеизложенного следует, что в процедуре рестартов (теорема 3.5) можно использовать вместо предписанного числа итераций $\bar{N}$ на каждом рестарте какой-нибудь критерий останова. В частности, дожи-даться, когда норма (или квадрат нормы) градиента (а в общем случае, когда минимум достигается не в точке экстремума, – норма градиентного отображения) уменьшиться вдвое. С таким критерием останова нет необходимости делать предписанного числа итераций на каждом рестарте. Однако пока не известен способ рассуждений, который позволял бы показать, что такая процедура сохраняет при перенесении оптимальность оценок (метод, работающий оптимально не в сильно выпуклом случае, порождает оптимальный метод и в сильно выпуклом случае). Впрочем, имеются различные эффективные на практике способы более раннего выхода с каждого рестарта (см., например, [139]), позволяющие (в случае задач безусловной оптимизации с евклидовой прокс-структурой и отсутствием композитного члена) ускорить описанную выше конструкцию на порядок.

## 3.4. Универсальный метод подобных треугольников

В ряде приложений значение константы $L$, необходимой МПТ для работы (см. формулу (3.16)), неизвестно. Однако, как следует из формул (3.10), (3.18), знание константы $L$ необязательно, если разрешается на одной итерации запрашивать значение функции в нескольких точках. Опишем соответствующий адаптивный вариант МПТ (3.14), (3.16), (3.4) (АМПТ).

Положим $A_0 = \alpha_0 = 1/L_0^0$ – параметр метода (считаем здесь и везде в дальнейшем $L_0^0 \le L$, иначе во всех приводимых далее оценках следует полагать $L := \max\left\{L_0^0, L\right\}$),

$$k = 0,\; j_0 = 0;\; x^0 = u^0 = \arg\min_{x\in Q} \varphi_0(x).$$

До тех пор пока

$$f\left(x^0\right) > f\left(y^0\right) + \left\langle \nabla f\left(y^0\right), x^0 - y^0 \right\rangle + \frac{L_0^{j_0}}{2}\left\|x^0 - y^0\right\|^2,$$

где

$$x^0 := u^0 := \arg\min_{x\in Q} \varphi_0(x),\; (A_0 :=)\,\alpha_0 := \frac{1}{L_0^{j_0}},$$

выполнять

$$j_0 := j_0 + 1;\; L_0^{j_0} := 2^{j_0} L_0^0.$$

4. $L_{k+1}^0 = L_k^{j_k}/2$, $j_{k+1} = 0$.

5. $$\begin{cases} \alpha_{k+1} := \dfrac{1 + A_k\tilde{\mu}}{2L_{k+1}^{j_{k+1}}} + \sqrt{\dfrac{1 + A_k\tilde{\mu}}{4\left(L_{k+1}^{j_{k+1}}\right)^2} + \dfrac{A_k \cdot (1 + A_k\tilde{\mu})}{L_{k+1}^{j_{k+1}}}},\ A_{k+1} := A_k + \alpha_{k+1}; \\ y^{k+1} := \dfrac{\alpha_{k+1}u^k + A_k x^k}{A_{k+1}},\ u^{k+1} := \arg\min_{x\in Q} \varphi_{k+1}(x),\ x^{k+1} := \dfrac{\alpha_{k+1}u^{k+1} + A_k x^k}{A_{k+1}}. \end{cases} \quad (*)$$

До тех пор пока

$$f\left(y^{k+1}\right) + \left\langle \nabla f\left(y^{k+1}\right), x^{k+1} - y^{k+1}\right\rangle + \frac{L_{k+1}^{j_{k+1}}}{2}\left\|x^{k+1} - y^{k+1}\right\|^2 < f\left(x^{k+1}\right),$$

выполнять

$$j_{k+1} := j_{k+1} + 1;\ L_{k+1}^{j_{k+1}} = 2^{j_{k+1}} L_{k+1}^0;\ (*).$$

6. Если не выполнен критерий останова, то $k := k+1$ и **go to** 1.

В качестве критерия останова, например, можно брать условие

$$\left\| x^{k+1} - \arg\min_{x\in Q}\left\{ \left\langle \nabla f\left(x^{k+1}\right), x - x^{k+1}\right\rangle + \frac{L_{k+1}^{j_{k+1}}}{2}\left\|x - x^{k+1}\right\|^2 \right\}\right\| \le \tilde{\varepsilon}.$$

Здесь и везде в дальнейшем под «до тех пор пока … выполнять …» подразумевается, что после каждого $j_{k+1} := j_{k+1} + 1$ при следующей проверке условия выхода из этого цикла должным образом меняется не только $L_{k+1}^{j_{k+1}}$, но и $x^{k+1}$, $y^{k+1}$, также входящие в это условие.

**Теорема 3.6.** *Пусть справедливы предположения* 3.1, 3.2. *Тогда АМПТ для задачи* (3.1) *сходится согласно оценке*

$$F\left(x^N\right) - \min_{x\in Q} F(x) \le \min\left\{ \frac{8LR^2}{(N+1)^2}, 2LR^2 \exp\left( -\frac{N}{2}\sqrt{\frac{\mu}{2L\tilde{\omega}_n}} \right) \right\}. \quad (3.26)$$

*При этом среднее число вычислений значения функции на одной итерации будет* $\approx 4$, *а градиента функции* $\approx 2$.

**Доказательство.** Нетривиальным ввиду оценки (3.21) и свойства, что все $L_k^{j_k} \le 2L$, представляется только последняя часть формулировки теоремы. Докажем именно её. Оценим общее число обращений за значениями функции (аналогично получается оценка общего числа обращений за значением градиента функции):

$$\sum_{k=1}^{N} 2\left(j_k+1\right)=\sum_{k=1}^{N} 2\left[\left(j_k-1\right)+2\right]=\sum_{k=1}^{N} 2\left[\log_2\left(\frac{L_k^{j_k}}{L_{k-1}^{j_{k-1}}}\right)+2\right]=$$

$$=4N+\log_2\left(\frac{L_N^{j_N}}{L_0^0}\right)\le 4N+\log_2\left(\frac{2L}{L_0^0}\right).$$

Деля обе части на $N$, получим в правой части приблизительно 4. ■

В действительности, оценка (3.26) оказывается, как правило, сильно завышенной, поскольку метод адаптивно настраивается на константу Липшица градиента $L$ на данном участке своего пребывания, а константа $L$, входящая в оценку (3.26), соответствует (согласно предположению 3.1) самому плохому случаю (самому плохому участку).

Предположим теперь, что по каким-то причинам невозможно получить точные значения функции и градиента. Тогда соотношение (аналогичное неравенство выписывается и при $k=0$, см. начало этого раздела)

$$f\left(y^{k+1}\right)+\left\langle\nabla f\left(y^{k+1}\right),x^{k+1}-y^{k+1}\right\rangle+\frac{L_{k+1}^{j_{k+1}}}{2}\left\|x^{k+1}-y^{k+1}\right\|^2\ge f\left(x^{k+1}\right)$$

может не выполниться не при каком $L_{k+1}^{j_{k+1}}$. Допустим, однако, что при этом имеет место

$$f\left(y^{k+1}\right)+\left\langle\nabla f\left(y^{k+1}\right),x^{k+1}-y^{k+1}\right\rangle+\frac{L}{2}\left\|x^{k+1}-y^{k+1}\right\|^2+\frac{\alpha_{k+1}}{A_{k+1}}\varepsilon\ge f\left(x^{k+1}\right).$$

Тогда заменим в АМПТ соответствующую часть шага 2 на

$$f\left(y^{k+1}\right)+\left\langle\nabla f\left(y^{k+1}\right),x^{k+1}-y^{k+1}\right\rangle+\frac{L_{k+1}^{j_{k+1}}}{2}\left\|x^{k+1}-y^{k+1}\right\|^2+\frac{\alpha_{k+1}}{A_{k+1}}\varepsilon\ge f\left(x^{k+1}\right). \quad (3.27)$$

**Теорема 3.7.** *Пусть справедливо предположение* 3.2 *и существует такое число* $L>0$, *что любого* $k=1,\dots,N$ *справедливо неравенство*

$$f\left(y^{k+1}\right)+\left\langle\nabla f\left(y^{k+1}\right),x^{k+1}-y^{k+1}\right\rangle+\frac{L}{2}\left\|x^{k+1}-y^{k+1}\right\|^2+\frac{\alpha_{k+1}}{A_{k+1}}\varepsilon\ge f\left(x^{k+1}\right). \quad (3.28)$$

*Тогда АМПТ с* (3.27) *для задачи* (3.1) *сходится согласно оценке*

$$F\left(x^N\right)-\min_{x\in Q}F\left(x\right)\le\min\left\{\frac{8LR^2}{\left(N+1\right)^2},2LR^2\exp\left(-\frac{N}{2}\sqrt{\frac{\mu}{2L\tilde{\omega}_n}}\right)\right\}+\varepsilon. \quad (3.29)$$

*При этом среднее число вычислений значения функции на одной итерации будет* $\approx 4$, *а градиента функции* $\approx 2$

**Доказательство.** Ключевым элементом в доказательстве является следующее уточнение леммы 4:

$$A_kF\left(x^k\right)\le\varphi_k^*+A_k\varepsilon, \quad (3.30)$$

из которого будет следовать формула (3.29). Чтобы доказать (3.30), будем рассуждать по индукции. База индукции $k=0$ очевидна. Итак, по предположению индукции

$$A_k F\left(x^k\right) - A_k \varepsilon \le \varphi_k^* = \varphi_k\left(u^k\right),$$

поэтому

$$\varphi_k\left(u^{k+1}\right) \ge \varphi_k\left(u^k\right) + \frac{1+A_k\tilde{\mu}}{2}\left\|u^{k+1}-u^k\right\|^2 \ge A_k F\left(x^k\right) - A_k\varepsilon + \frac{1+A_k\tilde{\mu}}{2}\left\|u^{k+1}-u^k\right\|^2.$$

Отсюда

$$\varphi_{k+1}^* \ge A_{k+1}\cdot\left(\frac{A_k}{A_{k+1}}h\left(x^k\right) + \frac{\alpha_{k+1}}{A_{k+1}}h\left(u^{k+1}\right)\right) + A_{k+1}f\left(y^k\right) - A_k\varepsilon$$

$$+\left\langle \nabla f\left(y^{k+1}\right), \alpha_{k+1}\cdot\left(u^{k+1}-y^{k+1}\right) + A_k\cdot\left(x^k - y^{k+1}\right)\right\rangle + \frac{\left(1+A_k\tilde{\mu}\right)}{2}\left\|u^{k+1}-u^k\right\|^2.$$

Следовательно,

$$\varphi_{k+1}^* + A_k\varepsilon \ge A_{k+1}\times$$

$$\times\left[f\left(y^{k+1}\right) + \left\langle \nabla f\left(y^{k+1}\right), x^{k+1}-y^{k+1}\right\rangle + \frac{A_{k+1}\cdot\left(1+A_k\tilde{\mu}\right)}{2\alpha_{k+1}^2}\left\|x^{k+1}-y^{k+1}\right\|^2 + h\left(x^{k+1}\right)\right].$$

Отсюда и из условия (3.27), $A_{k+1}\cdot\left(1+A_k\tilde{\mu}\right)/\alpha_{k+1}^2 = L_{k+1}^{j_{k+1}}$ (с учетом (3.28)) получаем

$$\varphi_{k+1}^* + A_{k+1}\varepsilon = \varphi_{k+1}^* + A_k\varepsilon + \alpha_{k+1}\varepsilon \ge A_{k+1}F\left(x^{k+1}\right),\ L_{k+1}^{j_{k+1}} \le 2L. \ \blacksquare$$

В действительности, выше установлено более сильное утверждение – в оценке (3.29) можно улучшить константу $L$.

**Теорема 3.8.** *Пусть справедливо предположение* 3.2. *Тогда АМПТ с* (3.27) *для задачи* (3.1) *сходится согласно оценке*

$$F\left(x^N\right) - \min_{x\in Q} F\left(x\right) \le \frac{R^2}{A_N} + \varepsilon \le \min\left\{\frac{4LR^2}{\left(N+1\right)^2}, LR^2\exp\left(-\frac{N}{2}\sqrt{\frac{\mu}{L\tilde{\omega}_n}}\right)\right\} + \varepsilon, \quad (3.31)$$

*где* $L = \max_{k=0,\dots,N} L_k^{j_k}$. *При этом среднее число вычислений значения функции на одной итерации будет* $\approx 4$, *а градиента функции* $\approx 2$.

Попробуем «сыграть» на условии (3.27), искусственно вводя неточность.

**Лемма 3.5** (см. [132]). *Пусть*

$$\left\|\nabla f\left(y\right) - \nabla f\left(x\right)\right\|_* \le L_\nu \left\|y-x\right\|^\nu \quad (3.32)$$

*при некотором* $\nu \in [0,1]$. *Тогда*

$$f(y)+\langle \nabla f(y), x-y\rangle + \frac{L}{2}\|x-y\|^2 + \delta \ge f(x), \; L = L_\nu \left[\frac{L_\nu}{2\delta}\frac{1-\nu}{1+\nu}\right]^{\frac{1-\nu}{1+\nu}}. \qquad (3.33)$$

Основным результатом данного раздела является описание и последующая оценка скорости сходимости нового варианта универсального быстрого (ускоренного) градиентного метода на базе МПТ (УМПТ).

Положим

$$A_0 = \alpha_0 = 1/L_0^0, \; k=0, \; j_0 = 0; \; x^0 = u^0 = \arg\min_{x\in Q} \varphi_0(x).$$

До тех пор пока

$$f(x^0) > f(y^0) + \langle \nabla f(y^0), x^0 - y^0\rangle + \frac{L_0^{j_0}}{2}\|x^0 - y^0\|^2 + \frac{\alpha_0}{2A_0}\varepsilon,$$

где

$$x^0 := u^0 := \arg\min_{x\in Q}\varphi_0(x), \; (A_0 :=)\alpha_0 := \frac{1}{L_0^{j_0}},$$

выполнять

$$j_0 := j_0 + 1; \; L_0^{j_0} := 2^{j_0} L_0^0.$$

**Универсальный метод подобных треугольников**

7. $L_{k+1}^0 = L_k^{j_k}/2, \; j_{k+1} = 0$.

8. $$\begin{cases} \alpha_{k+1} := \dfrac{1 + A_k\tilde{\mu}}{2L_{k+1}^{j_{k+1}}} + \sqrt{\dfrac{1 + A_k\tilde{\mu}}{4\left(L_{k+1}^{j_{k+1}}\right)^2} + \dfrac{A_k\cdot(1+A_k\tilde{\mu})}{L_{k+1}^{j_{k+1}}}}, \; A_{k+1} := A_k + \alpha_{k+1}; \\ y^{k+1} := \dfrac{\alpha_{k+1}u^k + A_k x^k}{A_{k+1}}, \; u^{k+1} := \arg\min_{x\in Q}\varphi_{k+1}(x), \; x^{k+1} := \dfrac{\alpha_{k+1}u^{k+1} + A_k x^k}{A_{k+1}}. \end{cases} \qquad (*)$$

До тех пор пока

$$f(y^{k+1}) + \langle \nabla f(y^{k+1}), x^{k+1} - y^{k+1}\rangle + \frac{L_{k+1}^{j_{k+1}}}{2}\|x^{k+1} - y^{k+1}\|^2 + \frac{\alpha_{k+1}}{2A_{k+1}}\varepsilon < f(x^{k+1}),$$

выполнять

$$j_{k+1} := j_{k+1} + 1; \; L_{k+1}^{j_{k+1}} = 2^{j_{k+1}} L_{k+1}^0; \; (*).$$

3. Если не выполнен критерий останова, то $k := k+1$ и **go to** 1.

**Теорема 3.9.** *Пусть выполняется условие* (3.32) *хотя бы для* $\nu = 0$ *и справедливо предположение* 3.2 *с* $\mu \ge 0$ (*допускается брать* $\mu = 0$). *Тогда УМПТ для задачи* (3.1) *сходится согласно оценке*

$$F\left(x^N\right)-\min_{x\in Q}F\left(x\right)\le\varepsilon ,$$

$$N\approx\min\left\{\inf_{\nu\in[0,1]}\left(\frac{L_\nu\cdot(16R)^{1+\nu}}{\varepsilon}\right)^{\frac{2}{1+3\nu}},\right.$$

$$\left.\inf_{\nu\in[0,1]}\left\{\left(\frac{8L_\nu^{\frac{2}{1+\nu}}\tilde{\omega}_n}{\mu\varepsilon^{\frac{1-\nu}{1+\nu}}}\right)^{\frac{1+\nu}{1+3\nu}}\ln^{\frac{2+2\nu}{1+3\nu}}\left(\frac{16L_\nu^{\frac{4+6\nu}{1+\nu}}R^2}{\left(\mu/\tilde{\omega}_n\right)^{\frac{1+\nu}{1+3\nu}}\varepsilon^{\frac{5+7\nu}{2+6\nu}}}\right)\right\}\right\}. \tag{3.34}$$

*При этом среднее число вычислений значения функции на одной итерации будет* $\approx 4$, *а градиента функции* $\approx 2$.

**Доказательство.** Рассмотрим два случая, когда $\mu\ge 0$ – мало: $\mu\ll\varepsilon/\left(2R^2\right)$, $\mu$ – велико: $\mu\gg\varepsilon/\left(2R^2\right)$, см. формулу (3.23).

В первом случае будем считать, что

$$A_{k+1}/\alpha_{k+1}^2\approx A_{k+1}\cdot\left(1+A_k\tilde{\mu}\right)/\alpha_{k+1}^2=L_{k+1}^{j_{k+1}},\text{ т. е. }\varepsilon\frac{\alpha_{k+1}}{2A_{k+1}}\approx\frac{\varepsilon}{2}\sqrt{\frac{1}{L_{k+1}^{j_{k+1}}A_{k+1}}}, \tag{3.35}$$

а во втором случае

$$A_{k+1}^2\tilde{\mu}/\alpha_{k+1}^2\approx A_{k+1}\cdot\left(1+A_k\tilde{\mu}\right)/\alpha_{k+1}^2=L_{k+1}^{j_{k+1}},\text{ т. е. }\varepsilon\frac{\alpha_{k+1}}{2A_{k+1}}\approx\frac{\varepsilon}{2}\sqrt{\frac{\tilde{\mu}}{L_{k+1}^{j_{k+1}}}}. \tag{3.36}$$

Из формулы (3.31) (см. теорему 3.8) имеем

$$\frac{R^2}{A_N}+\frac{\varepsilon}{2}\approx\varepsilon ,$$

т. е. $A_N\approx 2R^2/\varepsilon$, а также ( $L=\max_{k=0,\dots,N}L_k^{j_k}$ )

$$N^2\approx\frac{8LR^2}{\varepsilon}\text{ (в первом случае)}, \tag{3.37}$$

$$N^2\approx 4\frac{L}{\tilde{\mu}}\ln^2\left(\frac{2LR^2}{\varepsilon}\right)\text{ (во втором случае)}. \tag{3.38}$$

Из формул (3.33), (3.35) – (3.38) имеем, что в первом случае

$$L\le 2L_\nu\left[\frac{L_\nu}{2\frac{\varepsilon}{2}\sqrt{\frac{1}{LA_N}}}\frac{1-\nu}{1+\nu}\right]^{\frac{1-\nu}{1+\nu}}\le 2L_\nu\left[\frac{L_\nu\sqrt{LA_N}}{\varepsilon}\frac{1-\nu}{1+\nu}\right]^{\frac{1-\nu}{1+\nu}}\le 2L_\nu^{\frac{2}{1+\nu}}\left[\frac{N}{2\varepsilon}\frac{1-\nu}{1+\nu}\right]^{\frac{1-\nu}{1+\nu}}, \tag{3.39}$$

а во втором случае

$$L \le 2L_\nu \left[ \frac{L_\nu}{2\frac{\varepsilon}{2}\sqrt{\frac{\tilde{\mu}}{L}}} \frac{1-\nu}{1+\nu} \right]^{\frac{1-\nu}{1+\nu}} \le 2L_\nu \left[ \frac{L_\nu\sqrt{L/\tilde{\mu}}}{\varepsilon} \frac{1-\nu}{1+\nu} \right]^{\frac{1-\nu}{1+\nu}} \le 2L_\nu^{\frac{2}{1+\nu}} \left[ \frac{N}{2\varepsilon} \frac{1-\nu}{1+\nu} \right]^{\frac{1-\nu}{1+\nu}}. \qquad (3.40)$$

Подставляя (3.39) в (3.37), а (3.40) в (3.38) и учитывая, что параметр $\nu \in [0,1]$ можно выбирать произвольно (допускается, что $L_\nu = \infty$ при некоторых $\nu$, важно, чтобы существовало хотя бы одно значение $\nu$, при котором $L_\nu < \infty$; по условию $L_0 < \infty$), получим соответственно

$$N^2 \approx \frac{16L_\nu^{\frac{2}{1+\nu}} \left[ \frac{N}{2\varepsilon} \frac{1-\nu}{1+\nu} \right]^{\frac{1-\nu}{1+\nu}} R^2}{\varepsilon} \Rightarrow N^{\frac{1+3\nu}{1+\nu}} \approx \frac{16L_\nu^{\frac{2}{1+\nu}} R^2}{\varepsilon^{\frac{2}{1+\nu}}} \Rightarrow$$

$$\Rightarrow N \approx \inf_{\nu\in[0,1]} \left( \frac{L_\nu \cdot (16R)^{1+\nu}}{\varepsilon} \right)^{\frac{2}{1+3\nu}}, \qquad (3.41)$$

$$N^2 \approx \frac{8L_\nu^{\frac{2}{1+\nu}} \left[ \frac{N}{2\varepsilon} \frac{1-\nu}{1+\nu} \right]^{\frac{1-\nu}{1+\nu}}}{\tilde{\mu}} \ln^2\left( \frac{2L_\nu^2 R^2 N}{\varepsilon^{3/2}} \right) \Rightarrow N^{\frac{1+3\nu}{1+\nu}} \approx \frac{8L_\nu^{\frac{2}{1+\nu}}}{\tilde{\mu}\varepsilon^{\frac{1-\nu}{1+\nu}}} \ln^2\left( \frac{2L_\nu^2 R^2 N}{\varepsilon^{3/2}} \right)$$

$$\Rightarrow N \approx \left( \frac{8L_\nu^{\frac{2}{1+\nu}}}{\tilde{\mu}\varepsilon^{\frac{1-\nu}{1+\nu}}} \right)^{\frac{1+\nu}{1+3\nu}} \ln^{\frac{2+2\nu}{1+3\nu}}\left( \frac{2L_\nu^2 R^2 N}{\varepsilon^{3/2}} \right) \Rightarrow$$

$$\Rightarrow N \approx \inf_{\nu\in[0,1]} \left\{ \left( \frac{8L_\nu^{\frac{2}{1+\nu}}}{\tilde{\mu}\varepsilon^{\frac{1-\nu}{1+\nu}}} \right)^{\frac{1+\nu}{1+3\nu}} \ln^{\frac{2+2\nu}{1+3\nu}}\left( \frac{16L_\nu^{\frac{4+6\nu}{1+\nu}} R^2}{\tilde{\mu}^{\frac{1+\nu}{1+3\nu}} \varepsilon^{\frac{5+7\nu}{2+6\nu}}} \right) \right\}. \qquad (3.42)$$

Из формул (3.41), (3.42) получаем оценку (3.34). ■

Более аккуратные рассуждения позволяют в несколько раз уменьшить константы, входящие в оценку (3.34). Оценка (3.34) согласуется с нижними оценками для соответствующих классов задач [103].

# ЛИТЕРАТУРА

**Гасников** Александр Владимирович
**Гасникова** Евгения Владимировна

# МОДЕЛИ
# РАВНОВЕСНОГО РАСПРЕДЕЛЕНИЯ
# ТРАНСПОРТНЫХ ПОТОКОВ
# В БОЛЬШИХ СЕТЯХ